\documentclass[11pt]{amsart}

\allowdisplaybreaks[1]
\def\activecomma{\mathchar"602C \penalty0\hskip0pt\relax}
\begingroup\catcode`\,=\active
\gdef\enablecommabreak{\mathcode`\,="8000 \let,\activecomma}
\endgroup
\everymath{\enablecommabreak}

\usepackage[margin=1in]{geometry}
\usepackage{graphicx}
\usepackage{tikz}
\usetikzlibrary{positioning}
\usepackage{longtable}
\usepackage{array}
\usepackage{enumerate}
\usepackage{amsmath,amssymb}
\usepackage[hidelinks]{hyperref}
\usepackage{hyperxmp}
\hypersetup{
  pdftitle={Divisible additive codes},
  pdfauthor={Sascha Kurz},
  pdfsubject={A code is Delta-divisible if all its codeword weights are divisible by Delta. We study the possible lengths of Delta-divisible additive codes, i.e., F\string_q-linear codes over F\string_{q\string^h}, in their geometric description as multisets of (h-1)-spaces of a projective space over F\string_q. We determine the base lengths, i.e., a minimal set of lengths from which all possible lengths arise as sums, for all Delta dividing q\string^h, for Delta=q\string^{h+1} with h>=2, and for several parameters with Delta=q\string^{h+2} and Delta=q\string^{h+3}. We also treat divisibility constants that are not powers of the characteristic. Over F\string_2 and F\string_3 we show that p\string^r-divisible multisets of points of a given cardinality exist and are all sums of p\string^i-fold (r-i)-spaces if and only if the cardinality satisfies an explicit condition on its digits in a suitable number system. Finally, we classify the types of partitions of 8-divisible binary point sets into points, lines, and planes for small cardinalities.},
  pdfkeywords={additive codes, divisible codes, projective systems, field reduction, lengths of codes, vector space partitions},
  pdfcopyright={Creative Commons Lizenz CC BY 4.0},
  pdflicenseurl={https://creativecommons.org}
}
\theoremstyle{plain}
\newtheorem{Theorem}{Theorem}[section]
\newtheorem{Proposition}[Theorem]{Proposition}

\newtheorem{Corollary}[Theorem]{Corollary}
\newtheorem{Lemma}[Theorem]{Lemma}

\theoremstyle{definition}
\newtheorem{Definition}[Theorem]{Definition}

\newtheorem{Example}[Theorem]{Example}

\theoremstyle{remark}
\newtheorem{Remark}[Theorem]{Remark}

\newcommand{\cB}{\mathcal{B}}
\newcommand{\cC}{\mathcal{C}}
\newcommand{\cM}{\mathcal{M}}

\newcommand{\cK}{\mathcal{K}}
\newcommand{\cN}{\mathcal{N}}
\newcommand{\cP}{\mathcal{P}}
\newcommand{\cS}{\mathcal{S}}
\newcommand{\cG}{\mathcal{G}}
\newcommand{\cT}{\mathcal{T}}
\newcommand{\cX}{\mathcal{X}}

\newcommand{\cF}{\mathcal{F}}
\newcommand{\Cone}{\operatorname{Cone}}
\newcommand{\PG}{\operatorname{PG}}
\newcommand{\N}{\mathbb{N}}
\newcommand{\Z}{\mathbb{Z}}
\newcommand{\Q}{\mathbb{Q}}

\newcommand{\F}{\mathbb{F}}
\newcommand{\sg}[1]{\langle #1\rangle_{\N_0}}
\newcommand{\wt}{\operatorname{wt}}
\makeatletter
\newcommand{\itemlab}[2]{\item[(#1)]\protected@edef\@currentlabel{\theTheorem(#1)}\label{#2}}
\makeatother
\newcommand{\qbin}[3]{\genfrac{[}{]}{0pt}{}{#1}{#2}_{#3}}

\makeatletter
\def\section{\@startsection{section}{1}%
  \z@{.7\linespacing\@plus\linespacing}{.5\linespacing}%
  {\normalfont\scshape\centering}}
\def\subsection{\@startsection{subsection}{2}%
  \z@{.5\linespacing\@plus.7\linespacing}{-.5em}%
  {\normalfont\bfseries\boldmath}}
\makeatother
\begin{document}

\title{Divisible additive codes}
\author{Sascha Kurz}
\date{}

\begin{center}
  {\large\bfseries Divisible additive codes\par}
  \bigskip
  {Sascha Kurz\par}
  \smallskip
  {\small Mathematisches Institut, Universit\"at Bayreuth, Universit\"atsstra{\ss}e 30, 95447 Bayreuth, Germany\par
  \texttt{sascha.kurz@uni-bayreuth.de}\par}
\end{center}
\bigskip

\begin{center}
\begin{minipage}{0.88\textwidth}
\small\noindent
{\scshape Abstract.} A code is $\Delta$-divisible if all its codeword weights are divisible by $\Delta$. We study the possible lengths of $\Delta$-divisible additive codes, i.e., $\F_q$-linear codes over $\F_{q^h}$, in their geometric description as multisets of $(h-1)$-spaces of a projective space over $\F_q$. We determine the base lengths, i.e., a minimal set of lengths from which all possible lengths arise as sums, for all $\Delta$ dividing $q^h$, for $\Delta=q^{h+1}$ with $h\ge 2$, and for several parameters with $\Delta=q^{h+2}$ and $\Delta=q^{h+3}$. We also treat divisibility constants that are not powers of the characteristic. Over $\F_2$ and $\F_3$ we show that $p^r$-divisible multisets of points of a given cardinality exist and are all sums of $p^i$-fold $(r-i)$-spaces if and only if the cardinality satisfies an explicit condition on its digits in a suitable number system. Finally, we classify the types of partitions of $8$-divisible binary point sets into points, lines, and planes for small cardinalities.


\noindent
\textbf{MSC classification (2020):} Primary 94B05, secondary 51E20, 05B25, 94B65.
\end{minipage}
\end{center}
\bigskip

\section{Introduction}
A code $C$ is called $\Delta$-divisible if the weight $\wt(c)$ of every codeword $c\in C$ is divisible by $\Delta$. The study of $\Delta$-divisible linear codes was initiated by Ward \cite{ward1981divisible}. 
They are a central tool for packing and covering problems for subspaces in projective spaces, see \cite{kurz2021divisible} for a survey. For instance, restrictions on the possible lengths of $\Delta$-divisible linear codes 
improve the Johnson bound for constant-dimension codes \cite{kiermaier2020lengths}. For projective codes, i.e., codes without repeated columns, they give the currently best known upper bounds on the maximum size of partial spreads, 
see \cite[Section 9]{kurz2021divisible}, and they can be used to rule out the existence of certain vector space partitions, see, e.g., \cite{heden2012survey,honold2018partial,kurz2022vector}.

The possible lengths of $\Delta$-divisible linear codes over $\F_q$ are completely known. If $\Delta$ is coprime to $q$, then every $\Delta$-divisible code is the $\Delta$-fold repetition of a smaller code 
\cite[Theorem 1]{ward1981divisible}. For $\Delta=q^r$ the lengths were determined in \cite{kiermaier2020lengths}, and the remaining case that $\Delta$ is a power of the characteristic $p$ of $\F_q$ was settled in \cite{kurz2024lengths}. 
In both cases there is a finite list of explicitly described $\Delta$-divisible codes, namely certain simplex codes and their repetitions, such that the possible lengths are exactly the non-negative integer combinations of 
their lengths. Refined questions impose additional restrictions on the codes, e.g., on the maximum column multiplicity \cite{korner2024lengths}, and several applications need full classifications of divisible codes with prescribed properties.

Additive codes have recently been studied more intensively, among other reasons because additive codes over $\F_4$ describe stabilizer quantum codes \cite{calderbank1998quantum}, see also \cite{blokhuis2004small} for small additive quaternary codes. An additive code over $\F_{q^h}$ is a subset of $\F_{q^h}^n$ that is linear over $\F_q$. 
Geometrically, it corresponds to a multiset $\cM$ of $(h-1)$-spaces in a projective space $\PG(v-1,q)$, whose members play the role of the points of a linear code, and the weight of a codeword is the length $n=\#\cM$ minus the number of members contained in the corresponding hyperplane, see \cite{ball2025griesmer}. Multisets of subspaces also occur in vector space partitions, whose existence is restricted by divisibility conditions, and divisible multisets of points that can be partitioned into subspaces arise for generalized Hamming weights of additive codes \cite{dhaeseleer2025generalized}. 
The Griesmer bound and the Solomon--Stiffler construction generalize to additive codes. For sufficiently large minimum distance the bound can always be attained \cite{kurz2024additive}, which yields infinitely many series 
of additive codes with better parameters than all linear codes. Divisibility of additive codes has so far been considered mainly in the direction of the Gleason--Pierce--Ward theorem \cite{kim2009generalized}, and we are not aware of results on the possible lengths of $\Delta$-divisible additive codes. The aim of this paper is to initiate this study, with the possible lengths as the main question.

Throughout, a $\Delta$-divisible additive code over $\F_{q^h}$ is identified with a multiset $\cM$ of subspaces of $\PG(v-1,q)$ of dimension at most $h$, and $\cM$ is \emph{faithful} if all its members have dimension exactly $h$, i.e., are $(h-1)$-spaces. Members of smaller dimension are allowed only in result (4) below. The \emph{length} of the code is $n=\#\cM$. The set of lengths of $\Delta$-divisible faithful additive codes is closed under addition, and its inclusion-minimal set of generators is denoted by $\cN_q^{\Delta,h}$. Its elements are the \emph{base lengths}: every length is a non-negative integer combination of them, and no proper subset has this property. The notions of canonical multisets and of the digit rule, which occur only in result (5), are made precise in Section~\ref{sec_isomorphism_types}. For the introduction it suffices to know that the digit rule is a condition on the digits of the cardinality $n$ in a mixed number system (Theorem~\ref{thm_digit_rule}).

\subsection*{Main results}
The core of the paper determines base lengths.
\begin{enumerate}
  \item[(1)] \emph{$\Delta$ dividing $q^h$ or $\Delta=q^{h+1}$} (Theorems~\ref{thm_delta_divides_qh}, \ref{thm_C}, and~\ref{thm_N_qh1}). If $\Delta>1$ divides $q^h$, then $\cN_q^{\Delta,h}=\{\Delta,q^h+1\}$. For $h\ge 2$ and $\Delta=q^{h+1}$ the lengths of faithful codes are exactly the integers satisfying the necessary conditions of Lemmas~\ref{lemma_exclusion} and~\ref{lemma_counting}. For $h\ge 3$ the base lengths are $\cN_q^{q^{h+1},h}=\cG_q^h$, and for $h=2$ we have $\cN_q^{q^3,2}=\cG_q^2\backslash\{q^4+q^3+q^2+1\}$, where
  $$
    \cG_q^h:=\left\{q^{h+1},\,q^{h+1}+q,\,q^{2h}+q^h+1\right\}\cup\left\{q^{h+2}+q^{h+1}+c,\,q^{h+2}+q^{h+1}+q^2+c\,:\,1\le c\le q-1\right\}.
  $$
  Over $\F_4$ the sets $\cN_2^{2^l,2}$ are determined for all $l\ge 1$ (Theorems~\ref{thm_D} and~\ref{thm_N_explicit}).
  \item[(2)] \emph{$\Delta=q^{h+2}$ and $\Delta=q^{h+3}$} (Theorem~\ref{thm_apery_qhj} and Corollary~\ref{cor_qh3}). For $\Delta=q^{h+j}$ with $1\le j\le q$ and $h$ large enough, the semigroup given by the necessary conditions and its minimal generators are determined explicitly. For the sets of Table~\ref{table_N_sets} with status ``computer'' and for $(q,h,\Delta)=(2,6,256)$, every minimal generator is known to be a length, except $1297$ and $1299$ for $(2,6,256)$, whose status is open. For the other parameters, Remark~\ref{rem_N_qh3} lists the minimal generators for which no construction is known, see Subsection~\ref{sub_qh23}.
  \item[(3)] \emph{$\Delta$ not a power of $p$} (Theorems~\ref{thm_finite_exceptions}, \ref{thm_digit_sum}, and~\ref{thm_N_t1}). For $\Delta=q^l[h]_q$ with $l\ge 0$ and $(h,l)\neq(1,0)$ all but finitely many positive integers are lengths. A refined counting bound together with a digit sum condition determines $\cN_q^{\Delta,h}$ for many parameters.
  \item[(4)] \emph{Mixed dimensions} (Theorem~\ref{thm_mixed_qpower}). For $\Delta=q^l$ with $l\ge 1$ and $h\ge 2$ and members of mixed dimensions, $n$ is the cardinality of a $\Delta$-divisible projective $h$-system if and only if $n\equiv 0\pmod\Delta$ or $n>q\rho$, where $\rho\in\{0,\dots,\Delta-1\}$ is the residue of $n$ modulo $\Delta$.
  \item[(5)] \emph{The digit rule} (Theorems~\ref{thm_digits_only_if}, \ref{thm_digits_p2}, and~\ref{thm_p3_induction}). It is necessary for every prime $p$, sufficient for $p\in\{2,3\}$, and not sufficient for $p\ge 5$.
  \item[(6)] \emph{Partitions into points, lines, and planes} (Theorem~\ref{thm_types_8div_small}). We determine which types $3^{m_3}2^{m_2}1^{m_1}$, i.e., partitions into $m_3$ planes, $m_2$ lines, and $m_1$ points, occur for $8$-divisible binary point sets of cardinality at most $51$ and give the status for $60\le n\le 70$ (Table~\ref{table_types_60_70}). Cardinalities between $52$ and $59$ do not occur \cite{honold2019lengths}.
\end{enumerate}
A key tool for results (1)--(3) is Lemma~\ref{lemma_exclusion}: a faithful $\cM$ is $\Delta$-divisible if and only if the multiset of points obtained by replacing each member by its $[h]_q$ points is $(\Delta q^{h-1})$-divisible. This transfers all restrictions on the lengths of divisible multisets of points to additive codes.

\subsection*{Open problems}
The following questions are left open.
\begin{itemize}
  \item Are $1297$ and $1299$ lengths for $(q,h,\Delta)=(2,6,256)$? Which members of the families of minimal generators for $\Delta=q^{h+2}$ and $q^{h+3}$ in Remark~\ref{rem_N_qh3} are lengths?
  \item Find a digit rule for $p\ge 5$ and for non-prime fields.
  \item Determine the types of $8$-divisible binary point sets for the cardinalities $60\le n\le 70$ with $n\ne 63$ for which Table~\ref{table_types_60_70} lists open types.
\end{itemize}

The remaining part of the paper is structured as follows. After collecting preliminaries and structural results in Section~\ref{sec_preliminaries}, we turn to canonical multisets of points in Section~\ref{sec_isomorphism_types}, where the digit rule is proven. Section~\ref{sec_length_divisible} contains the necessary conditions on the lengths and the results for divisibility constants that are not powers of $p$ and for members of different dimensions. Section~\ref{sec_length_powers}, the longest part of the paper, treats divisibility constants that are powers of $p$ by partitions into subspaces and explicit constructions. In Section~\ref{sec_partitions} we classify the partitions of $8$-divisible binary point sets into points, lines, and planes. Longer proofs, the table of known sets $\cN_q^{\Delta,h}$, and the computer verifications are deferred to the appendices.

\subsection*{Reading guide}
\begin{itemize}
  \item Section~\ref{sec_preliminaries} fixes the language and is needed throughout.
  \item The digit rule of Section~\ref{sec_isomorphism_types} can be read on its own. Section~\ref{sec_length_powers} takes only the notation for the $S_q(a,b)$-adic expansion from it.
  \item Sections~\ref{sec_length_divisible} and~\ref{sec_length_powers} belong together. The first one collects the necessary conditions and the results for divisibility constants that are not powers of $p$. The second one shows that these conditions are sharp for $\Delta=q^l$ with small $l-h$ in the cases of Table~\ref{table_N_sets} and for $(q,h,\Delta)=(2,6,256)$, except possibly for $1297$ and $1299$ in the latter case.
  \item Subsection~\ref{sub_partitions} on types and $h$-partitions is also used in Subsection~\ref{sub_nonpower}.
  \item Section~\ref{sec_partitions} uses Lemmas~\ref{lemma_mixed_counting} and~\ref{lemma_exchange} of Subsection~\ref{sub_mixed} and continues the partitions into subspaces of Section~\ref{sec_length_powers}, but its results are otherwise self-contained.
  \item Readers interested only in the base lengths $\cN_q^{\Delta,h}$ can read Sections~\ref{sec_preliminaries}, \ref{sec_length_divisible}, and~\ref{sec_length_powers}, together with Subsection~\ref{sub_lengths_points} for the notation of the $S_q(a,b)$-adic expansion.
\end{itemize}

Proofs that are deferred to the appendices are preceded by a sketch in the main text. Appendix~\ref{app_comp} lists the computations.
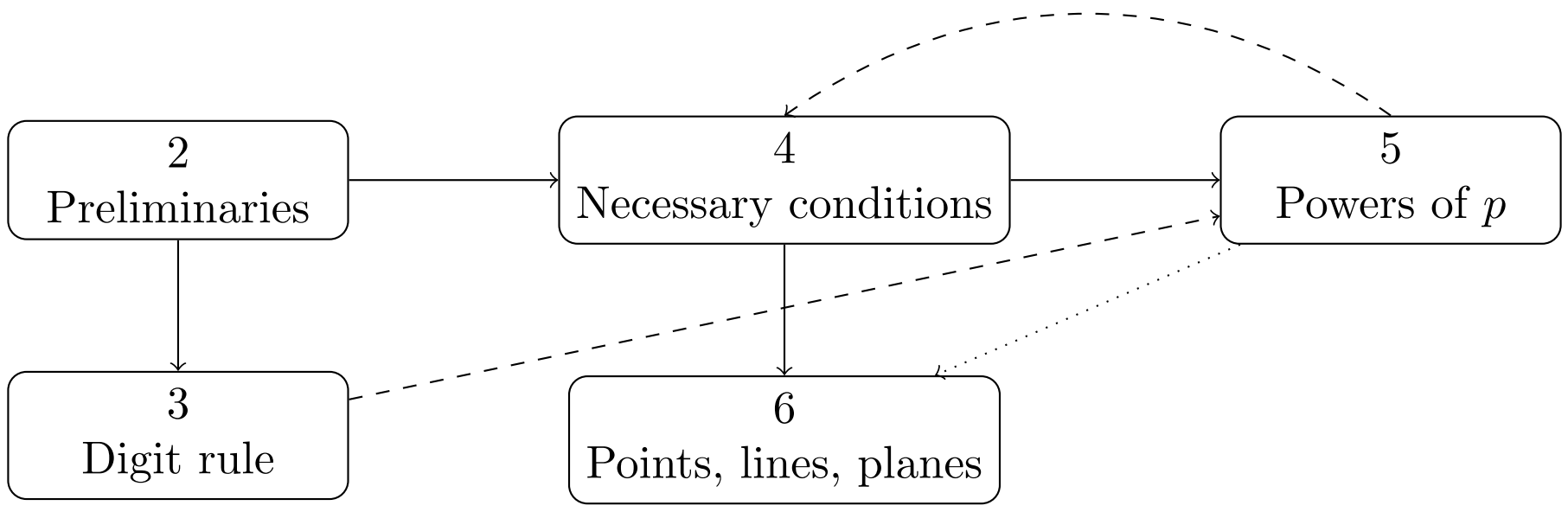
\begin{figure}[ht]
  \centering
  \begin{tikzpicture}[node distance=1.0cm and 1.6cm, every node/.style={font=\small}, box/.style={draw,rounded corners,minimum height=0.8cm,minimum width=2.6cm,align=center}]
    \node[box] (s2) {\ref*{sec_preliminaries}\\ Preliminaries};
    \node[box,right=of s2] (s4) {\ref*{sec_length_divisible}\\ Necessary conditions};
    \node[box,right=of s4] (s5) {\ref*{sec_length_powers}\\ Powers of $p$};
    \node[box,below=of s4] (s6) {\ref*{sec_partitions}\\ Points, lines, planes};
    \node[box,below=of s2] (s3) {\ref*{sec_isomorphism_types}\\ Digit rule};
    \draw[->] (s2) -- (s4);
    \draw[->] (s4) -- (s5);
    \draw[->] (s4) -- (s6);
    \draw[->,dotted] (s5) -- (s6);
    \draw[->,dashed] (s5.north) to[bend right=35] (s4.north);
    \draw[->,dashed] (s3) -- (s5);
    \draw[->] (s2) -- (s3);
  \end{tikzpicture}
  \caption{Dependencies between the sections. An arrow from $A$ to $B$ means that $B$ uses results of $A$. Dashed arrows: only notation or notions are used. Dotted arrow: only the setting is continued.}
  \end{figure}

\section{Preliminaries}
\label{sec_preliminaries}
This section provides the tools for the rest of the paper. After fixing the geometric language for additive codes, we show that for divisibility constants coprime to $q$ every divisible multiset of points is a repetition. We then collect an incidence count that holds for members of arbitrary dimensions, show that residuals inherit divisibility, and use this to characterize the Desarguesian line spreads of $\PG(5,q)$.

A linear $[n,k]_q$ code $C$ is a $k$-dimensional subspace of $\F_q^n$, and an $[n,k,d]_q$ code additionally has minimum Hamming distance $d$. A subset $C\subseteq\F_Q^n$ is called
\emph{additive} if $c+c'\in C$ for all $c,c'\in C$. Every additive code is linear over some subfield of $\F_Q$, see, e.g., \cite{ball2025griesmer}. Writing this subfield as $\F_q$ and $Q=q^h$,
an \emph{additive $[n,r/h,d]_q^h$ code} is a subset $C\subseteq\F_{q^h}^n$ that is linear over $\F_q$, has cardinality $\#C=q^r$, and has minimum Hamming distance $d$. The number
$r/h\in\tfrac{\N_0}{h}\subset\Q$ is called the \emph{dimension} of $C$. A linear $[n,k,d]_q$ code is an additive $[n,k,d]_q^1$ code. An additive code $C$ is \emph{$\Delta$-divisible}
if $\Delta\in\N_{>0}$ divides the Hamming weight $\wt(c)$ of every codeword $c\in C$.

By $\PG(v-1,q)$ we denote the projective geometry of $\F_q^v$, which has projective dimension $v-1$. An $(h-1)$-space in $\PG(v-1,q)$
is an $h$-dimensional subspace of $\F_q^v$. We call $0$-, $1$-, $2$-, $3$-, and $(v-2)$-spaces points, lines, planes, solids, and hyperplanes, respectively. By $[i]_q:=\tfrac{q^i-1}{q-1}$ we denote
the number of points of an $(i-1)$-space. A \emph{projective $h$-system} in $\PG(v-1,q)$ is a multiset $\cM$ of subspaces, each of which is an $(i-1)$-space for some $i\le h$, where $i$
may differ between members. It is \emph{spanning} if no hyperplane contains all elements of $\cM$. An additive $[n,r/h,d]_q^h$ code $C$ is in one-to-one
correspondence with a spanning projective $h$-system $\cM$ in $\PG(r-1,q)$ such that $\#\cM=n$ and each hyperplane contains at most $n-d$ elements of $\cM$, with equality for some hyperplane, see \cite{ball2025griesmer}.
Recall that $\cM$ is called faithful if all members of $\cM$ have dimension exactly $h$, and $C$ is called faithful if $\cM$ is. To simplify notation we exclude zero columns, i.e., every member of $\cM$ has dimension
at least $1$. Thus every coordinate is in the support of some codeword. For additive codes the dual minimum distance is at least $2$ if and only if $\cM$ is faithful. A projective $1$-system is just a multiset of points,
the well-known geometric description of a linear code. A multiset of points or a projective $h$-system is called $\Delta$-divisible if the corresponding additive code is.
In a list of weights $w_1^{m_1}w_2^{m_2}\cdots$ the exponent $m_i$ is the number of hyperplanes $H$ with $n-\cM(H)=w_i$, so that there are $(q-1)m_i$ codewords of weight $w_i$. The only exception is Table~\ref{table_few_weight}, where the exponents count codewords.

For a multiset of points $\cM$ in $\PG(v-1,q)$ let $\cM(P)$ denote the multiplicity of a point $P$. We extend this notion additively to each subspace $U$ by
$\cM(U):=\sum_{P\le U\,:\,\dim(P)=1} \cM(P)$. For a subspace $Y$ let $\cM|_Y$ denote the restriction of $\cM$ to the points of $Y$, a multiset of points of $Y$. Note that $\#\cM=\cM(V)$ for the ambient space $V$. Since the weight of the codeword belonging to $H$ is $\#\cM-\cM(H)$, $\cM$ is $\Delta$-divisible if and only if $\cM(H)\equiv \#\cM\pmod\Delta$
for every hyperplane $H$, provided $v\ge 2$. In $\PG(0,q)$ we call a multiset of points $\cM$ $\Delta$-divisible if and only if $\Delta$ divides $\#\cM$. We start with a geometric proof of the classic
result \cite[Theorem 1]{ward1981divisible} (see also, e.g., \cite[Theorem 3.2]{kurz2021divisible}) that each $\Delta$-divisible linear code over $\F_q$ is a $\Delta$-fold repetition of a smaller code if $\Delta$ is coprime to $q$.
\begin{Proposition}
  \label{thm_coprime}
  Let $\cM$ be a $\Delta$-divisible multiset of points in $\PG(v-1,q)$, where $v\ge 2$. If $\gcd(\Delta,q)=1$, then we have $\cM(P)\equiv 0\pmod \Delta$ for every point $P$. 
\end{Proposition}     
\begin{proof}
  Let $P$ be a point and $L$ a line through $P$. By Lemma~\ref{lemma_restriction_points} below, whose proof does not use this proposition, the restriction $\cM|_L$ is
  $\Delta/\gcd(\Delta,q^{v-2})=\Delta$-divisible, i.e., $\cM(P')\equiv\cM(L)\pmod\Delta$ for each of the $q+1$ points $P'$ of $L$. Summing over them gives
  $\cM(L)\equiv(q+1)\cM(L)\pmod\Delta$. Since $\gcd(\Delta,q)=1$, this implies $\cM(L)\equiv 0$ and thus $\cM(P)\equiv 0\pmod\Delta$.
\end{proof}

By $\qbin{v}{h}{q}$ we denote the number of $(h-1)$-spaces in $\PG(v-1,q)$. We have $\qbin{v}{h}{q}=\prod_{i=1}^{h} \tfrac{[v+1-i]_q}{[i]_q}$.

\medskip
\noindent\emph{Incidence counting.}\ Let $\cM$ be a projective $h$-system in $\PG(v-1,q)$, i.e., a multiset of subspaces, which we call \emph{members}, where a member $U$ has dimension $k_U\in\{1,\dots,h\}$, and let $n=\#\cM$ and $V=\F_q^v$.
For a subspace $X$ let $\widetilde{\cM}(X):=\#\{U\in\cM\,:\,U\le X\}$, so that $\cM$ is $\Delta$-divisible if and only if $n-\widetilde{\cM}(H)\equiv 0\pmod\Delta$ for every hyperplane $H$.
More generally, we say that $\cM$ is \emph{$\Delta$-divisible with respect to codimension $i$} if $n-\widetilde{\cM}(W)\equiv 0\pmod\Delta$ for every subspace $W$ of codimension $i$.
For a subspace $S$, the \emph{residual} of $\cM$ with respect to $S$ is the multiset of members contained in $S$, considered as a projective $h$-system in $S$.

\begin{Lemma}
  \label{lemma_incidence}
  Let $\cM$ be a $\Delta$-divisible projective $h$-system in $\PG(v-1,q)$, where $v\ge 2$. Then $\Delta$ divides $\sum_{U\in\cM}q^{v-k_U}[k_U]_q$. If $\cM$ is faithful, i.e., all members have dimension $h$, then
  $\Delta$ divides $q^{v-h}[h]_q\cdot n$, and $\Delta'/\gcd(\Delta',[h]_q)$ divides $n$ for the largest divisor $\Delta'$ of $\Delta$ that is coprime to $q$.
\end{Lemma}
\begin{proof}
  Each of the $[v]_q$ hyperplanes $H$ satisfies $n-\widetilde{\cM}(H)\equiv 0\pmod\Delta$. A member $U$ is contained in exactly $[v-k_U]_q$ hyperplanes, so it contributes $[v]_q-[v-k_U]_q=q^{v-k_U}[k_U]_q$ to
  $\sum_H\left(n-\widetilde{\cM}(H)\right)$. For faithful $\cM$ we obtain $\Delta\mid q^{v-h}[h]_qn$, so that $\Delta'\mid[h]_qn$. Dividing by $g:=\gcd(\Delta',[h]_q)$ gives $\Delta'/g\mid\left([h]_q/g\right)n$,
  and $\Delta'/g$ is coprime to $[h]_q/g$.
\end{proof}
For $\gcd(\Delta,q)=1$ the second part of Lemma~\ref{lemma_incidence} says that $\Delta/\gcd(\Delta,[h]_q)$ divides $n$ for faithful $\cM$, which for $h=1$ is a consequence of \cite[Theorem 1]{ward1981divisible}.
The faithfulness is needed. In $\PG(1,2)$ let $\cM$ consist of three copies of the full space $V$ and the three points.
Each hyperplane, i.e., point, contains exactly one member. Hence $\cM$ is $5$-divisible with $\#\cM=6$, while $5/\gcd(5,[2]_2)=5$ does not divide $6$. For this example
the first statement gives $5\mid 3\cdot[2]_2+3\cdot 2=15$.

For sets of points the incidences with hyperplanes satisfy the following standard equations. Lemma~\ref{lemma_hyperplane_moments} is used in Appendix~\ref{app_p3}.
\begin{Lemma}
  \label{lemma_hyperplane_moments}
  Let $\cS$ be a set of $n$ points of $\PG(v-1,q)$ with $v\ge 3$, let $t$ be the number of $3$-subsets of $\cS$ that are contained in a line (for $q=2$ the number of lines contained in $\cS$), and for a hyperplane $H$ let $\cS(H):=\#(\cS\cap H)$ and $w(H):=n-\cS(H)$. Then, summing over all hyperplanes $H$,
  \begin{align*}
    &\sum_H 1=[v]_q,\quad \sum_H\cS(H)=n[v-1]_q,\quad \sum_H\binom{\cS(H)}{2}=\binom{n}{2}[v-2]_q,\\
    &\sum_H\binom{\cS(H)}{3}=\binom{n}{3}[v-3]_q+tq^{v-3},\quad \sum_Hw(H)=nq^{v-1},\quad \sum_H\binom{w(H)}{2}=\binom{n}{2}(q-1)q^{v-2},\\
    &\sum_H\binom{w(H)}{3}=\left(\binom{n}{3}-t\right)(q-1)^2q^{v-3}+t(q-2)q^{v-2}.
  \end{align*}
\end{Lemma}
\begin{proof}
  We count pairs of a hyperplane $H$ and a $k$-subset $T$ of $\cS$ with $k\le 3$ such that $T\subseteq H$ resp.\ $T\cap H=\emptyset$. A point lies in $[v-1]_q$ hyperplanes and outside of $q^{v-1}$. Two distinct points lie in $[v-2]_q$ common hyperplanes
  and outside of $[v]_q-2[v-1]_q+[v-2]_q=(q-1)q^{v-2}$. Three points that are not on a line span a plane, so that they lie in $[v-3]_q$ common hyperplanes and outside of $[v]_q-3[v-1]_q+3[v-2]_q-[v-3]_q=(q-1)^2q^{v-3}$.
  Three points on a line lie in $[v-2]_q=[v-3]_q+q^{v-3}$ common hyperplanes. A hyperplane that does not contain the line meets it in exactly one point, so it avoids all three points if and only if this point is one of the other $q-2$ points of the line, which happens for $(q-2)q^{v-2}$ hyperplanes. For $q=2$ no hyperplane avoids all three points.
\end{proof}

For $\gcd(\Delta,q)=1$ the divisibility of a projective $h$-system can be read off from the numbers of members through the points.
\begin{Corollary}
  \label{cor_coprime}
  Let $\cM$ be a projective $h$-system in $\PG(v-1,q)$ and $\gcd(\Delta,q)=1$. Then $\cM$ is $\Delta$-divisible if and only if $\sum_{U\ni P}q^{h-k_U}\equiv 0\pmod\Delta$ for every point $P$.
  In particular, a faithful $\cM$ is $\Delta$-divisible if and only if every point is contained in a multiple of $\Delta$ members of $\cM$.
\end{Corollary}
\begin{proof}
  Let $\chi_U$ denote the characteristic function of the point set of $U$, regarded as a multiset of points, and let $\cM''$ be the multiset of points $\sum_{U}q^{h-k_U}\chi_U$. A member $U\not\le H$ contributes $q^{h-k_U}\cdot q^{k_U-1}=q^{h-1}$ to $\cM''(V)-\cM''(H)$, hence
  $\cM''(V)-\cM''(H)=q^{h-1}\left(n-\widetilde{\cM}(H)\right)$. So $\cM$ is $\Delta$-divisible if and only if $\cM''$ is, and Proposition~\ref{thm_coprime} applies (its converse is trivial).
\end{proof}
The weights $q^{h-k_U}$ cannot be dropped: in the $5$-divisible example in $\PG(1,2)$ above every point is contained in $4$ members.

Let $m_\cM(U)$ denote the multiplicity of $U$ as a member of $\cM$. If $\gcd(\Delta,q)=1$ and a member $U$ is disjoint from all members different from $U$, then $m_\cM(U)\equiv 0\pmod\Delta$ by Corollary~\ref{cor_coprime}, since the points of $U$ are
contained in exactly $m_\cM(U)$ members. Without this disjointness assumption the multiplicities $m_\cM(U)$ of the members need not be divisible by $\Delta$, only the numbers of
members through the points are. The set $\cK$ of all lines in $\PG(3,q)$ satisfies $n-\widetilde{\cK}(H)=q^2[3]_q$ for every hyperplane $H$. Thus $\cK$ is
$[3]_q$-divisible with $\gcd([3]_q,q)=1$, but $m_\cK(U)=1\not\equiv 0\pmod{[3]_q}$ for every line $U$. Similarly, the seven lines of a plane in $\PG(3,2)$ are $3$-divisible
(even $6$-divisible, the values of $n-\widetilde{\cM}(H)$ are $0$ and $6$) and all multiplicities are $1$, while every point of the plane is contained in exactly $3$ of the lines.

\medskip

\noindent\emph{Residuals and line spreads.}\ In coding theory language, the residual code of a $\Delta$-divisible linear code over $\F_q$ with respect to some codeword is $\Delta/q$-divisible
if $q$ divides $\Delta$, see \cite[Lemma 3.4]{kurz2021divisible}. For multisets of points, divisibility is inherited by restrictions in the following slightly more general form.
\begin{Lemma}
  \label{lemma_restriction_points}
  Let $\cN$ be a $\Delta$-divisible multiset of points in $\PG(m-1,Q)$ and $Y$ a subspace of codimension $r$. Then $\cN|_Y$ is $\Delta/\gcd(\Delta,Q^r)$-divisible.
\end{Lemma}
\begin{proof}
  It suffices to consider $r=1$, since $\Delta'/\gcd(\Delta',Q)$ with $\Delta'=\Delta/\gcd(\Delta,Q^{r-1})$ equals $\Delta/\gcd(\Delta,Q^r)$ (compare the exponents of each prime).
  For $m=1$ there is nothing to show. Let $n:=\#\cN$, let $K$ be a hyperplane of $Y$, and let $Y=H_0,H_1,\dots,H_Q$ be the hyperplanes through $K$. Then $\sum_{i=0}^Q\left(\cN(H_i)-\cN(K)\right)=n-\cN(K)$ and
  $\cN(H_i)-\cN(K)\equiv n-\cN(K)\pmod\Delta$ for all $i$. Thus, $Q\cdot\left(n-\cN(K)\right)\equiv 0\pmod\Delta$, which gives $\cN(Y)-\cN(K)\equiv n-\cN(K)\equiv 0\pmod{\Delta/\gcd(\Delta,Q)}$.
\end{proof}

For $Q=q^h$ every multiset of points $\cN$ in $\PG(m-1,Q)$, i.e., every linear code over $\F_Q$ in geometric form, can be regarded as an additive code over $\F_Q$.
As such it is described by the multiset $\cM$ of $(h-1)$-spaces in $\PG(mh-1,q)$ obtained by \emph{field reduction}, where each point $\langle x\rangle_{\F_Q}$
of $V=\F_Q^m$ is regarded as an $h$-dimensional $\F_q$-subspace. For an $\F_q$-subspace $S$ of $V$ let $S_Q:=\{x\in S\,:\,\F_Q x\subseteq S\}$
denote the largest $\F_Q$-subspace contained in $S$.

\begin{Proposition}
  Let $\cN$ be a $\Delta$-divisible multiset of points in $\PG(m-1,Q)$, $Q=q^h$, and $\cM$ its field reduction. Let $S\neq\{0\}$ be an $\F_q$-subspace of $V=\F_Q^m$ of $\F_q$-codimension $r$
  and let $r^\star$ be the $\F_Q$-codimension of $S_Q$. Then the residual of $\cM$ with respect to $S$ is $\Delta/\gcd(\Delta,Q^{r^\star})$-divisible. We also have $\lceil r/h\rceil\le r^\star\le r$.
  In particular, $\cM$ is $\Delta$-divisible and the residuals with respect to hyperplanes are $\Delta/\gcd(\Delta,Q)$-divisible.
\end{Proposition}
\begin{proof}
  First we show that for an $\F_Q$-subspace $Y$ and an $\F_q$-hyperplane $Z$ of $Y$ the subspace $Z_Q$ is an $\F_Q$-hyperplane of $Y$. We have $Z=\ker(\varphi)$ for a non-zero $\F_q$-linear
  map $\varphi\colon Y\to\F_q$, and $\varphi=\operatorname{Tr}\circ\psi$ for a unique non-zero $\F_Q$-linear map $\psi\colon Y\to\F_Q$, where $\operatorname{Tr}$ is the trace of $\F_Q/\F_q$.
  Since the trace form is non-degenerate, $Z_Q=\{y\,:\,\operatorname{Tr}(\lambda\psi(y))=0\ \forall\lambda\in\F_Q\}=\ker(\psi)$.

  A member $\langle x\rangle_{\F_Q}$ is contained in an $\F_q$-subspace $X$ if and only if it is contained in $X_Q$. So $\widetilde{\cM}(X)=\cN(X_Q)$. Let $H'$ be an $\F_q$-hyperplane of $S$.
  Then $H'_Q=(H'\cap S_Q)_Q$ and $H'\cap S_Q$ is either $S_Q$ or an $\F_q$-hyperplane of $S_Q$. Therefore $H'_Q$ is either $S_Q$ or an $\F_Q$-hyperplane of $S_Q$ and
  $\widetilde{\cM}(S)-\widetilde{\cM}(H')=\cN(S_Q)-\cN(H'_Q)$ is divisible by $\Delta/\gcd(\Delta,Q^{r^\star})$ by Lemma~\ref{lemma_restriction_points}.

  Since $S$ contains the common kernel of $r$ $\F_Q$-linear forms (choose $\psi_1,\dots,\psi_r$ with $S=\bigcap\ker(\operatorname{Tr}\circ\psi_i)$), we have $r^\star\le r$, and
  $hr^\star\ge r$ follows from $S_Q\le S$. For a hyperplane $S=H$ we have $r^\star=r=1$ by the first paragraph, and $\cM$ is $\Delta$-divisible since $n-\widetilde{\cM}(H)=n-\cN(H_Q)$.
\end{proof}

For additive codes that do not arise by field reduction, residuals need not inherit divisibility, not even for spreads, see Proposition~\ref{prop_line_spreads_PG52} below.
An $(h-1)$-spread $\cS$ of $\PG(3h-1,q)$ is $q^{2h}$-divisible, since every hyperplane contains exactly $q^h+1$ elements. An $(h-1)$-spread of $\PG(hm-1,q)$ is called
\emph{Desarguesian} if it arises, up to a collineation, by field reduction from the set of all points of $\PG(m-1,q^h)$. For instance, a Desarguesian line spread of $\PG(5,q)$
is the field reduction of the point set of $\PG(2,q^2)$. A spread is called \emph{normal} (or geometric) if the span of any two of its elements is partitioned by elements of the spread.

\begin{Proposition}
  \label{prop_line_spreads_PG52}
  Let $\cS$ be a line spread of $\PG(5,q)$. Recall that $\cS$ is $q^4$-divisible. The following are equivalent:
  \begin{enumerate}
    \item[(i)] $\cS$ is Desarguesian.
    \item[(ii)] Every solid contains either $1$ or $q^2+1$ lines of $\cS$.
    \item[(iii)] $\cS$ is $q^2$-divisible with respect to codimension $2$, i.e., with respect to solids.
    \item[(iv)] The residuals of $\cS$ with respect to all hyperplanes are $q^2$-divisible.
  \end{enumerate}
  For $q=2$ these conditions are also equivalent to: (v) $\cS$ is $\Delta$-divisible with respect to solids for some $\Delta>1$, and to: (vi) the residuals of $\cS$ with respect to all hyperplanes are $\Delta$-divisible for some $\Delta>1$.
\end{Proposition}
\begin{proof}
  Let $x_K$ denote the number of lines of $\cS$ contained in a solid $K$. Then $0\le x_K\le q^2+1$. Each hyperplane contains $q^2+1$ lines of $\cS$ and $n=\#\cS=q^4+q^2+1\equiv 1\pmod{q^2}$.
  Since every solid is contained in a hyperplane, (iii) and (iv) are both equivalent to $x_K\equiv 1\pmod{q^2}$ for all $K$, which is equivalent to (ii).
  If (ii) holds, then the solid spanned by two lines of $\cS$ contains $q^2+1\ge 2$ lines of $\cS$, which partition its points. Hence $\cS$ is normal. Normal $(t-1)$-spreads of $\PG(rt-1,q)$
  with $r\ge 3$ are Desarguesian by a classical result of Barlotti and Cofman \cite{barlotti1974finite}, see also \cite{bader2011desarguesian} for a self-contained proof. Conversely, let $\cS$ be Desarguesian, i.e., obtained by field reduction from $\F_{q^2}^3$, and $K$ a solid.
  With $Q=q^2$ and $\lambda\in\F_{Q}\setminus\F_q$ we have $K_Q=K\cap\lambda K$, which has $\F_q$-dimension at least $4+4-6=2$. Therefore $K_Q$ is an $\F_{q^2}$-point or an $\F_{q^2}$-line, and $x_K\in\{1,q^2+1\}$.

  Now let $q=2$. Double counting gives $\sum_K 1=\qbin{6}{4}{2}=651$, $\sum_K x_K=21\cdot\qbin{4}{2}{2}=735$, and $\sum_K x_K(x_K-1)=21\cdot 20=420$, since two lines of $\cS$ span a unique solid.
  Hence $F(\alpha,\beta):=\sum_K(x_K-\alpha)(x_K-\beta)=21\left(55-35(\alpha+\beta)+31\alpha\beta\right)$, which vanishes for $0\le\alpha\le\beta\le 5$ only if $(\alpha,\beta)=(1,5)$.
  Conditions (v) and (vi) mean $x_K\equiv 21\pmod\Delta$ and $x_K\equiv 5\pmod\Delta$, respectively, for all $K$, and we may assume that $\Delta$ is a prime. If $\Delta$ is odd, then
  all $x_K$ lie in $\{\alpha,\beta\}$ for some $0\le\alpha\le\beta\le 5$ with $\Delta\mid\beta-\alpha$. Then $F(\alpha,\beta)=0$ gives $(\alpha,\beta)=(1,5)$, so that $\Delta\mid 4$, which is impossible for odd $\Delta$. If $\Delta=2$, then
  all $x_K$ are odd. So $(x_K-1)(x_K-5)\le 0$ with equality if and only if $x_K\in\{1,5\}$, and $F(1,5)=0$ gives (ii).
\end{proof}
In particular, among the $131\,044$ isomorphism types of line spreads of $\PG(5,2)$ \cite{mateva2009line} only the Desarguesian one is $\Delta$-divisible with respect to solids for some $\Delta>1$.
For a non-Desarguesian line spread of $\PG(5,2)$ there is no $\Delta>1$ such that all its residuals with respect to hyperplanes are $\Delta$-divisible.

\section{Canonical divisible multisets of points and the digit rule}
\label{sec_isomorphism_types}
Let $q=p^e$ and $\Delta=p^{ae-b}$ with $a,b\in\N_0$, $a\ge 1$, and $b\le e-1$. As recalled in Subsection~\ref{sub_lengths_points} below, the $a$-spaces and, for $1\le i\le a$, the
$p^{ie-b}$-fold $(a-i)$-spaces are $\Delta$-divisible multisets of points over $\F_q$, and so is every sum of them. We call such a sum \emph{canonical}.
The question of this section is for which cardinalities $n$ all $\Delta$-divisible multisets of points over $\F_q$ of cardinality $n$ are canonical. For $q=p\in\{2,3\}$
the answer is given by the following digit rule.
\begin{Theorem}[Digit rule]
  \label{thm_digit_rule}
  Let $p\in\{2,3\}$, $q=p$, $\Delta=p^r$ with $r\in\N_{>0}$, and $s_i:=p^i\cdot[r+1-i]_p$ for $0\le i\le r$. Write $n=\sum_{i=0}^r c_is_i$ with $c_0,\dots,c_{r-1}\in\{0,\dots,p-1\}$ and $c_r\in\Z$. Then there exists a $\Delta$-divisible multiset 
  of points over $\F_p$ of cardinality $n$ and all of them are canonical if and only if $c_r\ge 0$ and $c_i+c_{i+1}\le p-1$ for all $0\le i\le r$, where $c_{r+1}:=0$.
\end{Theorem}
The theorem follows by combining Theorem~\ref{thm_digits_only_if} (necessity, all primes), Theorem~\ref{thm_digits_p2} ($p=2$), and Theorem~\ref{thm_p3_induction} ($p=3$, using the computer
classifications for $\Delta\in\{3,9\}$). For $p\ge 5$ the condition is not sufficient, see Proposition~\ref{prop_triangle} in Appendix~\ref{app_p5}, which is why the theorem is stated only for $p\in\{2,3\}$.
For non-prime fields the condition would have to be modified, since the digits of the $S_q(a,b)$-adic expansion in Subsection~\ref{sub_lengths_points} no longer all range over $\{0,\dots,p-1\}$.

For $p=2$ the condition means that all digits, including $c_r$, lie in $\{0,1\}$ and that no two consecutive digits are equal to $1$. The resulting cardinalities are $\{2,3\}$ for $\Delta=2$, $\{4,6,7,11\}$ for $\Delta=4$, and
$\{8,12,14,15,22,23,27\}$ for $\Delta=8$.
For $p=3$ the rule gives $\{3,4,6,7,8\}$ for $\Delta=3$ and $\{9,12,13,18,21,22,24,25,26,31,34,35,44\}$ for $\Delta=9$.

For $\Delta=2$ the cardinalities $2$ and $3$ admit only canonical examples, whereas for cardinality four there is also the affine plane $\chi_E-\chi_L$ with a line $L$ in a plane $E$,
and for cardinality five the projective base (frame). The $2$-divisible multisets of points over $\F_2$ with cardinality at most seven and the $4$-divisible ones with cardinality at most $13$ are classified
in \cite{korner2024lengths}, in accordance with Theorem~\ref{thm_digit_rule}. Moreover, \cite[Proposition 3]{korner2024lengths} together with \cite[Theorem 1]{kiermaier2020lengths} shows for all prime powers $q$ and all $r,n\in\N_0$
with $n\le [r+1]_q$ that a $q^r$-divisible multiset of points over $\F_q$ of cardinality $n$ either does not exist or is canonical.

For $q=3$ the possible cardinalities of $3$-divisible multisets of points are $\{3,4\}\cup\{n\in\N_0:n\ge 6\}$. Using the software package \texttt{LinCode} \cite{bouyukliev2021computer} we have enumerated
all $3$-divisible $[n,k]_3$ codes with $n\le 13$ and all $9$-divisible ones with $n\le 44$, see Tables~\ref{table_3_divisible} and~\ref{table_9_divisible} in Appendix~\ref{app_sec3}.
All examples are canonical for the cardinalities that satisfy the digit rule of Theorem~\ref{thm_digit_rule}. For instance, each $9$-divisible multiset of points over $\F_3$ with cardinality $44$ can be written as $\cM=\chi_E+\chi_{E'}+9\chi_P+9\chi_{P'}$ for two planes $E$, $E'$ and
two points $P$, $P'$. For cardinality $43$, which violates the digit rule, the non-canonical examples include two that span $\PG(6,3)$ and have large automorphism groups.
Non-canonical examples of cardinality $q^2$ are the cylinders $\sum_{j=1}^q\chi_{L_j}-q\chi_P$ for pairwise different lines through a point $P$, which are conjectured to be the only $q$-divisible sets of $q^2$ points over prime fields with
maximum point multiplicity one, see, e.g., \cite{kurz2021generalization}. Examples for $q=3$ and $n\in\{9,10,11\}$ are described in Appendix~\ref{app_sec3}.

\subsection{Canonical multisets and basic tools}\label{sub_lengths_points}
In this subsection we recall the possible lengths of $\Delta$-divisible codes over $\F_q$ when $\Delta$ is a power of the characteristic $p$ of $\F_q$, together with
the constructions behind them, which give the canonical multisets. For $q=p^e$ and $\Delta=p^{ae-b}$ as at the beginning of this section we use the elementary facts on $\Delta$-divisible multisets of points
(embeddings into larger ambient spaces, sums, multiples, and the divisibility of $\chi_U$ for subspaces $U$) and the notation of \cite[Section~1 and Section~2]{kurz2024lengths}, see also \cite{kiermaier2020lengths} for proofs.
The cardinalities of the $a$-space and the $p^{ie-b}$-fold $(a-i)$-spaces, $1\le i\le a$, are $s_q(a,b,0)=[a+1]_q$ and $s_q(a,b,i)=p^{ie-b}[a+1-i]_q$, see \cite[Equations (1) and (2)]{kurz2024lengths}. For $q=p$, $a=r$, and $b=0$ we have
$s_p(r,0,i)=s_i$ for $0\le i\le r$, with $s_i$ as in Theorem~\ref{thm_digit_rule}. For each $c_0,\dots,c_a\in\N_0$ there exists a
$p^{ae-b}$-divisible multiset of points of cardinality
\begin{equation}\label{eq_s_q_adic}
n=\sum_{i=0}^a c_i\cdot s_q(a, b,i)
\end{equation}
in $\PG(v-1,q)$ for sufficiently large dimension $v$ of the ambient space.\footnote{$v\ge a+1$ is always sufficient. If $0\le i\le a$ is the smallest index with $c_i>0$, then $v\ge a+1-i$ 
is sufficient.} The set $S_q(a,b):=\{s_q(a,b,0),\dots,s_q(a,b,a)\}$ is the base of a positional system, and each integer $n$ has a unique $S_q(a,b)$-adic expansion as in Equation~(\ref{eq_s_q_adic})
with $c_0 \in\left\{0,\dots, p^{e-b}-1\right\}$, $c_1,\dots,c_{a-1}\in \left\{0,\dots, q-1\right\}$, and leading coefficient $c_a \in\Z$, see \cite[Equation (3)]{kurz2024lengths}. By \cite[Theorem~2 and Lemma~3]{kurz2024lengths}
there exists a $p^{ae-b}$-divisible linear code of length $n$ over $\F_q$ if and only if the leading coefficient $c_a$ of the $S_q(a,b)$-adic expansion of $n$ is non-negative, and the
construction in the proof of \cite[Lemma~3]{kurz2024lengths} uses the full expansion. The $c_i$-fold use of an $(a-i)$-space in this construction can be replaced by $c_i$ arbitrary $(a-i)$-spaces, which gives the multisets
\begin{equation}
  \cM=\sum_{j=1}^{c_0} \chi_{S_{a+1}^j}\,+\,\sum_{i=1}^a\sum_{j=1}^{c_i} p^{ie-b}\cdot \chi_{S_{a+1-i}^j}\label{eq_canonical},
\end{equation}   
where $\dim\!\left(S_{a+1-i}^j\right)=a+1-i$ for all $0\le i\le a$.
These are exactly the canonical multisets defined at the beginning of this section.

A \emph{piece} is a multiset $w\cdot\chi_S$, where $S$ is a subspace and $w\in\N_{>0}$. For $q=p$ prime and $\Delta=p^r$ the pieces in Equation~(\ref{eq_canonical}), which we call \emph{canonical pieces}, are the
$p^i\cdot\chi_S$ with $\dim(S)=r+1-i$ and $0\le i\le r$, and their cardinalities are the numbers $s_i=p^i[r+1-i]_p$ of Theorem~\ref{thm_digit_rule}. So a multiset of points is
canonical if and only if it is a (possibly empty) sum of such pieces. An integer linear combination of $\Delta$-divisible multisets of points is again a $\Delta$-divisible multiset
of points whenever it is non-negative, since $\#\cM-\cM(H)$ is additive in $\cM$.

\begin{Lemma}[E.g.~{\cite[Lemma 19]{honold2018partial}}, {\cite[Definition 5.7 and Exercise 5.8]{kurz2021divisible}}]
  \label{lemma_cone}
  Let $U$ and $V$ be subspaces of $\F_q^v$ with $U\cap V=\{0\}$ and $t:=\dim(V)\ge 0$. Let $\cM'$ be a multiset of points of $U$ with
  $\#\cM'=s$ and let $z\in\N_0$. Then
  $$
    \Cone_V(\cM',z):=\sum_{P} \cM'(P)\cdot\left(\chi_{\langle V,P\rangle}-\chi_V\right)+z\cdot\chi_V,
  $$
  where $P$ runs over the points of $U$, is a multiset of points of cardinality $q^ts+z[t]_q$. For $t=0$ we have $\Cone_V(\cM',z)=\cM'$. A point $X\le V$ has multiplicity $z$, and a point $X\le V+U$ with
  $X\not\le V$ has multiplicity $\cM'(P)$, where $P$ is the unique point of $U$ with $X\le\langle V,P\rangle$. Let $\Delta'\in\N_{>0}$. If $\cM'$ is $\Delta'$-divisible (as a multiset of points of $U$) and
  $$
    z\equiv -(q-1)s\pmod {q\Delta'},
  $$
  then $\Cone_V(\cM',z)$ is $q^t\Delta'$-divisible.
\end{Lemma}
\begin{proof}
  Put $W:=V+U=V\oplus U$ and $\cC:=\Cone_V(\cM',z)$. For a point $X\le W$ with $X\not\le V$ the projection along $V$ maps $X$ to a point $P$ of $U$, and
  $\langle V,P\rangle$ is the only subspace of the form $\langle V,P'\rangle$ that contains $X$. This gives the stated multiplicities, which are
  non-negative. Since $\#(\langle V,P\rangle\setminus V)=[t+1]_q-[t]_q=q^t$ we have $\#\cC=q^ts+z[t]_q$.

  As the support of $\cC$ is contained in $W$, a hyperplane of the ambient space either contains $W$, so that $\#\cC-\cC(H)=0$, or meets $W$ in a hyperplane of $W$,
  and every hyperplane of $W$ arises in this way. So it suffices to consider the hyperplanes $H$ of $W$. If $V\le H$, then $H=V\oplus H_U$ for the hyperplane
  $H_U:=H\cap U$ of $U$, and every hyperplane $H_U$ of $U$ arises in this way (if $\dim(U)=1$, then $H_U=\{0\}$ and $H=V$). For a point $P$ of $U$ we have
  $\langle V,P\rangle\le H$ if $P\le H_U$ and $\langle V,P\rangle\cap H=V$ otherwise. Hence
  $$
    \#\cC-\cC(H)=q^t\left(s-\cM'(H_U)\right).
  $$
  If $V\not\le H$, then $\dim(V\cap H)=t-1$ and for every point $P$ of $U$ the $(t-1)$-space $\langle V,P\rangle\cap H$ contains $V\cap H$, so it contains
  $[t]_q-[t-1]_q=q^{t-1}$ points outside of $V$. This gives $\cC(H)=q^{t-1}s+z[t-1]_q$ and
  $$
    \#\cC-\cC(H)=q^{t-1}(q-1)s+z\left([t]_q-[t-1]_q\right)=q^{t-1}\big((q-1)s+z\big).
  $$
  Hence $\cC$ is $q^t\Delta'$-divisible if $\Delta'$ divides $s-\cM'(H_U)$ for all hyperplanes $H_U$ of $U$, i.e., if $\cM'$ is $\Delta'$-divisible, and $q\Delta'$ divides $(q-1)s+z$.
\end{proof}

\begin{Lemma}
  \label{lemma_direct_sum}
  Let $U_1$, $U_2$ be subspaces of $\F_q^v$ with $U_1\cap U_2=\{0\}$ and let $\cM_1$, $\cM_2$ be multisets of points whose supports are contained in $U_1$ and $U_2$,
  respectively. If $\cM:=\cM_1+\cM_2=\sum_{j} w_j\chi_{S_j}$ is a sum of pieces, then each $S_j$ is contained in $U_1$ or in $U_2$, and $\cM_l=\sum_{j\,:\,S_j\le U_l} w_j\chi_{S_j}$
  for $l=1,2$. In particular, if $\cM_1+\cM_2$ is canonical, then $\cM_1$ and $\cM_2$ are canonical.
\end{Lemma}
\begin{proof}
  Since $w_j>0$, every point of $S_j$ has positive multiplicity in $\cM$, so $S_j=(S_j\cap U_1)\cup(S_j\cap U_2)$. A vector space is not the union of two proper subspaces, hence
  $S_j\le U_1$ or $S_j\le U_2$. Since the support of $\cM_2$ meets $U_1$ trivially, we have
  $\cM_1=\cM|_{U_1}=\sum_{j\,:\,S_j\le U_1} w_j\chi_{S_j}$, and similarly for $\cM_2$. The pieces of $\cM_l$ are pieces of $\cM$, hence $\cM_l$ is canonical if $\cM$ is.
\end{proof}

\begin{Lemma}
  \label{lemma_cone_noncanonical}
  Let $p$ be a prime, $r\ge 1$, $0\le i\le r-1$, $t:=r-1-i$, and $b_1,b_2\in\{0,\dots,p-1\}$ with $b_2\ge 1$. Let $U$, $V$ be subspaces with $U\cap V=\{0\}$
  and $\dim(V)=t$, and let $\cB$ be a $p$-divisible multiset of points of $U$ with $\#\cB=b_1(p+1)+b_2p$ and all multiplicities at most $p-1$.
  Then $\cM:=\Cone_V\!\left(p^i\cB,z\right)$ with $z:=p^i(b_1+b_2p)$ is a $p^r$-divisible multiset of points of cardinality $b_1s_i+b_2s_{i+1}$ that is not canonical.
\end{Lemma}
\begin{proof}
  The multiset $p^i\cB$ is $p^{i+1}$-divisible. Since $-(p-1)\left(b_1(p+1)+b_2p\right)=-b_1(p^2-1)-b_2p^2+b_2p\equiv b_1+b_2p\pmod{p^2}$, we have
  $z\equiv -(p-1)p^i\#\cB\pmod{p^{i+2}}$, and $\cM$ is $p^{t+i+1}=p^r$-divisible by Lemma~\ref{lemma_cone} with $\Delta'=p^{i+1}$.
  Using $[t+2]_p=p^{t+1}+p^t+[t]_p$ and $[t+1]_p=p^t+[t]_p$ we compute
  $$
    b_1s_i+b_2s_{i+1}=b_1p^i[t+2]_p+b_2p^{i+1}[t+1]_p=p^{t+i}\left(b_1(p+1)+b_2p\right)+z[t]_p=\#\cM.
  $$
  Assume $\cM=\sum_{j=1}^J p^{e_j}\chi_{S_j}$ with $\dim(S_j)=r+1-e_j\ge 1$. By Lemma~\ref{lemma_cone} the points of $V$ have multiplicity $z<p^{i+2}$ and all other
  points have multiplicity at most $p^i(p-1)<p^{i+1}$. Hence $e_j\le i+1$. If $e_j=i+1$, then all points of $S_j$ have multiplicity at least $p^{i+1}$ and thus lie in $V$,
  so that $\dim(S_j)=t+1>\dim(V)$, a contradiction. So $e_j\le i$ for all $j$.
  Counting points with $(p-1)p^e[r+1-e]_p=p^{r+1}-p^e$ gives $X:=\sum_j p^{e_j}=(J-b_1-b_2)p^{r+1}+(b_1+b_2p)p^i$. Since $0<(b_1+b_2p)p^i<p^{r+1}$,
  $X>0$ implies $J\ge b_1+b_2$, and $X\le Jp^i$ implies $(J-b_1-b_2)\left(p^{r+1}-p^i\right)\le -b_2(p-1)p^i<0$, a contradiction.
\end{proof}

\subsection{Base examples and the necessary condition}

\begin{Proposition}[E.g.~for $p=2$ {\cite[Proposition 15]{korner2024lengths}}]
  \label{prop_base_examples}
  Let $p$ be a prime and $1\le b\le p-1$. Let $\ell$ be a line, $Q_0,\dots,Q_b$ pairwise different points of $\ell$, and $M_0,\dots,M_b$ pairwise disjoint lines
  with $M_j\cap\ell=Q_j$, e.g., $M_j=\langle Q_j,X_j\rangle$ for points $X_0,\dots,X_b$ that span a $b$-space disjoint from $\ell$. Then
  $$
    \cB(p,b):=(p-1)\chi_\ell+\sum_{j=0}^b\left(\chi_{M_j}-p\chi_{Q_j}\right)
  $$
  is a $p$-divisible multiset of points over $\F_p$ with $\#\cB(p,b)=p^2+b$ and all multiplicities at most $p-1$.
\end{Proposition}
\begin{proof}
  The points $Q_j$ have multiplicity $(p-1)+1-p=0$, the other points of $\ell$ have multiplicity $p-1$, the points of $M_j\setminus\{Q_j\}$ have multiplicity $1$,
  and all other points have multiplicity $0$. Lines and $p$-fold points are $p$-divisible, and so is every integer linear combination of them with non-negative multiplicities. Finally,
  $\#\cB(p,b)=(p-1)(p+1)+(b+1)(p+1)-(b+1)p=p^2+b$.
\end{proof}
For $p=2$ the five points of $\cB(2,1)$ form a projective base (frame) of $\PG(3,2)$.

For $p=3$ the examples $\cB(3,1)$ and $\cB(3,2)$ have cardinalities $10$ and $11$. Projective examples of these cardinalities are the ovoid of $\PG(3,3)$ and the
point set of the dual of the ternary Golay code.

\begin{Theorem}
  \label{thm_digits_only_if}
  Let $p$ be a prime, $r\ge 1$, $\Delta=p^r$, and $n=\sum_{i=0}^r c_is_i$ with $c_0,\dots,c_{r-1}\in\{0,\dots,p-1\}$ and $c_r\in\Z$. Set $c_{r+1}:=0$.
  \begin{itemize}
    \item[(a)] If $c_r<0$, then there is no $p^r$-divisible multiset of points over $\F_p$ with cardinality $n$.
    \item[(b)] If $c_r\ge 0$ and $c_i+c_{i+1}\ge p$ for some $0\le i\le r$, then there exists a $p^r$-divisible multiset of points over $\F_p$ with cardinality $n$ that is not canonical.
  \end{itemize}
\end{Theorem}
\begin{proof}
  Part (a) is the main result of \cite{kurz2024lengths}. For (b) we construct a non-canonical $p^r$-divisible multiset $\cM_1$ and a canonical $p^r$-divisible multiset $\cM_2$
  with $\#\cM_1+\#\cM_2=n$. Embedding them into subspaces $U_1$, $U_2$ of a sufficiently large ambient space with $U_1\cap U_2=\{0\}$, the multiset $\cM_1+\cM_2$ is $p^r$-divisible (sums of $p^r$-divisible multisets are $p^r$-divisible, and embeddings preserve the divisibility, see \cite{kiermaier2020lengths}), and it is not canonical by Lemma~\ref{lemma_direct_sum}, since otherwise $\cM_1$ would be canonical.

  First assume that $c_r\ge p$. Let $\cM_1:=\chi_S-\chi_{S'}$, where $S$ is an $(r+1)$-space and $S'$ is a hyperplane of $S$, i.e., $\cM_1$ is an affine space of cardinality
  $p^{r+1}=p\cdot s_r$. It is $p^r$-divisible since $\chi_S$ is $p^{r+1}$-divisible and $\chi_{S'}$ is $p^r$-divisible. All multiplicities of $\cM_1$ are at most one.
  So a decomposition into canonical pieces can only use pairwise disjoint $r$-spaces of multiplicity one. This would imply that $[r+1]_p$ divides $p^{r+1}$, which is impossible since $[r+1]_p>1$
  and $\gcd\left([r+1]_p,p\right)=1$. So $\cM_1$ is not canonical. Let $\cM_2$ be the canonical multiset (\ref{eq_canonical}) with coefficients $c_0,\dots,c_{r-1},c_r-p$.

  Otherwise we have $c_r\le p-1$ and the index $i$ with $c_i+c_{i+1}\ge p$ satisfies $i\le r-1$ and $c_i\ge 1$ (if $c_i=0$, then $c_{i+1}\ge p$, which is excluded since $c_{i+1}\le p-1$ by the digit range for $i+1<r$ and by $c_r\le p-1$ for $i+1=r$).
  Let $\cM_1$ be the multiset of Lemma~\ref{lemma_cone_noncanonical} with $b_1=c_i$, $b_2=p-c_i$, and $\cB=\cB(p,c_i)$ of Proposition~\ref{prop_base_examples}, where $\#\cB=p^2+c_i=c_i(p+1)+(p-c_i)p$.
  It is non-canonical, $p^r$-divisible, and has cardinality $c_is_i+(p-c_i)s_{i+1}$. Since $c_i+c_{i+1}-p\ge 0$ we can let $\cM_2$ be the canonical multiset with coefficients
  $c_k$ for $k\notin\{i,i+1\}$, $0$ at position $i$, and $c_i+c_{i+1}-p$ at position $i+1$.
\end{proof}
The condition $c_r\ge 0$ and $c_i+c_{i+1}\le p-1$ for all $0\le i\le r$ of Theorem~\ref{thm_digit_rule} is necessary for every prime $p$.

\subsection{Sufficiency of the digit rule for \texorpdfstring{$p\in\{2,3\}$}{p in {2,3}}}\label{sub_p3}
\noindent\emph{The case $p=2$.}\ The converse of Theorem~\ref{thm_digits_only_if}(b) holds for $p=2$. The proof is the same as for $p=3$ below. It uses the defect of a hyperplane (Lemma~\ref{lemma_p3_defect}), the residue function modulo $p$ (Lemma~\ref{lemma_p3_residue}), and the minimum weight and the second weight of the binary Reed--Muller codes. The case $c_0=0$ is reduced by Lemma~\ref{lemma_p3_div3}, which shows that all multiplicities are divisible by $p$, and the case $c_0=1$ is Lemma~\ref{thm_p3_c01}.

\begin{Theorem}
  \label{thm_digits_p2}
  The digit rule, Theorem~\ref{thm_digit_rule}, holds for $p=2$. More precisely, let $r\ge 1$ and $n=\sum_{i=0}^rc_is_i$ with $c_0,\dots,c_r\in\{0,1\}$ such that no two consecutive $c_i$ are equal to one.
  Then every $2^r$-divisible multiset of points $\cM$ over $\F_2$ with $\#\cM=n$ is of the form $\cM=\sum_{i=0}^r c_i2^i\chi_{F_i}$, where $F_i$ is an $(r-i)$-space for every $i$ with $c_i=1$.
  These $F_i$ are uniquely determined: a point $P$ lies in $F_i$ if and only if the binary digit of $\cM(P)$ at position $i$ (counted from $0$) equals one.
\end{Theorem}
\begin{proof}
  By Theorem~\ref{thm_p2_star} in Appendix~\ref{app_p2} the multiset $\cM$ is canonical, and by Lemma~\ref{lemma_digits_unique} each decomposition of $\cM$ into canonical pieces consists of one piece $2^i\chi_{F_i}$
  with $\dim(F_i)=r+1-i$ for every $i$ with $c_i=1$. Since the $2^i$ with $c_i=1$ are distinct powers of two, $\cM(P)=\sum_{i\,:\,P\le F_i}c_i2^i$ is the binary expansion of $\cM(P)$.
  Together with Theorem~\ref{thm_digits_only_if} this proves Theorem~\ref{thm_digit_rule} for $p=2$.
\end{proof}

\begin{Remark}
  The proof of Theorem~\ref{thm_p2_star} uses the main result of \cite{kurz2024lengths} only through Lemma~\ref{lemma_p3_defect}(a), i.e., only for $2^j$-divisible multisets with $j<r$.
\end{Remark}

\medskip
\noindent\emph{The cases $p=2$ and $p=3$.}\ Throughout this part let $p\in\{2,3\}$, $r\ge 1$, and $s_i=p^i[r+1-i]_p$ for $0\le i\le r$. We call $(c_0,\dots,c_r)\in\{0,\dots,p-1\}^{r+1}$ a \emph{digit-good pattern} (of level $r$) if $c_i+c_{i+1}\le p-1$ for all $0\le i<r$,
and we call $n$ \emph{digit-good} if $n=\sum_{i=0}^r c_is_i$ for a digit-good pattern. These are exactly the cardinalities of the ``if'' part of Theorem~\ref{thm_digit_rule}.
By $(\star_r)$ we denote the statement that every $p^r$-divisible multiset of points over $\F_p$ with digit-good cardinality is canonical. The statement $(\star_0)$ holds trivially,
since every multiset of points is a sum of points.

For $p=3$ the proof of $(\star_r)$ for $r\ge 3$ is by induction on $r$ and is given in Appendix~\ref{app_p3}. It uses two tools. The first is the defect $d(H)\ge 0$ of a hyperplane $H$
(Lemma~\ref{lemma_p3_defect}). If $d(H)=0$, then $\cM|_H$ is $p^{r-1}$-divisible with digit pattern $(c_0,\dots,c_{r-1})$. The second is the residue function
$x\mapsto\cM(\langle x\rangle)\bmod p$, which lies in a generalized Reed--Muller code (Lemma~\ref{lemma_p3_residue}). The minimum weight and the second weight of these codes
\cite{delsarte1970generalized,rolland2010second} control the set $\mathcal{R}(\cM)$ of points whose multiplicity is not divisible by $p$.
If $c_0=0$, then all multiplicities are divisible by $p$ and $(\star_{r-1})$ applies to $\cM/p$ (Proposition~\ref{prop_p3_reduction}). It remains to treat the cases $c_0=1$ and $c_0=2$.

\begin{Lemma}
  \label{thm_p3_c01}
  Let $p\in\{2,3\}$ and $r\ge 1$, where $r\ge 3$ if $p=3$. Assume $(\star_{r-1})$ and let $(c_0,\dots,c_r)$ be a digit-good pattern with $c_0=1$. Then every $p^r$-divisible multiset of points over $\F_p$ with
  cardinality $n=\sum_{i=0}^rc_is_i$ is canonical.
\end{Lemma}
\begin{proof}[Proof sketch, see Appendix~\ref{app_p2}]
  Using $(\star_{r-1})$ on the hyperplanes of defect $0$, Lemma~\ref{lemma_p3_count} gives $\#\mathcal{R}(\cM)<\left((p+1)p^r+p-3\right)/(p-1)$. By the second weight of the generalized Reed--Muller codes the residue function then has minimum weight, i.e., $\#\mathcal{R}(\cM)=[r+1]_p$, and Proposition~\ref{prop_p3_reduction}(b) applies.
\end{proof}

\begin{Lemma}
  \label{thm_p3_c02}
  Let $r\ge 3$, assume $(\star_{r-1})$, and let $(c_0,\dots,c_r)$ be a digit-good pattern with $c_0=2$. Then every $3^r$-divisible multiset of points over $\F_3$ with
  cardinality $n=\sum_{i=0}^rc_is_i$ is canonical.
\end{Lemma}
\begin{proof}[Proof sketch, see Appendix~\ref{app_p3_case}]
  Induction on the dimension of the span of the support. Either the residue function is $2\chi_F$ for an $r$-space $F$ and Lemma~\ref{thm_p3_c01} applies to $\cM-\chi_F$, or the projection of $\cM$ from a point of multiplicity $1$ is analyzed. For $r=4$ and $n=638$ the general estimates do not suffice, and Lemma~\ref{lemma_p3_172}, which also uses the third moment of the hyperplane intersection numbers, closes this case.
\end{proof}

\begin{Theorem}
  \label{thm_p3_induction}
  Let $r\ge 3$ and assume $(\star_{r-1})$. Then every $3^r$-divisible multiset of points over $\F_3$ with digit-good cardinality $n$ is canonical, i.e., $(\star_r)$ holds.
  Consequently, using the computer classifications for $\Delta\in\{3,9\}$ in Table~\ref{table_3_divisible} and Table~\ref{table_9_divisible}, $(\star_r)$ holds for all $r\ge 1$,
  and Theorem~\ref{thm_digit_rule} holds for $p=3$.
\end{Theorem}
\begin{proof}
  Let $n=\sum_{i=0}^rc_is_i$ for a digit-good pattern $(c_0,\dots,c_r)$. The first statement follows from Proposition~\ref{prop_p3_reduction}(a) if $c_0=0$, from Lemma~\ref{thm_p3_c01}
  if $c_0=1$, and from Lemma~\ref{thm_p3_c02} if $c_0=2$. The classifications of Table~\ref{table_3_divisible} and Table~\ref{table_9_divisible} (computer results) state that
  $(\star_1)$ and $(\star_2)$ hold. So induction on $r$ gives $(\star_r)$ for all $r\ge 1$. Together with Theorem~\ref{thm_digits_only_if} this proves Theorem~\ref{thm_digit_rule} for $p=3$.
\end{proof}
For $r=3$ we have $(s_0,s_1,s_2,s_3)=(40,39,36,27)$. Here the theorem covers $30$ non-zero digit-good cardinalities between $27$ and $152=2s_0+2s_2$. Among them, $13$ have $c_0=0$,
$11$ have $c_0=1$, and the six cardinalities $80$, $107$, $116$, $134$, $143$, $152$ have $c_0=2$.

\section{Necessary conditions and general results on the lengths}
\label{sec_length_divisible}
This section and the next one study the possible lengths of divisible faithful additive codes, encoded by the sets $\cN_q^{\Delta,h}$ of base lengths. Subsection~\ref{sub_necessary} gives necessary conditions on the length and settles the case that $\Delta$ divides $q^h$. Divisibility constants that are not powers of $p$ are treated in Subsection~\ref{sub_nonpower}, and projective systems with members of different dimensions in Subsection~\ref{sub_mixed}. Appendix~\ref{app_N_sets} lists known sets $\cN_q^{\Delta,h}$. Results that are only sketched in the main text are proved in full in Appendix~\ref{app_sec4}.
\subsection{Necessary conditions and basic lengths}\label{sub_necessary}

The sum of two $\Delta$-divisible projective $h$-systems over $\F_q$ with cardinalities $n$ and $n'$ is a $\Delta$-divisible projective $h$-system over $\F_q$ with cardinality $n+n'$.
Hence the set of lengths of faithful $\Delta$-divisible $\F_q$-linear codes over $\F_{q^h}$ is additively closed. Every additively closed set of non-negative integers is generated by finitely many of its elements, since the elements divided by their greatest common divisor form a numerical semigroup. A generating set that is minimal with respect to inclusion is unique, which leads to the following definition.
\begin{Definition}
  The set $\cN_q^{\Delta,h}$ is the inclusion-minimal generating set of the set of lengths of faithful $\Delta$-divisible $\F_q$-linear codes over $\F_{q^h}$. Its elements are the \emph{base lengths}.
\end{Definition}

Replacing each member by its points reduces the divisibility of multisets of subspaces to that of multisets of points.
\begin{Lemma}[E.g.~{\cite[Lemma 4]{kurz2024additive}}]
  \label{lemma_exclusion}
  Let $\cM$ be a multiset of $(h-1)$-spaces in $\PG(v-1,q)$ and $\cM'$ be the multiset of points arising by replacing each $(h-1)$-space by its $[h]_q$ points. Then $\#\cM'=[h]_q\cdot \#\cM$ and 
  $\cM$ is $\Delta$-divisible if and only if $\cM'$ is $\Delta\cdot q^{h-1}$-divisible.
\end{Lemma}
\begin{proof}
  For a hyperplane $H$ each $(h-1)$-space contained in $H$ contributes $[h]_q$ points to $\cM'(H)$ and each other one contributes $[h-1]_q$ points. Hence
  $\cM'(V)-\cM'(H)=\left([h]_q-[h-1]_q\right)\cdot\left(\#\cM-\widetilde{\cM}(H)\right)=q^{h-1}\cdot\left(\#\cM-\widetilde{\cM}(H)\right)$.
\end{proof}
For a set $X\subseteq\N_0$ let $\sg{X}$ denote the set of non-negative integer combinations of its elements. As an example let $q=h=2$. The $2$-fold line, the $4$-fold line, and a line spread of $\PG(3,2)$ give the lengths $2$, $5$ for $\Delta=2$ and $4$, $5$ for $\Delta=4$, which are base lengths, i.e., $\{2,5\}\subseteq\cN_2^{2,2}$ and $\{4,5\}\subseteq\cN_2^{4,2}$. The positive integers outside $\sg{2,5}$ and $\sg{4,5}$ are $1,3$ and $1,2,3,6,7,11$, respectively. Lemma~\ref{lemma_exclusion} excludes all of them, since there is no $4$-divisible multiset of points of cardinality $3$ or $9$ and no $8$-divisible multiset of points with cardinality in $\{3,6,9,18,21,33\}$. Hence $\cN_2^{2,2}=\{2,5\}$ and $\cN_2^{4,2}=\{4,5\}$, which is also a special case of Theorem~\ref{thm_delta_divides_qh} below.
\begin{Lemma}
  \label{lemma_counting}
  Let $\cM$ be a $\Delta$-divisible multiset of $n$ $(h-1)$-spaces in $\PG(v-1,q)$, where $v>h$, and $n\equiv\rho\pmod\Delta$ with $0\le\rho<\Delta$. Then $n\cdot[v-h]_q\ge \rho\cdot[v]_q$.
  In particular, $\rho\ge 1$ implies $n>\rho\cdot q^h$.
\end{Lemma}
\begin{proof}
  For each hyperplane $H$ the number $\widetilde{\cM}(H)$ is non-negative and congruent to $\#\cM\equiv\rho$ modulo $\Delta$. So $\widetilde{\cM}(H)\ge\rho$. Each element of $\cM$ is contained in
  $[v-h]_q$ of the $[v]_q$ hyperplanes. Summing over all hyperplanes gives $n[v-h]_q=\sum_H\widetilde{\cM}(H)\ge\rho[v]_q=\rho q^h[v-h]_q+\rho[h]_q$, i.e., $\left(n-\rho q^h\right)[v-h]_q\ge\rho[h]_q>0$ for $\rho\ge 1$.
\end{proof}

\begin{Theorem}
  \label{thm_delta_divides_qh}
  Let $\Delta>1$ be a divisor of $q^h$. Then there exists a $\Delta$-divisible faithful additive code of length $n$ over $\F_{q^h}$ that is $\F_q$-linear if and only if $n=a\Delta+b\left(q^h+1\right)$ for non-negative integers $a$, $b$. 
  In other words, $\cN_q^{\Delta,h}=\left\{\Delta,q^h+1\right\}$.
\end{Theorem}
\begin{proof}
  \emph{Constructions.} A $\Delta$-fold $(h-1)$-space has length $\Delta$, and every hyperplane contains it or misses it, so all weights are $0$ or $\Delta$. Let $\cS$ be an $(h-1)$-spread of $\PG(2h-1,q)$, i.e., a set of $q^h+1$ pairwise disjoint $(h-1)$-spaces, which partition the points since $(q^h+1)[h]_q=[2h]_q$. A hyperplane $H$ contains $x$ elements of $\cS$ and meets each other element in an $(h-2)$-space. Counting the $[2h-1]_q=[h]_q+q^h[h-1]_q$ points of $H$ gives $[2h-1]_q=x[h]_q+(q^h+1-x)[h-1]_q$, hence $x=1$. So every hyperplane contains exactly one element of $\cS$, all non-zero weights equal $q^h$, and $\cS$ is $\Delta$-divisible. Since sums of $\Delta$-divisible projective $h$-systems are $\Delta$-divisible, all lengths $a\Delta+b(q^h+1)$ occur.

  \emph{Necessity.} Let $\cM$ be a $\Delta$-divisible multiset of $n$ $(h-1)$-spaces. If $n\equiv 0\pmod\Delta$, then $n=\frac{n}{\Delta}\cdot\Delta$. So let $n\equiv\rho\pmod\Delta$ with $1\le\rho<\Delta$. Embedding the ambient space into a larger one does not change the weights, since a hyperplane of the larger space contains the old space or meets it in a hyperplane. We may assume $v>h$. Lemma~\ref{lemma_counting} gives $n=\rho q^h+s$ with $s\ge 1$. Since $\Delta$ divides $q^h$, we have $s\equiv\rho\pmod\Delta$, and $1\le\rho<\Delta$ implies $s\ge\rho$. Hence $s-\rho$ is a non-negative multiple of $\Delta$, and $n=\frac{s-\rho}{\Delta}\cdot\Delta+\rho\,(q^h+1)$.

  \emph{Minimality.} Neither number is a non-negative combination of the other one, since $\Delta\le q^h<q^h+1$ and $q^h+1\equiv 1\not\equiv 0\pmod\Delta$ for $\Delta>1$. So $\cN_q^{\Delta,h}=\{\Delta,q^h+1\}$.
\end{proof}
\begin{Example}
  For $q=2$, $h=4$, and $\Delta=8$ the lengths are the numbers $8a+17b$, and the largest integer that is not a length is $8\cdot 17-8-17=111$. For example, $n=9$ has $\rho=1$, so Lemma~\ref{lemma_counting} requires $n>2^4=16$, and the smallest admissible $n\equiv 1\pmod 8$ is $17$, realized by a spread of $\PG(7,2)$. Similarly, $n\equiv 2\pmod 8$ forces $n>32$, so $34=2\cdot 17$ is the smallest admissible length with residue $2$.
\end{Example}

\begin{Corollary}
  \label{cor_3_lines}
  For $q=3$ and $h=2$ we have $\cN_3^{3,2}=\{3,10\}$. So no $3$-divisible multiset of $n$ lines exists in $\PG(v-1,3)$ for $n\in\{1,2,4,5,7,8,11,14,17\}$, that is, for all $n\notin\sg{3,10}$.
\end{Corollary}
In particular, no $3$-divisible multiset of $11$ lines exists in $\PG(v-1,3)$, although $9$-divisible multisets of $44$ points exist in $\PG(v-1,3)$. None of them is $2$-partitionable (see Subsection~\ref{sub_partitions}), so that the divisibility of the point multiset alone does not suffice in general.

\subsection{Divisibility constants that are not powers of \texorpdfstring{$p$}{p}}\label{sub_nonpower}
We now consider divisibility constants $\Delta$ that are not powers of $p$, starting with $\Delta=q^l[h]_q$.

\begin{Theorem}
  \label{thm_finite_exceptions}
  Let $h\ge 1$ and $l\ge 0$ be integers with $(h,l)\neq(1,0)$ and $\Delta:=q^l\cdot[h]_q$. Let
  $$
    F_{q,h,l}:=\left\{\begin{array}{ll} q[h]_q^2-q[h]_q-1 & \text{if } l=0,\\
    q^{l-1}\left(q^2[h]_q^2-q[h]_q-1\right)+\left(q^{l-1}-1\right)[hm]_q & \text{if } l\ge 1,\end{array}\right.
  $$
  where $m:=\left\lceil l/h\right\rceil+1$. Then for every integer $n>F_{q,h,l}$ there exists a $\Delta$-divisible faithful additive code over $\F_{q^h}$ of length $n$ that is $\F_q$-linear.
  Moreover, there exists an integer $n\ge(\Delta-1)\left(q^h-1\right)$ such that no such code of length $n$ exists.
\end{Theorem}
\begin{proof}[Proof sketch, see Appendix~\ref{app_sec4_necessary}]
  Suitable multiples of the set of all $(h-1)$-spaces of a subspace $\PG(h,q)$, of a $\Delta$-fold $(h-1)$-space, and of a Desarguesian $(h-1)$-spread of $\PG(hm-1,q)$ are $\Delta$-divisible, and the bound follows from the Frobenius number of two coprime integers. The lower bound follows from Lemma~\ref{lemma_counting}.
\end{proof}
\begin{Remark}
  For $(q,h,l)=(2,2,0)$ the bound $F_{2,2,0}=11$ is sharp by Theorem~\ref{thm_N_t1}, which gives $\cN_2^{3,2}=\{3,7\}$.
  In general, the largest non-length is at least $(\Delta-1)\left(q^h-1\right)$, and at least the largest integer excluded by Lemma~\ref{lemma_exclusion} or Lemma~\ref{lemma_counting}.
  It is at most the Frobenius number of the semigroup generated by the lengths of the three constructions in the proof of Theorem~\ref{thm_finite_exceptions} and the lengths of $[h]_q$-fold copies of $q^l$-divisible $\F_{q^h}$-linear codes.
  This Frobenius number can be smaller than $F_{q,h,l}$. For $(q,h,l)=(2,2,1)$, i.e., $\Delta=6$, it is the Frobenius number $23$ of $\sg{6,7,15}$, while $F_{2,2,1}=29$.
\end{Remark}
\begin{Remark}
  \label{rem_finite_iff}
  For a general divisibility constant $\Delta=p^a\cdot m'$ with $\gcd(m',p)=1$ the set of non-lengths is finite if and only if $m'$ divides $[h]_q$. Indeed, if $m'\mid[h]_q$, then $\Delta$ divides $q^a[h]_q$ and Theorem~\ref{thm_finite_exceptions} applies (the excluded case $(h,a)=(1,0)$ is $\Delta=1$, where every $n$ is a length).
  Conversely, double counting gives $\sum_{H}\left(n-\widetilde{\cM}(H)\right)=n\left([v]_q-[v-h]_q\right)=nq^{v-h}[h]_q$, where $H$ runs over all hyperplanes. Hence $m'$ divides $n[h]_q$ for every $\Delta$-divisible $\cM$ of cardinality $n$,
  i.e., $m'/\gcd(m',[h]_q)$ divides $n$. If $m'\nmid[h]_q$, then $m'/\gcd(m',[h]_q)>1$, and infinitely many integers are non-lengths.
\end{Remark}

For divisibility constants that are not powers of $p$ the counting argument of Lemma~\ref{lemma_counting} can be sharpened by counting only the hyperplanes that avoid a fixed point.
We state the identity for subspaces of arbitrary dimensions, since it is used again in Subsection~\ref{sub_mixed}.
\begin{Lemma}
  \label{lemma_point_identity}
  Let $\cM$ be a multiset of subspaces of $V=\F_q^v$ of dimensions $k_U\in\{1,\dots,h\}$, $U\in\cM$, let $N:=\sum_{U\in\cM}q^{h-k_U}$, $\cM'':=\sum_{U\in\cM}q^{h-k_U}\chi_U$, and $\rho\in\Z$. Then every point $P$ satisfies
  $$
    \sum_{H\,:\,P\not\le H}\left(\widetilde{\cM}(H)-\rho\right)=q^{v-h-1}\left(N-\rho q^h-\cM''(P)\right),
  $$
  where $H$ runs over the hyperplanes. If all members have dimension $h$, then $N=\#\cM$ and $\cM''$ is the point multiset $\cM'$ of $\cM$.
\end{Lemma}
\begin{proof}
  Exactly $q^{v-1}$ hyperplanes do not contain $P$. If $P\not\le U$, then exactly $q^{v-k_U-1}$ of them contain $U$, namely the preimages of the hyperplanes of $V/U$ that do not contain the point $P+U$, and if $P\le U$, then none.
  So, the left hand side equals $\sum_{U\,:\,P\not\le U}q^{v-k_U-1}-\rho q^{v-1}$, which is the right hand side.
\end{proof}
\begin{Lemma}
  \label{lemma_moment}
  Let $\Delta=p^a\cdot m'$ with $\gcd(m',p)=1$, let $\cM$ be a $\Delta$-divisible multiset of $n$ $(h-1)$-spaces in $\PG(v-1,q)$, and let $n\equiv\rho\pmod\Delta$ with $1\le\rho\le\Delta-1$. Then
  $n\ge\rho q^h+m'$. If $n=\rho q^h+m'$, then $[h]_q n/m'=[v']_q$ for some integer $v'$ and the point multiset of $\cM$ is $m'\cdot\chi_W$ for a $v'$-dimensional subspace $W$.
\end{Lemma}
\begin{proof}
  The point multiset $\cM'$ of $\cM$ is $\Delta q^{h-1}$-divisible by Lemma~\ref{lemma_exclusion}. All its multiplicities are divisible by $m'$ by Proposition~\ref{thm_coprime}. Since $\widetilde{\cM}(H)\equiv\rho\pmod\Delta$,
  we have $\widetilde{\cM}(H)\ge\rho$ for every hyperplane $H$, and Lemma~\ref{lemma_point_identity} gives $\cM'(P)\le n-\rho q^h$ for every point $P$. Choosing $P$ with $\cM'(P)>0$ gives $n-\rho q^h\ge m'$.
  Now let $n=\rho q^h+m'$ and without loss of generality let the members of $\cM$ span $V$. Then $\cM'(P)\in\{0,m'\}$ for every point $P$, and for $\cM'(P)=m'$ the identity gives $\widetilde{\cM}(H)=\rho$ for all hyperplanes $H$ with $P\not\le H$.
  Since the support of $\cM'$ spans $V$, every hyperplane misses one of its points. Therefore $\widetilde{\cM}(H)=\rho$ for all $H$. The identity then gives $\cM'(P)=m'$ for all points $P$, i.e., $\cM'=m'\chi_V$ and $[h]_qn=\#\cM'=m'[v]_q$.
\end{proof}
For $m'=1$ the inequality is the bound $n>\rho q^h$ of Lemma~\ref{lemma_counting}.
If $\Delta$ divides $[h]_q$, then $a=0$, $m'=\Delta$, and $q^h\equiv 1\pmod\Delta$. So the smallest length in the residue class of $\rho$ allowed by Lemma~\ref{lemma_counting} is $\rho q^h+\Delta$,
and Lemma~\ref{lemma_moment} excludes it unless $[h]_q\rho/\Delta=[u]_q$ for some integer $u\ge 1$.

Further examples arise as orbits of a Singer cycle.
\begin{Lemma}
  \label{lemma_singer}
  Let $v\ge h$, let $\sigma$ be a Singer cycle of $\PG(v-1,q)$, i.e., a collineation of order $[v]_q$ that permutes the points regularly, and let $U$ be an $(h-1)$-space. Then the multiset
  $\left\{\sigma^i(U)\,:\,0\le i<[v]_q\right\}$ of $[v]_q$ $(h-1)$-spaces covers every point exactly $[h]_q$ times. Hence it is $q^{v-h}[h]_q$-divisible.
  More generally, let $S_1<S_2<\dots<S_v$ be subspaces with $\dim(S_j)=j$ and $S_0:=\{0\}$. For $0\le b<h$ and $h<i\le v$ the multiset of points
  $[h-b]_q\chi_{S_i}+\left([i-b]_q-[h-b]_q\right)\chi_{S_b}$ is a sum of characteristic functions of $(h-1)$-spaces, i.e., it is $h$-partitionable in the sense of Subsection~\ref{sub_partitions}.
\end{Lemma}
\begin{proof}
  A point $P$ lies in $\sigma^i(U)$ if and only if $\sigma^{-i}(P)\in U$, and since $\sigma$ acts regularly on the points, this happens for exactly $[h]_q$ values of $i$. So, the point multiset is $[h]_q\chi_V$ and all non-zero weights equal
  $[h]_q\left([v]_q-[v-1]_q\right)/q^{h-1}=q^{v-h}[h]_q$. For the second statement apply the first one in $S_i/S_b$ to $(h-b-1)$-spaces and take preimages.
\end{proof}
Lemma~\ref{lemma_singer} gives a $q^l[h]_q$-divisible faithful additive code of length $[h+l]_q$, which is coprime to $q$. If it replaces the $[h]_q$-fold Desarguesian $(h-1)$-spread of $\PG(hm-1,q)$ in the proof of Theorem~\ref{thm_finite_exceptions},
then every $n>q^{l-1}\left(q^2[h]_q^2-q[h]_q-1\right)+\left(q^{l-1}-1\right)[h+l]_q$ is a length of a $q^l[h]_q$-divisible code. Since $hm\ge h+l$ with equality if and only if $h\mid l$, this improves $F_{q,h,l}$ for $l\ge 2$ unless $h$ divides $l$.

Let $\mathcal{A}_q^{\Delta,h}$ be the set of non-negative integers $n$ that pass the necessary conditions of Lemma~\ref{lemma_exclusion}, Lemma~\ref{lemma_counting}, and Lemma~\ref{lemma_moment}, so that every length of a $\Delta$-divisible faithful $\F_q$-linear code over $\F_{q^h}$ lies in $\mathcal{A}_q^{\Delta,h}$.
For the condition of Lemma~\ref{lemma_exclusion} we use that the lengths of $\Delta q^{h-1}$-divisible multisets of points are the $m'$-fold multiples of the lengths of $p^aq^{h-1}$-divisible multisets of points by Proposition~\ref{thm_coprime}.
The minimal generators of $\mathcal{A}_q^{\Delta,h}$ can be computed exactly.
\begin{Proposition}
  \label{prop_N_coprime}
  \begin{enumerate}
    \item[(a)] $\cN_2^{6,2}=\{6,7,15\}$, $\cN_3^{2,2}=\{2,13\}$, and $\cN_3^{6,2}=\{6,13,20\}$.
    \item[(b)] We have
               \begin{eqnarray*}
                 \cN_2^{12,2}&=&\{12,14,15,31,37\},\\
                 \cN_2^{24,2}&=&\{24,28,30,31,63,74,75,77\},\\
                 \cN_2^{48,2}&=&\{48,56,60,62,63,127,148,150,151,154,155,157\},\\
                 \cN_3^{18,2}&=&\{18,20,39,73,121\}.
               \end{eqnarray*}
  \end{enumerate}
\end{Proposition}
\begin{proof}[Proof sketch, see Appendix~\ref{app_sec4_necessary}]
  Lemma~\ref{lemma_exclusion}, Lemma~\ref{lemma_counting}, and Lemma~\ref{lemma_moment} give lower bounds in each residue class, in (b) by a computer calculation, and the minimal generators are realized by multiples of $(h-1)$-spaces, line spreads, Singer orbits (Lemma~\ref{lemma_singer}), and $h$-partitionable types (Subsection~\ref{sub_partitions}) found by integer linear programming.
\end{proof}

The hyperplanes that contain more than the minimum number of members carry a divisible structure of their own.
\begin{Lemma}
  \label{lemma_dual}
  Let $\Delta=p^a\cdot m'$ with $\gcd(m',p)=1$ and $m'\mid[h]_q$, and let $t:=[h]_q/m'$. Let $\cM$ be a $\Delta$-divisible multiset of $n$ $(h-1)$-spaces whose members span $V=\F_q^v$, let $n\equiv\rho\pmod\Delta$ with
  $1\le\rho\le\Delta-1$, and let $r:=v-h-1$. Then $s:=\left(n-\rho q^h\right)/m'$ is a positive integer, and $m_P:=\cM'(P)/m'$ is an integer with $0\le m_P\le s$ for every point $P$, where $\cM'$ is the point multiset of $\cM$.
  For every hyperplane $H$ let $j_H:=\left(\widetilde{\cM}(H)-\rho\right)/\Delta$, and let $\mathcal{D}$ be the multiset of points of the dual space in which the point $H$ has multiplicity $j_H$. Then
  $$
    \sum_{H\,:\,P\not\le H}j_H=\frac{q^r}{p^a}\left(s-m_P\right)
  $$
  for every point $P$. So $p^a\mathcal{D}$ is a $q^r$-divisible multiset of points of cardinality $s[r+1]_q-\rho t$ whose weights lie in $\left\{0,q^r,2q^r,\dots,sq^r\right\}$.
\end{Lemma}
\begin{proof}
  Since $\rho\neq 0$ we have $v>h$, i.e., $r\ge 0$. The $j_H$ are non-negative integers since $\widetilde{\cM}(H)\equiv\rho\pmod\Delta$ and $0\le\rho<\Delta$, and the $m_P$ are integers by Lemma~\ref{lemma_exclusion} and Proposition~\ref{thm_coprime}.
  Since $m'$ divides $n-\rho$ and $q^h-1=(q-1)[h]_q$, it divides $n-\rho q^h$. Hence $s$ is an integer. Lemma~\ref{lemma_point_identity} gives $\Delta\sum_{H\,:\,P\not\le H}j_H=q^r\left(n-\rho q^h-m'm_P\right)=q^rm'\left(s-m_P\right)$,
  which is the stated identity and implies $m_P\le s$. Moreover, $s\ge 1$ by Lemma~\ref{lemma_moment}. Since each member lies in $[r+1]_q$ hyperplanes and $[v]_q=q^h[r+1]_q+[h]_q$, the cardinality of $p^a\mathcal{D}$ is
  $p^a\left(n[r+1]_q-\rho[v]_q\right)/\Delta=\left(sm'[r+1]_q-\rho[h]_q\right)/m'=s[r+1]_q-\rho t$. The hyperplanes of the dual space are the sets
  of hyperplanes $H$ through a point $P$, and the weight of such a hyperplane with respect to $p^a\mathcal{D}$ is $q^r\left(s-m_P\right)$.
\end{proof}
For the multiset of $22$ lines consisting of the $7$ lines of a plane $E$ of $\PG(3,2)$ and three copies of a line spread, where $(q,h,\Delta)=(2,2,6)$, $\rho=4$, and $s=2$, we have $m_P=2$ for $P\le E$ and $m_P=1$ otherwise, and
$\mathcal{D}$ consists of the point $E$ of the dual space with multiplicity $1$.

Let $\sigma_q(x)$ denote the sum of the digits of a non-negative integer $x$ in base $q$, i.e., the minimum number of powers of $q$ with sum $x$.
\begin{Theorem}
  \label{thm_digit_sum}
  Let $\Delta=p^a\cdot m'$ with $\gcd(m',p)=1$ and $m'\mid[h]_q$, let $t:=[h]_q/m'$, and let $\cM$ be a $\Delta$-divisible multiset of $n$ $(h-1)$-spaces with $n\equiv\rho\pmod\Delta$ and $1\le\rho\le\Delta-1$.
  Then $s:=\left(n-\rho q^h\right)/m'$ is a positive integer with
  $$
    \sigma_q\!\left(\frac{n-\rho}{m'}\right)\le s.
  $$
\end{Theorem}
\begin{proof}[Proof sketch, see Appendix~\ref{app_sec4_necessary}]
  Let $r:=v-h-1$, where $v$ is the dimension of the span of the members. By Lemma~\ref{lemma_dual} and \cite{kiermaier2020lengths} the number $s[r+1]_q-\rho t$ is a non-negative integer combination of the numbers $q^i[r+1-i]_q$, and a rearrangement of the corresponding powers of $q$ writes $s+(q-1)\rho t$ as a sum of at most $s$ powers of $q$.
\end{proof}
In the situation of Theorem~\ref{thm_digit_sum}, let $\mathcal{A}'^{\Delta,h}_q\subseteq\mathcal{A}_q^{\Delta,h}$ be the set of integers $n$ that additionally satisfy the condition of Theorem~\ref{thm_digit_sum}, i.e., $\rho=0$ or $\sigma_q\!\left((n-\rho)/m'\right)\le s$ for $s=(n-\rho q^h)/m'$.
\begin{Lemma}
  The set $\mathcal{A}'^{\Delta,h}_q$ is closed under addition.
\end{Lemma}
\begin{proof}
  Let $n_1,n_2$ belong to $\mathcal{A}'^{\Delta,h}_q$, with $s_i$, $\rho_i$, $x_i=(n_i-\rho_i)/m'$, and let $\rho\le\rho_1+\rho_2$ be the residue of $n_1+n_2$. First we show $n_1+n_2\in\mathcal{A}_q^{\Delta,h}$.
  The lengths of $p^aq^{h-1}$-divisible multisets of points are closed under addition, so the condition of Lemma~\ref{lemma_exclusion} holds. From Lemma~\ref{lemma_counting} only the inequality $n>\rho q^h$ for $\rho\ge 1$ is used, and it follows from $n\ge\rho q^h+m'$ of Lemma~\ref{lemma_moment}. If $\rho=0$, then there is nothing to show.
  If exactly one of $\rho_1,\rho_2$ is zero, say $\rho_1=0$, then $\rho=\rho_2$ and $n_1+n_2-\rho q^h>n_2-\rho_2q^h\ge m'$ for $n_1>0$. If both are non-zero, then $n_1+n_2-\rho q^h\ge n_1+n_2-\left(\rho_1+\rho_2\right)q^h\ge 2m'$.
  In both cases $n_1+n_2>\rho q^h+m'$, so the equality case of Lemma~\ref{lemma_moment} does not occur. It remains to check the condition on the digit sums, which is void for $\rho=0$.
  If $\rho_1=0$, then $n_1=m'k$ for some $k\ge 1$, $\rho=\rho_2$, $s=s_2+k$, $x=x_2+k$, and $\sigma_q(x)\le\sigma_q(x_2)+\sigma_q(k)\le s_2+k=s$. The case $\rho_2=0$ is symmetric.
  If $\rho_1,\rho_2\ge 1$ and $\rho=\rho_1+\rho_2<\Delta$, then $s=s_1+s_2$, $x=x_1+x_2$, and $\sigma_q(x)\le\sigma_q(x_1)+\sigma_q(x_2)\le s$.
  If $\rho=\rho_1+\rho_2-\Delta\ge 1$, then $s=s_1+s_2+p^aq^h$, and $\Delta=p^am'$ gives $x=x_1+x_2+p^a$, so $\sigma_q(x)\le s_1+s_2+p^a\le s$.
\end{proof}

For $s=1$ the condition says that $(n-\rho)/m'$ is a power of $q$, i.e., $\rho[h]_q/m'=[u]_q$ for some integer $u\ge 1$, which is the condition of Lemma~\ref{lemma_moment} for $n=\rho q^h+m'$. On the other hand, if $(n-\rho)/m'$ is a sum of $s$ powers $q^{i}$, then
$s[r+1]_q-\rho t$ is the length of a $q^r$-divisible multiset of points for all $r\ge\max i$. So Theorem~\ref{thm_digit_sum} extracts everything that the lengths of $q^r$-divisible multisets of points can say about the cardinality of $p^a\mathcal{D}$ in Lemma~\ref{lemma_dual}.
For $(q,h,\Delta)=(2,4,3)$ and $n=38$, that is, $\rho=2$ and $s=2$, we have $(n-\rho)/m'=12=2^3+2^2$, and $38$ is indeed a length, see Proposition~\ref{prop_25_38}.

For $t=1$, i.e., $m'=[h]_q$, and $p^a\mid q$, Theorem~\ref{thm_digit_sum} is sharp, and the Singer orbits of Lemma~\ref{lemma_singer} realize all lengths.
\begin{Theorem}
  \label{thm_N_t1}
  Let $h\ge 2$, let $p$ be the characteristic of $\F_q$, let $a\ge 0$ with $p^a\mid q$, and let $\Delta=p^a[h]_q$. An integer $n\ge 1$ with $n\equiv\rho\pmod\Delta$, $1\le\rho\le\Delta-1$, is the length of a
  $\Delta$-divisible faithful additive code over $\F_{q^h}$ that is $\F_q$-linear if and only if $s:=\left(n-\rho q^h\right)/[h]_q$ satisfies $s\ge 1$ and $\sigma_q\!\left((n-\rho)/[h]_q\right)\le s$.
  Moreover, $\cN_q^{[h]_q,h}=\left\{[h]_q,[h+1]_q,\dots,[2h-1]_q\right\}$, and $\cN_q^{\Delta,h}=\left\{\Delta,[h+1]_q,[h+2]_q,\dots,[2h]_q\right\}$ if $a\ge 1$.
\end{Theorem}
\begin{proof}[Proof sketch, see Appendix~\ref{app_sec4_necessary}]
  Necessity is Theorem~\ref{thm_digit_sum}. Conversely, a representation $(n-\rho)/[h]_q=\sum_ic_iq^i$ with $\sum_ic_i=s$ and $p^a\mid c_0$ gives $n=\tfrac{c_0}{p^a}\Delta+\sum_{i\ge 1}c_i[h+i]_q$, and $[h+i]_q$ is realized by Lemma~\ref{lemma_singer}.
\end{proof}
Theorem~\ref{thm_N_t1} gives $\cN_2^{3,2}=\{3,7\}$ and $\cN_2^{6,2}=\{6,7,15\}$, see Proposition~\ref{prop_N_coprime}.(a), and the sets marked with $^\ast$ in Proposition~\ref{prop_N_coprime_dual} below.
For $p^a\nmid q$ the Singer orbits do not suffice, e.g., $37\in\cN_2^{12,2}$.

\begin{Proposition}
  \label{prop_N_coprime_dual}
  The sets $\cN_q^{\Delta,h}$ of the following table hold. The sets marked with $^\ast$ are the cases of Theorem~\ref{thm_N_t1}.
  \begin{center}
  \begin{tabular}{rrrl}
    \hline
    $q$ & $h$ & $\Delta$ & $\cN_q^{\Delta,h}$\\
    \hline
    2 & 3 & 7$^\ast$ & $\{7,15,31\}$\\
    2 & 3 & 14$^\ast$ & $\{14,15,31,63\}$\\
    2 & 3 & 28 & $\{28,30,31,63,85,127\}$\\
    2 & 4 & 3 & $\{3,25,38\}$\\
    2 & 4 & 5 & $\{5,21\}$\\
    2 & 4 & 6 & $\{6,25,38,51\}$\\
    2 & 4 & 12 & $\{12,25,38,51\}$\\
    3 & 2 & 4$^\ast$ & $\{4,13\}$\\
    3 & 2 & 12$^\ast$ & $\{12,13,40\}$\\
    3 & 2 & 36 & $\{36,39,40,121,145,146\}$\\
    4 & 2 & 5$^\ast$ & $\{5,21\}$\\
    4 & 2 & 10$^\ast$ & $\{10,21,85\}$\\
    4 & 2 & 20$^\ast$ & $\{20,21,85\}$\\
    4 & 2 & 40 & $\{40,42,85,201,341\}$\\
    5 & 2 & 2 & $\{2,31\}$\\
    5 & 2 & 3 & $\{3,31\}$\\
    5 & 2 & 6$^\ast$ & $\{6,31\}$\\
    \hline
  \end{tabular}
  \end{center}
\end{Proposition}
\begin{proof}[Proof sketch, see Appendix~\ref{app_sec4_necessary}]
  An exact computer calculation shows that Lemma~\ref{lemma_exclusion}, Lemma~\ref{lemma_counting}, Lemma~\ref{lemma_moment}, and Theorem~\ref{thm_digit_sum} leave exactly the semigroups generated by the stated sets, and the stated lengths are realized by $\Delta$-fold $(h-1)$-spaces, Lemma~\ref{lemma_singer}, spreads, field reduction, $h$-partitionable types, and Proposition~\ref{prop_25_38}.
\end{proof}

The constructions for $25$ and $38$ are most easily described in the dual space.
\begin{Lemma}
  \label{lemma_dual_cover}
  Let $\cM$ be a multiset of $n$ $(h-1)$-spaces in $\PG(v-1,q)$, $v>h$, and let $\cM^\perp$ be the multiset of the $(v-h-1)$-spaces $U^\perp$, $U\in\cM$, where $\perp$ is taken with respect to a non-degenerate symmetric
  bilinear form. Then $\widetilde{\cM}(H)$ is the number of members of $\cM^\perp$ that contain the point $H^\perp$. Hence $\cM$ is $\Delta$-divisible if and only if every point of $\PG(v-1,q)$ is contained in a number of members of $\cM^\perp$
  that is congruent to $n$ modulo $\Delta$.
\end{Lemma}
\begin{proof}
  We have $U\le H$ if and only if $H^\perp\le U^\perp$.
\end{proof}

\begin{Proposition}
  \label{prop_line_cone}
  Let $h\ge 2$, let $\ell$ be a line of $\PG(h+1,q)$, and let $\Delta\in\N_{>0}$, and let $\rho'\ge 1$, $\mu\ge 0$, and $L\ge 0$ be integers such that the quotient geometry $\PG(h+1,q)/\ell\cong\PG(h-1,q)$ contains $L$ lines that cover every point at most $\rho'$ times.
  Let $C:=\rho'[h]_q-(q+1)L$ and let $c_X\ge 0$, for the points $X$ of $\ell$, be integers with $\sum_Xc_X=C$. Then there is a multiset $\cK$ of $n:=\mu+q^2L+qC$ lines of $\PG(h+1,q)$ such that every point outside $\ell$
  is contained in exactly $\rho'$ members of $\cK$ and every point $X$ of $\ell$ in exactly $\mu+qc_X$ members. If $\mu+qc_X\equiv\rho'\equiv n\pmod\Delta$ for all points $X$ of $\ell$, then $\cK^\perp$ is a $\Delta$-divisible
  multiset of $n$ $(h-1)$-spaces in $\PG(h+1,q)$.
\end{Proposition}
\begin{proof}[Proof sketch, see Appendix~\ref{app_sec4_necessary}]
  For each of the $L$ lines of the quotient take $q^2$ lines of a line spread of the corresponding solid through $\ell$, complete the coverage in the planes through $\ell$ by pencils of lines through points of $\ell$, and add $\mu$ copies of $\ell$. The last statement is Lemma~\ref{lemma_dual_cover}.
\end{proof}
For $(q,h)=(2,4)$, $\rho'=1$, $L=3$ (three lines of a line spread of $\PG(3,2)$), $c_X=6$ for one point $X$ of $\ell$, and $\mu=1$ we obtain $n=25$, where $X$ is covered $13$ times and all other points once. This gives a $12$-divisible multiset of $25$ solids in $\PG(5,2)$.
For $(q,h)=(3,4)$, $\rho'=1$, $L=8$ (eight pairwise disjoint lines of $\PG(3,3)$), $c_X=8$ for one point $X$, and $\mu=1$ we obtain $n=97$, where $X$ is covered $25$ times. This gives a $24$-divisible multiset of $97$ solids in $\PG(5,3)$ with weights $72^1\,96^{363}$.
The same multiset is $\Delta$-divisible for the divisors $\Delta\in\{4,8,12\}$ of $24$, and $97\equiv 1\pmod\Delta$ for $\Delta\in\{4,8,12,24\}$. In fact, $97\in\cN_3^{\Delta,4}$ for these four values of $\Delta$. By Lemma~\ref{lemma_counting} every length $n\not\equiv 0\pmod\Delta$ satisfies $n>3^4=81$. So if $97$ were the sum of two smaller positive lengths, then one summand would be a multiple of $\Delta$ and the other one a length $n\equiv 1\pmod\Delta$ with $81<n<97$.
Theorem~\ref{thm_digit_sum} excludes these candidates. For $\Delta=4$ we have $t=10$, and $n=85,89,93$ give $s=1,2,3$ with $\sigma_3(21)=3$, $\sigma_3(22)=4$, and $\sigma_3(23)=5$, respectively. For $\Delta=12$ only $n=85$ remains, with the same values of $s$ and $t$.
For $\Delta=8$ only $n=89$ remains, where $t=5$, $s=1$, and $\sigma_3(11)=3$. For $\Delta=24$ there is no candidate.

\begin{Proposition}
  \label{prop_25_38}
  There exist a $12$-divisible multiset of $25$ solids in $\PG(5,2)$ with weights $12^1\,24^{62}$ and a $12$-divisible multiset of $38$ solids in $\PG(6,2)$ with weights $24^1\,36^{126}$.
  So, $25$ and $38$ are lengths of $\Delta$-divisible faithful additive codes over $\F_{16}$ that are $\F_2$-linear for $\Delta\in\{3,6,12\}$.
\end{Proposition}
\begin{proof}[Proof sketch, see Appendix~\ref{app_sec4_necessary}]
  For $n=25$ this is the example after Proposition~\ref{prop_line_cone}. For $n=38$ we give $38$ planes of $\PG(6,2)$ such that one point lies on $14$ and every other point on $2$ of them, using a group of order $42$, and apply Lemma~\ref{lemma_dual_cover}.
\end{proof}

\subsection{Projective systems of mixed dimensions}\label{sub_mixed}
Let $C\subseteq\F_{q^{k_1}}\times\dots\times\F_{q^{k_n}}$ be an $\F_q$-linear code over a mixed alphabet, given as the image of $x\mapsto\left(\pi_1(x),\dots,\pi_n(x)\right)$ for $\F_q$-linear
maps $\pi_i\colon\F_q^v\to\F_q^{k_i}$. Replacing the $i$th alphabet by a smaller one if necessary, we may assume that each $\pi_i$ is surjective. Then $U_i:=\ker(\pi_i)^\perp$ has dimension $k_i$,
and for $x\neq 0$ and $H=x^\perp$ we have $\pi_i(x)=0$ if and only if $U_i\le H$. So the codeword of $x$ has Hamming weight $n-\widetilde{\cM}(H)$, where $\cM=\{U_1,\dots,U_n\}$ is regarded as a projective $h$-system
for any $h\ge\max_i k_i$, and $\Delta$-divisibility of $\cM$ is $\Delta$-divisibility of $C$ with respect to the Hamming weight. Throughout this subsection
$\cM$ is a projective $h$-system in $\PG(v-1,q)$, $v\ge 2$, whose members $U$ have dimensions $k_U\in\{1,\dots,h\}$, and $n=\#\cM$.

Lemma~\ref{lemma_incidence} and Corollary~\ref{cor_coprime} were stated for such $\cM$. The weighted multiset of points from the proof of Corollary~\ref{cor_coprime} also generalizes Lemma~\ref{lemma_exclusion}.

Let $\cM''=\sum_{U\in\cM}q^{h-k_U}\chi_U$. Then $\#\cM''=\sum_{U\in\cM}q^{h-k_U}[k_U]_q$ and $\cM''(V)-\cM''(H)=q^{h-1}\left(n-\widetilde{\cM}(H)\right)$ for every hyperplane $H$, since a member $U\not\le H$ contributes $q^{h-k_U}\cdot q^{k_U-1}=q^{h-1}$ to the left hand side, as in the proof of Corollary~\ref{cor_coprime}, and a member $U\le H$ contributes $0$.
So $\cM$ is $\Delta$-divisible if and only if $\cM''$ is $\Delta q^{h-1}$-divisible.

\begin{Lemma}
  \label{lemma_mixed_counting}
  Let $\cM$ be $\Delta$-divisible with $n\equiv\rho\pmod\Delta$, $1\le\rho<\Delta$, and $N:=\sum_{U\in\cM}q^{h-k_U}$. Then $[v]_q\left(N-\rho q^h\right)\ge\#\cM''>0$. Since $N\le q^{h-1}n$,
  this implies $n>q\rho$ and $n\ge\Delta$.
\end{Lemma}
\begin{proof}
  Every hyperplane $H$ satisfies $\widetilde{\cM}(H)\equiv\rho\pmod\Delta$ and hence $\widetilde{\cM}(H)\ge\rho$. Thus, Lemma~\ref{lemma_point_identity} gives $N-\rho q^h\ge\cM''(P)$ for every point $P$.
  Summing over the $[v]_q$ points gives the first inequality. Finally $q^{h-k_U}\le q^{h-1}$, and if $n<\Delta$, then $\rho=n$, contradicting $n>q\rho$.
\end{proof}
For faithful $\cM$ we have $N=n$ and $\#\cM''=[h]_qn$, and since $[v]_q-[h]_q=q^h[v-h]_q$, Lemma~\ref{lemma_mixed_counting} reduces to Lemma~\ref{lemma_counting}. Constructions for mixed dimensions are based on the following exchange.

\begin{Lemma}
  \label{lemma_exchange}
  Let $U$ be a $k$-dimensional subspace and $\cX$ a multiset of subspaces of $U$ such that $\#\cX-\widetilde{\cX}(K)=w$ for every hyperplane $K$ of $U$. If $\cM$ contains $w$ copies of $U$, then
  replacing them by $\cX$ does not change $n-\widetilde{\cM}(H)$ for any hyperplane $H$. Examples are the $[k]_q$ points of $U$ with $w=q^{k-1}$ and a $(j-1)$-spread of $U$ with $w=q^{k-j}$.
\end{Lemma}
\begin{proof}
  If $U\le H$, then both the $w$ copies of $U$ and the elements of $\cX$ are contained in $H$. Otherwise $K=H\cap U$ is a hyperplane of $U$ and both contribute $w$.
\end{proof}
Let $\Delta\in\N_{>0}$, $2\le k\le h$, and let $c,m\ge 0$ be integers with $cq^{k-1}\le m\Delta$. Applying Lemma~\ref{lemma_exchange} $c$ times to $m\Delta$ copies of a $(k-1)$-space $U$ gives a $\Delta$-divisible projective $h$-system
of cardinality $m\Delta+c[k-1]_q$. For $k=2$ these are $m\Delta-qc$ copies of a line $L$ and $c$ copies of each point of $L$.
\begin{Corollary}
  \label{cor_mixed_finite}
  Let $h\ge 2$ and $\Delta\ge 2$. Every integer $n\ge q\Delta$ and every multiple of $\Delta$ is the cardinality of a $\Delta$-divisible projective $h$-system. If $\Delta\ge q-1$, then there
  is none of cardinality $(q-1)\Delta-1$.
\end{Corollary}
\begin{proof}
  Write $n=m\Delta+c$ with $0\le c<\Delta$ and $m\ge q$. Then $qc\le m\Delta$, and the construction above with $k=2$ gives cardinality $n$. For $n=(q-1)\Delta-1$ we have $\rho=\Delta-1$ and $n\le q\rho$, which is excluded by Lemma~\ref{lemma_mixed_counting}.
\end{proof}
So in contrast to the faithful case, see Remark~\ref{rem_finite_iff}, the set of non-lengths is finite for every $\Delta$ if mixed dimensions are allowed. For $\Delta=q^l$ the necessary
condition of Lemma~\ref{lemma_mixed_counting} is also sufficient.
\begin{Theorem}
  \label{thm_mixed_qpower}
  Let $h\ge 2$, $l\ge 1$, and $\Delta=q^l$. There exists a $\Delta$-divisible projective $h$-system of cardinality $n\ge 1$ if and only if $n\equiv 0\pmod\Delta$ or $n>q\rho$, where
  $\rho\in\{0,\dots,\Delta-1\}$ and $n\equiv\rho\pmod\Delta$. Equivalently, $n$ is a non-negative integer combination of $q^l,q^l+1,\dots,[l+1]_q$. The examples can be chosen in $\PG(l,q)$
  and consist of lines and points only. The largest non-length is $(q-1)q^l-1$.
\end{Theorem}
\begin{proof}[Proof sketch, see Appendix~\ref{app_sec4_mixed}]
  Necessity is Lemma~\ref{lemma_mixed_counting}. The set of admissible $n$ is additively closed and generated by $q^l,\dots,[l+1]_q$, and these cardinalities are realized in $\PG(l,q)$ by applying Lemma~\ref{lemma_exchange} to copies of a line and to the set of all points of $\PG(l,q)$.
\end{proof}
For $q=2$ Theorem~\ref{thm_mixed_qpower} says that all $n\ge 2^l$ occur. For example, for $\Delta=8$ the faithful lengths over $\F_4$ are $\sg{8,10,21,25}$ and the
lengths of $8$-divisible binary linear codes are $\sg{8,12,14,15}$. So the lengths $n\in\{9,11,13,17,19\}$ need mixed dimensions, and none of them is a sum of a faithful length and a points-only length, which for $17$ and $19$ does not follow from their size.
Examples for $n=9$ and $n=13$ are six copies of a line together with its three points, and, for a solid $T$, a plane $E<T$, and a point $P\in E$, the four lines of $E$ not through $P$ together with the nine points of $T$ not covered by these lines.

Also for $q\ge 3$ points alone do not suffice. For example, for $(q,\Delta)=(3,27)$ the lengths $28,\dots,35$, $37$, and $38$ satisfy the conditions of Theorem~\ref{thm_mixed_qpower}
but are not cardinalities of $27$-divisible multisets of points. Other divisibility constants are discussed in Remark~\ref{rem_mixed_general_Delta} in Appendix~\ref{app_sec4_mixed}.

\noindent
\textbf{Subspace designs.} The multiset of blocks of a subspace design is a one-weight code, and the conditions of Lemmas~\ref{lemma_exclusion}, \ref{lemma_counting}, and~\ref{lemma_incidence} never exclude admissible design parameters, see Proposition~\ref{prop_designs} in Appendix~\ref{app_sec4_mixed}.

\section{Lengths for divisibility constants that are powers of \texorpdfstring{$p$}{p}}
\label{sec_length_powers}
For divisibility constants that are powers of $p$ a complete answer is out of reach in general, and most of our results concern small parameters. Partitions into lines determine the sets $\cN_2^{\Delta,2}$ over $\F_4$ in Subsection~\ref{sub_partitions}, and Subsections~\ref{sub_constructions} and~\ref{sub_qh23} contain constructions for $\Delta=q^{h+1}$, $q^{h+2}$, and $q^{h+3}$. The longer proofs are in Appendix~\ref{app_sec4}.

\subsection{Partitions into lines and the lengths over \texorpdfstring{$\F_4$}{F\textunderscore{}4}}\label{sub_partitions}

Since the possible lengths of divisible linear codes are known \cite{kiermaier2020lengths,kurz2024lengths}, the criterion of Lemma~\ref{lemma_exclusion} can be evaluated directly. The base examples
over $\F_4$ for $\cN_4^{2,1}=\left\{2,5\right\}$ are a double point and a line, and those for $\cN_4^{4,1}=\left\{4,5\right\}$ are a $4$-fold point and a line, i.e., $p^l$-fold copies of simplex codes.
Field reduction turns them into the base examples for $\cN_2^{2,2}=\{2,5\}$ and $\cN_2^{4,2}=\{4,5\}$ of Theorem~\ref{thm_delta_divides_qh}. Since field reduction maps $\F_4$-linear codes to additive codes,
every element of $\sg{8,10,21}=\sg{\cN_4^{8,1}}$ is a length of an $8$-divisible additive code over $\F_4$, and Lemma~\ref{lemma_exclusion} leaves $25$, $33$, $35$, and $43$ as the only further candidates. By Theorem~\ref{thm_N_explicit} below, $8$, $10$, and $21$ remain base lengths.

One strategy to construct faithful additive codes is to partition suitable multisets of points into $(h-1)$-spaces, as in the Solomon--Stiffler type construction for additive codes
in \cite{kurz2024additive}.

\begin{Definition}
  A multiset of points $\cM$ in $\PG(v-1,q)$ is $h$-partitionable if there exist $(h-1)$-spaces $H_1,\dots,H_l$ with $\cM=\sum_{i=1}^{l} \chi_{H_i}$, see also \cite[Definition 5.1]{dhaeseleer2025generalized}.
  The multiset $\{H_1,\dots,H_l\}$ is then called an $h$-partition of $\cM$.
\end{Definition}
If $\cM$ is $h$-partitionable, then $\cM$ is $q^{h-1}$-divisible and $\#\cM$ is divisible by $[h]_q$.
By Lemma~\ref{lemma_exclusion}, every $h$-partition of a $q^r$-divisible $h$-partitionable multiset $\cM$ of points is $q^{r-h+1}$-divisible. So the divisibility of an $h$-partition is determined by $\cM$, and the difficulty lies in deciding whether a $q^r$-divisible multiset of points is $h$-partitionable at all,
see Corollary~\ref{cor_3_lines}.

\begin{Definition}[E.g.~{\cite[Definition 5.2]{dhaeseleer2025generalized}}]
  Let $S_1<S_2<\dots<S_v$ be a chamber of $\PG(v-1,q)$, i.e., $\dim(S_i)=i$ for all $1\le i\le v$. We say that a multiset of points $\cM$ in $\PG(v-1,q)$ is of type $\sum_{i=1}^v \varepsilon_i[i]$
  if $\cM=\sum_{i=1}^v\varepsilon_i\cdot\chi_{S_i}$, where $\varepsilon_1,\dots,\varepsilon_v\in\mathbb{Z}$.
\end{Definition}
Since $\mathrm{GL}(v,q)$ acts transitively on chambers, the type determines $\cM$ up to isomorphism. Note that $[i]$ in a type is a formal symbol for $\chi_{S_i}$, while $[i]_q=\#\chi_{S_i}$.
\begin{Lemma}[E.g.~{\cite[Lemmas 29 and 44]{kurz2024additive}}]
  \label{lemma_type_divisible}
  Let $\cM$ be a multiset of type $\sum_{i=1}^v \varepsilon_i[i]$ in $\PG(v-1,q)$ and for a hyperplane $H$ let $j(H)$ be the smallest index with $S_{j(H)}\not\le H$. Then
  $$
    \#\cM-\cM(H)=\sum_{i=j(H)}^{v}\varepsilon_i q^{i-1},
  $$
  and exactly $q^{v-j}$ hyperplanes $H$ satisfy $j(H)=j$. In particular, $\cM$ is $q^{r-1}$-divisible if and only if $q^{r-i}$ divides $\varepsilon_i$ for all $1\le i\le r$.
\end{Lemma}
\begin{proof}
  For $i<j(H)$ we have $S_i\le H$, and for $i\ge j(H)$ the intersection $S_i\cap H$ is a hyperplane of $S_i$. The difference $\#\chi_{S_i}-\chi_{S_i}(H)$ equals $0$ for $i<j(H)$ and $q^{i-1}$ for $i\ge j(H)$,
  which gives the formula. The hyperplanes with $j(H)=j$ are those that contain $S_{j-1}$ but not $S_j$, and there are $[v-j+1]_q-[v-j]_q=q^{v-j}$ of them. In particular, every $j\in\{1,\dots,v\}$ occurs.
  If $q^{r-i}$ divides $\varepsilon_i$ for $1\le i\le r$, then every summand $\varepsilon_iq^{i-1}$ is divisible by $q^{r-1}$. Conversely, the values of $\#\cM-\cM(H)$ for $j(H)=j$ and $j(H)=j+1$
  differ by $\varepsilon_jq^{j-1}$, and the value for $j(H)=v$ is $\varepsilon_vq^{v-1}$. Hence $q^{r-1}$-divisibility forces $q^{r-1}\mid\varepsilon_jq^{j-1}$, i.e., $q^{r-j}\mid\varepsilon_j$, for $1\le j\le r$.
\end{proof}
\begin{Remark}
  \label{rem_type_weights}
  If $\cM$ is the point multiset of a multiset $\cM'$ of $(h-1)$-spaces, then by Lemmas~\ref{lemma_exclusion} and~\ref{lemma_type_divisible} the weight of a hyperplane $H$ for $\cM'$ is $q^{1-h}\sum_{i\ge j(H)}\varepsilon_iq^{i-1}$.
  So, a multiset of $(h-1)$-spaces whose point multiset has a type with $s$ non-zero coefficients has at most $s$ distinct non-zero weights. This explains why all constructions of
  Subsection~\ref{sub_constructions} are few-weight codes. E.g.\ the type $c[2h+1]+(q-c)q^h[h+1]+q^2[2h-1]$ of Theorem~\ref{thm_plus_q2} gives the weights
  $cq^{h+1}$ for $j(H)\ge 2h$, $(c+1)q^{h+1}$ for $h+2\le j(H)\le 2h-1$, and $(q+1)q^{h+1}$ for $j(H)\le h+1$, attained by $q+1$, $q^2[h-2]_q$, and $q^h[h+1]_q$ hyperplanes, respectively.
\end{Remark}

Proposition~\ref{thm_partition} below gives a $16$-divisible, $2$-partitionable multiset of $75$ points in $\PG(4,2)$ of type $[5]+4[3]+16[1]$ (the case $r=5$, $S=\{1,3\}$).
So we have $25\in \cN_2^{8,2}$ and $50=2\cdot 25\in \cN_2^{16,2}$. Theorem~\ref{thm_N_explicit} below gives $\cN_2^{8,2}=\left\{8,10,21,25\right\}$, $\cN_2^{16,2}=\left\{16,20,21,50\right\}$,
and $\cN_2^{32,2}=\left\{32,40,42,85,100,101,105,121\right\}$, where the lengths $101$, $105$, and $121$ come from Corollary~\ref{cor_conj_partition} below.

The following lemma lists basic $h$-partitionable types. Proposition~\ref{thm_partition} uses parts~(a) and~(c) for $q=h=2$, and Subsection~\ref{sub_constructions} uses them for general $q$ and $h$.
\begin{Lemma}[E.g.~for part~(a) {\cite[Lemma 3]{honold2018partial}}, {\cite[Lemma 10]{kurz2024additive}}]
  \label{lemma_moves}
  Let $h\ge 2$ and let $S_1<S_2<\dots<S_v$ be a chamber of $\PG(v-1,q)$.
  \begin{enumerate}
    \item[(a)] (Lifted MRD codes) For $2h\le i\le v$ the multiset $\chi_{S_i}-\chi_{S_{i-h}}$ is $h$-partitionable by $q^{i-h}$ $(h-1)$-spaces.
    \item[(b)] (Cones) For $h\le j\le v$ the $[j-h+1]_q$ $(h-1)$-spaces in $S_j$ that contain $S_{h-1}$ have the point multiset $\chi_{S_j}+q[j-h]_q\cdot\chi_{S_{h-1}}$.
    \item[(c)] (Pencils in $S_{h+1}$) Let $\alpha,\beta\in\N_0$ such that $q$ divides $\beta$, $\alpha\ge \beta/q$, and $[h]_q$ divides $\alpha-\beta/q$. Then
               $\alpha\cdot\chi_{S_{h+1}}+\beta\cdot\chi_{S_{h-1}}$ is $h$-partitionable. More precisely, let $z:=\left(\alpha-\beta/q\right)/[h]_q$ and take each of the
               $q+1$ $(h-1)$-spaces of $S_{h+1}$ through $S_{h-1}$ exactly $z+\beta/q$ times and each other $(h-1)$-space of $S_{h+1}$ exactly $z$ times.
  \end{enumerate}
\end{Lemma}
\begin{proof}
  (a) Let $m:=i-h\ge h$, write $S_i=W\oplus S_{m}$ with $\dim(W)=h$, identify $S_{m}$ with $\F_{q^m}$, and fix an injective $\F_q$-linear map $\iota\colon W\to\F_{q^m}$.
  For $a\in\F_{q^m}$ let $X_a:=\{w+a\iota(w)\,:\,w\in W\}$. Each $X_a$ is an $h$-dimensional subspace with $X_a\cap S_m=\{0\}$, and $X_a\cap X_b=\{0\}$ for $a\neq b$ since $(a-b)\iota(w)=0$ implies $w=0$.
  The $q^m$ subspaces $X_a$ contain $q^m\left(q^h-1\right)$ non-zero vectors, which is the number of vectors in $S_i\backslash S_m$. Hence they partition the points of $S_i$ outside $S_m$.
  (b) The $(h-1)$-spaces of $S_j$ containing $S_{h-1}$ correspond to the points of $S_j/S_{h-1}$, so that there are $[j-h+1]_q$ of them, every point of $S_{h-1}$ is covered $[j-h+1]_q=1+q[j-h]_q$ times, and every
  other point $P$ of $S_j$ is covered exactly once, namely by $S_{h-1}+P$.
  (c) Let $P$ be a point of $S_{h+1}$. If $P\notin S_{h-1}$, then $P$ lies in exactly one of the $q+1$ $(h-1)$-spaces through $S_{h-1}$ and in $[h]_q-1$ of the others. It is covered
  $\left(z+\beta/q\right)+([h]_q-1)z=[h]_qz+\beta/q=\alpha$ times. If $P\in S_{h-1}$, then $P$ lies in all $q+1$ $(h-1)$-spaces through $S_{h-1}$ and in $[h]_q-q-1$ of the others. It is covered $(q+1)\left(z+\beta/q\right)+([h]_q-q-1)z=[h]_qz+\beta+\beta/q=\alpha+\beta$ times.
\end{proof}

\begin{Proposition}
  \label{thm_partition}
  Let $r\ge 2$, $S\subseteq\{1,\dots,r-1\}$, and $\cT(r,S)$ be the multiset of points of type $[r]+\sum_{s\in S}2^{r-s}\cdot[s]$ in $\PG(r-1,2)$. Then $\cT(r,S)$ is $2^{r-1}$-divisible and $\#\cT(r,S)=2^r-1+\sum_{s\in S}\left(2^r-2^{r-s}\right)$.
  The multiset $\cT(r,S)$ is $2$-partitionable if and only if $3$ divides $\#\cT(r,S)$, which is equivalent to the following condition on the number $o$ of odd elements of $S$: $o\equiv 2\pmod 3$ if $r$ is odd, and $o\equiv 0\pmod 3$ if $r$ is even.
\end{Proposition}
\begin{proof}
  The divisibility follows from Lemma~\ref{lemma_type_divisible}. The cardinality is clear, and modulo $3$ we have $\#\cT(r,S)\equiv 1+o$ for odd $r$ and $\#\cT(r,S)\equiv 2o$ for even $r$, since $2^{j}\equiv(-1)^j$.

  Clearly, $3\mid\#\cT(r,S)$ is necessary for $2$-partitionability. For the converse we use that all subsequent modifications change the cardinality by multiples of $3$. First, we replace $[i]$ by $[i-2]$ for all $i\ge 4$ using Lemma~\ref{lemma_moves}.(a) with $q=h=2$.
  We end up with lines $[2]$ and a residual multiset of type $B\cdot[1]+A\cdot[3]$ with $B=\varepsilon_1$ and $A=\sum_{i\ge 3\text{ odd}}\varepsilon_i$, where $\varepsilon_i$ denotes the coefficient of $[i]$ in the type.
  \begin{itemize}
    \item If $1\notin S$, then $B=0$ and $3\mid A$. So Lemma~\ref{lemma_moves}.(c) with $q=h=2$, $\alpha=A$, and $\beta=0$ applies (if $A=0$, e.g., for $r=2$, nothing is left to partition).
    \item Now let $1\in S$, i.e., $\varepsilon_1=2^{r-1}$. The set of all lines through $S_1$ in $S_r$ is of type $[r]+\left(2^{r-1}-2\right)[1]$, since every point
          other than $S_1$ is contained in exactly one such line. Replacing $[r]$ by these lines, we are left with a multiset of type $2\cdot[1]+\sum_{1<i<r}\varepsilon_i[i]$. After the previous reduction we end up with $B=2$ and $A=\sum_{3\le i<r,\,i\text{ odd}}\varepsilon_i$.
          Since the cardinality $B+7A\equiv 2+A$ is divisible by $3$, we have $A\ge 1$ (in particular $r\ge 4$). Thus, $B/2=1\le A$ and $3$ divides $A-B/2$, so that again Lemma~\ref{lemma_moves}.(c) with $q=h=2$, $\alpha=A$, and $\beta=B$ applies.
  \end{itemize}
\end{proof}
\begin{Corollary}
  \label{cor_conj_partition}
  Let $r\ge 5$ be an odd integer and $S$ be a $2$-subset of the positive odd integers less than $r$. Then there exists a $2^{r-1}$-divisible multiset of points in $\PG(r-1,2)$ of type 
  $[r]+\sum_{s\in S}2^{r-s}\cdot [s]$ that is $2$-partitionable. There exists a $2^{r-2}$-divisible additive code over $\F_4$ of length $\left(2^r-1+\sum_{s\in S}\left(2^r-2^{r-s}\right)\right)/3$.
\end{Corollary}
\begin{proof}
  We have $o=2$. Hence Proposition~\ref{thm_partition} applies, and Lemma~\ref{lemma_exclusion} gives the divisibility of the additive code.
\end{proof}
For $r=7$ Corollary~\ref{cor_conj_partition} gives the lengths $101$, $105$, and $121$.
Let $\mathcal{L}_q^r$ denote the set of cardinalities of $q^r$-divisible multisets of points over $\F_q$ (including $0$). By \cite[Theorem 1]{kiermaier2020lengths}, $\mathcal{L}_2^m$ is the semigroup generated by the numbers $2^{m+1-j}\left(2^j-1\right)$ with $1\le j\le m+1$.

\begin{Theorem}
  \label{thm_D}
  Let $l\ge 1$. There exists a $2^l$-divisible faithful additive code over $\F_4$ of length $n$, i.e., a $2^l$-divisible multiset of $n$ lines in some $\PG(v-1,2)$, if and only if $3n\in \mathcal{L}_2^{l+1}$.
\end{Theorem}
\begin{proof}
  Necessity follows from Lemma~\ref{lemma_exclusion}, since replacing each line by its three points gives a $2^{l+1}$-divisible multiset of $3n$ points, and $3n\in \mathcal{L}_2^{l+1}$ by \cite{kiermaier2020lengths}. For sufficiency set $G_j:=2^{l+2-j}\left(2^j-1\right)=2^{l+2}-2^{l+2-j}$ for $1\le j\le l+2$, so that $\mathcal{L}_2^{l+1}$ is the semigroup generated by $G_1,\dots,G_{l+2}$. We use the following examples.
  \begin{enumerate}
    \item[(a)] For $1\le t\le (l+3)/2$ let $X_t:=2^{\max\{0,l+2-2t\}}\cdot(4^t-1)/3$. Field reduction of the $\F_4$-linear simplex code of dimension $t$, which is $4^{t-1}$-divisible, gives a line spread of $\PG(2t-1,2)$ that is $2^{2t-2}$-divisible. Taking it
          $2^{\max\{0,l+2-2t\}}$ times gives a $2^l$-divisible multiset of $X_t$ lines. Note $3X_t=G_{2t}$ for $2t\le l+2$ and $3X_t=2^{l+3}-1$ for $2t=l+3$.
    \item[(b)] For odd integers $1\le c<b<r\le l+2$ let $Y_{r,b,c}:=2^{l+2-r}\cdot\#\cT(r,\{b,c\})/3$. By Proposition~\ref{thm_partition} the multiset $\cT(r,\{b,c\})$ is $2^{r-1}$-divisible and $2$-partitionable. So Lemma~\ref{lemma_exclusion}
          gives a $2^{r-2}$-divisible multiset of lines. Taking it $2^{l+2-r}$ times gives a $2^l$-divisible multiset of $Y_{r,b,c}$ lines. Note $3Y_{r,b,c}=G_r+G_b+G_c$.
  \end{enumerate}
  Now let $3n=\sum_{j}e_jG_j$ with $e_j\in\N_0$. For odd $b\ge 3$ we have $2G_b=2G_1+G_{b-1}$ and $G_{b-1}=3X_{(b-1)/2}$. Replacing two copies of $G_b$ in this way does not change the number of odd generators, so that we may assume that every odd index $b\ge 3$ occurs at most once.
  For even $j$ we have $3\mid G_j$, and for odd $j$ we have $G_j\equiv(-1)^{l+1}\not\equiv 0\pmod 3$. Thus, the number $\sum_{j\text{ odd}}e_j$ of odd generators is divisible by $3$. We form triples of pairwise different odd indices $\ge 3$ until at most two such indices remain,
  complete them by two resp.\ one copies of $G_1$, and group the remaining copies of $G_1$ into triples. It suffices to show that $G_{2t}/3$ and the following sums are non-negative integer combinations of the numbers in (a) and (b), where $G_j=2^{l+2}-2^{l+2-j}$.
  The first statement is (a). The remaining triples are:
  \begin{itemize}
    \item $a>b>c$ odd: $(G_a+G_b+G_c)/3=Y_{a,b,c}$.
    \item $a>b=c=1$: $G_a+2G_1=2^{l+3}-2^{l+2-a}$, which equals $2G_{a+1}=6X_{(a+1)/2}$ for $a\le l+1$ and $3X_{(l+3)/2}$ for $a=l+2$.
    \item $a=b=c=1$: $(3G_1)/3=G_1=2X_1$.
  \end{itemize}
  So $n$ is such a combination. Since sums of multisets of lines in different ambient spaces can be realized in a common ambient space, the length $n$ is attained.
\end{proof}

\begin{Theorem}
  \label{thm_N_explicit}
  For every $l\ge 1$ we have
  \begin{multline*}
    \cN_2^{2^l,2}=\left\{2^{\max\{0,l+2-2t\}}\cdot\frac{4^t-1}{3}\,:\,1\le t\le \frac{l+3}{2}\right\}\cup\\
    \left\{2^{l+2-r}\cdot\frac{3\cdot 2^r-1-2^{r-b}-2^{r-c}}{3}\,:\,1\le c<b<r\le l+2,\ r,b,c\text{ odd}\right\}
  \end{multline*}
  and $\#\cN_2^{2^l,2}=\left\lfloor\frac{l+3}{2}\right\rfloor+\binom{\lfloor (l+1)/2\rfloor+1}{3}$.
\end{Theorem}
\begin{proof}[Proof sketch, see Appendix~\ref{app_sec4_lines}]
  The stated numbers generate the semigroup by the proof of Theorem~\ref{thm_D}. Their indecomposability follows by comparing the number of generators of $\mathcal{L}_2^{l+1}$ in a representation of $3n$ with a binary digit sum.
\end{proof}

So $\cN_2^{2^l,2}$ consists of the lengths of the $\F_4$-linear simplex codes of dimension $t\le (l+3)/2$, taken $2^{\max\{0,l+2-2t\}}$ times, and of the examples of
Corollary~\ref{cor_conj_partition} with $r\le l+2$, taken $2^{l+2-r}$ times. In particular, the examples of Corollary~\ref{cor_conj_partition} are base examples (choose $l=r-2$).

\subsection{Constructions for \texorpdfstring{$\Delta=q^{h+1}$}{Delta=q\textasciicircum{}(h+1)}}\label{sub_constructions}
For $h\ge 2$ let
\begin{equation}
  \label{eq_S_qh}
  \cG_q^h:=\left\{q^{h+1},\,q^{h+1}+q,\,q^{2h}+q^h+1\right\}\cup\left\{q^{h+2}+q^{h+1}+c,\,q^{h+2}+q^{h+1}+q^2+c\,:\,1\le c\le q-1\right\}.
\end{equation}
In this subsection we give families of constructions that realize all elements of $\cG_q^h$ as lengths of $q^{h+1}$-divisible faithful additive codes over $\F_{q^h}$.
Theorem~\ref{thm_N_qh1} at the end of this subsection shows that $\cG_q^h$ determines $\cN_q^{q^{h+1},h}$. The constructions are also the main source of the other known sets $\cN_q^{\Delta,h}$,
which are collected in Table~\ref{table_N_sets} in Appendix~\ref{app_N_sets}.

Besides Lemma~\ref{lemma_moves} we use the following $h$-partitionable types.

\begin{Lemma}
  \label{lemma_moves_general}
  Let $h\ge 2$ and let $S_1<S_2<\dots<S_v$ be a chamber of $\PG(v-1,q)$, $S_0:=\{0\}$, and $0\le a\le h-1$.
  \begin{enumerate}
    \item[(d)] (Intervals) For $h\le b\le v$ the $(h-1)$-spaces $X$ with $S_a\le X\le S_b$ have the point multiset
               $\genfrac{[}{]}{0pt}{}{b-a-1}{h-a-1}_q\chi_{S_b}+\left(\genfrac{[}{]}{0pt}{}{b-a}{h-a}_q-\genfrac{[}{]}{0pt}{}{b-a-1}{h-a-1}_q\right)\chi_{S_a}$.
    \item[(e)] (Cones over lifted MRD codes) For $m\ge h$ and $i:=m+h-a\le v$ the multiset $\chi_{S_i}-\chi_{S_m}+q^{m-a}\chi_{S_a}$ is $h$-partitionable by $q^{m-a}$ $(h-1)$-spaces containing $S_a$.
    \item[(f)] (Cones over spreads) If $h-a$ divides $i-a$, then $\chi_{S_i}+\frac{q^{i-a}-q^{h-a}}{q^{h-a}-1}\chi_{S_a}$ is $h$-partitionable.
    \item[(g)] (Parallel classes) Let $h\ge 3$, $u\ge 1$ with $h+1+u\le v$, $k\ge 1$, $N:=k[u]_q$, and $m_0\ge q\left\lceil N/q^{h-1}\right\rceil$. Then
               $k\chi_{S_{h+1+u}}+\left(m_0[h]_q-k[u+1]_q\right)\chi_{S_{h+1}}$ is $h$-partitionable. For $h\ge 4$ the same holds if $m_0\ge (q+1)\left\lceil N/\left(q^2[h-2]_q\right)\right\rceil$.
  \end{enumerate}
\end{Lemma}
\begin{proof}
  (d) The $(h-1)$-spaces $X$ with $S_a\le X\le S_b$ correspond to the $(h-a)$-dimensional subspaces of $S_b/S_a$. A point of $S_a$ lies in all of them, and a point $P$ of $S_b$ outside $S_a$ lies in those
  that contain $S_a+P$. (e) Apply Lemma~\ref{lemma_moves}.(a) in $S_i/S_a$ to $S_m/S_a$ (for $h-a=1$ take all points of $S_i/S_a$ outside $S_m/S_a$) and take preimages. Each preimage contains $S_a$, the points
  outside $S_m$ are covered once and the points of $S_m$ outside $S_a$ are not covered. (f) Take the preimages of the elements of a Desarguesian $(h-a-1)$-spread of $S_i/S_a$.
  Part~(g) is proved in Appendix~\ref{app_sec4_qh1}.
\end{proof}
Lemma~\ref{lemma_moves}.(b) is the case $a=h-1$ of part~(d). Since sums of $h$-partitionable multisets of points are $h$-partitionable, every non-negative integer combination of the
types in Lemmas~\ref{lemma_moves} and~\ref{lemma_moves_general}, all with respect to the same chamber, is an $h$-partitionable type, and its divisibility can be read off from Lemma~\ref{lemma_type_divisible}.
All constructions of this subsection are of this form.

For example, Theorem~\ref{thm_C} combines (b) for $j\in\{h+1,2h+1\}$ with (c), and Proposition~\ref{prop_ring_family} combines (a) with (d) for $(a,b)=(0,h+1)$. The types of
Theorem~\ref{thm_plus_q2} combines (a) with (g) for $(k,u,m_0)=\left(q^2,h-2,(q-1)(q-c)+q\right)$. Finding a shortest such type with a prescribed residue
of the cardinality modulo $\Delta$ is a small integer linear program. We use it for the results for $\Delta=q^{h+2}$ below.

\begin{Theorem}
  \label{thm_C}
  Let $h\ge 2$ and $1\le c\le q-1$. Let $A<B$ be subspaces of $V=\F_q^{2h+1}$ with $\dim(A)=h-1$ and $\dim(B)=h+1$. Let $\cM_c$ be the multiset of $(h-1)$-spaces consisting of
  \begin{itemize}
    \item $c$ copies of each of the $[h+2]_q$ $(h-1)$-spaces of $V$ containing $A$,
    \item $c$ copies of each of the $[h+1]_q-q-1$ $(h-1)$-spaces of $B$ not containing $A$, and
    \item $y_c:=q^{h+1}-cq[h]_q$ further copies of each of the $q+1$ $(h-1)$-spaces between $A$ and $B$.
  \end{itemize}
  Then $\cM_c$ is $q^{h+1}$-divisible and $\#\cM_c=q^{h+2}+q^{h+1}+c$. So, $q^{h+2}+q^{h+1}+c$ is the length of a $q^{h+1}$-divisible faithful additive code over $\F_{q^h}$.
\end{Theorem}
\begin{proof}
  Since $cq[h]_q\le (q-1)q[h]_q=q^{h+1}-q$, we have $y_c\ge q\ge c$. By Lemma~\ref{lemma_moves}.(b) with $j=2h+1$ the $c$ copies of the $(h-1)$-spaces of $V$ containing $A$ have the point multiset
  $c\left(\chi_V+q[h+1]_q\chi_A\right)$. The remaining elements of $\cM_c$ lie in $B$: the $q+1$ $(h-1)$-spaces between $A$ and $B$ are taken $y_c$ times and the other $(h-1)$-spaces of $B$ are taken $c$ times.
  By Lemma~\ref{lemma_moves}.(c) with $\alpha=[h]_qc+y_c-c$ and $\beta=q\left(y_c-c\right)$ their point multiset is $\alpha\chi_B+\beta\chi_A$.
  Since $\alpha=q^{h+1}-c\left(q^h-1\right)-c=(q-c)q^h$ and $cq[h+1]_q+\beta=cq[h+1]_q+q^{h+2}-cq^2[h]_q-cq=q^{h+2}$, the point multiset of $\cM_c$ is
  $$
    \cM_c'=c\cdot\chi_V+(q-c)q^h\cdot\chi_B+q^{h+2}\cdot\chi_A,
  $$
  i.e., it has type $c[2h+1]+(q-c)q^h[h+1]+q^{h+2}[h-1]$.
  Since $q^{h}$ divides $(q-c)q^h$ and $q^{h+2}$ divides $q^{h+2}$, $\cM_c'$ is $q^{2h}$-divisible by Lemma~\ref{lemma_type_divisible} (with $r=2h+1$). By Lemma~\ref{lemma_exclusion} the multiset $\cM_c$ is
  $q^{2h}/q^{h-1}=q^{h+1}$-divisible. Finally,
  $\#\cM_c=\#\cM_c'/[h]_q$ and $c[2h+1]_q+(q-c)q^h[h+1]_q+q^{h+2}[h-1]_q=[h]_q\left(q^{h+2}+q^{h+1}+c\right)$, which follows from $[2h+1]_q=q^{h+1}[h]_q+[h+1]_q$ and $[h+1]_q=q^h+[h]_q$.
\end{proof}
By Remark~\ref{rem_type_weights} the non-zero weights of $\cM_c$ are $cq^{h+1}$, $q^{h+2}$, and $q^{h+2}+q^{h+1}$. For $c=1$ we obtain, e.g., $25\in\cN_2^{8,2}$,
$49\in\cN_2^{16,3}$, $97\in\cN_2^{32,4}$, $109\in\cN_3^{27,2}$, and $321\in\cN_4^{64,2}$.

\begin{Proposition}
  \label{prop_ring_family}
  Let $h\ge 2$ and $c,d\in\N_0$ with $c\ge 1$ and $[h]_q\mid c+dq^h$. Then the multiset of points $c\cdot\chi_V+dq^h\cdot\chi_B$ of type $c[2h+1]+dq^h[h+1]$ in $V=\F_q^{2h+1}$, where $B=S_{h+1}$, is $h$-partitionable: take $c$ times the lifted MRD partition of
  $V\backslash B$ from Lemma~\ref{lemma_moves}.(a) and $(c+dq^h)/[h]_q$ times all $(h-1)$-spaces of $B$. The resulting multiset of $(h-1)$-spaces is $q^{h+1}$-divisible. For $c+d=[h]_q$ its cardinality is $q^h[h+1]_q+c$.
\end{Proposition}
\begin{proof}
  Every point of $B$ lies in $[h]_q$ $(h-1)$-spaces of $B$, which gives the type. It is $q^{2h}$-divisible by Lemma~\ref{lemma_type_divisible}. Hence Lemma~\ref{lemma_exclusion} gives $q^{h+1}$-divisibility.
  For $c+d=[h]_q$ the cardinality is $cq^{h+1}+\left(q^h-c(q-1)\right)[h+1]_q=q^h[h+1]_q+c$.
\end{proof}
For $h=2$ this gives the lengths $q^4+q^3+q^2+c$, $1\le c\le q+1$, e.g., $119\in\cN_3^{27,2}$ and $338,339\in\cN_4^{64,2}$. The lengths $q^4+q^3+q^2+1=\left(q^4+q^2+1\right)+q^3$, $q^4+q^3+q^2+q=(q+1)\left(q^3+q\right)$, and $[5]_q=\left(q^4+q^2+1\right)+\left(q^3+q\right)$
are decomposable. So only $2\le c\le q-1$ yields new elements of $\cN_q^{q^3,2}$.

\begin{Theorem}
  \label{thm_plus_q2}
  Let $h\ge 3$ and $1\le c\le q-1$. Then the multiset of points
  $$
    \cT:=c\cdot\chi_V+(q-c)q^h\cdot\chi_B+q^2\cdot\chi_D
  $$
  of type $c[2h+1]+(q-c)q^h[h+1]+q^2[2h-1]$ in $V=\F_q^{2h+1}$, where $B<D$, $\dim(B)=h+1$, $\dim(D)=2h-1$, is $h$-partitionable. The corresponding multiset of $(h-1)$-spaces is $q^{h+1}$-divisible and has cardinality
  $q^{h+2}+q^{h+1}+q^2+c$. So, there exists a $q^{h+1}$-divisible faithful additive code over $\F_{q^h}$ of length $q^{h+2}+q^{h+1}+q^2+c$.
\end{Theorem}
\begin{proof}
  The type is $q^{2h}$-divisible by Lemma~\ref{lemma_type_divisible}, since $q^2$ divides the coefficient of $[2h-1]$ and $q^h$ divides the coefficient of $[h+1]$. Hence the divisibility follows from Lemma~\ref{lemma_exclusion}.
  The cardinality follows from $[2h-1]_q-q^h[h-1]_q=[h]_q$ and Theorem~\ref{thm_C}. We have $\cT=c\left(\chi_V-\chi_B\right)+\left(q^2\chi_D+\left((q-c)q^h+c\right)\chi_B\right)$.
  The first summand is $c$ times the lifted MRD partition of Lemma~\ref{lemma_moves}.(a) with $i=2h+1$. For the second summand apply Lemma~\ref{lemma_moves_general}.(g) with $(k,u,m_0)=\left(q^2,h-2,(q-1)(q-c)+q\right)$,
  which gives $m_0[h]_q-q^2[h-1]_q=(q-c)q^h+c$. Here $N=q^2[h-2]_q$, so that $m_0\ge q=q\lceil N/q^{2}\rceil$ for $h=3$ and $m_0\ge q+1=(q+1)\lceil N/(q^2[h-2]_q)\rceil$ for $h\ge 4$, since $(q-1)(q-c)\ge 1$.
\end{proof}
For $h=3$ the construction uses a point $\langle u\rangle$ outside $B$ and a line spread of $B\cong\PG(3,q)$ that consists of $q^2$ lines $L$ and one further line $L_\infty$.
The construction then consists of $c$ lifted MRD partitions of $V\backslash B$ and, for each spread line $L\neq L_\infty$, the $q^2$ planes $L+\langle u+w\rangle$ with $w$ in a complement of $L$ in $B$.
In addition, it takes every plane of $B$ through $L_\infty$ exactly $(q-1)(q-c)+q$ times and every other plane of $B$ exactly $(q-1)(q-c)$ times.
This gives, e.g., $53\in\cN_2^{16,3}$ and $334,335\in\cN_3^{81,3}$. For $h=2$ the type of Theorem~\ref{thm_plus_q2} coincides with the one of Proposition~\ref{prop_ring_family}.

\begin{Theorem}
  \label{thm_N_qh1}
  Let $h\ge 2$ and let $\cG_q^h$ be the set defined in~(\ref{eq_S_qh}).
  \begin{enumerate}
    \item[(a)] For every $h\ge 2$ every element of $\cG_q^h$ is the length of a $q^{h+1}$-divisible faithful additive code over $\F_{q^h}$.
    \item[(b)] For every $h\ge 2$ and every prime power $q$ a positive integer $n$ is the length of a $q^{h+1}$-divisible faithful additive code over $\F_{q^h}$ if and only if $n\in\mathcal{A}_q^h(h+1)$ (see Subsection~\ref{sub_qh23}),
               i.e., if and only if $n$ satisfies the necessary conditions of Lemma~\ref{lemma_exclusion} and Lemma~\ref{lemma_counting}, and this is the case if and only if $n$ is a non-negative integer combination of the elements of $\cG_q^h$.
               Moreover, $\cN_q^{q^{h+1},h}=\cG_q^h$ for $h\ge 3$ and $\cN_q^{q^3,2}=\cG_q^2\backslash\left\{q^4+q^3+q^2+1\right\}$. For example, $\cN_3^{243,4}=\cG_3^4=\{243,246,973,974,982,983,6643\}$.  \end{enumerate}
\end{Theorem}
\begin{proof}
  (a) The first three elements are realized by a $q^{h+1}$-fold $(h-1)$-space, by $q$ copies of a Desarguesian $(h-1)$-spread of $\PG(2h-1,q)$, and by the field reduction of all points of $\PG(2,q^h)$, respectively.
  The others follow from Theorem~\ref{thm_C}, Proposition~\ref{prop_ring_family} (for $h=2$), and Theorem~\ref{thm_plus_q2} (for $h\ge 3$).

  Part~(b) is proved after Corollary~\ref{cor_gens_qhj}.
\end{proof}
\begin{Example}
  For $q=2$ the set $\cG_2^h$ is
  $$
    \left\{2^{h+1},\,2^{h+1}+2,\,2^{h+2}+2^{h+1}+1,\,2^{h+2}+2^{h+1}+5,\,2^{2h}+2^h+1\right\}.
  $$
  For $h=3$ Theorem~\ref{thm_N_qh1}.(b) gives $\cN_2^{16,3}=\{16,18,49,53,73\}$, so the lengths of $16$-divisible faithful additive codes over $\F_8$ are the non-negative combinations of these five numbers. The number $17$ is not a length, although it satisfies the counting condition: for $\rho=1$ we need $n>2^3=8$. Lemma~\ref{lemma_exclusion} requires $[3]_2\cdot 17=119$ to be the cardinality of a $64$-divisible multiset of points over $\F_2$. The cardinalities of such multisets are the combinations of $64,96,112,120,124,126,127$, and $119$ is not one of them. In contrast, $18$ satisfies both conditions, since $7\cdot 18=126$ is one of these generators and $18>2\cdot 8$, so that it is a length by Theorem~\ref{thm_N_qh1}.(b).
\end{Example}

\subsection{The cases \texorpdfstring{$\Delta=q^{h+2}$}{Delta=q\textasciicircum{}(h+2)} and \texorpdfstring{$\Delta=q^{h+3}$}{Delta=q\textasciicircum{}(h+3)}}\label{sub_qh23}
The set $\cG_q^h$ is a special case of a general pattern, which we now describe for all divisibility constants $\Delta=q^l$ with $l\ge h$ simultaneously.
Which of the minimal generators for $\Delta=q^{h+2}$ and $\Delta=q^{h+3}$ are realized by the constructions of this subsection is recorded in Remark~\ref{rem_N_qh3} in Appendix~\ref{app_N_sets}.
Recall from \cite[Theorem 1]{kiermaier2020lengths} that the set $\mathcal{L}_q^r$ of cardinalities of $q^r$-divisible multisets of points consists of the integers whose $S_q(r,0)$-adic expansion has a non-negative
leading coefficient, where $s_q(r,0,i)=q^i[r+1-i]_q$ for $0\le i\le r$, see Equation~(\ref{eq_s_q_adic}). For $l\ge h$ let $\mathcal{A}_q^{h}(l)$ be the set of non-negative integers $n$ such that
\begin{itemize}
  \item[(N1)] $[h]_q\cdot n\in\mathcal{L}_q^{l+h-1}$ and
  \item[(N2)] $n\equiv\rho\pmod{q^l}$ with $1\le\rho<q^l$ implies $n>\rho q^h$.
\end{itemize}
By Lemma~\ref{lemma_exclusion} and Lemma~\ref{lemma_counting} every length of a $q^l$-divisible faithful additive code over $\F_{q^h}$ lies in $\mathcal{A}_q^{h}(l)$, and $\mathcal{A}_q^h=\mathcal{A}_q^{h}(h+1)$. In other words, $\mathcal{A}_q^{h}(l)$ is the set of integers that pass the necessary conditions of Lemma~\ref{lemma_exclusion} and Lemma~\ref{lemma_counting} for $\Delta=q^l$. It contains the set $\mathcal{A}_q^{q^l,h}$ of Subsection~\ref{sub_nonpower}, which in addition uses Lemma~\ref{lemma_moment}.
The set $\mathcal{A}_q^{h}(l)$ is closed under addition. For (N1) this follows from taking sums of multisets of points. For (N2) note that $n_i\equiv\rho_i\pmod{q^l}$, $n_i\ge 1$, and $\rho_1+\rho_2\ge 1$ imply $n_1+n_2>(\rho_1+\rho_2)q^h$, and the residue of $n_1+n_2$ is at most $\rho_1+\rho_2$. Since $\rho q^h\ge\rho$, (N2) excludes all $1\le n<q^l$, so that $q^l$ is the smallest positive element of $\mathcal{A}_q^{h}(l)$. For $1\le\rho<q^l$ let $m_\rho^{(l)}$ denote the smallest element of $\mathcal{A}_q^{h}(l)$
that is congruent to $\rho$ modulo $q^l$. So $\left\{m_\rho^{(l)}\right\}$ is the Ap\'ery set of $\mathcal{A}_q^{h}(l)$ with respect to $q^l$. The smallest integer $k$ with $\rho+kq^l>\rho q^h$ is
$$
  \kappa_l(\rho):=\left\lfloor\frac{\rho\left(q^h-1\right)}{q^l}\right\rfloor+1,
$$
so that (N2) alone gives $m_\rho^{(l)}\ge\rho+\kappa_l(\rho)q^l$. Finally, let $v_q(\rho)$ be the largest integer $i$ with $q^i\mid\rho$.

\begin{Lemma}
  \label{lemma_divide_q}
  Let $l\ge h+1$. Then $\mathcal{A}_q^{h}(l)\cap q\Z=q\cdot\mathcal{A}_q^{h}(l-1)$. In particular, $m_{q\rho'}^{(l)}=q\cdot m_{\rho'}^{(l-1)}$ for all $1\le\rho'<q^{l-1}$.
\end{Lemma}
\begin{proof}
  Let $r\ge 2$ and $x\in\N_0$ with $q\mid x$. Then the $S_q(r,0)$-adic expansion of $x$ has $c_0=0$, and since $s_q(r,0,i)=q\cdot s_q(r-1,0,i-1)$ for $1\le i\le r$, the remaining digits form the $S_q(r-1,0)$-adic expansion of
  $x/q$ with the same leading coefficient. So, $x\in\mathcal{L}_q^r$ if and only if $x/q\in\mathcal{L}_q^{r-1}$. Since $[h]_q\equiv 1\pmod q$, we can apply this to $x=[h]_q n$ with $r=l+h-1$. So (N1) for $n=qn'$ and $l$ is
  equivalent to (N1) for $n'$ and $l-1$. Moreover, $qn'\equiv q\rho'\pmod{q^l}$ if and only if $n'\equiv\rho'\pmod{q^{l-1}}$, and $qn'>q\rho'q^h$ if and only if $n'>\rho'q^h$. The conditions (N2) correspond to each other as well.
\end{proof}

For $l=h$ the proof of Theorem~\ref{thm_delta_divides_qh} gives $m_\rho^{(h)}=\rho\left(q^h+1\right)=\rho+\kappa_h(\rho)q^h$.

\begin{Theorem}
  \label{thm_apery_qhj}
  Let $j$ and $h$ be integers with $1\le j\le q$ and $h\ge j+2$, or $j=1$ and $h=2$, and let $l=h+j$. Then for all $1\le\rho<q^{l}$ we have
  \begin{equation}
    \label{eq_apery_qhj}
    m_\rho^{(l)}=\rho+\left(\kappa_l(\rho)+q\cdot\max\left\{0,\,j-v_q(\rho)-\left\lfloor\rho/q^{j+1}\right\rfloor\right\}\right)q^l.
  \end{equation}
\end{Theorem}
\begin{proof}[Proof sketch, see Appendix~\ref{app_sec4_apery}]
  The proof is by induction on $j$. Residues divisible by $q$ are reduced to $j-1$ by Lemma~\ref{lemma_divide_q}, and for $q\nmid\rho$ the right hand side of Equation~(\ref{eq_apery_qhj}) is written as a sum of explicit elements of $\mathcal{A}_q^{h}(l)$, while an explicit $S_q(l+h-1,0)$-adic expansion shows that (N1) fails for all smaller elements of the residue class.
\end{proof}

\begin{Corollary}
  \label{cor_gens_qhj}
  Let $1\le j\le q$ and $h\ge j+2$, where $h\ge j+3$ if $j=q$. Then the minimal generating set of the semigroup $\mathcal{A}_q^{h}(h+j)$ is
  \begin{eqnarray*}
    \mathcal{G}_q^{h,j}&:=&\left\{q^{h+j},q^{h+j}+q^j,q^{2h}+q^h+1\right\}\\
    & &\cup\left\{\left((j-i)q+1\right)q^{h+j}+q^iu+fq^{j+1}\,:\,(i,u,f)\in I_q^j\right\},
  \end{eqnarray*}
  where $I_q^j$ is the set of triples of integers with $0\le i\le j-1$, $1\le u<q^{j-i}$, $q\nmid u$, and $0\le f\le j-i$. This set has $3+(q-1)\sum_{t=1}^{j}(t+1)q^{t-1}$ elements. In particular, $\mathcal{G}_q^{h,1}=\cG_q^h$ for $h\ge 3$.
  For $j=2$, i.e., $\Delta=q^{h+2}$, and $h\ge 4$, where $h\ge 5$ if $q=2$, the minimal generators are $q^{h+2}$, $q^{h+2}+q^2$, $q^{2h}+q^h+1$, $q^{h+3}+q^{h+2}+cq+\eta q^3$, and $(2q+1)q^{h+2}+c+dq+fq^3$, where $1\le c\le q-1$, $\eta\in\{0,1\}$, $0\le d\le q-1$, and $f\in\{0,1,2\}$.
\end{Corollary}
\begin{proof}[Proof sketch, see Appendix~\ref{app_sec4_apery}]
  The proof of Theorem~\ref{thm_apery_qhj} writes every $m_\rho^{(h+j)}$ as a sum of elements of $\mathcal{G}_q^{h,j}$, and minimality follows by induction on $j$ from Lemma~\ref{lemma_divide_q} and the lower bounds of Theorem~\ref{thm_apery_qhj}.
\end{proof}
\begin{proof}[Proof of Theorem~\ref{thm_N_qh1}.(b)]
  By Lemma~\ref{lemma_exclusion} and Lemma~\ref{lemma_counting} every length lies in $\mathcal{A}_q^h(h+1)$. The proof of Corollary~\ref{cor_gens_qhj} that $\mathcal{G}_q^{h,1}=\cG_q^h$ generates $\mathcal{A}_q^h(h+1)$ uses only Lemma~\ref{lemma_divide_q} and Theorem~\ref{thm_apery_qhj}, which includes the case $(j,h)=(1,2)$.
  So, $\cG_q^h$ generates $\mathcal{A}_q^h(h+1)$ for all $h\ge 2$, and by (a) the lengths are exactly the positive elements of $\langle\cG_q^h\rangle=\mathcal{A}_q^h(h+1)$. For $h\ge 3$ Corollary~\ref{cor_gens_qhj} with $j=1$ gives $\cN_q^{q^{h+1},h}=\cG_q^h$.

  Let $h=2$. The element $q^4+q^3+q^2+1=\left(q^4+q^2+1\right)+q^3$ is decomposable. For $q^3$, $q^3+q$, and $q^4+q^2+1$ the argument of the proof of Corollary~\ref{cor_gens_qhj} applies unchanged.
  So let $x=y+z$ with $x=(q+1)q^3+\rho_x$, $\rho_x\in\{c,q^2+c\}$ and $1\le c\le q-1$, and $y,z$ as there, i.e., $k_y+k_z=q+1$, $k_y,k_z\ge 1$, $\rho_y+\rho_z=\rho_x$, and $q\nmid\rho_y$.
  If $\rho_y<q^2$, then Equation~(\ref{eq_apery_qhj}) gives $k_y\ge\kappa_3(\rho_y)+q\ge q+1$, which is impossible. Otherwise $\rho_x=q^2+c$ and $\rho_y=q^2+c_y$ with $1\le c_y\le c$, and $\rho_z=c-c_y$.
  Equation~(\ref{eq_apery_qhj}) gives $k_y\ge\kappa_3(q^2+c_y)$, which equals $q$ for $c_y=1$ and $q+1$ for $c_y\ge 2$. Hence $c_y=1$, $k_z=1$, and $\rho_z=c-1$. For $c=1$ this is the decomposition of $x=q^4+q^3+q^2+1$ above.
  For $c\ge 2$ we have $q\nmid\rho_z<q^2$, so that $k_z\ge q+1$ as before, a contradiction.
\end{proof}

\begin{Remark}
  Theorem~\ref{thm_apery_qhj} explains the shape of $\cG_q^h$ and of the generators in Remark~\ref{rem_N_qh3}. The counting condition (N2) determines $m_\rho^{(l)}$ up to the penalty term $q\cdot\max\{\dots\}$ in
  Equation~(\ref{eq_apery_qhj}), which is non-zero only for the residues $\rho<\left(j-v_q(\rho)\right)q^{j+1}$, and the point condition (N1) is needed exactly for these residues. For $j=1$, i.e., for the lengths of $q^{h+1}$-divisible faithful additive codes in Theorem~\ref{thm_N_qh1}.(b), these are the residues $\rho<q^2$ with $q\nmid\rho$, while all other lower bounds already follow from (N2).

  Apart from the three lengths $q^{h+j}$, $q^j\left(q^h+1\right)$, and $[3]_{q^h}$ of the $q^{h+j}$-fold $(h-1)$-space, $q^j$ copies of a Desarguesian spread, and the field reduction of $\PG(2,q^h)$,
  the elements of $\mathcal{G}_q^{h,j}$ come in $j$ layers indexed by $0\le i\le j-1$. The elements of layer $0$ are $(jq+1)q^{h+j}+u+fq^{j+1}$, i.e., the residues $u+fq^{j+1}$ with $q\nmid u<q^j$
  and $0\le f\le j$, all with the same quotient $jq+1$. For $j=1$ these are the lengths of Theorem~\ref{thm_C} ($f=0$) and Theorem~\ref{thm_plus_q2} ($f=1$).
  Layer $i$ consists of $q^i$ times the elements of layer $0$ for $\Delta=q^{h+j-i}$, which are inherited via $q$-fold repetition (Lemma~\ref{lemma_divide_q}). The hypothesis $j\le q$ cannot be dropped, and for $j=q$ and $h=j+2$ the set $\mathcal{G}_q^{h,j}$ is not minimal, see Remark~\ref{rem_apery_nonminimal}.
\end{Remark}
The elements of $\cN_4^{32,2}=\{32,34,161,169,273\}$ are realized as follows. For $161$ and $169$ let $V=\F_4^5$ and let $P$ be a point and $E\ni P$ be a plane.
\begin{enumerate}[(a)]
  \item The $85$ lines through $P$, together with $12$ copies of each of the $5$ lines of $E$ through $P$ and one copy of each of the $16$ other lines of $E$, have the point multiset $\chi_V+16\chi_E+128\chi_P$ by Lemma~\ref{lemma_moves}.(b),(c) (with $\alpha=16$, $\beta=44$).
  \item A lifted MRD partition of $V\backslash E$ into $64$ lines together with $5$ copies of each of the $21$ lines of $E$ has the point multiset $\chi_V+24\chi_E$.
  \item Both point multisets are $128$-divisible, since $\#\chi_V-\chi_V(H)=256$, $\#\chi_E-\chi_E(H)\in\{0,16\}$, and $\#\chi_P-\chi_P(H)\in\{0,1\}$ for every hyperplane $H$. Hence Lemma~\ref{lemma_exclusion} gives $32$-divisible multisets of $161$ and $169$ lines, respectively.
\end{enumerate}
The lengths $32$, $34$, and $273$ are realized by a $32$-fold line, two copies of a line spread of $\PG(3,4)$, and the field reduction of $\PG(2,16)$. Here the necessary condition of Lemma~\ref{lemma_exclusion} uses the lengths of $2^7$-divisible multisets of points over $\F_4$.

For $\Delta=q^{h+j}$ the analog of Theorem~\ref{thm_C} is a family of types in $\F_q^{2h+j}$ that is described by the following lemma. The basic building blocks are cones with vertex $A=S_{h-1}$, which put too much weight on $A$.
The numbers of the following units replace cones and put less weight on $A$:
\begin{itemize}
  \item $\mu_k$ is the number of lifted MRD partitions of $S_{2h+j-k}\backslash S_{h+j-k}$, for $0\le k\le j-1$,
  \item $\sigma$ is the number of spread units in $F=S_{2h+1}$,
  \item $\tau$ and $\nu$ are the numbers of units of Lemma~\ref{lemma_moves_general}.(g) in $D=S_{2h-1}$ and in $E=S_{h+2}$, respectively.
\end{itemize}
The remaining two numbers are determined by the others: $y_{AB}$ is the number of copies of each $(h-1)$-space between $A$ and $B=S_{h+1}$, and $y_B$ is the number of copies of each hyperplane of $B$.
They balance the coefficients of $\chi_B$ and $\chi_A$.
\begin{Lemma}
  \label{lemma_cone_family}
  Let $h\ge 2$, $j\ge 1$, $0\le g\le j$, let $S_1<\dots<S_{2h+j}=V=\F_q^{2h+j}$ be a chamber, $A:=S_{h-1}$, $B:=S_{h+1}$, $E:=S_{h+2}$, $D:=S_{2h-1}$, and $F:=S_{2h+1}$. Let $c_0,\dots,c_{j-1}\in\{0,\dots,q-1\}$ with $c_0\ge 1$,
  $\rho:=\sum_{i=0}^{j-1}c_iq^i$. Let $\mu_0,\dots,\mu_{j-1},\sigma,\tau,\nu$ be non-negative integers with $\mu_k\le c_kq^k$ for $k\le j-2$, $\mu_{j-1}+\sigma\le c_{j-1}q^{j-1}$, $\tau\le gq^{j-1}$,
  $\nu\le q\left((c_{j-1}+1)q^{j-1}-\rho\right)$, and $\tau=\nu=0$ if $h=2$, and $\nu=0$ if $j=1$. If
  \begin{eqnarray*}
    y_B&:=&\rho+gq^{j+1}-\sum_{k=0}^{j-1}q^{j-k}\mu_k-\sigma-(q+1)\tau-q\nu\quad\text{and}\\
    y_{AB}&:=&q^h\left(q^j-\rho\right)+\sum_{k=0}^{j-1}\mu_k-\left([h]_q+q\right)\tau-\left([h]_q-1\right)\nu-[h]_qy_B
  \end{eqnarray*}
  are non-negative, then the multiset of points of type
  \begin{eqnarray*}
    \cT_j&:=&\sum_{k=0}^{j-1}c_kq^k[2h+j-k]+(q-c_0)q^h[h+j]+\sum_{k=1}^{j-1}(q-1-c_k)q^{h+k}[h+j-k]\\
    & &+gq^{j+1}[2h-1]+(j-g)q^{h+j+1}[h-1]
  \end{eqnarray*}
  in $V$ is $h$-partitionable. The corresponding multiset of $(h-1)$-spaces is $q^{h+j}$-divisible and has cardinality $(jq+1)q^{h+j}+\rho+gq^{j+1}$.
\end{Lemma}
\begin{proof}[Proof sketch, see Appendix~\ref{app_sec4_qh2}]
  The partition combines the units listed above, and comparing the coefficients of $\chi_B$ and $\chi_A$ gives the definitions of $y_{AB}$ and $y_B$.
\end{proof}

\begin{Corollary}
  \label{cor_qh3}
  \begin{enumerate}
    \item[(a)] For all $j\ge 1$, $h\ge 2$, and $1\le\rho<q^j$ with $q\nmid\rho$ the number $(jq+1)q^{h+j}+\rho$ is the length of a $q^{h+j}$-divisible faithful additive code over $\F_{q^h}$.
    \item[(b)] Let $h\ge 3$, $0\le a\le q-1$, and $1\le b\le q-1$. Then the multiset of points of type
                $$
                  b[2h+2]+aq[2h+1]+(q-b)q^h[h+2]+(q-1-a)q^{h+1}[h+1]+q^3[2h-1]+q^{h+3}[h-1]
                $$
                in $V=\F_q^{2h+2}$ is $h$-partitionable. Consequently, $2q^{h+3}+q^{h+2}+q^3+aq+b$ is the length of a $q^{h+2}$-divisible faithful additive code over $\F_{q^h}$, and
                $2q^{h+4}+q^{h+3}+q^4+aq^2+bq$ is the length of a $q^{h+3}$-divisible faithful additive code over $\F_{q^h}$.
    \item[(c)] Let $1\le c\le q-1$ and $0\le d,e\le q-1$. Then $(3q+1)q^{h+3}+c+dq+eq^2+gq^4$ is the length of a $q^{h+3}$-divisible faithful additive code over $\F_{q^h}$ for $g=1$ and all $h\ge 4$,
               and for $g=2$, all $h\ge 4$, and $q\in\{3,4,5,7,8,9\}$, provided that $e\le q-2$ or $c\ge 3$, and $(q,c,d,e)\neq(3,1,0,1)$.
  \end{enumerate}
\end{Corollary}
\begin{proof}[Proof sketch, see Appendix~\ref{app_sec4_qh2}]
  All parts are special cases of Lemma~\ref{lemma_cone_family}, with $g=0$ in~(a), $(j,g)=(2,1)$ in~(b), and $j=3$ in~(c). For $g=2$ in~(c) the parameters
  are found by a computer search.
\end{proof}
For $j=1$ part~(a) gives the lengths of Theorem~\ref{thm_C}. For $j=2$ and $\rho=b+aq$ the construction of part~(a) consists of $b$ copies of each $(h-1)$-space
through $A$, $a(q-1)$ copies of each such space in $S_{2h+1}$, $a$ lifted MRD partitions of $S_{2h+1}\backslash B$, $(q-b)q^h$ copies of each of the $q^2$ spaces $X$
with $A<X<S_{h+2}$ and $X\not\le B$, $y_{AB}$ further copies of each space between $A$ and $B$, and $b$ copies of each hyperplane of $B$.
By Remark~\ref{rem_type_weights} the non-zero weights of these codes are $bq^{h+2}$, $(a+b)q^{h+2}$, $(q+a)q^{h+2}$, $(2q-1)q^{h+2}$, and $(2q+1)q^{h+2}$, where coinciding values (for $a=0$ or $a=q-1$) are merged. For $q=3$ and $h=2$ part~(a) with $j=2$ gives the six lengths $568,569,571,572,574,575\in\cN_3^{81,2}$, which were first found by integer linear programming.
For $q=2$ it gives $2^{h+4}+2^{h+2}+1$ and $2^{h+4}+2^{h+2}+3$, e.g., $161,163\in\cN_2^{32,3}$.

\begin{Example}
  For $q=2$, $h=3$, and $\Delta=32$ the necessary conditions of Lemmas~\ref{lemma_exclusion} and~\ref{lemma_counting} read as follows. By Lemma~\ref{lemma_exclusion} the number $7n$ has to be the cardinality of a $128$-divisible multiset of points over $\F_2$, i.e., a combination of $128,192,224,240,248,252,254,255$. By Lemma~\ref{lemma_counting}, $n\equiv\rho\pmod{32}$ with $1\le\rho<32$ implies $n>8\rho$. For $n=33,\dots,36$ the conditions give the following.
  \begin{center}
    \begin{tabular}{c|c|c|c}
      $n$ & $n>8\rho$ & $7n$ & combination of the generators\\
      \hline
      $33$ & yes & $231$ & no\\
      $34$ & yes & $238$ & no\\
      $35$ & yes & $245$ & no\\
      $36$ & yes & $252$ & yes, a generator
    \end{tabular}
  \end{center}
  A computation of the semigroup $\mathcal{A}_2^3(5)$ of all integers satisfying both conditions gives the minimal generators $32$, $36$, $73$, $98$, $106$, $107$, $161$, $163$, and all of them are lengths, see Proposition~\ref{prop_N_qh2} in Appendix~\ref{app_N_sets}.
\end{Example}

For small parameters the sets are determined completely, for instance $\cN_2^{32,3}=\{32,36,73,98,106,107,161,163\}$, and $\cN_q^{q^4,2}$ for $q\in\{3,4,5,7\}$, see Proposition~\ref{prop_N_qh2} in Appendix~\ref{app_N_sets}.
The same method, i.e., a computation of the minimal generators of the semigroup of the necessary conditions and a computer search for non-negative integer combinations of the types of
Lemma~\ref{lemma_moves} and Lemma~\ref{lemma_moves_general} realizing each of them, gives the further sets $\cN_q^{\Delta,h}$ with status ``computer'' in Table~\ref{table_N_sets}, see Appendix~\ref{app_N_sets}.
In this way all minimal generators of the semigroup of necessary conditions are realized, with the exceptions of the following table.
\begin{center}
\small
\begin{tabular}{lll}
  \hline
  $(q,h,\Delta)$ & minimal generators not realized by this method & realized by\\
  \hline
  $(2,4,64)$ & $339$ & Proposition~\ref{thm_qh2_f2}\\
  $(2,5,128)$ & $657,659$ & Proposition~\ref{prop_qh2_packing}\\
  $(2,6,256)$ & $1297,1299$ & open\\
  $(3,4,729)$ & $5158,5159,5161,5162,5164,5165$ & Proposition~\ref{thm_qh2_f2}\\
  $(4,3,1024)$ & $9359$ & Proposition~\ref{thm_qh2_f2}\\
  \hline
\end{tabular}
\end{center}
All of them have the form $2q^{h+3}+q^{h+2}+2q^3+aq+b$. The constructions of Propositions~\ref{thm_qh2_f2} and~\ref{prop_qh2_packing} in Appendix~\ref{app_sec4_qh2} realize all of them except $1297$ and $1299$.

In contrast to Proposition~\ref{thm_qh2_f2}, the construction of Corollary~\ref{cor_qh3}.(b) works with the chamber types of Lemma~\ref{lemma_moves} and Lemma~\ref{lemma_moves_general} directly. For $(q,h)$ equal to $(2,4)$, $(2,5)$, and $(3,4)$ it gives the lengths $329,331\in\cN_2^{64,4}$, $649,651\in\cN_2^{128,5}$, and
$5131,5132,5134,5135,5137,5138\in\cN_3^{729,4}$, which were also first found by integer linear programming, see Table~\ref{table_N_sets}.

\begin{Remark}
  By Remark~\ref{rem_type_weights} all constructions of this subsection are few-weight codes. Table~\ref{table_few_weight} in Appendix~\ref{app_N_sets} lists the parameters $[n,r/h,d]_q^h$ and the weight
  enumerators of some of the resulting additive codes, where $A_w=(q-1)\cdot\#\left\{H\,:\,n-\widetilde{\cM}(H)=w\right\}$ is the number of codewords of weight $w$. The codes of
  Theorems~\ref{thm_C} and~\ref{thm_plus_q2} have three weights, and the codes of Corollary~\ref{cor_qh3}.(a) with $j=2$ have four or five weights.
  The codes of Proposition~\ref{prop_ring_family} with $c+d=[h]_q$ and $c<[h]_q$ have the two weights
  $cq^{h+1}$ and $[h]_qq^{h+1}$ with $A_{cq^{h+1}}=(q-1)[h]_q$ and $A_{[h]_qq^{h+1}}=(q-1)q^h[h+1]_q$. No code in Table~\ref{table_few_weight} is equivalent to an $\F_{q^h}$-linear
  code, since either $h$ does not divide $r$ or, for Corollary~\ref{cor_qh3}.(a) with $j=2$ and $h=2$, some $A_w$ is not divisible by $q^h-1$.
\end{Remark}

For $q=2$ and $6\le h\le 8$ the minimal generators listed in Remark~\ref{rem_N_qh3} in Appendix~\ref{app_N_sets} include the numbers $5\cdot 2^{h+3}+1+2d+4e+32g$ with $d,e,g\in\{0,1\}$. For $g=0$ the following
proposition realizes them by a type that contains $S_{h-2}$ instead of $S_{h-1}$.
\begin{Proposition}
  \label{thm_qh3_q2}
  Let $h\ge 4$ and $d,e\in\{0,1\}$. Then the multiset of points of type
  $$
    [2h+3]+2d[2h+2]+4e[2h+1]+2^h[h+3]+(1-d)2^{h+1}[h+2]+(1-e)2^{h+2}[h+1]+2^{h+5}[h-2]
  $$
  in $V=\F_2^{2h+3}$ is $h$-partitionable. Thus $5\cdot 2^{h+3}+1+2d+4e$ is the length of a $2^{h+3}$-divisible faithful additive code over $\F_{2^h}$.
\end{Proposition}
\begin{proof}[Proof sketch, see Appendix~\ref{app_sec4_qh3}]
  The partition combines telescoping unions of the cones over $S_{h-2}$ from Lemma~\ref{lemma_moves_general}.(e), lifted MRD partitions, the subspaces of Lemma~\ref{lemma_moves_general}.(d) with $(a,b)=(h-2,h+1)$, and hyperplanes of $S_{h+1}$.
\end{proof}

\section{Partitions of divisible sets of points}
\label{sec_partitions}

In Section~\ref{sec_length_powers} faithful additive codes over $\F_{q^h}$ were constructed by partitioning divisible multisets of points into $(h-1)$-spaces, see Subsection~\ref{sub_partitions}. In this section we study
partitions of small divisible sets of points into subspaces of possibly different dimensions.

We say that a set of points $\cS$ in $\PG(v-1,q)$ is of \emph{partition type} $k^{m_k}\dots 2^{m_2}1^{m_1}$ if it can be partitioned into $m_i$ $(i-1)$-spaces for $1\le i\le k$, where factors with $m_i=0$ may be omitted. If not stated otherwise, the set is also assumed to be $q^k$-divisible for the index $k$ of the leading factor in the notation, also if $m_k=0$, so that $2^11^5$ stands for a $4$-divisible set over $\F_2$. Higher divisibility is stated explicitly.

\noindent\emph{Counting bound.}\ A $q^r$-divisible multiset of points $\cS=\sum_{j=1}^m\chi_{U_j}$, where $U_j$ has dimension $k_j$, corresponds by Lemma~\ref{lemma_exchange} to the
projective system with $q^{k_j-1}$ copies of each $U_j$, which is $q^r$-divisible, since $\#\cS-\cS(H)=\sum_{j\,:\,U_j\not\le H}q^{k_j-1}$. So this section deals with
mixed projective systems with these multiplicities. Lemma~\ref{lemma_mixed_counting} with $h=\max_jk_j$, $N=q^{h-1}m$, and $\#\cM''=q^{h-1}\#\cS$ gives the following bound.
\begin{Corollary}
  Let $\cS=\sum_{j=1}^m\chi_{U_j}$ be a $q^r$-divisible multiset of points in $\PG(v-1,q)$, where $U_j$ has dimension $k_j$, and $\rho\equiv\sum_{j=1}^mq^{k_j-1}\pmod{q^r}$ with $1\le\rho<q^r$.
  Then $[v]_q\left(m-q\rho\right)\ge\#\cS$. In particular $m>q\rho$.
\end{Corollary}
For $q=2$ this excludes the $4$-divisible type $2^11^5$ ($m=6$, $\rho=3$), the $8$-divisible types $2^11^{13}$ and $2^91^5$ ($m=14$, $\rho=7$), and the $8$-divisible types $3^72^41^3$ and $3^72^5$ ($n=64$, with $m=14$, $\rho=7$ and $m=12$, $\rho=6$, respectively).
This gives short proofs of Lemma~\ref{lemma_n16}, of parts of Lemmas~\ref{lemma_exclusion_1} and~\ref{lemma_n32}, and of two exclusions in Table~\ref{table_types_60_70}. No realized type in
Section~\ref{sec_partitions} violates it. Among the types whose $4$-divisible part $2^b1^c$ is realizable, it excludes for $60\le n\le 70$ only these two types, and for $n\le 51$ further types
with many planes, e.g., $3^7$ for $n=49$.

The cardinality $\sum_{i=1}^k m_i[i]_q$ of a set of type $k^{m_k}\dots 1^{m_1}$ has to be the length of a projective $q^k$-divisible code over $\F_q$. Replacing one of the $(i-1)$-spaces by a vector space partition of $\PG(i-1,q)$
turns a realizable type into another realizable type, which we call \emph{derived} from the first one. For example, if type $3^{m_3}2^{m_2}1^{m_1}$ with $m_3>0$ can be realized in $\PG(v-1,q)$, then so can type $3^{m_3-1}2^{m_2+1}1^{m_1+q^2}$.

\begin{Lemma}[E.g.~{\cite[Theorem 13]{honold2018partial}}, {\cite[Proposition 3]{kurz2022vector}}]
  \label{lemma_proj_lengths}
  Let $n$ be the length of a non-trivial projective linear code.
  \begin{enumerate}
    \itemlab{a}{lemma_2_div_gf2} If the code is binary and $2$-divisible, then $n\ge 3$.
    \itemlab{b}{lemma_4_div_gf2} If the code is binary and $4$-divisible, then $n\in\{7,8\}$ or $n\ge 14$.
    \itemlab{c}{lemma_8_div_gf2} If the code is binary and $8$-divisible, then $n\in\{15, 16, 30, 31, 32, 45, 46, 47, 48, 49, 50, 51\}$ or $n\ge 60$ \cite{honold2019lengths}.
    \itemlab{d}{lemma_2_div_gf4} If the code is $2$-divisible over $\F_4$, then $n\ge 5$ \cite[Proposition 18]{kurz2024lengths}.
    \itemlab{e}{lemma_4_div_gf4} If the code is $4$-divisible over $\F_4$, then $n\in\{5, 10, 15, 16, 17\}$ or $n\ge 20$ \cite{kurz2021divisible}.
    \itemlab{f}{lemma_2_div_gf8} If the code is $2$-divisible over $\F_8$, then $n\in\{9,10\}$ or $n\ge 12$ \cite[Proposition 25]{kurz2024lengths}.
  \end{enumerate}
  All listed lengths are attained.
\end{Lemma}
The lemmas determining the types of $4$- and $8$-divisible binary sets of small cardinality are proven in Appendix~\ref{app_sec5}. In the $4$-divisible case they combine into the following statement.
\begin{Corollary}
  \label{cor_types_2_1}
  Let $m_1,m_2\in\N_0$ and $n=m_1+3m_2$. A set of points of type $2^{m_2}1^{m_1}$ exists in $\PG(v-1,2)$ for some $v$ if and only if $n=0$ or $n\in\{7,8\}\cup\{14,15,16,\dots\}$, $m_1\notin\{1,2\}$, and $(m_1,m_2)\notin\{(5,1),(5,3),(5,4)\}$.
\end{Corollary}
\begin{proof}
  The necessity follows from Lemma~\ref{lemma_4_div_gf2} (cardinality), the remark after Lemma~\ref{lemma_2_div_lines} ($m_1\notin\{1,2\}$), Lemma~\ref{lemma_exclusion_1} (for $(5,1)$ and $(5,3)$), and Lemma~\ref{lemma_exclusion_2} (for $(5,4)$). 
  For the sufficiency we use Proposition~\ref{prop_types_2_1_large} for $m_1\ge 14$ and the lemmas of Appendix~\ref{app_sec5_4div} for $m_1\le 13$, i.e., Lemma~\ref{lemma_2_div_lines} for $m_1=0$, Lemma~\ref{lemma_2_div_lines_p4} for $m_1=4$, and Lemmas~\ref{lemma_types_2_13t}, \ref{lemma_types_2_143t}, \ref{lemma_types_2_15}, \ref{lemma_types_2_18}, and~\ref{lemma_types_2_111} for the remaining residues of $m_1$.
\end{proof}
In particular, the bound $m_1\ge 14$ in Proposition~\ref{prop_types_2_1_large} is sharp.

\begin{Theorem}
  \label{thm_types_8div_small}
  Let $\cS$ be an $8$-divisible set of points in $\PG(v-1,2)$ with cardinality $n\le 51$, so that $n\in\{15,16,30,31,32,45,\dots,51\}$ by Lemma~\ref{lemma_8_div_gf2}. Then $\cS$ can be of each type in the
  second column of Table~\ref{table_types_small} and of each type derived from them, but not of any type in the third column. Every type $3^{m_3}2^{m_2}1^{m_1}$ with $7m_3+3m_2+m_1=n$ falls into one of two classes. Either it is derived from a type in the second column, or replacing some of its parts by vector space partitions yields a type in the third column. A type occurs if and only if it is derived from a type in the second column.
\end{Theorem}
\begin{proof}
  For each $n$ this is the lemma listed in the last column of Table~\ref{table_types_small}, see Appendix~\ref{app_sec5_8div}. Several of the non-existence statements rest on the enumeration of the projective $8$-divisible binary codes in Table~\ref{table_enumeration_8_div} and on small integer linear programs, see Appendix~\ref{app_sec5_ilp}. The additional statement follows by enumerating all types with $7m_3+3m_2+m_1=n$ for the cardinalities $n$ of Table~\ref{table_types_small}.
\end{proof}
\begin{table}[htp]
  \begin{center}
  \small
  \begin{tabular}{rp{5.4cm}p{5.4cm}r}
    \hline
    $n$ & realizable types & excluded types & Lemma\\
    \hline
    15 & $3^11^8$, $2^5$ & $3^21^1$, $3^12^11^5$, $3^12^21^2$ & \ref{lemma_n15}\\
    16 & $1^{16}$ & $2^11^{13}$ & \ref{lemma_n16}\\
    30 & $3^21^{16}$, $3^12^51^8$, $2^{10}$ & $3^31^9$, $3^22^11^{13}$, $3^12^61^5$ & \ref{lemma_n30}\\
    31 & $3^12^8$, $2^91^4$ & $3^21^{17}$, $2^{10}1^1$ & \ref{lemma_n31}\\
    32 & $2^81^8$ & $3^11^{25}$, $2^91^5$ & \ref{lemma_n32}\\
    45 & $3^32^8$, $3^22^91^4$, $3^12^{10}1^8$, $2^{15}$ & $3^41^{17}$, $3^22^{10}1^1$, $3^12^{11}1^5$ & \ref{lemma_n45}\\
    46 & $3^22^81^8$, $3^12^{13}$, $2^{14}1^4$ & $3^31^{25}$, $3^22^91^5$, $2^{15}1^1$ & \ref{lemma_n46}\\
    47 & $3^22^91^6$, $3^12^{12}1^4$, $2^{14}1^5$ & $3^31^{26}$, $3^22^{10}1^3$, $3^12^{13}1^1$, $2^{15}1^2$ & \ref{lemma_n47}\\
    48 & $3^12^{11}1^8$, $2^{16}$ & $3^21^{34}$, $3^12^{12}1^5$ & \ref{lemma_n48}\\
    49 & $2^{14}1^7$ & $3^11^{42}$, $2^{15}1^4$ & \ref{lemma_n49}\\
    50 & $2^{10}1^{20}$ & $3^11^{43}$, $2^{11}1^{17}$ & \ref{lemma_n50}\\
    51 & $2^{17}$ & $3^11^{44}$ & \ref{lemma_n51}\\
    \hline
  \end{tabular}
  \caption{Partition types $3^{m_3}2^{m_2}1^{m_1}$ of $8$-divisible sets of points in $\PG(v-1,2)$ with $n\le 51$, see Theorem~\ref{thm_types_8div_small}.}
  \label{table_types_small}
  \end{center}
\end{table}

For $8$-divisible sets of $n\ge 60$ points in $\PG(v-1,2)$ we use the following two constructions.
\begin{Lemma}[E.g.~{\cite[Lemma 19]{honold2018partial}}]
  \label{lemma_cone_minus_apex}
  Let $\cX$ be a $4$-divisible set of $m$ points in a hyperplane $H_0$ of $\PG(v-1,2)$, where $m\equiv 0\pmod 8$, and let $P\notin H_0$ be a point.
  Then $\cC:=\cX\cup\left(P+\cX\right)$, where $P+x$ denotes the third point on the line through $P$ and $x$, is an $8$-divisible set of $2m$ points.
  If $\cX$ can be partitioned into $\alpha$ planes, $\beta$ lines, and $\gamma$ points, then
  $\cC$ has type $3^\alpha 2^\beta 1^{m+\gamma}$.
\end{Lemma}
\begin{proof}
  For a hyperplane $H$ with $P\in H$ we have $P+x\in H$ if and only if $x\in H$. Hence $2m-\cC(H)=2\left(m-\cX(H\cap H_0)\right)\equiv 0\pmod 8$. For $P\notin H$ exactly one of
  $x$ and $P+x$ is contained in $H$, which gives $2m-\cC(H)=m\equiv 0\pmod 8$. For the type keep the planes and lines of $\cX$ and use the $m+\gamma$ points of
  $\left(P+\cX\right)\cup\left(\cX\setminus(\text{planes}\cup\text{lines})\right)$ as single points.
\end{proof}
For $m=32$, two disjoint planes and the field reduction of a hyperoval in $\PG(2,4)$ (type $2^6$) form a $4$-divisible set $\cX$ of type $3^22^6$. Lemma~\ref{lemma_cone_minus_apex} and the replacement of lines by
three points show that the types $3^22^\beta1^{50-3\beta}$ with $0\le\beta\le 6$ can be realized for $n=64$.

\begin{Lemma}
  \label{lemma_plane_gluing}
  Let $\cS_1$ and $\cS_2$ be $8$-divisible sets of points in $\PG(v-1,2)$ with $\cS_1\cap\cS_2=E$ for a plane $E$. Then the symmetric difference
  $\cS_1\triangle\cS_2$ is an $8$-divisible set of $\#\cS_1+\#\cS_2-14$ points. If $\cS_i$ has type $3^{a_i}2^{b_i}1^{c_i}$ such that $E$ is one of the planes of the
  partition for $i=1,2$, then $\cS_1\triangle\cS_2$ has type $3^{a_1+a_2-2}2^{b_1+b_2}1^{c_1+c_2}$. In particular, if type $3^a2^b1^c$ with $a\ge 1$ can be realized
  for cardinality $n$, then type $3^{a-1}2^b1^{c+8}$ can be realized for cardinality $n+1$.
\end{Lemma}
\begin{proof}
  We have $\chi_{\cS_1\triangle\cS_2}=\chi_{\cS_1}+\chi_{\cS_2}-2\chi_E$, and $2\chi_E$ is $8$-divisible since every hyperplane contains $3$ or $7$ points of $E$.
  The parts of the two partitions other than $E$ partition $(\cS_1\backslash E)\cup(\cS_2\backslash E)=\cS_1\triangle\cS_2$. Given two realizations, we can embed them
  into a common ambient space such that the two copies of $E$ coincide and $\langle\cS_1\rangle\cap\langle\cS_2\rangle=E$, which implies $\cS_1\cap\cS_2=E$.
  For the last statement choose $\cS_2$ as a solid containing $E$, partitioned into $E$ and eight points.
\end{proof}
Applied to a plane spread of $\PG(5,2)$ and $t$ solids, Lemma~\ref{lemma_plane_gluing} gives type $3^{9-t}1^{8t}$ for $n=63+t$, $1\le t\le 9$.
Lemma~\ref{lemma_plane_gluing} with $\cS_2=\PG(4,2)$ of type $3^12^8$ increases the cardinality by $31-14=17$. For example, type $3^12^{12}1^4$ for $n=47$ yields type $2^{20}1^4$ for $n=64$.

Further types were found by a computer search, and Table~\ref{table_types_60_70} in Appendix~\ref{app_types_60_70} summarizes the status of all types $3^{m_3}2^{m_2}1^{m_1}$ for $60\le n\le 70$.

\begin{Proposition}[E.g.~for part~(b) {\cite[Lemma 11.2]{kurz2021divisible}}]
  Let $\cS$ be a set of points of type $2^{m_2}1^{m_1}$ in $\PG(v-1,4)$, i.e., a $16$-divisible set of points that is partitioned into $m_2$ lines and $m_1$ points.
  \begin{itemize}
    \item[(a)] $m_1=0$ or $m_1\in\{5,10,15,16,17\}\cup\{20,21,\dots\}$.
    \item[(b)] If $m_1=5$, then the five points form a line. In particular, type $2^{m_2}1^5$ exists if and only if type $2^{m_2+1}$ exists.
    \item[(c)] All types $2^{m}1^{85-5m}$ with $0\le m\le 17$, the types $1^{21}$ and $2^11^{16}$, and type $2^m$ for every length $m$ of a projective $4$-divisible
               linear code over $\F_{16}$ can be realized, e.g., $m\in\{17,21,52,65\}$. If types $2^{a_2}1^{a_1}$ and $2^{b_2}1^{b_1}$ are realizable, then so is type $2^{a_2+b_2}1^{a_1+b_1}$,
               by taking two realizations whose spans intersect trivially.
  \end{itemize}
\end{Proposition}
\begin{proof}
  (a) Lines over $\F_4$ are $4$-divisible. So the $m_1$ points form a $4$-divisible set, and Lemma~\ref{lemma_4_div_gf4} applies.
  (b) Let $\cF$ be a $4$-divisible set of five points over $\F_4$. Every hyperplane contains $1$ or $5$ points of $\cF$. A hyperplane containing two points of $\cF$
  contains all of them. Thus $\cF$ is contained in the intersection of all hyperplanes through two of its points, which is a line.
  (c) $\PG(3,4)$ is $64$-divisible, and we can take $m$ lines of a line spread of $\PG(3,4)$. Moreover, $\PG(2,4)$ is $16$-divisible. For the field reduction of a
  multiset $\cM$ of $m$ points over $\F_{16}$ and an $\F_4$-hyperplane $H$ containing the $\F_{16}$-hyperplane $K$ we have $\cS(H)=m+4\cM(K)$. Hence $\cS$ is
  $16$-divisible if and only if $\cM$ is $4$-divisible. Examples over $\F_{16}$: a line ($17$ points), a Baer subplane $\PG(2,4)$ ($21$ points, line intersections $1$ and $5$),
  a Denniston maximal arc of degree $4$ in $\PG(2,16)$ ($52$ points, line intersections $0$ and $4$), and the Hermitian unital in $\PG(2,16)$ ($65$ points, line
  intersections $1$ and $5$).
\end{proof}

\appendix
\section{The digit rule for \texorpdfstring{$p\in\{2,3\}$}{p in {2,3}}}\label{app_p3}
In this appendix we prove Theorem~\ref{thm_p2_star} and Lemmas~\ref{thm_p3_c01} and~\ref{thm_p3_c02}, and we show that the digit rule fails for $p\ge 5$ (Subsection~\ref{app_p5}). Except in Subsection~\ref{app_p5}, let $p\in\{2,3\}$, $r\ge 0$, and $s_i=p^i[r+1-i]_p$ for $0\le i\le r$.
As in Subsection~\ref{sub_p3}, but now also for $r=0$, we call $(c_0,\dots,c_r)\in\{0,\dots,p-1\}^{r+1}$ a \emph{digit-good pattern} (of level $r$) if $c_i+c_{i+1}\le p-1$ for all $0\le i<r$,
and $n$ \emph{digit-good} if $n=\sum_{i=0}^r c_is_i$ for a digit-good pattern. For $p=2$ this means $c_ic_{i+1}=0$ for all $0\le i<r$. By $(\star_r)$ we denote the statement that every
$p^r$-divisible multiset of points over $\F_p$ with digit-good cardinality is canonical, and by $(\star_r,v)$ the same statement for multisets of points of $\PG(v-1,p)$.
For $n=\sum_{i=0}^r c_is_i$ we write $N:=\sum_{i=0}^r c_ip^i$ and $C:=\sum_{i=0}^r c_i$. Since $(p-1)s_i=p^{r+1}-p^i$, we have $(p-1)n=p^{r+1}C-N$. If $c_0,\dots,c_r\in\{0,\dots,p-1\}$, then $N\le p^{r+1}-1$.
An easy induction on $r$ gives $N\le p^{r+2}/(p+1)$ for every digit-good pattern. If $c_r=p-1$, then $c_{r-1}=0$ and $N\le(p-1)p^r+p^r/(p+1)=p^{r+2}/(p+1)$.
If $c_r\le p-2$, then $N\le(p-2)p^r+p^{r+1}/(p+1)<p^{r+2}/(p+1)$. For $p=3$ this is the bound $N\le\tfrac{9}{4}3^r$. Hence $(p-1)n/(p^{r+1}-N)\le p^{r+1}C/(p^{r+1}-N)\le(p+1)C$. Moreover $2C\le(p-1)(r+2)$,
since $c_{2j}+c_{2j+1}\le p-1$ for all $j$ (with $c_{r+1}=0$). For $p=3$ this gives $C\le r+2$.

\subsection{Tools for \texorpdfstring{$p\in\{2,3\}$}{p in {2,3}}}
\begin{Lemma}
  \label{lemma_digits_unique}
  Let $n=\sum_{i=0}^r c_is_i$ with $c_0,\dots,c_r\in\{0,\dots,p-1\}$. If $n=\sum_{i=0}^r e_is_i$ with non-negative integers $e_i$, then $e_i=c_i$ for all $i$. In particular, every decomposition of a
  canonical multiset of points of cardinality $n$ into canonical pieces contains exactly $c_i$ pieces $p^i\chi_S$ with $\dim(S)=r+1-i$ for every $0\le i\le r$.
\end{Lemma}
\begin{proof}
  Let $E:=\sum_ie_i$ and $M:=\sum_ie_ip^i$. We have $p^{r+1}E-M=(p-1)n=p^{r+1}C-N$, i.e., $M=N+p^{r+1}(E-C)$. As $M\ge 0$ and $0\le N<p^{r+1}$, we get $E\ge C$.
  Starting from $(e_0,\dots,e_r,0)$ we repeatedly replace $p$ units at a position $i\le r$ by one unit at position $i+1$. This keeps $\sum_ie_ip^i=M$, lowers the sum of the entries by $p-1$, and stops
  after some $t\ge 0$ steps at $(c_0,\dots,c_r,E-C)$, since $M\equiv N\pmod{p^{r+1}}$. Comparing the sums of the entries gives $E-(p-1)t=C+(E-C)$. So $t=0$ and $e_i=c_i$ for all $i$.
\end{proof}

\begin{Lemma}[E.g.~for part~(a) {\cite[Theorem 2]{kiermaier2020lengths}}]
  \label{lemma_p3_defect}
  Let $r\ge 1$ and $\cM$ be a $p^r$-divisible multiset of points of $\PG(v-1,p)$, $v\ge 2$, of cardinality $n=\sum_{i=0}^rc_is_i$ with $c_0,\dots,c_r\in\{0,\dots,p-1\}$. For a hyperplane $H$ let $d(H):=C-(n-\cM(H))/p^r\in\Z$.
  \begin{itemize}
    \item[(a)] $d(H)\ge 0$ for every hyperplane $H$.
    \item[(b)] $\sum_H d(H)=\sum_{i=0}^r c_i\left(p^{v-1-r+i}-1\right)/(p-1)$, where $H$ runs over all hyperplanes.
    \item[(c)] $\sum_{H\ni P}d(H)=C[v-1]_p-(n-\cM(P))p^{v-2-r}$ for every point $P$.
    \item[(d)] If $d(H)=0$, then $\cM|_H$ is a $p^{r-1}$-divisible multiset of points of $H$ with cardinality $\sum_{i=0}^{r-1}c_ip^i[r-i]_p$.
  \end{itemize}
\end{Lemma}
\begin{proof}
  Put $s_i':=p^i[r-i]_p$ for $0\le i\le r-1$. Then $s_i=p^r+s_i'$, $n=p^rC+\sum_{i<r}c_is_i'$, and $p^r=ps_{r-1}'$. Hence
  $$
    \cM(H)=n-p^r\left(C-d(H)\right)=\sum_{i=0}^{r-2}c_is_i'+\left(c_{r-1}+p\,d(H)\right)s_{r-1}'.
  $$
  Since $\cM|_H$ is $p^{r-1}$-divisible \cite[Lemma 3.4]{kurz2021divisible}, Theorem~\ref{thm_digits_only_if}(a) for $\Delta=p^{r-1}$ gives $c_{r-1}+p\,d(H)\ge 0$ if $r\ge 2$. For $r=1$ the displayed
  identity reads $\cM(H)=c_0+p\,d(H)$, which is non-negative. As $c_{r-1}\le p-1$ we obtain (a), and (d) is the displayed identity for $d(H)=0$.

  Every point lies outside of $p^{v-1}$ hyperplanes, and
  exactly $[v-1]_p-[v-2]_p=p^{v-2}$ of the hyperplanes through a point $P$ do not contain a given point $Q\neq P$. So $\sum_H(n-\cM(H))=p^{v-1}n$ and $\sum_{H\ni P}(n-\cM(H))=p^{v-2}(n-\cM(P))$.
  This gives (c), and (b) follows from $(p-1)\left([v]_pC-p^{v-1-r}n\right)=p^{v-1-r}N-C$.
\end{proof}

For $R\ge 0$ let $\operatorname{GRM}(R,v)$ denote the generalized Reed--Muller code of length $p^v$ over $\F_p$, i.e., the set of functions $\F_p^v\to\F_p$ given by polynomials of degree at most $R$ in which
every variable occurs with exponent at most $p-1$. Every function $\F_p^v\to\F_p$ has a unique such reduced polynomial. For $p=2$ these are the binary Reed--Muller codes. We use the following classical facts, see
\cite[Chapter 13]{macwilliams1977theory} for $p=2$ and \cite{delsarte1970generalized} and \cite[Chapter 5]{assmus1992designs} for $p=3$.
\begin{itemize}
  \item[(G1)] For $0\le k\le v$ the $\F_p$-span of the characteristic functions of the $k$-dimensional flats of $\F_p^v$ is $\operatorname{GRM}((p-1)(v-k),v)$, and the dual code of $\operatorname{GRM}(R,v)$ is $\operatorname{GRM}((p-1)v-1-R,v)$.
  \item[(G2)] For $0\le a\le v$ the minimum weight of $\operatorname{GRM}((p-1)a,v)$ is $p^{v-a}$, and the codewords of weight $p^{v-a}$ are the non-zero multiples of the characteristic functions of the
              $(v-a)$-dimensional flats.
  \item[(G3)] For $2\le a\le v-2$ every codeword of $\operatorname{GRM}((p-1)a,v)$ with weight larger than $p^{v-a}$ has weight at least $(p+1)p^{v-a-1}$.
              For $p=3$ this is the case $q=3$ of \cite[Theorem 3.8]{rolland2010second}, where the second weight of $\operatorname{GRM}(a(q-1),v)$ over $\F_q$ is shown to be
              $2(q-1)q^{v-a-1}$ for $q\ge 3$, $v\ge 3$, and $2\le a\le v-1$. For $p=2$ the second weight $3\cdot 2^{v-a-1}$ follows from the classification of the codewords
              of weight less than twice the minimum weight due to Kasami and Tokura \cite{kasami1970weight}, see also \cite[Introduction]{rolland2010second}.
\end{itemize}

\begin{Lemma}
  \label{lemma_p3_residue}
  Let $v\ge r+1$ and $\cM$ be a $p^r$-divisible multiset of points of $\PG(v-1,p)$ with cardinality $n$. Let $g:\F_p^v\to\F_p$ be given by $g(0)\equiv n\pmod p$ and $g(x)\equiv\cM(\langle x\rangle)\pmod p$ for $x\neq 0$, and
  let $\mathcal{R}(\cM)$ be the set of points $P$ with $\cM(P)\not\equiv 0\pmod p$. Then $g\in\operatorname{GRM}((p-1)(v-r-1),v)$. Consequently:
  \begin{itemize}
    \item[(a)] If $p\mid n$, then $\mathcal{R}(\cM)=\emptyset$ or $\#\mathcal{R}(\cM)\ge[r+1]_p+1$.
    \item[(b)] If $p\nmid n$, then $\#\mathcal{R}(\cM)\ge[r+1]_p$, with equality if and only if there is an $r$-space $F$ (i.e., $\dim(F)=r+1$) such that $\mathcal{R}(\cM)$ is the set of points of $F$ and $\cM(P)\equiv n\pmod p$ for all points $P\le F$.
  \end{itemize}
\end{Lemma}
\begin{proof}
  Let $A$ be a flat of dimension $v-r$.

  If $0\notin A$, then $r\ge 1$ and $A=a+K$ with a subspace $K$ of dimension $v-r$ and $a\notin K$. The restriction of $\cM$ to the subspace $W:=K+\langle a\rangle$ of codimension $r-1$ is
  $p$-divisible \cite[Lemma 3.4]{kurz2021divisible}, $K$ is a hyperplane of $W$, and every point of $W$ outside of $K$ has exactly one representative in $A$. Hence $\sum_{x\in A}g(x)\equiv\cM(W)-\cM(K)\equiv 0\pmod p$.

  If $0\in A$, then $A$ is a subspace of codimension $r$. Choose subspaces $\F_p^v=U_0>U_1>\dots>U_r=A$ with $\dim(U_k)=v-k$. As $\cM|_{U_k}$ is $p^{r-k}$-divisible we have $\cM(U_{k+1})\equiv\cM(U_k)\pmod p$ for
  $0\le k<r$. Hence $\cM(A)\equiv n\pmod p$. Since every point has $p-1$ representatives, $\sum_{x\in A}g(x)\equiv n+(p-1)\cM(A)\equiv pn\equiv 0\pmod p$.

  Therefore $g$ is orthogonal to all characteristic functions of
  $(v-r)$-dimensional flats, and (G1) gives $g\in\operatorname{GRM}((p-1)(v-r)-1,v)$. This proves the claim for $p=2$. For $p=3$ we have $g(-x)=g(x)$ for all $x$. As the reduced polynomial of $g$ is unique, only monomials of even degree occur in it,
  and thus $g\in\operatorname{GRM}(2(v-r-1),v)$.

  The weight of $g$ is $(p-1)\#\mathcal{R}(\cM)$ if $p\mid n$ and $(p-1)\#\mathcal{R}(\cM)+1$ otherwise, and by (G2) every non-zero codeword has weight at least $p^{r+1}$. In case (a) this gives
  $(p-1)\#\mathcal{R}(\cM)\ge p^{r+1}>(p-1)[r+1]_p$ if $\mathcal{R}(\cM)\neq\emptyset$. In case (b) it gives $\#\mathcal{R}(\cM)\ge (p^{r+1}-1)/(p-1)=[r+1]_p$. If equality holds, then (G2) gives $g=c\cdot\chi_A$ for an
  $(r+1)$-dimensional flat $A$ and $c\in\F_p\setminus\{0\}$. As $g(0)\neq 0$, $A$ is a subspace $F$ and $c\equiv n\pmod p$. The converse is obvious.
\end{proof}

\begin{Lemma}
  \label{lemma_p3_div3}
  Let $\cM$ be a $p^r$-divisible multiset of points over $\F_p$ with cardinality $n=\sum_{i=0}^rc_is_i$, where $c_0=0$. Assume that $(c_0,\dots,c_r)$ is digit-good. Then all multiplicities of $\cM$ are divisible by $p$.
\end{Lemma}
\begin{proof}
  Induction on $r$. If $C=0$, then $\cM=0$. So let $C\ge 1$, and embed $\cM$ into $\PG(v-1,p)$ with $v\ge r+2$. Let $\mathcal{Z}$ be the set of hyperplanes $H$ with $d(H)=0$, see Lemma~\ref{lemma_p3_defect}.
  For $H\in\mathcal{Z}$ the multiset $\cM|_H$ is $p^{r-1}$-divisible with cardinality $\sum_{i=1}^{r-1}c_ip^i[r-i]_p$ by Lemma~\ref{lemma_p3_defect}(d). Its pattern $(0,c_1,\dots,c_{r-1})$ satisfies the assumption of the lemma
  for $r-1$, and by induction all multiplicities of $\cM|_H$ are divisible by $p$. So no point of $\mathcal{R}(\cM)$ lies in a hyperplane of $\mathcal{Z}$.

  Let $P\in\mathcal{R}(\cM)$. Then $d(H)\ge 1$ for all $[v-1]_p$
  hyperplanes $H$ through $P$. Lemma~\ref{lemma_p3_defect}(c) gives $(n-\cM(P))p^{v-2-r}\le (C-1)[v-1]_p$, and hence $\cM(P)\ge n-(C-1)p^{r+1}/(p-1)=(p^{r+1}-N)/(p-1)$, where $N<p^{r+1}$. Summing over $\mathcal{R}(\cM)$ gives
  $$
    \#\mathcal{R}(\cM)\le\frac{(p-1)n}{p^{r+1}-N}=\frac{p^{r+1}C-N}{p^{r+1}-N}.
  $$
  As $c_0=0$, we have $2C\le(p-1)(r+1)$. With $N\le p^{r+2}/(p+1)$ the right hand side is at most $(p+1)C\le(p^2-1)(r+1)/2$, which is smaller than $[r+1]_p+1$ if $p=2$
  or $r\ge 2$. For $p=3$ and $r=1$ the patterns $(0,1)$ and $(0,2)$ give the values $1$ and $4$.
  So $\#\mathcal{R}(\cM)<[r+1]_p+1$. Since $p\mid n$, Lemma~\ref{lemma_p3_residue}(a) gives $\mathcal{R}(\cM)=\emptyset$.
\end{proof}

\begin{Proposition}
  \label{prop_p3_reduction}
  Let $r\ge 1$, $v\ge 1$, and assume $(\star_{r-1},v)$. Let $\cM$ be a $p^r$-divisible multiset of points of $\PG(v-1,p)$ with digit-good cardinality $n=\sum_{i=0}^rc_is_i$.
  \begin{itemize}
    \item[(a)] If $c_0=0$, then $\cM=p\cM'$ for a $p^{r-1}$-divisible multiset $\cM'$ with digit-good cardinality, and $\cM$ is canonical.
    \item[(b)] If $c_0\ge 1$, then $\cM$ is canonical if and only if there are $r$-spaces $F_1,\dots,F_{c_0}$ such that $\cM(P)\equiv\#\{j\,:\,P\le F_j\}\pmod p$ for all points $P$. For $c_0=1$ this is the case if and only if
               $\#\mathcal{R}(\cM)\le[r+1]_p$.
  \end{itemize}
  In particular, $(\star_{r-1})$ implies that all $p^r$-divisible multisets of points over $\F_p$ with digit-good cardinality $n\equiv 0\pmod p$ are canonical.
\end{Proposition}
\begin{proof}
  (a) By Lemma~\ref{lemma_p3_div3} we have $\cM=p\cM'$. The multiset $\cM'$ is $p^{r-1}$-divisible and has cardinality $n/p=\sum_{i=0}^{r-1}c_{i+1}p^i[r-i]_p$, i.e., pattern $(c_1,\dots,c_r)$, which is digit-good.
  By $(\star_{r-1},v)$ it is a sum of pieces $p^i\chi_S$ with $\dim(S)=r-i$. Hence $\cM$ is a sum of pieces $p^{i+1}\chi_S$ with $\dim(S)=r+1-(i+1)$.

  (b) If $\cM$ is canonical, then Lemma~\ref{lemma_digits_unique} shows that a decomposition contains exactly $c_0$ pieces $\chi_{F_j}$ with $\dim(F_j)=r+1$, while all other pieces have multiplicities divisible by $p$.
  Conversely, let $\cM'':=\cM-\sum_{j=1}^{c_0}\chi_{F_j}$. If a point $P$ lies in $t\ge 1$ of the $F_j$, then $t\le c_0\le p-1$ and $\cM(P)\equiv t\pmod p$, which implies $\cM(P)\ge t$. Hence $\cM''$ is a $p^r$-divisible multiset of points
  with all multiplicities divisible by $p$ and cardinality $n-c_0s_0$ with digit-good pattern $(0,c_1,\dots,c_r)$. As in (a), $\cM''/p$ is canonical by $(\star_{r-1},v)$, and so is $\cM$. For $c_0=1$ embed $\cM$ into a space of dimension at least $r+1$, which does not change canonicity. Then the
  last statement follows from Lemma~\ref{lemma_p3_residue}(b).
\end{proof}

\begin{Lemma}
  \label{lemma_p3_count}
  Let $v\ge r+2$ and let $\cM$ be a $p^r$-divisible multiset of points of $\PG(v-1,p)$ with cardinality $n=\sum_{i=0}^r c_is_i$, where $c_0,\dots,c_r\in\{0,\dots,p-1\}$ and $c_0\neq 0$.
  Let $d(H)$ be as in Lemma~\ref{lemma_p3_defect}, $\mathcal{R}:=\mathcal{R}(\cM)$, and $X:=p^{v-2-r}$. Let $\mathcal{B}$ be a set of hyperplanes $H$ with $d(H)=0$ such that
  $\#(\mathcal{R}\cap H)=[r]_p$ for every hyperplane $H\notin\mathcal{B}$ with $d(H)=0$. Then $\#(\mathcal{R}\cap H)\ge[r]_p$ for every hyperplane $H$ and
  $$
    \left(\#\mathcal{R}-p[r]_p\right)\cdot\left(p^{r+1}-N\right)\le (p-1)\cM(\mathcal{R})+\frac{(p-1)\#\mathcal{B}}{X}\left(\#\mathcal{R}-[r]_p\right),
  $$
  where $\cM(\mathcal{R}):=\sum_{P\in\mathcal{R}}\cM(P)\le n$.
\end{Lemma}
\begin{proof}
  Let $g$ be the function of Lemma~\ref{lemma_p3_residue}. Since $n\equiv c_0\not\equiv 0\pmod p$, we have $g(0)\neq 0$. For a hyperplane $H$ the restriction $g|_H$ is given by a polynomial of degree at most
  $(p-1)(v-r-1)$ in $v-1$ variables. So (G2) gives $1+(p-1)\#(\mathcal{R}\cap H)\ge p^{(v-1)-(v-r-1)}=p^r$, i.e., $\#(\mathcal{R}\cap H)\ge[r]_p$. In particular $\#\mathcal{R}\ge[r]_p$.

  For a hyperplane $H$ let $e(H):=\#(\mathcal{R}\cap H)-[r]_p\ge 0$. If $d(H)\ge 1$, then $(1-d(H))e(H)\le 0$. Otherwise $d(H)=0$ by Lemma~\ref{lemma_p3_defect}(a). Then $e(H)=0$ for
  $H\notin\mathcal{B}$ by assumption, and $e(H)\le\#\mathcal{R}-[r]_p$ for $H\in\mathcal{B}$. Summing $(1-d(H))e(H)$ over all $[v]_p$ hyperplanes and using
  $\sum_H(1-d(H))\#(\mathcal{R}\cap H)=\sum_{P\in\mathcal{R}}\sum_{H\ni P}(1-d(H))$ we obtain
  $$
    \sum_{P\in\mathcal{R}}\left([v-1]_p-\sum_{H\ni P}d(H)\right)\le[r]_p\left([v]_p-\sum_H d(H)\right)+\#\mathcal{B}\left(\#\mathcal{R}-[r]_p\right).
  $$
  By Lemma~\ref{lemma_p3_defect}(b),(c) and $(p-1)n=p^{r+1}C-N$ we have
  \begin{align*}
    [v-1]_p-\sum_{H\ni P}d(H)&=\tfrac{1}{p-1}\left(X\left(p^{r+1}-N\right)+C-1\right)-X\cM(P)\quad\text{and}\\
    [v]_p-\sum_Hd(H)&=\tfrac{1}{p-1}\left(pX\left(p^{r+1}-N\right)+C-1\right).
  \end{align*}
  Inserting this and multiplying by $p-1$ we get
  $$
    \left(\#\mathcal{R}-p[r]_p\right)X\left(p^{r+1}-N\right)\le (p-1)X\cM(\mathcal{R})+(p-1)\#\mathcal{B}\left(\#\mathcal{R}-[r]_p\right)-(C-1)\left(\#\mathcal{R}-[r]_p\right).
  $$
  The claim follows from $C\ge 1$ and $\#\mathcal{R}\ge[r]_p$.
\end{proof}

\subsection{The case \texorpdfstring{$c_0=1$}{c\textunderscore{}0=1} and the case \texorpdfstring{$p=2$}{p=2}}\label{app_p2}
\begin{proof}[Proof of Lemma~\ref{thm_p3_c01}]
  Let $K_p:=\left((p+1)p^r+p-3\right)/(p-1)$, i.e., $K_2=3\cdot 2^r-1$ and $K_3=2\cdot 3^r$. We first show
  \begin{equation}
    \label{eq_p3_dagger}
    p[r]_p+\frac{(p-1)n}{p^{r+1}-N}<K_p
  \end{equation}
  for every digit-good pattern with $c_0\ge 1$, where $r\ge 1$ if $p=2$ and $r\ge 3$ if $p=3$. We have $K_p-p[r]_p=(p^r+2p-3)/(p-1)$ and
  $(p-1)n/(p^{r+1}-N)\le(p+1)C\le(p^2-1)(r+2)/2$. This gives (\ref{eq_p3_dagger}) for $r\ge 3$ if $p=2$ and for $r\ge 4$ if $p=3$. For $p=2$ and $r\in\{1,2\}$
  the patterns $(1,0)$, $(1,0,0)$, and $(1,0,1)$ give $(p-1)n/(p^{r+1}-N)\in\{1,1,11/3\}$. For $p=3$ and $r=3$ the maximum of $2n/(81-N)$ over the $17$ digit-good
  patterns with $c_0\ge 1$ is $266/23<15$, attained for $n=133$.

  Embedding $\cM$ into a larger space does not change canonicity. We may therefore assume that $\cM$ is a multiset of points of $\PG(v-1,p)$ with $v\ge r+3$.
  Let $H$ be a hyperplane with $d(H)=0$. By Lemma~\ref{lemma_p3_defect}(d) and $(\star_{r-1})$ the restriction $\cM|_H$ is canonical with pattern $(1,c_1,\dots,c_{r-1})$.
  By Lemma~\ref{lemma_digits_unique} it is the sum of $\chi_T$ for an $(r-1)$-space $T$ and of pieces whose multiplicities are divisible by $p$. Hence
  $\mathcal{R}(\cM)\cap H$ is the set of points of $T$ and $\#(\mathcal{R}(\cM)\cap H)=[r]_p$.

  Lemma~\ref{lemma_p3_count} with $\mathcal{B}=\emptyset$ and (\ref{eq_p3_dagger}) give $\#\mathcal{R}(\cM)<K_p$. As $K_p$ is an integer, the function
  $g\in\operatorname{GRM}((p-1)(v-r-1),v)$ of Lemma~\ref{lemma_p3_residue} has weight $(p-1)\#\mathcal{R}(\cM)+1\le(p-1)(K_p-1)+1=(p+1)p^r-1$.
  As $2\le v-r-1\le v-2$, (G3) shows that $g$ has the minimum weight $p^{r+1}$, i.e., $\#\mathcal{R}(\cM)=[r+1]_p$, and Proposition~\ref{prop_p3_reduction}(b)
  shows that $\cM$ is canonical.
\end{proof}

\begin{Theorem}
  \label{thm_p2_star}
  For $p=2$ the statement $(\star_r)$ holds for all $r\ge 0$.
\end{Theorem}
\begin{proof}
  The statement $(\star_0)$ holds trivially. Let $r\ge 1$ and assume $(\star_{r-1})$. For $p=2$ every digit-good pattern has $c_0\in\{0,1\}$. Proposition~\ref{prop_p3_reduction}(a)
  covers $c_0=0$, and Lemma~\ref{thm_p3_c01} covers $c_0=1$. Hence $(\star_r)$ holds, and induction on $r$ gives the claim.
\end{proof}

\subsection{The case \texorpdfstring{$p=3$}{p=3} with \texorpdfstring{$c_0=2$}{c\textunderscore{}0=2}}\label{app_p3_case}
In this subsection let $p=3$.

\begin{Lemma}
  \label{lemma_p3_172}
  Let $w\ge 6$ and let $\cM$ be a $27$-divisible multiset of points of $\PG(w-1,3)$ with cardinality $172=s_0+2s_1+2s_3$, where $s_i=3^i[4-i]_3$, i.e., with digits $(1,2,0,2)$.
  Let $d(H)$ be as in Lemma~\ref{lemma_p3_defect} (with $r=3$), let $\mathcal{R}:=\mathcal{R}(\cM)$, and let $\mathcal{B}$ be a set of at most $[w-5]_3$ hyperplanes such that $\#(\mathcal{R}\cap H)=13$
  for every hyperplane $H\notin\mathcal{B}$ with $d(H)=0$. Then $\#\mathcal{R}\notin\{54,55,56,57,58\}$.
\end{Lemma}
\begin{proof}
  \emph{Overview.} By the weights of generalized Reed--Muller codes, the intersection numbers $h_H=\#(\mathcal{R}\cap H)$ lie in a small set $A_x$. The first three binomial moments of $h_H$ over all hyperplanes are determined by counting (the third only as a lower bound), and Lemma~\ref{lemma_p3_defect}(c) bounds $\sum_Hh_Hd(H)$ from above.
  A polynomial $\psi$ of degree at most $3$ with $\psi(13)=0$ and $\psi(h_H)\le h_Hd(H)+\nu\cdot\mathbf 1_{\mathcal{B}}(H)$ for every hyperplane then gives a linear inequality in $X=3^{w-5}$, which is violated for $x\in\{54,\dots,58\}$.

  Assume $x:=\#\mathcal{R}\in\{54,\dots,58\}$, let $X:=3^{w-5}$, and for a hyperplane $H$ let $h_H:=\#(\mathcal{R}\cap H)$. By Lemma~\ref{lemma_p3_count} (with $r=3$) we have $h_H\ge 13$.
  Let $g\in\operatorname{GRM}(2(w-4),w)$ be the function of Lemma~\ref{lemma_p3_residue} and $a\notin H$. In affine coordinates on the coset $a+H$ the restriction of $g$ is given by a polynomial of degree at most $2(w-4)$
  in $w-1$ variables, and every point outside of $H$ has exactly one representative in $a+H$. So $x-h_H$ is the weight of a codeword of $\operatorname{GRM}(2(w-4),w-1)$, and (G2), (G3) with $2\le w-4\le w-3$ give
  $x-h_H\in\{0,27\}$ or $x-h_H\ge 36$. Hence $h_H\in A_x:=\{13,\dots,x-36\}\cup\{x-27,x\}$ for every hyperplane $H$.

  By Lemma~\ref{lemma_hyperplane_moments} with $(v,q)=(w,3)$ we have
  $$
    \sum_H 1=[w]_3,\quad \sum_H h_H=x[w-1]_3,\quad \sum_H\binom{h_H}{2}=\binom{x}{2}[w-2]_3,\quad \sum_H\binom{h_H}{3}\ge\binom{x}{3}[w-3]_3.
  $$
  Lemma~\ref{lemma_p3_defect}(c) with $C=5$ and $n=172$ gives
  $$
    \sum_H h_Hd(H)=\sum_{P\in\mathcal{R}}\sum_{H\ni P}d(H)=5x[w-1]_3-X\left(172x-\cM(\mathcal{R})\right)\le 5x[w-1]_3-172X(x-1).
  $$
  Let $\varphi$ and $\nu$ be as in the following table, and put $\psi(h):=h-13\varphi(h)/\varphi(13)$.
  $$
    \begin{array}{rlrrr}
      x & \varphi(h) & \nu & \alpha & \beta\\
      \hline
      54 & (h-18)(h-27) & 27 & 5861/70 & 14823/70\\
      55 & (h-19)(h-28) & 28 & 182/5 & 971/5\\
      56 & (h-20)(h-29)(70-h) & 29 & 423/76 & 10667/76\\
      57 & (h-20)(h-21)(63-h) & 21 & 1543/350 & 100137/700\\
      58 & (h-21)(h-22)(66-h) & 22 & 2165/53 & 7693/53
    \end{array}
  $$
  In all cases $\varphi(13)>0$ and $\psi(13)=0$. No element of $A_x$ lies strictly between the two smallest roots of $\varphi$, and in the cubic cases all elements of $A_x$ are smaller than the third root $70$, $63$, or $66$, respectively. So $\varphi(h)\ge 0$ and $\psi(h)\le h$ for all $h\in A_x$.

  A direct computation shows $\nu=\max\{\psi(h)\,:\,h\in A_x\setminus\{13\}\}$. We claim $\psi(h_H)\le h_Hd(H)+\nu\cdot\mathbf 1_{\mathcal{B}}(H)$ for every hyperplane $H$, where $\mathbf 1_{\mathcal{B}}(H)$ is $1$ for $H\in\mathcal{B}$ and $0$ otherwise. If $d(H)\ge 1$, then $\psi(h_H)\le h_H\le h_Hd(H)$.
  If $d(H)=0$ and $H\notin\mathcal{B}$, then $h_H=13$ and $\psi(h_H)=0$. If $d(H)=0$ and $H\in\mathcal{B}$, then $\psi(h_H)\le\max\{\psi(13),\nu\}=\nu$.

  Write $\psi(h)=e_0+e_1h+e_2\binom{h}{2}+e_3\binom{h}{3}$, where $e_3>0$ if $\varphi$ is cubic and $e_3=0$ otherwise. Summing over all hyperplanes and using $\#\mathcal{B}\le[w-5]_3$ we obtain
  $$
    e_0[w]_3+e_1x[w-1]_3+e_2\binom{x}{2}[w-2]_3+e_3\binom{x}{3}[w-3]_3\le 5x[w-1]_3-172X(x-1)+\nu[w-5]_3.
  $$
  Since $[w-k]_3=(3^{5-k}X-1)/2$, the left hand side minus the right hand side equals $\alpha X+\beta$ with $\alpha,\beta$ as in the table, computed in exact rational arithmetic. As $\alpha,\beta>0$,
  this is a contradiction.
\end{proof}

\begin{proof}[Proof of Lemma~\ref{thm_p3_c02}]
  \emph{Overview.} After the numerical preparation by the inequalities (\ref{eq_p3_C1}) and (\ref{eq_p3_C2}) we use induction on the dimension of the span of the support. If no hyperplane of defect $0$ contains a point of multiplicity $1$ (Case~1), then as for Lemma~\ref{thm_p3_c01} the residue function is $2\chi_F$ for an $r$-space $F$, and $\cM-\chi_F$ has the digit pattern $(1,0,c_2,\dots,c_r)$, so that Lemma~\ref{thm_p3_c01} applies. Otherwise (Case~2) we project from such a point $Y$. Steps~1 and~2 show that, modulo $3$ and after dividing by $3$, the projection is described by an $r$-space $F'$ and an $(r-1)$-space $G$. Step~3 lifts $G$ to an $r$-space $\hat G$ and shows that $\cM-\chi_{\hat G}$ again has the pattern $(1,0,c_2,\dots,c_r)$. The single pattern $(2,0,2,0,2)$ with $r=4$, where the estimates do not suffice, is settled by Lemma~\ref{lemma_p3_172}.

  Note that $c_1=0$. We put $N'':=(N+1)/3$ and $n'':=(n-1-s_0)/3$. We need inequality (\ref{eq_p3_C1}) for all digit-good patterns with $c_0=2$, and (\ref{eq_p3_C2}) if $c_2=2$:
  \begin{eqnarray}
    3[r]_3+\frac{2n-2}{3^{r+1}-N-2}&<&2\cdot 3^r,\label{eq_p3_C1}\\
    \left(2\cdot 3^{r-1}-3[r-1]_3\right)\left(3^r-N''\right)-2\cdot 3^{r-1}+[r-1]_3&\ge&2n''.\label{eq_p3_C2}
  \end{eqnarray}
  For $r\ge 4$ inequality (\ref{eq_p3_C1}) follows as (\ref{eq_p3_dagger}) from $(3^{r+1}C-N-2)/(3^{r+1}-N-2)\le 3^{r+1}C/(\tfrac{3}{4}3^r-2)<\tfrac{17}{4}(r+2)<(3^r+3)/2$. For (\ref{eq_p3_C2}) we use
  $N''=1+\sum_{i\ge 2}c_i3^{i-1}\le 1+\tfrac{3}{4}3^r$, since $(c_2,\dots,c_r)$ is a digit-good pattern of level $r-2$, and $2n''=3^r(C-1)-N''<3^r(r+1)$. So the left hand side of
  (\ref{eq_p3_C2}) is at least $\tfrac{1}{2}(3^{r-1}+3)(\tfrac{1}{4}3^r-1)-2\cdot 3^{r-1}>3^r(r+1)$ for $r\ge 5$. For $r\in\{3,4\}$ the two inequalities were checked directly for all digit-good patterns with $c_0=2$, resp.\ $c_0=c_2=2$, using
  exact rational arithmetic. The only failure is (\ref{eq_p3_C2}) for $r=4$ and $(c_0,\dots,c_4)=(2,0,2,0,2)$, and this case is treated separately at the end of Step 2.

  We use induction on $\dim(V)$, where $V$ is the span of the support of $\cM$, i.e., we assume the statement for all multisets whose support spans a subspace of smaller dimension.
  We consider $\cM$ as a multiset of points of $\PG(v-1,3)\supseteq V$ with $v\ge r+4$. All hyperplanes below are hyperplanes of the ambient space or of its quotients, and $\mathcal{Z}$ denotes the set
  of hyperplanes $H$ with $d(H)=0$. For $H\in\mathcal{Z}$ Lemma~\ref{lemma_p3_defect}(d), $(\star_{r-1})$, and Lemma~\ref{lemma_digits_unique} give
  \begin{equation}
    \label{eq_p3_restr}
    \cM|_H=\chi_{T_1^H}+\chi_{T_2^H}+\mathcal{K}_H,
  \end{equation}
  where $T_1^H$, $T_2^H$ are $(r-1)$-spaces and $\mathcal{K}_H$ is a sum of at most two pieces $3^i\chi_S$ with $\dim(S)=r-i$ for each $2\le i\le r-1$. The support of $\mathcal{K}_H$ has at most
  $2\sum_{j=1}^{r-2}[j]_3<[r-1]_3$ points, while $T_1^H\neq T_2^H$ implies $\#(T_1^H\setminus T_2^H)\ge 3^{r-1}>[r-1]_3$. So if $T_1^H\neq T_2^H$, then $H$ contains a point of multiplicity $1$.

  \emph{Case 1: no hyperplane of $\mathcal{Z}$ contains a point of multiplicity $1$.} Then $T_1^H=T_2^H$ and $\mathcal{R}(\cM)\cap H=T_1^H$ for all $H\in\mathcal{Z}$. As in the proof of Lemma~\ref{thm_p3_c01},
  Lemma~\ref{lemma_p3_count} with $\mathcal{B}=\emptyset$, (\ref{eq_p3_dagger}), and (G3) show that $g$ has weight $3^{r+1}$. By (G2) and $g(0)=2$ we have $g=2\chi_F$ for an $r$-space $F$. Thus $\cM(P)\ge 2$
  for all points $P\le F$, and $\cM-\chi_F$ is a $3^r$-divisible multiset of points with digit-good pattern $(1,0,c_2,\dots,c_r)$. It is canonical by Lemma~\ref{thm_p3_c01}, and so is $\cM$.

  \emph{Case 2: some $H_0\in\mathcal{Z}$ contains a point $Y$ with $\cM(Y)=1$.} Let $\cM'$ be the projection of $\cM$ from $Y$ (compare, e.g., \cite[Lemma 5]{korner2024lengths}), i.e., the multiset of points of the quotient $\F_3^v/Y$ given by $\cM'(L/Y)=\cM(L)-\cM(Y)$
  for the lines $L\ge Y$. The hyperplanes of $\F_3^v/Y$ are the $H/Y$ with $Y\le H$, and $\#\cM'-\cM'(H/Y)=n-\cM(H)$. So $\cM'$ is $3^r$-divisible with cardinality
  $n-1=s_0+s_1+\sum_{i\ge 2}c_is_i$ (note $s_0-1=s_1$), that is, with digits $c'=(1,1,c_2,\dots,c_r)\in\{0,1,2\}^{r+1}$, with the same digit sum $C$, with $N+2$ in place of $N$, and defect $d(H)$ at $H/Y$.

  For a subspace $S$ the projection of $\chi_S$ from $Y$ is
  $3\chi_{S/Y}$ if $Y\le S$ and $\chi_{(S+Y)/Y}$ otherwise. Let $H\in\mathcal{Z}$ with $Y\le H$. Since $\cM(Y)=1$, (\ref{eq_p3_restr}) shows that $Y$ lies in exactly one of $T_1^H$, $T_2^H$, say in $T_1^H$, and not in the support
  of $\mathcal{K}_H$. Hence
  \begin{equation}
    \label{eq_p3_proj}
    \cM'|_{H/Y}=3\chi_{T_1^H/Y}+\chi_{(T_2^H+Y)/Y}+\overline{\mathcal{K}}_H,
  \end{equation}
  where $\overline{\mathcal{K}}_H$, the projection of $\mathcal{K}_H$, is a sum of pieces $3^i\chi_{(S+Y)/Y}$ with $i\ge 2$.

  \emph{Step 1: There is an $r$-space $F'\le\F_3^v/Y$ with $\cM'(P)\equiv\chi_{F'}(P)\pmod 3$ for all points $P$.} By (\ref{eq_p3_proj}) the set $\mathcal{R}(\cM')\cap(H/Y)$ consists of the points of the $(r-1)$-space $(T_2^H+Y)/Y$ for every hyperplane $H/Y$ with
  $d(H)=0$. Lemma~\ref{lemma_p3_count} (applied to $\cM'$ with $\mathcal{B}=\emptyset$) and (\ref{eq_p3_C1}) give $\#\mathcal{R}(\cM')<2\cdot 3^r$, and (G3), (G2) show that the residue function of $\cM'$ is $\chi_{F'}$
  for an $r$-space $F'$.

  \emph{Step 2: There is an $(r-1)$-space $G\le\F_3^v/Y$ such that $\cM'':=\tfrac{1}{3}(\cM'-\chi_{F'})$ satisfies $\cM''(P)\equiv\chi_G(P)\pmod 3$ for all points $P$.} Since $\cM'(P)\ge 1$ for all points $P\le F'$, Step 1 shows that $\cM''$ is a multiset of points. It is
  $3^{r-1}$-divisible with cardinality $(n-1-s_0)/3=\sum_{j=0}^{r-1}c''_j3^j[r-j]_3$, where $c''=(1,c_2,\dots,c_r)$. For $H\in\mathcal{Z}$ with $Y\le H$ a comparison of (\ref{eq_p3_proj}) with Step 1 gives
  $F'\cap(H/Y)=(T_2^H+Y)/Y$, so that $F'\not\le H/Y$ and
  \begin{equation}
    \label{eq_p3_proj2}
    \cM''|_{H/Y}=\chi_{T_1^H/Y}+\tfrac{1}{3}\overline{\mathcal{K}}_H.
  \end{equation}
  If $c_2\le 1$, then $c''$ is a digit-good pattern of level $r-1$, and the claim follows from $(\star_{r-1})$ and Lemma~\ref{lemma_digits_unique}. Let $c_2=2$ and let $d''$ be the defect function of $\cM''$
  (Lemma~\ref{lemma_p3_defect} with $r-1$ instead of $r$). From $\#\cM''-\cM''(H/Y)=\left(n-\cM(H)-s_0+\#(F'\cap H/Y)\right)/3$ we get $d''(H/Y)=d(H)$ if $F'\not\le H/Y$ and $d''(H/Y)=d(H)-1$ if $F'\le H/Y$.
  So a hyperplane $H/Y$ with $d''(H/Y)=0$ either satisfies $H\in\mathcal{Z}$, and then $\mathcal{R}(\cM'')\cap(H/Y)$ is the $(r-2)$-space $T_1^H/Y$ by (\ref{eq_p3_proj2}), or it contains $F'$. Let $\mathcal{B}$ be the set of the latter
  hyperplanes. Then $\#\mathcal{B}\le[v-r-2]_3<\tfrac{1}{2}3^{v-r-2}$, and Lemma~\ref{lemma_p3_count} (with $r-1$ instead of $r$ in a space of dimension $v-1$, so that $X=3^{v-r-2}$) gives
  $(x-3[r-1]_3)(3^r-N'')<2n''+x-[r-1]_3$ for $x=\#\mathcal{R}(\cM'')$ if $x>[r-1]_3$. By (\ref{eq_p3_C2}), and since $3^r-N''\ge 1$ makes the difference of both sides non-decreasing in $x$, this excludes $x\ge 2\cdot 3^{r-1}$, except for $r=4$ and $(c_0,\dots,c_4)=(2,0,2,0,2)$.

  In this exceptional case we have $n''=172$ and $N''=61$, and the inequality reads $20(x-39)<331+x$, i.e., $x\le 58$. Moreover, $\cM''$ has digits $(1,2,0,2)$, the ambient space $\F_3^v/Y$ has dimension $w:=v-1\ge 7$, and $\#\mathcal{B}\le[w-5]_3$. Therefore Lemma~\ref{lemma_p3_172} excludes $54\le x\le 58$.

  So in all cases $x<2\cdot 3^{r-1}$, and (G3), (G2) give the claim.

  \emph{Step 3: Conclusion.} Since $\cM$ is supported in $V$, the supports of $\cM'$ and $\cM''$, and hence $G$, are contained in $V/Y$. Let $\hat G\le V$ be the $r$-space with $Y\le\hat G$ and $\hat G/Y=G$.
  For every point $P$ of $G$ we have $\cM'(P)=\chi_{F'}(P)+3\cM''(P)\equiv\chi_{F'}(P)+3\pmod 9$. Assume that $\cM(Z)=0$ for a point $Z\le\hat G$, let $L:=Y+Z$, and let $\cM^Z$ be the projection of $\cM$ from $Z$.
  It is a $3^r$-divisible multiset of points with cardinality $n$ whose support spans $V/Z$. By induction $\cM^Z$ is canonical. Since $c_1=0$, its pieces have multiplicity $1$ or a multiplicity divisible by $9$,
  and at most two of them have multiplicity $1$. So all multiplicities of $\cM^Z$ are congruent to $0$, $1$, or $2$ modulo $9$. However,
  $\cM^Z(L/Z)=\cM(L)-\cM(Z)=\cM(Y)+\cM'(L/Y)\equiv 1+\chi_{F'}(L/Y)+3\pmod 9$ is congruent to $4$ or $5$. Hence $\cM(P)\ge 1$ for all points $P\le\hat G$. So $\cM-\chi_{\hat G}$ is a $3^r$-divisible multiset of points
  with digit-good pattern $(1,0,c_2,\dots,c_r)$, which is canonical by Lemma~\ref{thm_p3_c01}. Thus $\cM$ is canonical.
\end{proof}

\begin{Remark}
  The arguments of Lemma~\ref{thm_p3_c01} also work for $r\in\{1,2\}$, except for the pattern $(1,0,2)$. The construction of Proposition~\ref{prop_triangle} has no analogue that is relevant for $p=3$: with $\gamma=(p+1)/2=2$ it gives the cardinalities $2s_i+s_{i+1}$, which violate the digit condition.
\end{Remark}

\subsection{Failure of the digit rule for \texorpdfstring{$p\ge 5$}{p>= 5}}\label{app_p5}

For a point $b$ of $E=\PG(2,p)$ let $b^\perp$ denote the line $\{x\,:\,b\cdot x=0\}$. Then $P\le b^\perp$ if and only if $b\le P^\perp$. A set of points is a blocking set if it meets every line.

\begin{Proposition}
  \label{prop_triangle}
  Let $p\ge 5$ be a prime, $\gamma:=(p+1)/2$, $r\ge 1$, and $0\le i\le r-1$. Then there exists a $p^r$-divisible multiset of points over $\F_p$ of cardinality $\gamma s_i+s_{i+1}$
  that is not canonical. Consequently, if $n=\sum_{k=0}^r c_ks_k$ with $c_0,\dots,c_{r-1}\in\{0,\dots,p-1\}$ and $c_r\ge 0$ such that $c_i\ge\gamma$ and $c_{i+1}\ge 1$ for some $0\le i\le r-1$,
  then there exists a $p^r$-divisible multiset of points over $\F_p$ of cardinality $n$ that is not canonical. Since $\gamma+1\le p-1$, the condition of Theorem~\ref{thm_digits_only_if} is
  not sufficient for $p\ge 5$. E.g.\ for $p=5$ and $r=1$ there is a non-canonical $5$-divisible multiset of $23=3s_0+s_1$ points in $\PG(2,5)$.
\end{Proposition}
\begin{proof}
  Let $Q\subseteq\F_p$ be the set of squares including $0$, thus $\#Q=\gamma$, and let
  $$
    \cT:=\left\{(0,1,-\sigma),\,(-\sigma,0,1),\,(1,-\sigma,0)\,:\,\sigma\in Q\right\}
  $$
  be the projective triangle in $E$. The three families are pairwise disjoint (compare the positions of the zero coordinates). Hence $\#\cT=3\gamma$. We show that $\cT$ is a blocking set.
  Let $L$ be the line $u_1x_1+u_2x_2+u_3x_3=0$. If $u_1=0$, $u_2=0$, or $u_3=0$, then $L$ contains $(1,0,0)$, $(0,1,0)$, or $(0,0,1)$, respectively, which are points of $\cT$.
  If $u_1u_2u_3\neq 0$, then $L$ contains $(0,1,-\sigma)$, $(-\sigma,0,1)$, and $(1,-\sigma,0)$ for $\sigma=u_2/u_3$, $\sigma=u_3/u_1$, and $\sigma=u_1/u_2$, respectively. The product of these three
  values is $1$. Since a product of three non-squares is a non-square, one of them is a square and $L$ meets $\cT$. The line $x_1=0$ contains exactly the $\gamma+1$ points
  $(0,1,-\sigma)$, $\sigma\in Q$, and $(0,0,1)$ of $\cT$, and similarly for $x_2=0$ and $x_3=0$. Every other line meets each of these three lines, which cover $\cT$, in one point.
  So every line contains at most $\max\{\gamma+1,3\}=\gamma+1<p+1$ points of $\cT$.

  Let $\cM_0:=\sum_{b\in\cT}\chi_{b^\perp}-\chi_E$. For a point $P$ we have $\cM_0(P)=\#\left(\cT\cap P^\perp\right)-1$, hence $0\le\cM_0(P)\le\gamma<p$ since $\cT$ is a blocking set.
  For a line $L=a^\perp$ we have $\cM_0(L)=\sum_{b\in\cT}\#\left(b^\perp\cap a^\perp\right)-(p+1)=\#\cT+p\cdot\#(\cT\cap\{a\})-(p+1)\equiv\#\cT-1\pmod p$. Thus $\cM_0$ is $p$-divisible. Moreover
  $\#\cM_0=3\gamma(p+1)-(p^2+p+1)=\gamma(p+1)+p$.

  The multiset $\cM_0$ is $p$-divisible, has all multiplicities at most $\gamma\le p-1$, and $\#\cM_0=\gamma(p+1)+p$. So Lemma~\ref{lemma_cone_noncanonical} with $\cB=\cM_0$, $b_1=\gamma$, and $b_2=1$
  gives a non-canonical $p^r$-divisible multiset $\cM$ of cardinality $\gamma s_i+s_{i+1}$.

  For the second statement let $\cM_2$ be the canonical multiset with coefficients $c_k$ for $k\notin\{i,i+1\}$, $c_i-\gamma$ at position $i$, and $c_{i+1}-1$ at position $i+1$.
  As in the proof of Theorem~\ref{thm_digits_only_if}, Lemma~\ref{lemma_direct_sum} shows that $\cM+\cM_2$, embedded with independent supports, is a non-canonical $p^r$-divisible multiset of
  cardinality $n$.
\end{proof}

The bound $c_i\ge\gamma$ cannot be improved by examples in the plane with $r=1$ and $c_1=1$. To see this, let $\cM$ be a $p$-divisible multiset of points of $E=\PG(2,p)$ with $\#\cM=c(p+1)+p$ and $0\le c\le p-1$.
Then $\cM(L)\equiv c\pmod p$ for every line $L$, so $y_L:=(\cM(L)-c)/p$ is a non-negative integer. Counting the points on the $p+1$ lines through a point $P$ gives $\sum_{L\ni P}y_L=1+\cM(P)$, i.e.\
$\cM=\sum_L y_L\chi_L-\chi_E$ with $\sum_L y_L=c+p+1$. If the lines with $y_L>0$ contain all lines through some point $R$, then $\cM=p\chi_R+\sum_{k=1}^c\chi_{L_k}$ is canonical. Otherwise the points dual to these lines form a blocking set
that contains no line, and thus a non-trivial minimal blocking set, which has at least $3(p+1)/2$ points \cite{blokhuis1994size}. It has at most $\sum_Ly_L=c+p+1$ points, which is less than $3(p+1)/2$ for $c<(p+1)/2$. So $c<(p+1)/2$ implies that $\cM$ is canonical.

\section{Proofs for Sections~\ref{sec_length_divisible} and~\ref{sec_length_powers}}\label{app_sec4}
This appendix contains the longer technical proofs of Sections~\ref{sec_length_divisible} and~\ref{sec_length_powers}.
\subsection{Proofs for Subsections~\ref{sub_necessary} and~\ref{sub_nonpower}}\label{app_sec4_necessary}
\begin{proof}[Proof of Theorem~\ref{thm_finite_exceptions}]
  We use the following three $\Delta$-divisible building blocks.
  \begin{itemize}
    \item Let $\cM_1$ be the multiset of all $(h-1)$-spaces of an $h$-dimensional projective space $\PG(h,q)$ (as a subspace of the ambient space), taken $q^{l-1}$ times if $l\ge 1$ and once if $l=0$.
          A hyperplane of the ambient space either contains $\PG(h,q)$, or it meets $\PG(h,q)$ in exactly one of the $[h+1]_q$ elements. Thus all non-zero weights of $\cM_1$ equal $q^{\max\{l-1,0\}}\left([h+1]_q-1\right)=q^{\max\{l,1\}}[h]_q$,
          which is divisible by $\Delta$. We have $\#\cM_1=q^{l-1}[h+1]_q$ for $l\ge 1$ and $\#\cM_1=[h+1]_q$ for $l=0$.
    \item Let $\cM_2$ be the $[h]_q$-fold copy of the field reduction of all points of $\PG(m-1,q^h)$, i.e., of a Desarguesian $(h-1)$-spread of $\PG(hm-1,q)$. Every hyperplane $H$ of $\PG(hm-1,q)$ contains a unique hyperplane of
          $\PG(m-1,q^h)$ (as point set). Hence $H$ contains exactly $[m-1]_{q^h}$ spread elements, and all weights of $\cM_2$ equal $[h]_q\left([m]_{q^h}-[m-1]_{q^h}\right)=[h]_q\,q^{h(m-1)}$, which is divisible by $\Delta$ since
          $h(m-1)\ge l$. We have $\#\cM_2=[h]_q[m]_{q^h}=[hm]_q$, which is coprime to $q$.
    \item Let $\cM_3$ be a $\Delta$-fold $(h-1)$-space, with $\#\cM_3=\Delta$.
  \end{itemize}

  For $l=0$ the integers $[h]_q=\#\cM_3$ and $[h+1]_q=q[h]_q+1=\#\cM_1$ are coprime. Hence every integer larger than $[h]_q[h+1]_q-[h]_q-[h+1]_q=q[h]_q^2-q[h]_q-1$ is a non-negative integer combination of them.
  For $l\ge 1$ let $d:=q^{l-1}$, $a:=[h+1]_q$, and $b:=q[h]_q=a-1$. Then $\#\cM_1=da$, $\#\cM_3=db$, and $\gcd(a,b)=1$. Let $n>F_{q,h,l}$. Since $\gcd(\#\cM_2,d)=1$ there is $t\in\{0,\dots,d-1\}$ with
  $t\cdot\#\cM_2\equiv n\pmod d$. Then $n-t\cdot\#\cM_2$ is divisible by $d$ and $\left(n-t\cdot\#\cM_2\right)/d\ge\left(n-(d-1)[hm]_q\right)/d>ab-a-b$. Therefore $\left(n-t\cdot\#\cM_2\right)/d=xa+yb$ for some $x,y\in\N_0$.
  So $n=x\cdot\#\cM_1+y\cdot\#\cM_3+t\cdot\#\cM_2$, and the union of $x$ copies of $\cM_1$, $y$ copies of $\cM_3$, and $t$ copies of $\cM_2$ (embedded in a common ambient space) is $\Delta$-divisible of length $n$.
  Note that $ab-a-b=q^2[h]_q^2-q[h]_q-1$.

  For the lower bound let $n_0$ be the largest integer with $n_0\le(\Delta-1)q^h$ and $n_0\equiv-1\pmod\Delta$. Then $n_0\ge(\Delta-1)q^h-(\Delta-1)$, and since $\Delta\ge 2$ we have $n_0\ge 1$ and $n_0\equiv\Delta-1\not\equiv 0\pmod\Delta$.
  Embedding a $\Delta$-divisible multiset of $n_0$ $(h-1)$-spaces into an ambient space of dimension larger than $h$, which keeps all non-zero weights, Lemma~\ref{lemma_counting} (with $\rho=\Delta-1$) gives $n_0>(\Delta-1)q^h$, a contradiction. So, no $\Delta$-divisible faithful additive code of length $n_0$ exists.
\end{proof}

\begin{proof}[Proof of Proposition~\ref{prop_N_coprime}]
  (a) For $(q,h,\Delta)=(2,2,6)$ this is the case $a=1$ of Theorem~\ref{thm_N_t1}. For $(q,h,\Delta)=(3,2,2)$ every odd length satisfies $n>9$, and $n=11$ is excluded by Lemma~\ref{lemma_moment} since $4\cdot 11/2=22$ is not of the form $[v']_3$. The lengths $2$
  and $13$ are realized by a double line and the $13$ lines of $\PG(2,3)$ (all weights $12$). For $(q,h,\Delta)=(3,2,6)$ Lemma~\ref{lemma_counting} gives the lower bounds $13$, $20$, $33$, $40$, and $47$ for
  $\rho=1,\dots,5$, and $47=5\cdot 9+2$ is excluded by Lemma~\ref{lemma_moment} since $4\cdot 47/2=94$ is not of the form $[v']_3$. Since $33=13+20$, $40=20+20$, and $53=13+20+20$, every length lies in $\sg{6,13,20}$,
  and $6$, $13$, $20$ are realized by a $6$-fold line, the $13$ lines of $\PG(2,3)$, and two copies of a line spread of $\PG(3,3)$ (weights $6$, $12$, $18$).

  (b) By an exact computer calculation the given sets are the minimal generating sets of $\mathcal{A}_q^{\Delta,h}$. For $\Delta=3\cdot 2^l$ the lengths $[i]_2\cdot 2^{\max\{0,l+2-i\}}$ with $2\le i\le l+3$ are realized by $2^{\max\{0,l+2-i\}}$ copies of the multiset of Lemma~\ref{lemma_singer} in $\PG(i-1,2)$, whose weights are $3\cdot 2^{i-2}$.
  For $\Delta=18$ the lengths $18$, $20$, $39$, and $121=[5]_3$ are realized by an $18$-fold line, two copies of a line spread of $\PG(3,3)$, three copies of the lines of $\PG(2,3)$, and Lemma~\ref{lemma_singer}. The remaining lengths are realized by multisets of points of the types
$3[4]+6[3]+24[1]$ ($n=37$), $6[4]+12[3]+48[1]$ ($74$), $3[5]+12[3]+48[1]$ ($75$), $3[5]+6[4]+48[1]$ ($77$), $12[4]+24[3]+96[1]$ ($148$), $6[5]+24[3]+96[1]$ ($150$), $3[6]+24[3]+96[1]$ ($151$), $6[5]+12[4]+96[1]$ ($154$), $3[6]+12[4]+96[1]$ ($155$), $3[6]+6[5]+96[1]$ ($157$), and $4[4]+6[3]+54[1]$ ($73$ for $q=3$).
The decompositions use the types $C_j:=[j]+q[j-2]_q[1]$ (the cone of Lemma~\ref{lemma_moves}.(b) for $h=2$), $A_i:=[2]_q[i]$ (Lemma~\ref{lemma_singer} with $b=0$), $D_4:=7[4]$ (all lines of $S_4$ for $q=2$), $F_i:=[i]$ (a line spread of $S_i$ for even $i$), and $X:=[5]-[4]+8[1]$ (Lemma~\ref{lemma_moves_general}.(e) with $(q,a,m)=(2,1,4)$):
$37$: $A_3+3C_3+3C_4$, $74$: twice the one for $37$, $75$: $3A_3+3C_3+3C_5$, $77$: $2F_4+A_4+C_4+3C_5$, $148$: four times the one for $37$, $150$: twice the one for $75$,
$151$: $7A_3+3C_3+3C_6$, $154$: $D_4+A_4+2C_4+6C_5$, $155$: $D_4+F_4+A_4+C_4+3C_6$, $157$: $F_4+A_5+F_6+X+2C_5+2C_6$, and for $q=3$, $73$: $A_3+2C_3+4C_4$.
  These types are $h$-partitionable as non-negative integer combinations of the types of
  Lemma~\ref{lemma_moves}, Lemma~\ref{lemma_moves_general}, and Lemma~\ref{lemma_singer}, as the decompositions show. Their divisibility follows from Remark~\ref{rem_type_weights} and Lemma~\ref{lemma_exclusion}.
\end{proof}

\begin{proof}[Proof of Theorem~\ref{thm_digit_sum}]
  Without loss of generality the members of $\cM$ span $\F_q^v$, and let $r:=v-h-1$ and $x:=(n-\rho)/m'=s+(q-1)\rho t$. By Lemma~\ref{lemma_dual} there is a $q^r$-divisible multiset of points of cardinality $s[r+1]_q-\rho t$. Hence
  $s[r+1]_q-\rho t=\sum_{i=0}^{r}a_iq^i[r+1-i]_q$ for some non-negative integers $a_i$ by \cite{kiermaier2020lengths}. Let $A:=\sum_{i=0}^{r}a_i$. Multiplying by $q-1$ and using $(q-1)q^i[r+1-i]_q=q^{r+1}-q^i$ gives
  $$
    \sum_{i=0}^{r}a_iq^i=(A-s)q^{r+1}+x,
  $$
  i.e., $(A-s)q^{r+1}+x$ is the sum of a multiset $T$ of $A$ powers $q^i$ with $0\le i\le r$. If $A\le s$, then $x$ is the sum of the $A$ elements of $T$ and $s-A$ copies of $q^{r+1}$.
  If $A>s$, let $k:=A-s$, order the elements of $T$ as $q^{e_1}\ge\dots\ge q^{e_A}$. If the sum of the first $l-1$ elements is less than $kq^{r+1}$, then the sum of the first $l$ elements is at most $kq^{r+1}$, since $q^{e_l}$ divides both the former sum and $kq^{r+1}$.
  The sums of the first $l$ elements increase from $0$ for $l=0$ to $kq^{r+1}+x\ge kq^{r+1}$ for $l=A$. So the sum of the first $l$ elements equals $kq^{r+1}$ for some $l$, and $l\ge qk$ since all elements of $T$ are at most $q^r$. So, $x$ is the sum of the remaining $A-l\le s+k-qk\le s$ elements of $T$.
  In both cases $x$ is the sum of at most $s$ powers of $q$, i.e., $\sigma_q(x)\le s$.
\end{proof}

\begin{proof}[Proof of Theorem~\ref{thm_N_t1}]
  Since $q^h\equiv 1\pmod{[h]_q}$, $s$ is an integer. The conditions are necessary by Theorem~\ref{thm_digit_sum} with $m'=[h]_q$ and $t=1$. Conversely, let $x:=(n-\rho)/[h]_q=s+(q-1)\rho$. Since $p^a$ divides $q$, we have
  $q^h\equiv 0$ and $[h]_q\equiv 1\pmod{p^a}$. Hence $n\equiv\rho\pmod{p^a}$ gives $s\equiv\rho$ and $x\equiv q\rho\equiv 0\pmod{p^a}$. The base $q$ expansion writes $x$ as a sum of $\sigma_q(x)$ powers of $q$, and the
  number $x\bmod q$ of summands equal to $1$ is divisible by $p^a$. Replacing a summand $q^i$ with $i\ge 1$ by $q$ summands $q^{i-1}$ increases the number of summands by $q-1$ and keeps the number of summands equal to $1$ divisible by $p^a$.
  Since $\sigma_q(x)\le s\le x$ and $s\equiv x\pmod{q-1}$, we obtain $x=\sum_{i\ge 0}c_iq^i$ with non-negative integers $c_i$, $\sum_ic_i=s$, and $p^a\mid c_0$. Then $\sum_ic_i[i]_q=(x-s)/(q-1)=\rho$, and
  $[h+i]_q=[i]_q+q^i[h]_q$ gives
  $$
    n=\rho q^h+[h]_qs=\rho+[h]_qx=\frac{c_0}{p^a}\cdot\Delta+\sum_{i\ge 1}c_i[h+i]_q.
  $$
  The summand $\Delta$ is realized by a $\Delta$-fold $(h-1)$-space and $[h+i]_q$ with $i\ge 1$ by Lemma~\ref{lemma_singer} in $\PG(h+i-1,q)$, whose weights $q^i[h]_q$ are divisible by $\Delta$. So, $n$ is a length.

  For $i\ge h+1$ the length $[h+i]_q=[i]_q+q^i[h]_q$ is the sum of the length $[i]_q$ and a multiple of $\Delta$, and $[2h]_q=\left(q^h+1\right)[h]_q$ is a multiple of $\Delta$ if $a=0$. So, the stated sets generate all lengths.
  For $1\le i\le h$ we have $[h+i]_q=[i]_qq^h+[h]_q$, that is, $\rho=[i]_q$ (for $a=0$ and $i<h$, and for $a\ge 1$) and $s=1$. Therefore $[h+i]_q$ is the smallest element of its residue class allowed by Lemma~\ref{lemma_moment}.
  If $[h+i]_q=n_1+n_2$ for lengths $n_1,n_2\ge 1$ with residues $\rho_1,\rho_2$, then $\rho_1=0$ would give a smaller length $n_2$ in the residue class of $[h+i]_q$, and $\rho_1,\rho_2\ge 1$ would give
  $n_1+n_2\ge(\rho_1+\rho_2)q^h+2[h]_q>[i]_qq^h+[h]_q$, since $\rho_1+\rho_2\equiv[i]_q\pmod\Delta$. So, no stated generator of the form $[h+i]_q$ is a sum of smaller lengths.
  Also $\Delta$ is not a sum of smaller lengths: $\rho_1=0$ would give $n_1\ge\Delta$, and $\rho_1,\rho_2\ge 1$ with $\rho_1+\rho_2\equiv 0\pmod\Delta$ gives $\rho_1+\rho_2\ge\Delta$ and $n_1+n_2\ge\Delta q^h+2[h]_q>\Delta$ by Lemma~\ref{lemma_moment}.
\end{proof}

\begin{proof}[Proof of Proposition~\ref{prop_N_coprime_dual}]
  For the eight parameters $(2,3,7)$, $(2,3,14)$, $(3,2,4)$, $(3,2,12)$, $(4,2,5)$, $(4,2,10)$, $(4,2,20)$, and $(5,2,6)$, which satisfy $t=1$ and $p^a\mid q$, the statement is Theorem~\ref{thm_N_t1}.
  For the remaining nine parameters Theorem~\ref{thm_digit_sum} excludes the following lengths $n$, where we state $\rho$, $s$, and $x=(n-\rho)/m'$ in base $q$, and $\sigma_q(x)>s$ in each case.
  \begin{center}
  \small
  \begin{tabular}{lll}
    \hline
    $(q,h,\Delta)$ & $t$ & $n\,(\rho,s,x)$\\
    \hline
    $(2,3,28)$ & $1$ & $222\,(26,2,11100_2)$\\
    $(2,4,3)$ & $5$ & $22\,(1,2,111_2)$\\
    $(2,4,6)$ & $5$ & $70\,(4,2,10110_2)$\\
    $(2,4,12)$ & $5$ & $166\,(10,2,110100_2)$\\
    $(2,4,5)$ & $3$ & $58\,(3,2,1011_2)$, $74\,(4,2,1110_2)$, $79\,(4,3,1111_2)$\\
    $(3,2,36)$ & $1$ & $323\,(35,2,2200_3)$\\
    $(4,2,40)$ & $1$ & $298\,(18,2,320_4)$, $554\,(34,2,1220_4)$, $639\,(39,3,1320_4)$\\
    $(5,2,2)$ & $3$ & $29\,(1,2,24_5)$\\
    $(5,2,3)$ & $2$ & $56\,(2,2,33_5)$, $59\,(2,3,34_5)$\\
    \hline
  \end{tabular}
  \end{center}
  By an exact computer calculation, for each of these nine parameters and each $\rho\in\{1,\dots,\Delta-1\}$ the smallest $n\equiv\rho\pmod\Delta$ that satisfies Lemma~\ref{lemma_exclusion}, Lemma~\ref{lemma_counting}, Lemma~\ref{lemma_moment},
  and Theorem~\ref{thm_digit_sum} is the smallest element of the residue class of $\rho$ in the semigroup generated by the stated set. Let $m(\rho)$ be this smallest element, so that every length $n\equiv\rho\pmod\Delta$ satisfies $n\ge m(\rho)$. The stated set is a minimal generating set of the semigroup it generates (an exact computation), so each stated generator $g\ne\Delta$ equals $m(\rho_g)$.
  If $g=n_1+n_2$ with lengths $n_1,n_2\ge 1$, then $n_i\ge m(\rho_i)$ and $m(\rho_g)\le m(\rho_1)+m(\rho_2)\le g$, hence $g=m(\rho_1)+m(\rho_2)$ with both summands in the semigroup and smaller than $g$ (if $\rho_1=0$, then $n_1\ge\Delta$ and $n_2<g$ lies in the residue class of $g$), a contradiction to the minimality of the generating set.
  Since the set of lengths is closed under addition and contains $\Delta$, it remains to realize the stated lengths. They are realized as follows.
  \begin{itemize}
    \item The $\Delta$-fold $(h-1)$-space has length $\Delta$.
    \item Lemma~\ref{lemma_singer} gives $[v]_q$ $(h-1)$-spaces with weights $q^{v-h}[h]_q$, which realizes $15$, $31$, $63$, and $127$ for $(q,h)=(2,3)$,
          $13$, $40$, and $121$ for $q=3$, $21$, $85$, and $341$ for $q=4$, and $31$ for $q=5$. For each stated set that contains such a length, $\Delta$ divides the weight $q^{v-h}[h]_q$.
    \item Two copies of $15$, three copies of $13$, and two copies of $21$ give $30$, $39$, and $42$, with weights $28$, $36$, and $40$.
    \item For $(q,h,\Delta)=(2,4,5)$ the field reduction of the $21$ lines of $\PG(2,4)$ is a set of $21$ solids of $\PG(5,2)$. Every hyperplane of $\PG(5,2)$ contains exactly one of them, and hence all weights equal $20$.
    \item The lengths $85$ for $(2,3,28)$, $145$ and $146$ for $(3,2,36)$, and $201$ for $(4,2,40)$ are realized by the types $7[5]+14[4]+56[2]$, $4[4]+24[3]+108[1]$, $8[4]+12[3]+108[1]$, and $5[4]+20[3]+160[1]$, respectively. They are $h$-partitionable as in the proof of
          Proposition~\ref{prop_N_coprime}.(b), and their divisibility follows from Remark~\ref{rem_type_weights}.
    \item The lengths $25$ and $38$ for $(2,4,\Delta)$ with $\Delta\in\{3,6,12\}$ are realized by Proposition~\ref{prop_25_38}, and $51$ for $\Delta\in\{6,12\}$ by three copies of a Desarguesian $3$-spread of $\PG(7,2)$, whose weights are $48$.
  \end{itemize}
\end{proof}

\begin{proof}[Proof of Proposition~\ref{prop_line_cone}]
  The planes and solids through $\ell$ correspond to the points and lines of the quotient geometry $\PG(h-1,q)$. For each of the $L$ given lines of the quotient take the $q^2$ lines different from $\ell$ of a line spread
  of the corresponding solid that contains $\ell$, which exists since the collineation group of $\PG(3,q)$ is transitive on lines. These $q^2$ lines partition the points of the solid outside $\ell$. So, a point outside $\ell$ in a plane
  $E\ge\ell$ is covered $\tau(E)\le\rho'$ times, where $\tau(E)$ is the number of the $L$ lines of the quotient through $E$, and no point of $\ell$ is covered. Since $\sum_E\left(\rho'-\tau(E)\right)=\rho'[h]_q-(q+1)L=C$, we can assign
  points $X$ of $\ell$ to the planes $E\ge\ell$, $\rho'-\tau(E)$ of them to $E$ (repetitions allowed), such that every $X$ is assigned $c_X$ times in total. For each assignment of $X$ to $E$ we take the $q$ lines through $X$ in $E$ that are different from $\ell$,
  which partition the points of $E$ outside $\ell$. Adding $\mu$ copies of $\ell$ gives $\cK$, and the last statement follows from Lemma~\ref{lemma_dual_cover}.
\end{proof}

\begin{proof}[Proof of Proposition~\ref{prop_25_38}]
  For $n=25$ let $\ell$ be a line of $\PG(5,2)$, $X$ a point of $\ell$, and $S_1,\dots,S_5$ the solids through $\ell$ that correspond to a line spread of the quotient $\PG(3,2)$, i.e., every point outside $\ell$
  lies in exactly one of them. The multiset consisting of $\ell$, the $12$ lines different from $\ell$ through $X$ in $S_1$ or $S_2$, and for $3\le i\le 5$ the four lines different from $\ell$ of a line spread of $S_i$ that contains $\ell$ is the
  example after Proposition~\ref{prop_line_cone}: $X$ lies on $13$ and every other point on exactly one of these $25$ lines. By Lemma~\ref{lemma_dual_cover} the $25$ orthogonal complements are solids, the hyperplane $X^\perp$ contains $13$ and every other
  hyperplane exactly one of them. Hence the weights are $12$ and $24$.

  For $n=38$, by Lemma~\ref{lemma_dual_cover} it suffices to give $38$ planes of $\PG(6,2)$ such that one point $P$ lies on $14$ and every other point on $2$ of them. Let $\alpha\in\F_8$ with $\alpha^3=\alpha+1$, identify $\F_2^7$ with
  $\F_8\times\F_8\times\F_2$, and let $P=(0,0,1)$ and $W=\F_8\times\F_8\times\{0\}$. For $\lambda\in\F_8$ let $E_\lambda=\{(\lambda y,y,0)\,:\,y\in\F_8\}$ and let $E_\infty=\F_8\times\{0\}\times\{0\}$, which are nine planes partitioning the points of $W$.
  Let $\Sigma$ be the set of the $21$ points of $E_0\cup E_1\cup E_\infty$ and $\Sigma'$ the set of the remaining $42$ points $(x,y,0)$ of $W$, i.e., $x,y\neq 0$ and $x\neq y$. For $c\in\F_8^*$ the set
  $T_c=\{(\alpha c,c,0),(\alpha^3c,\alpha c,0),(c,\alpha^3c,0)\}$ is a line since $\alpha+\alpha^3+1=0$, and its points lie in $E_\alpha$, $E_{\alpha^2}$, and $E_{\alpha^4}$. Let $T'_c$ be the image of $T_c$ under $(x,y,z)\mapsto(y,x,z)$, whose
  points lie in $E_{\alpha^6}$, $E_{\alpha^5}$, and $E_{\alpha^3}$. So, the $14$ lines $T_c$, $T'_c$ partition $\Sigma'$.

  Let $G$ be the group of the $42$ maps $(x,y,z)\mapsto\left(\lambda(ax+by),\lambda(cx+dy),z\right)$ with $\lambda\in\F_8^*$ and $\left(\begin{smallmatrix}a&b\\c&d\end{smallmatrix}\right)\in\mathrm{GL}(2,2)$. It fixes $P$ and $W$, the scalars fix every $E_\mu$, and $\mathrm{GL}(2,2)$
  permutes the $E_\mu$ as $\mathrm{PGL}(2,2)$ acts on $\F_8\cup\{\infty\}$. This action is transitive on $\{0,1,\infty\}$ and regular on $\F_8\backslash\F_2$, since a fixed point in $\F_8\backslash\F_2$ of a non-trivial element of $\mathrm{PGL}(2,2)$ would be a root of a
  polynomial over $\F_2$ of degree at most $2$. So, $\Sigma$ and $\Sigma'$ are orbits of $G$, and since $Q\mapsto Q+P$ commutes with $G$, also $\Sigma+P$ and $\Sigma'+P$ are orbits. Let
  \begin{eqnarray*}
    \pi&=&\left\langle(1,0,0),(\alpha^4,1,0),(\alpha^2,0,1)\right\rangle\\
    &=&\left\{(1,0,0),(\alpha^4,1,0),(\alpha^5,1,0),(\alpha^2,0,1),(\alpha^6,0,1),(\alpha,1,1),(\alpha^3,1,1)\right\},
  \end{eqnarray*}
  whose stabilizer in $G$ consists of the identity and $(x,y,z)\mapsto(x+y,y,z)$. So its orbit $G\pi$ consists of $21$ planes. The plane $\pi$ does not contain $P$ and meets $\Sigma$, $\Sigma'$, $\Sigma+P$, and $\Sigma'+P$ in $1$, $2$, $2$, and $2$ points, respectively.
  Double counting incidences between an orbit $O$ of points and $G\pi$ shows that every point of $O$ lies in $21\cdot\#(\pi\cap O)/\#O$ members of $G\pi$, that is, in $1$, $1$, $2$, and $1$ of them for $O=\Sigma,\Sigma',\Sigma+P,\Sigma'+P$, respectively.

  Now take the three planes $E_0$, $E_1$, $E_\infty$, which cover $\Sigma$ once, the $14$ planes $\langle P,T_c\rangle$ and $\langle P,T'_c\rangle$, which contain $P$ and cover $\Sigma'$ and $\Sigma'+P$ once, since $\langle P,T\rangle=\{P\}\cup T\cup(T+P)$, and the $21$ planes of $G\pi$.
  Then $P$ lies on $14$ and every other point on $2$ of the $38$ planes. By Lemma~\ref{lemma_dual_cover} the hyperplane $P^\perp$ contains $14$ and every other hyperplane $2$ of the $38$ orthogonal complements. Consequently, the weights are $24$ and $36$.
\end{proof}

The point multiset of the $38$ solids is $3\chi_V+3\chi_{P^\perp}$. In the variant of the proof for $n=25$ that uses the $6$ lines through a second point $Y$ of $\ell$ in $S_2$ instead of those through $X$ (weights $18^2\,24^{61}$), the point multiset of the $25$ solids is $3\chi_V+3\chi_{X^\perp}+3\chi_{Y^\perp}$.
Both point multisets are as predicted by Lemma~\ref{lemma_dual}. The variant for $n=25$ is $6$-divisible.

\subsection{Proofs for Subsection~\ref{sub_mixed}}\label{app_sec4_mixed}
\begin{proof}[Proof of Theorem~\ref{thm_mixed_qpower}]
  Necessity is Lemma~\ref{lemma_mixed_counting}. Let $\Lambda$ be the set of all $n\ge 1$ with $\rho=0$ or $n>q\rho$. If $n,n'\in\Lambda$ have residues $\rho,\rho'$, then the residue of $n+n'$ is at
  most $\rho+\rho'$ and $n+n'>q(\rho+\rho')$ unless $\rho=\rho'=0$. Hence $\Lambda$ is additively closed. An integer $q^l+\rho$ with $0\le\rho<q^l$ is in $\Lambda$ if and only if $(q-1)\rho<q^l$, i.e., if and only if
  $\rho\le[l]_q$. Conversely, let $n=m\Delta+\rho\in\Lambda$ with $m\ge 1$ and $\rho\ge 1$. Then $(q-1)\rho\le m\Delta-1$. For $m\ge q$ we have $\rho\le\Delta-1\le q[l]_q\le m[l]_q$. For $m<q$
  note that $(q-1)\rho$ and $m(\Delta-1)$ are multiples of $q-1$ and $m\Delta-1<m(\Delta-1)+q-1$. Hence $(q-1)\rho\le m(\Delta-1)$, i.e., $\rho\le m[l]_q$. In both cases
  $\rho=\rho_1+\dots+\rho_m$ with $0\le\rho_i\le[l]_q$ and $n=\sum_{i=1}^m(q^l+\rho_i)$. So, $\Lambda$ is generated by $q^l,\dots,q^l+[l]_q=[l+1]_q$.

  It remains to construct $\Delta$-divisible examples of each cardinality $q^l+\rho$, $0\le\rho\le[l]_q$. By the remark before Corollary~\ref{cor_mixed_finite} with $k=2$ and $m=1$ we obtain $q^l+c$
  for $0\le c\le q^{l-1}$. The set of all points of $\PG(l,q)$ is $q^l$-divisible of cardinality $[l+1]_q$. For $l\ge 3$ let $L_1,\dots,L_{q^{l-1}}$ be the lines of Lemma~\ref{lemma_moves}.(a) with $h=2$, $i=l+1$, and $S_{l-1}=W$, which partition the points of $\F_q^{l+1}\backslash W$ for a subspace $W$ of dimension $l-1$. Replacing the $q+1$ points of $L_1,\dots,L_j$ by $q$ copies of each of these lines does not change any weight by Lemma~\ref{lemma_exchange}. This yields
  $[l+1]_q-j$ for $0\le j\le q^{l-1}$. Since $[l+1]_q-q^{l-1}=q^l+[l-1]_q\le q^l+q^{l-1}+1$, all cardinalities between $q^l$ and $[l+1]_q$ occur. For $l\le 2$ we have
  $[l+1]_q\le q^l+q^{l-1}+1$ and the first two families suffice.

  The largest $n\notin\Lambda$ satisfies $\rho=\Delta-1$ and $m\Delta\le(q-1)(\Delta-1)$, that is, $m=q-2$. So $n=(q-1)q^l-1$.
\end{proof}

\begin{Remark}
  \label{rem_mixed_general_Delta}
  For $\Delta=p^c$ with $q=p^e$ and $e\nmid c$ write $\Delta=p^{ae-b}$ with $1\le b\le e-1$. Besides the family $m\Delta+c'$, $qc'\le m\Delta$, there are the points-only lengths
  $\left\langle S_q(a,b)\right\rangle$ and the exchanges of Lemma~\ref{lemma_exchange} in $p^{e-b}\chi_{\PG(a-1,q)}$, which give the cardinalities $p^{e-b}[a]_q-j$ for $0\le j\le p^{e-b}q^{a-2}$ if $a\ge 4$ and for
  $0\le j\le p^{e-b}$ if $a\in\{2,3\}$. For small $q$ and $\Delta$ these examples realize all $n\ge 1$ with $n\equiv 0\pmod\Delta$ or $n>q\rho$ for small $q$ and $\Delta$, and we conjecture the characterization of Theorem~\ref{thm_mixed_qpower} for all powers $\Delta$ of the characteristic.

  For $\Delta$ coprime to $q$ the condition $n>q\rho$ is not sufficient and the answer depends on $h$. For $(q,h,\Delta)=(2,2,3)$ Corollary~\ref{cor_coprime} says that every point $P$
  satisfies $d_P+2x_P\equiv 0\pmod 3$, where $d_P$ and $x_P$ are the numbers of line and point members through $P$. Summing over the points shows that $3$ divides the number $n_1$ of point members. For $n=5$
  the case $n_1=0$ is excluded by Lemma~\ref{lemma_counting}. If $n_1=3$, then every point on one of the two lines has $x_P\ge 1$ (since $d_P\in\{1,2\}$) and $x_P\ge 2$ if $d_P=2$. Hence $\sum_Px_P\ge 6$.
  So the lengths are $\sg{3,4}$ although $n=5>2\cdot 2$. The remaining statements follow from Lemma~\ref{lemma_mixed_counting}, the construction before Corollary~\ref{cor_mixed_finite}, and Corollary~\ref{cor_mixed_finite}.

  For $\Delta=5$ and $h=3$ a plane together with its seven points is $5$-divisible of cardinality $8$. For $h=2$ there is no $5$-divisible example of cardinality $8$. Here $d_P+2x_P\equiv 0\pmod 5$
  and $3n_2+2n_1\equiv 0\pmod 5$ give $n_1=n_2=4$. A point on a line member has $d_P\le 4$. Hence $x_P\ge f(d_P)$ with $f(1)=2$, $f(2)=4$, $f(3)=1$, $f(4)=3$, and $f(d)\ge d/3$ with equality
  only for $d=3$. Since $\sum_Pd_P=12$ and $\sum_Px_P=4$, every point on a line member lies on exactly three of the four lines. Then two of the lines share two points, so that all four lines coincide, a contradiction.
\end{Remark}

\begin{Remark}
  For unfaithful codes some results of Subsection~\ref{sub_necessary} change. Let $q=h=2$. By Subsection~\ref{sub_mixed} a multiset $\cS$ of points, regarded as a projective $2$-system, is $\Delta$-divisible if and only if $2\cS$ is $2\Delta$-divisible.
  For $\Delta=2$ the length $1$ is still impossible, since there is no $4$-divisible multiset of points of cardinality $2$ or $3$. However, the three points of a line $L$ form a $2$-divisible projective $2$-system, since the double line $2\chi_L$ is $4$-divisible.
  This gives a $2$-divisible additive code over $\F_4$ of length $3\notin\sg{2,5}$. Similarly, $4$-divisible multisets of points of cardinality $6$, $7$, and $11$ give $4$-divisible additive codes over $\F_4$ of these lengths, which are not contained in $\sg{4,5}$.
\end{Remark}

\medskip
\noindent\textbf{Subspace designs.} Let $\mathcal{B}$ be a $t$-$(v,k,\lambda)_q$ design with $t\ge 1$, i.e., a multiset of $(k-1)$-spaces, the blocks. Then every point lies in exactly $\lambda_1=\lambda\qbin{v-1}{t-1}{q}/\qbin{k-1}{t-1}{q}$ blocks,
so the multiset of points of $\mathcal{B}$ is $\lambda_1\chi_V$, and Lemma~\ref{lemma_exclusion} gives $\#\mathcal{B}-\widetilde{\mathcal{B}}(H)=\lambda_1q^{v-1}/q^{k-1}=\lambda_1q^{v-k}$ for every hyperplane $H$. So $\mathcal{B}$ is a one-weight code, and
it is $\Delta$-divisible if and only if $\Delta$ divides $\lambda_1q^{v-k}$. Every hyperplane contains $\lambda_1[v-k]_q/[k]_q$ blocks.
\begin{Proposition}
  \label{prop_designs}
  For admissible parameters of a $t$-$(v,k,\lambda)_q$ design, $t\ge 1$, the cardinality $b=\lambda_1[v]_q/[k]_q$ and $\Delta=\lambda_1q^{v-k}$ satisfy the necessary conditions of
  Lemma~\ref{lemma_exclusion}, Lemma~\ref{lemma_counting}, and Lemma~\ref{lemma_incidence}.
\end{Proposition}
\begin{proof}
  The multiset $\lambda_1\chi_V$ is a $\Delta q^{k-1}$-divisible multiset of $b[k]_q$ points. With $b-\Delta=\lambda_1[v-k]_q/[k]_q\ge 0$ and $b-\Delta\equiv b\pmod\Delta$ we have
  $\rho\le b-\Delta$ and $b[v-k]_q=(b-\Delta)[v]_q\ge\rho[v]_q$. Finally $\Delta\mid q^{v-k}[k]_qb=q^{v-k}\lambda_1[v]_q$.
\end{proof}
Consequently, the length results cannot exclude a design with admissible parameters. The same holds for the blocks contained in a subspace of codimension at most $t$, since their numbers are
determined by the parameters as well.

For the binary $q$-analog of the Fano plane, a $2$-$(7,3,1)_2$ design, we have $b=381$, $\lambda_1=21$, every hyperplane contains $45$ blocks, and every
subspace of codimension $2$ contains $5$ blocks. This design is $336$-divisible, and the $45$ blocks in a hyperplane $H$ form a one-weight code in $H$ with weight $40$. Three orbits of planes under a
Singer cycle of $\PG(6,2)$ form a $1$-$(7,3,21)_2$ design, i.e., a $336$-divisible multiset of $381$ planes (by Lemma~\ref{lemma_singer} each orbit covers every point $[3]_2=7$ times). The derived design in a point $P$ is a line spread of
$\PG(5,2)\cong V/P$. Therefore Proposition~\ref{prop_line_spreads_PG52} characterizes the points with a Desarguesian derived design by the $4$-divisibility of the hyperplane residuals of the derived design.

\subsection{Proofs for Subsections~\ref{sub_partitions} and~\ref{sub_constructions}}\label{app_sec4_lines}\label{app_sec4_qh1}
\begin{proof}[Proof of Theorem~\ref{thm_N_explicit}]
  By Theorem~\ref{thm_D} the set of lengths is the semigroup $\Lambda_l:=\{n\in\N_0\,:\,3n\in \mathcal{L}_2^{l+1}\}$, and $\cN_2^{2^l,2}$ is the set of its indecomposable elements, i.e., of the elements that are not the sum of two non-zero elements of $\Lambda_l$.
  The proof of Theorem~\ref{thm_D} shows that the stated set, consisting of the numbers $X_t$ and $Y_{r,b,c}$, generates $\Lambda_l$. So, it remains to show that each of these numbers is indecomposable.

  Let $M:=2^{l+1}$. Then $\mathcal{L}_2^{l+1}$ is generated by the numbers $G_j=2M-2^{l+2-j}$, $1\le j\le l+2$, of the proof of Theorem~\ref{thm_D}. A sum of $k$ generators has the form $2kM-T$, where $T$ is a sum of $k$ powers $2^e\le M$.
  Merging two equal summands $2^e<M$ into $2^{e+1}$ until all summands smaller than $M$ are distinct, we end with $\lfloor T/M\rfloor$ summands $M$ and the binary expansion of $T\bmod M$. Hence
  $k\ge f(T):=\lfloor T/M\rfloor+\wt(T\bmod M)$, where $\wt(u)$ is the number of ones in the binary expansion of $u$, with equality only if the summands of $T$ are already of this canonical form.
  For each stated number $x$ we have $3x=2k_0M-T_0$ with the following $k_0$ and canonical summands of $T_0$. For $2t\le l+2$ we have $3X_t=G_{2t}$, $k_0=1$, and $T_0=2^{l+2-2t}<M$.
  For $2t=l+3$ we have $3X_t=4M-1=2G_1+G_{l+2}$, $k_0=3$, and $T_0=M+M+1$. Finally $3Y_{r,b,c}=G_r+G_b+G_c$, $k_0=3$, and $T_0=2^{l+2-r}+2^{l+2-b}+2^{l+2-c}<2M$.
  In all cases $f(T_0)=k_0$, and for $k_0=3$ the indices of the three generators are odd.

  Assume $x=y+z$ with $y,z\in\Lambda_l\backslash\{0\}$. Combining representations of $3y$ and $3z$ writes $3x$ as a sum of $k\ge 2$ generators with $T=T_0+2(k-k_0)M\ge k$.
  For $k_0=3$ and $k=2$ we have $T=T_0-2M<2$, which is impossible. So, $k\ge k_0$, and then $f(T)=f(T_0)+2(k-k_0)=2k-k_0$ and $k\ge f(T)$ give $k=k_0=f(T)$. This excludes $k_0=1$.
  For $k_0=3$ the representation consists of the three generators above, and one of $3y$ and $3z$ is a single generator $G_j$ with odd $j$, which contradicts $3\nmid 2^{l+2-j}\left(2^j-1\right)=G_j$.

  The stated numbers are pairwise different, since $3X_t<2M$ for $2t\le l+2$, $3x\ge 4M-1$ for all other stated numbers $x$, and $3x$ determines $T_0$ and hence the generators $G_j$ in the representation of $3x$ above. Counting gives $\lfloor (l+3)/2\rfloor$ numbers $X_t$ and, since for $r=2i+1$ there are $\binom{i}{2}$ pairs of odd integers $1\le c<b<r$, $\sum_{i=2}^{K}\binom{i}{2}=\binom{K+1}{3}$ numbers $Y_{r,b,c}$, where $K=\lfloor (l+1)/2\rfloor$.
\end{proof}

\begin{Lemma}
  \label{lemma_quadrangle}
  Let $B=\F_q^{k}$ with $k\ge 3$, let $\pi$ be a $3$-dimensional subspace of the dual space $B^*$, and let $P$ be a point of $\pi$. Let $L_1,\dots,L_{q^2}$ be the $2$-dimensional subspaces of $\pi$ that
  do not contain $P$. Then there are vectors $b_1,\dots,b_{q^2}\in B$ such that the cosets $b_i+L_i^\perp$, $1\le i\le q^2$, partition $B$.
\end{Lemma}
\begin{proof}
  Let $W:=\pi^\perp$. Then $\dim(B/W)=3$, and every $L_i^\perp$ contains $W$ with $\dim\left(L_i^\perp/W\right)=1$. Since $P\not\le L_i$, the subspace $L_i^\perp$ is not contained in $P^\perp$. So, the $L_i^\perp/W$ are
  exactly the $q^2$ one-dimensional subspaces of $B/W$ that are not contained in the $2$-dimensional subspace $P^\perp/W$, and it suffices to partition $B/W$ into cosets of them. Choose coordinates
  $B/W=\F_{q^2}\times\F_q$, where $\F_{q^2}$ is regarded as a $2$-dimensional vector space over $\F_q$, such that $P^\perp/W=\F_{q^2}\times\{0\}$. Then the subspaces in question are $\langle(s,1)\rangle$ with $s\in\F_{q^2}$.
  Fix $\omega\in\F_{q^2}\backslash\F_q$. A vector $(x,\lambda)\in\F_{q^2}\times\F_q$ lies in the coset $(\omega s,0)+\langle(s,1)\rangle$ if and only if $x=(\omega+\lambda)s$, i.e., if and only if $s=x/(\omega+\lambda)$, which is well
  defined since $\omega+\lambda\neq 0$. So, these $q^2$ cosets partition $B/W$.
\end{proof}

\begin{proof}[Proof of Lemma~\ref{lemma_moves_general}.(g)]
  Let $B:=S_{h+1}$, write $S_{h+1+u}=B\oplus U$ with $\dim(U)=u$, and let $B^*$ be the dual space of $B$. The point multiset is obtained in two steps. Step~(1) takes, for each point $\langle x\rangle$ of $U$, $k$ units, each of which
  is a set of $q^2$ subspaces of dimension $h$ that cover every point of $\langle B,x\rangle\backslash B$ exactly once. Step~(2) takes hyperplanes of $B$. Let $N:=k[u]_q$ be the number of units. We consider two kinds of units.
  \begin{itemize}
    \item[(P)] For an $(h-1)$-space $K\le B$ and a complement $W_K$ of $K$ in $B$ let $\cP(x,K):=\left\{K+\langle x+w\rangle\,:\,w\in W_K\right\}$. These $q^2$ subspaces meet $B$ in $K$, and each point of $\langle B,x\rangle\backslash B$ is
          $\langle x+b\rangle$ with a unique $b=w+k'$, $w\in W_K$, $k'\in K$. A point $Q$ of $B$ is covered $q^2$ times if $Q\le K$ and not at all otherwise.
    \item[(Q)] For a $3$-dimensional subspace $\pi$ of $B^*$ and a point $P$ of $\pi$ let $L_1,\dots,L_{q^2}$ be the $2$-dimensional subspaces of $\pi$ not containing $P$ and $b_i+L_i^\perp$ the partition of $B$ from Lemma~\ref{lemma_quadrangle}.
          Let $\mathcal{Q}(x,\pi,P)$ be the set of the $q^2$ subspaces $L_i^\perp+\langle x+b_i\rangle$ of dimension $h$. Since the points of $L_i^\perp+\langle x+b_i\rangle$ outside $B$ are the points $\langle x+b_i+k'\rangle$ with $k'\in L_i^\perp$,
          each point of $\langle B,x\rangle\backslash B$ is covered exactly once, and a point $Q$ of $B$ is covered $\#\{i\,:\,Q\le L_i^\perp\}$ times.
  \end{itemize}
  \emph{Case (P).} Choose $(h-1)$-spaces $K_1,\dots,K_N$ of $B$ and assign them to the $N$ units. Let $t(H)$ be the number of $K_j\le H$ and take each hyperplane $H$ of $B$ exactly $m(H):=m_0-q\,t(H)$ times. A point $Q$ of $B$ lies in $[h]_q$ hyperplanes of $B$,
  and exactly $q+1$ of them contain $K_j$ and $Q$ if $Q\le K_j$, and exactly one otherwise. So, step~(2) covers $Q$ exactly $m_0[h]_q-q\left(N+q\deg(Q)\right)$ times, where $\deg(Q)$ is the number of $j$ with $Q\le K_j$, and together with step~(1), which covers $Q$ exactly $q^2\deg(Q)$ times,
  every point of $B$ is covered $m_0[h]_q-qN=\left(m_0[h]_q-k[u+1]_q\right)+k$ times, since $q[u]_q+1=[u+1]_q$. Every point of $S_{h+1+u}\backslash B$ is covered $k$ times. It remains to choose the $K_j$ with $q\,t(H)\le m_0$.
  In $B^*\cong\F_q^2\oplus\F_q^{h-1}$ Lemma~\ref{lemma_moves}.(a) gives $q^{h-1}$ pairwise disjoint $2$-dimensional subspaces $L_a$. Choose the $K_j$ among the $L_a^\perp$, each at most $\tau:=\lceil N/q^{h-1}\rceil$ times.
  Since $K\le H$ if and only if $H^\perp\le K^\perp$ and every point of $B^*$ lies on at most one $L_a$, we have $t(H)\le\tau$, and $q\tau\le m_0$ by assumption.

  \emph{Case (Q), $h\ge 4$.} Choose $N$ pairs $(\pi_j,P_j)$ and assign them to the $N$ units. Let $t(H)$ be the number of $j$ with $H^\perp\in\pi_j\backslash\{P_j\}$ and take each hyperplane $H$ of $B$ exactly $m(H):=m_0-t(H)$ times.
  A point of $\pi_j$ other than $P_j$ lies on exactly $q$ of the corresponding lines $L_i$, and $P_j$ lies on none. So, $q\,t(H)$ is the number of the $q^2N$ triples $(j,i)$ with $L_i^\perp\le H$. For each triple exactly $q+1$ hyperplanes of $B$ contain $L_i^\perp$ and $Q$
  if $Q\le L_i^\perp$, and exactly one otherwise. Since step~(1) covers $Q$ once for each triple with $Q\le L_i^\perp$, the steps~(1) and~(2) cover $Q$ together exactly $m_0[h]_q-qN$ times, as in case~(P).
  It remains to choose the pairs such that $t(H)\le m_0$. Write $B^*=E\oplus Z$ with $\dim(E)=3$ and $\dim(Z)=h-2\ge 2$, choose a point $P=\langle p\rangle$ of $E$ and a complement $A$ of
  $\langle p\rangle$ in $E$, identify $A$ with $\F_{q^2}$, and fix $\omega\in\F_{q^2}\backslash\F_q$. Let $\cK$ be a multiset of $[h-2]_q$ lines of $Z$ such that every point of $Z$ lies on exactly $q+1$ of them.
  (If $h-2$ is even, take $q+1$ copies of a line spread of $Z$. If $h-2$ is odd, take $q+1$ copies of a partition of the points of $Z$ outside a plane $E'$ into lines and add the $[3]_q$ lines of $E'$.
  Both partitions exist by repeated application of Lemma~\ref{lemma_moves}.(a) with $h=2$.) For each line $\langle z_1,z_2\rangle$ of $\cK$ and each $a\in A$ let $\pi:=\langle p,a+z_1,\omega a+z_2\rangle$. This gives
  $q^2[h-2]_q$ planes $\pi$, each used with the point $P$. Take each of these pairs at most $\lceil N/\left(q^2[h-2]_q\right)\rceil$ times.

  The projection of $\pi$ to $Z$ along $E$ is $\langle z_1,z_2\rangle$. Hence $\pi\cap E=P$. Now let $R=\langle x'+z\rangle$ be a point of $B^*$ with $x'\in E$ and $0\neq z\in Z$. If $x'+z\in\pi$, then
  $x'+z=\mu p+\alpha(a+z_1)+\beta(\omega a+z_2)$ with $\mu,\alpha,\beta\in\F_q$, i.e., $z=\alpha z_1+\beta z_2$ and $x'=\mu p+(\alpha+\beta\omega)a$. Given a line of $\cK$ through $\langle z\rangle$, the pair
  $(\alpha,\beta)\neq(0,0)$ is determined by $z$. So $\alpha+\beta\omega\neq 0$, and there is exactly one $a\in A$ with $R\le\pi$. Thus, $R$ lies in exactly $q+1$ of the planes, and no point of $E$
  other than $P$ lies in one of them. So $t(H)\le(q+1)\lceil N/\left(q^2[h-2]_q\right)\rceil\le m_0$.
\end{proof}

\subsection{The semigroups \texorpdfstring{$\mathcal{A}_q^{h}(h+j)$}{A\textunderscore{}q\textasciicircum{}(h)(h+j)}}\label{app_sec4_apery}
\begin{proof}[Proof of Theorem~\ref{thm_apery_qhj}]
  We use induction on $j$, where the case $j=0$, i.e., $l=h$, was stated above.
  Let $r=l+h-1=2h+j-1$ and $s_i:=s_q(r,0,i)=q^i[2h+j-i]_q$ for $0\le i\le r$. We have $[h]_qq^l=s_{h+j}\in\mathcal{L}_q^r$. Hence
  $n+q^l$ satisfies (N1) whenever $n$ does.

  If $\rho=q\rho'$, then Lemma~\ref{lemma_divide_q} and the induction hypothesis for $j-1$ give Equation~(\ref{eq_apery_qhj}), since $\kappa_{l-1}(\rho')=\kappa_l(\rho)$, $v_q(\rho')=v_q(\rho)-1$, and
  $\lfloor\rho'/q^j\rfloor=\lfloor\rho/q^{j+1}\rfloor$. So, let $q\nmid\rho$ and write $\rho=u+wq^j$ with $1\le u<q^j$, $u=\sum_{i=0}^{j-1}d_iq^i$, $0\le d_i\le q-1$, and $w=Fq+d$ with $0\le d\le q-1$. Then
  $F=\lfloor\rho/q^{j+1}\rfloor$ and $1\le d_0\le q-1$. The following elements lie in $\mathcal{A}_q^{h}(l)$:
  \begin{itemize}
    \item $g_{\mathrm{sp}}:=q^l+q^j$, which is the length of $q^j$ copies of a Desarguesian $(h-1)$-spread of $\PG(2h-1,q)$,
    \item $g_{\mathrm{pg}}:=q^{2h}+q^h+1$, which is the length of the field reduction of all points of $\PG(2,q^h)$ (all weights are $q^{2h}$ and $h\ge j$),
    \item the elements of $q\cdot\mathcal{A}_q^{h}(l-1)$ (Lemma~\ref{lemma_divide_q}),
    \item $g_f:=u+fq^{j+1}+(jq+1)q^l$ for $0\le f\le j$. For these, (N2) holds since $\kappa_l\left(u+fq^{j+1}\right)\le fq+1$, and (N1) follows from the following identity, applied to the residue $u+fq^{j+1}$, that is, with $d=0$ and $F=f$.
  \end{itemize}
  For $0\le F\le j$ and $n:=\rho+(jq+1+d)q^l$ we have
  \begin{equation}
    \label{eq_apery_identity}
    [h]_qn=\sum_{i=0}^{j-1}d_is_i+ds_j+Fs_{j+1}+(q-d_0)s_h+\sum_{i=1}^{j-1}(q-1-d_i)s_{h+i}+(j-F)s_{h+j+1}.
  \end{equation}
  To verify it, note that $s_i-s_{h+i}=q^i[h]_q$ for $0\le i\le j+1$, which follows from $[2h+j-i]_q=q^h[h+j-i]_q+[h]_q$. In particular, $ds_j=d[h]_q\left(q^j+q^l\right)$. Moreover,
  $qs_h+(q-1)\sum_{i=1}^{j-1}s_{h+i}+js_{h+j+1}=(jq+1)q^l[h]_q$, i.e., $q[h+j]_q+(q-1)\sum_{i=1}^{j-1}q^i[h+j-i]_q+jq^{j+1}[h-1]_q=(jq+1)q^j[h]_q$, which follows from $(q-1)q^i[h+j-i]_q=q^{h+j}-q^i$.
  The indices $0,\dots,j+1,h,\dots,h+j-1,h+j+1$ are at most $r$ and pairwise distinct, except $j+1=h$ for $(j,h)=(1,2)$.

  \textit{Claim:} If $F\le j-1$, then $[h]_q\left(n-q^l\right)\notin\mathcal{L}_q^r$. Let $x:=[h]_q\left(n-q^l\right)=[h]_qn-s_{h+j}$. Since $F\le q-1$, and $F=0$ for $(j,h)=(1,2)$, the coefficients $c_i$ of $s_i$, $i<h+j$, in Equation~(\ref{eq_apery_identity}) lie in
  $\{0,\dots,q-1\}$ and $x-\sum_{i<h+j}c_is_i=q^{h+j}t$ with $t=(j-F)q[h-1]_q-[h]_q$. Since $s_{h+j+i}=q^{h+j}s_q(h-1,0,i)$, the uniqueness of the $S_q(r,0)$-adic expansion shows that $x\in\mathcal{L}_q^r$ if and only if
  $t\in\mathcal{L}_q^{h-1}$. Every element of $\mathcal{L}_q^{h-1}$ is non-negative, and if it is congruent to $c'\in\{0,\dots,q-1\}$ modulo $q$, then it is at least $c'[h]_q$, since all other elements of
  $S_q(h-1,0)$ are divisible by $q$. Now $t=(j-F-1)q[h-1]_q-1$. For $F=j-1$ we have $t<0$. For $F\le j-2$ we have $1\le j-F-1\le q-1$. Hence $t\equiv q-1\pmod q$ and
  $t<(q-1)q[h-1]_q+q-1=(q-1)[h]_q$. So, $t\notin\mathcal{L}_q^{h-1}$, which proves the claim. Hence, (N1) fails for $\rho+kq^l$ whenever $k\le jq+d$.

  If $\rho<uq^h$, then $0<u/q^j-\rho/q^l<1$ and $\kappa_l(\rho)=\left\lfloor w+u/q^j-\rho/q^l\right\rfloor+1=w+1$.

  (a) If $F\le j-1$, then $\rho<jq^{j+1}\le uq^h$ (since $j\le q$ and $h\ge j+2$, or $j=1$). So the right hand side of Equation~(\ref{eq_apery_qhj}) is $\rho+(w+1+q(j-F))q^l=\rho+(jq+1+d)q^l=g_F+d\cdot g_{\mathrm{sp}}$.
  By the claim it is equal to $m_\rho^{(l)}$.

  (b) If $F\ge j$, then the right hand side of Equation~(\ref{eq_apery_qhj}) is $\rho+\kappa_l(\rho)q^l$, which is a lower bound by (N2). It suffices to show that this number lies in $\mathcal{A}_q^{h}(l)$.
  If $\rho<uq^h$, then $w\ge jq$ and $g_j+(w-jq)\cdot g_{\mathrm{sp}}=\rho+(w+1)q^l=\rho+\kappa_l(\rho)q^l$. So, let $\rho\ge uq^h$.
  Since $\rho-d_0q^h\ge(u-d_0)q^h\ge 0$ and $\rho-d_0q^h\equiv d_0\not\equiv 0\pmod q$, the number $\sigma':=\left(\rho-d_0\left(q^h+1\right)\right)/q$ is a non-negative integer with $\sigma'<q^{l-1}$. Suppose
  $\sigma'\ge 1$ and $j-1-v_q(\sigma')-\lfloor\sigma'/q^j\rfloor>0$. Then $v_q(\sigma')\le j-2$ and $\sigma'<(j-1)q^j$. Hence $u=d_0+\left(q\sigma'\bmod q^j\right)\ge d_0+q$ and
  $uq^h\ge d_0q^h+q^{h+1}>d_0q^h+d_0+q\sigma'=\rho$, since $d_0+q\sigma'<q+(q-1)q^{j+1}<q^{h+1}$. This contradicts $\rho\ge uq^h$. So, by induction,
  $q\cdot m_{\sigma'}^{(l-1)}=q\sigma'+\kappa_{l-1}(\sigma')q^l$ for $\sigma'\ge 1$. Since $q^{2h}=q^{h-j}q^l$ we obtain that $d_0\cdot g_{\mathrm{pg}}+q\cdot m_{\sigma'}^{(l-1)}\in\mathcal{A}_q^{h}(l)$ is congruent to $\rho$ and equals
  $\rho+\left(d_0q^{h-j}+\kappa_{l-1}(\sigma')\right)q^l$, where we set $\kappa_{l-1}(0):=0$ and omit the second summand for $\sigma'=0$. Finally,
  $\rho\left(q^h-1\right)/q^l=d_0q^{h-j}-d_0/q^l+\sigma'\left(q^h-1\right)/q^{l-1}$, and the fractional part of $\sigma'\left(q^h-1\right)/q^{l-1}$ is either $0$, which happens if and only if $\sigma'=0$, or at least $1/q^{l-1}>d_0/q^l$.
  So $\kappa_l(\rho)=d_0q^{h-j}+\kappa_{l-1}(\sigma')$.
\end{proof}

\begin{proof}[Proof of Corollary~\ref{cor_gens_qhj}]
  Let $l=h+j$. In the proof of Theorem~\ref{thm_apery_qhj} every $m_\rho^{(l)}$ with $q\nmid\rho$ is written as a sum of elements $g_f$, $g_{\mathrm{sp}}$, $g_{\mathrm{pg}}$, and an element of $q\cdot\mathcal{A}_q^{h}(l-1)$, and $m_{q\rho'}^{(l)}\in q\cdot\mathcal{A}_q^{h}(l-1)$.
  By induction on $j$, where $\mathcal{A}_q^{h}(h)=\left\langle q^h,q^h+1\right\rangle$ is generated by the (not minimal) set $\mathcal{G}_q^{h,0}=\left\{q^h,q^h+1,q^{2h}+q^h+1\right\}$, the semigroup $q\cdot\mathcal{A}_q^{h}(l-1)$ is generated by $q\cdot\mathcal{G}_q^{h,j-1}$, whose elements are $q^l$, $q^l+q^j$, $q\,g_{\mathrm{pg}}=g_{\mathrm{pg}}+\dots+g_{\mathrm{pg}}$, and
  the elements of $\mathcal{G}_q^{h,j}$ with $i\ge 1$. So, $\mathcal{G}_q^{h,j}$ generates $\mathcal{A}_q^{h}(l)$. It remains to show that no $x\in\mathcal{G}_q^{h,j}$ is the sum of two positive elements $y=\rho_y+k_yq^l$ and
  $z=\rho_z+k_zq^l$ of $\mathcal{A}_q^{h}(l)$, where $0\le\rho_y,\rho_z<q^l$ and $k_y,k_z\ge 0$. Since $q^l$ is the smallest positive element, this is clear for $q^l$ and $q^l+q^j<2q^l$. Moreover $k_y,k_z\ge 1$, since a positive element with $k=0$ violates (N2).
  Since $y>\rho_yq^h$ and $z>\rho_zq^h$ by (N2), we have $y+z>q^{l+h}>x$ if $\rho_y+\rho_z\ge q^l$, and $y+z-(\rho_y+\rho_z)q^h\ge 2>1=g_{\mathrm{pg}}-\left(q^h+1\right)q^h$ otherwise, where $g_{\mathrm{pg}}\equiv q^h+1\pmod{q^l}$.
  So, $x\neq g_{\mathrm{pg}}$, $\rho_y+\rho_z=\rho_x$, and $k_y+k_z=k_x$, where $x=\rho_x+k_xq^l$ with $\rho_x=q^iu+fq^{j+1}<(j+1)q^{j+1}$ and $k_x=(j-i)q+1$.
  If one of $\rho_y,\rho_z$, say $\rho_y$, is not divisible by $q$, then $\rho_y\le\rho_x<(j+1)q^{j+1}\le q^h$ (since $j+1\le q$ or $h\ge j+3$). In the notation of the proof
  of Theorem~\ref{thm_apery_qhj}, $F\le j$, $\kappa_l(\rho_y)=w+1$, and Equation~(\ref{eq_apery_qhj}) gives $k_y\ge jq+1+d\ge k_x$, while $k_z\ge 1$, which is impossible. Otherwise $q$ divides $\rho_x$, i.e., $i\ge 1$, and $x/q=y/q+z/q$ contradicts
  the induction hypothesis for $j-1$, since $x/q\in\mathcal{G}_q^{h,j-1}$ and $y/q,z/q\in\mathcal{A}_q^{h}(l-1)$ by Lemma~\ref{lemma_divide_q}. For the number of elements note that there are $(q-1)q^{t-1}$ choices for $u$ if $t=j-i$.
\end{proof}

\subsection{Constructions for \texorpdfstring{$\Delta=q^{h+2}$}{Delta=q\textasciicircum{}(h+2)}}\label{app_sec4_qh2}
This subsection contains the proof of Lemma~\ref{lemma_cone_family} and of Corollary~\ref{cor_qh3} and further constructions for $\Delta=q^{h+2}$.
\begin{proof}[Proof of Lemma~\ref{lemma_cone_family}]
  The coefficient of $[i]$ in $\cT_j$ is divisible by $q^{2h+j-i}$ for every $i$. So, $\cT_j$ is $q^{2h+j-1}$-divisible by Lemma~\ref{lemma_type_divisible} with $r=2h+j$, and the corresponding multiset of $(h-1)$-spaces
  is $q^{h+j}$-divisible by Lemma~\ref{lemma_exclusion}. For a type $\sum_i\varepsilon_i[i]$ we have $(q-1)\sum_i\varepsilon_i[i]_q=\sum_i\varepsilon_iq^i-\sum_i\varepsilon_i$. For $\cT_j$ we have $\sum_i\varepsilon_iq^i=(jq+1)q^{2h+j}$
  and, since $(q-c_0)q^h+\sum_{k=1}^{j-1}(q-1-c_k)q^{h+k}=q^h\left(q^j-\rho\right)$, $\sum_i\varepsilon_i=(jq+1)q^{h+j}-\left(\rho+gq^{j+1}\right)\left(q^h-1\right)$. Since $q^h-1=(q-1)[h]_q$, this gives
  $\#\cT_j=[h]_q\left((jq+1)q^{h+j}+\rho+gq^{j+1}\right)$.

  \emph{Overview.} Most of $\cT_j$ is produced by cones with vertex $A$. The $c_kq^k$ copies of all $(h-1)$-spaces of $S_{2h+j-k}$ through $A$ give $c_kq^k\chi_{S_{2h+j-k}}$ plus a multiple of $\chi_A$, and the telescoping families $\mathcal{D}_i$ of step~(4) below produce the terms $\chi_{S_{h+j-k}}$ for $k\le j-2$. These cones put too much weight on $A$. The following units replace cones and put less weight on $A$.
  \begin{itemize}
    \item $\mu_k$ lifted MRD partitions of $S_{2h+j-k}\backslash S_{h+j-k}$ (step~(1)),
    \item $\sigma$ spread units in $F$ (step~(2)),
    \item $\tau$ resp.\ $\nu$ units of Lemma~\ref{lemma_moves_general}.(g), which replace $q^2$ cones in $D$ resp.\ $q^{h-1}$ copies of $\mathcal{D}_{h+2}$ (steps~(3) and~(5)).
  \end{itemize}
  Finally $y_{AB}$ copies of the $(h-1)$-spaces between $A$ and $B$ and $y_B$ copies of the hyperplanes of $B$ (step~(6)) balance the coefficients of $\chi_B$ and $\chi_A$. The definitions of $y_{AB}$ and $y_B$ are exactly these two balance equations.

  For $i\ge h+1$ let $\mathcal{D}_i$ be the set of the $q^{i-h}$ $(h-1)$-spaces $X$ with $A<X\le S_i$ and $X\not\le S_{i-1}$, which has the point multiset $\chi_{S_i}-\chi_{S_{i-1}}+q^{i-h}\chi_A$ by
  Lemma~\ref{lemma_moves_general}.(e) with $(a,m)=(h-1,i-1)$. Recall that the $(h-1)$-spaces of $S_i$ containing $A$ have the point multiset $\chi_{S_i}+q[i-h]_q\chi_A$ (Lemma~\ref{lemma_moves}.(b)).
  Let $\rho_k:=\sum_{i=0}^kc_iq^i$ for $k\ge -1$, so that $\rho_{-1}=0$ and $\rho_{j-1}=\rho$. For $0\le k\le j-2$ let $U_k:=q^h\left(q^{k+1}-\rho_k\right)+\sum_{i=0}^k\mu_i$, and let $U_{-1}:=0$. Take
  \begin{enumerate}
    \item[(1)] for $0\le k\le j-1$, $\mu_k$ copies of a lifted MRD partition of the points of $S_{2h+j-k}\backslash S_{h+j-k}$ (Lemma~\ref{lemma_moves}.(a)) and $c_kq^k-\mu_k$ copies, for $k=j-1$ only
               $c_kq^k-\mu_k-\sigma$ copies, of each $(h-1)$-space of $S_{2h+j-k}$ containing $A$,
    \item[(2)] $\sigma$ copies of the union of a Desarguesian $h$-spread of $S_{2h}$ and the $q^{h+1}$ $(h-1)$-spaces given by Lemma~\ref{lemma_moves_general}.(e) for
               $(a,m,i)=(h-1,2h,2h+1)$, which has the point multiset $\chi_F+q^{h+1}\chi_A$,
    \item[(3)] $\tau$ copies of the $h$-partition of $q^2\chi_D+\left([h]_q+q\right)\chi_B$ given by Lemma~\ref{lemma_moves_general}.(g) with $(k,u,m_0)=\left(q^2,h-2,q+1\right)$,
               and $gq^{j+1}-\tau q^2$ copies of each $(h-1)$-space of $D$ containing $A$,
    \item[(4)] $U_k$ copies of $\mathcal{D}_{h+j-k}$ for $0\le k\le j-3$, and, for $j\ge 2$, $U_{j-2}-\nu q^{h-1}$ copies of $\mathcal{D}_{h+2}$,
    \item[(5)] $\nu$ copies of the $h$-partition of $q^{h-1}\chi_E+q[h-2]_q\chi_B$ given by Lemma~\ref{lemma_moves_general}.(g) with $(k,u,m_0)=\left(q^{h-1},1,q\right)$,
    \item[(6)] $y_{AB}$ copies of each of the $q+1$ $(h-1)$-spaces between $A$ and $B$, which have the point multiset $\chi_B+q\chi_A$, and $y_B$ copies of each hyperplane of $B$.
  \end{enumerate}
  In (3) part~(g) applies, since $N=q^2[h-2]_q$. Hence $m_0=q+1$ satisfies the condition for $h\ge 4$, and for $h=3$ we have $N=q^{h-1}$ and $m_0\ge q$. Moreover, $(q+1)[h]_q-q^2[h-1]_q=[h]_q+q$.
  In (5) it applies, since $q\left\lceil q^{h-1}/q^{h-1}\right\rceil=q$ and $q[h]_q-q^{h-1}[2]_q=q[h-2]_q$. All numbers of copies are non-negative, where $U_{j-2}\ge \nu q^{h-1}$ follows from the bound on $\nu$.

  We compare the point multiset with $\cT_j$ term by term, where terms belonging to coinciding subspaces are added. The coefficients of $\chi_{S_{2h+j-k}}$ for $0\le k\le j-1$ and of $\chi_D$ are $c_kq^k$ and $gq^{j+1}$.
  For $0\le k\le j-2$ the coefficient of $\chi_{S_{h+j-k}}$ is $U_k-U_{k-1}-\mu_k$, which is $(q-c_0)q^h$ for $k=0$ and $q^h\left(q^{k+1}-q^k-c_kq^k\right)=(q-1-c_k)q^{h+k}$ for $k\ge 1$.
  Since $q^{h-1}+q[h-2]_q=[h]_q-1$, the coefficient of $\chi_B$ is
  $-\mu_{j-1}-U_{j-2}+\left([h]_q+q\right)\tau+\left([h]_q-1\right)\nu+y_{AB}+[h]_qy_B=q^h\left(q^j-\rho\right)-q^h\left(q^{j-1}-\rho_{j-2}\right)=(q-1-c_{j-1})q^{h+j-1}$.
  For $j=1$ we have $U_{-1}=0$, so that the right-hand side is $q^h\left(q-\rho\right)=(q-c_0)q^h$, which is again the coefficient of $\chi_B$ in $\cT_1$ by the case $k=0$ above.
  The coefficient of $\chi_A$ is
  $x:=\sum_{k=0}^{j-1}\left(c_kq^k-\mu_k\right)q[h+j-k]_q+\sigma\left(q^{h+1}-q[h+1]_q\right)+\left(gq^{j+1}-\tau q^2\right)q[h-1]_q+\sum_{k=0}^{j-2}q^{j-k}U_k-\nu q^{h+1}+qy_{AB}$.
  After inserting $U_k$ and $y_{AB}$, increasing $\mu_i$, $\sigma$, $\tau$, $\nu$, or $y_B$ by one changes $x$ by $-q[h+j-i]_q+q[j-i]_q=-q^{j-i+1}[h]_q$ (since $\sum_{k=i}^{j-2}q^{j-k}+q=q[j-i]_q$), $q^{h+1}-q[h+1]_q=-q[h]_q$,
  $-q^3[h-1]_q-q\left([h]_q+q\right)=-q(q+1)[h]_q$, $-q^{h+1}-q\left([h]_q-1\right)=-q^2[h]_q$, and $-q[h]_q$, respectively. If all of them are zero, then
  $x=\sum_kc_kq^{k+1}[h+j-k]_q+gq^{j+1}\left([h]_q-1\right)+\sum_{k=0}^{j-1}q^{h+j-k}\left(q^{k+1}-\rho_k\right)$, and $[h+j-k]_q=[h]_q+q^h[j-k]_q$ and
  $\sum_{k=0}^{j-1}q^{h+j-k}\rho_k=\sum_{i=0}^{j-1}c_iq^{h+i+1}[j-i]_q$ give $x=jq^{h+j+1}+q[h]_q\rho+gq^{j+1}\left([h]_q-1\right)$. So in general
  \begin{eqnarray*}
    x&=&jq^{h+j+1}+gq^{j+1}\left([h]_q-1\right)+q[h]_q\left(\rho-\sum_{i=0}^{j-1}q^{j-i}\mu_i-\sigma-(q+1)\tau-q\nu-y_B\right)\\
    &=&jq^{h+j+1}+gq^{j+1}\left([h]_q-1-q[h]_q\right)\,=\,(j-g)q^{h+j+1}.
  \end{eqnarray*}
\end{proof}
For $j=1$ the lemma contains the constructions of Theorem~\ref{thm_C} (with $g=0$, $\mu_0=\sigma=\tau=0$, where $y_B=c_0$ copies of the hyperplanes of $B$ account for the difference between the numbers of copies) and, for $h\ge 3$, the types of Theorem~\ref{thm_plus_q2} (with $g=1$, $\mu_0=c_0$, $\tau=1$, $\sigma=0$).
For $j=2$ and $(c_0,c_1)=(b,a)$ the types $\cT_2$ with $g=0,1,2$ are the types of Corollary~\ref{cor_qh3}.(a), of Corollary~\ref{cor_qh3}.(b), and of Proposition~\ref{thm_qh2_f2}, and for $j=3$ and $(c_0,c_1,c_2)=(c,d,e)$ the type $\cT_3$ is the type of Corollary~\ref{cor_qh3}.(c).

\medskip
Propositions~\ref{thm_qh2_f2} and~\ref{prop_qh2_packing} below realize the exceptional lengths of the form $2q^{h+3}+q^{h+2}+2q^3+aq+b$ listed in Subsection~\ref{sub_qh23}, except $1297$ and $1299$. For a chamber $S_1<\dots<S_{2h+2}$ of $\F_q^{2h+2}$, both constructions describe the part of the partition
inside $S_{h+2}$, resp.\ $S_{h+1}$, in the dual space.

\begin{Proposition}
  \label{thm_qh2_f2}
  Let $h\in\{3,4\}$, $0\le a\le q-1$, and $1\le b\le q-1$. Then the multiset of points of type
  $$
    b[2h+2]+aq[2h+1]+(q-b)q^h[h+2]+(q-1-a)q^{h+1}[h+1]+2q^3[2h-1]
  $$
  in $V=\F_q^{2h+2}$ is $h$-partitionable. Thus $2q^{h+3}+q^{h+2}+2q^3+aq+b$ is the length of a $q^{h+2}$-divisible faithful additive code over $\F_{q^h}$.
\end{Proposition}
By Remark~\ref{rem_type_weights} the non-zero weights of these codes are $bq^{h+2}$, $(a+b)q^{h+2}$, $(a+b+2)q^{h+2}$ (only for $h=4$), $(q+a+2)q^{h+2}$, and $(2q+1)q^{h+2}$, attained by $1$, $q^2+q$, $q^3$, $q^h$, and $q^{h+1}[h+1]_q$ hyperplanes, respectively, where coinciding values (for $a=0$ or $a=q-1$) are merged.
For $(q,h,a,b)$ equal to $(2,4,1,1)$, $(4,3,3,3)$, and $(3,4,a,b)$ Proposition~\ref{thm_qh2_f2} gives the lengths $339\in\cN_2^{64,4}$, $9359\in\cN_4^{1024,3}$, and $5158,5159,5161,5162,5164,5165\in\cN_3^{729,4}$.
For $h\ge 5$ the dual construction of Proposition~\ref{thm_qh2_f2} imposes much stronger conditions, and for $h\ge 8$ and $b=q-1$ it is impossible, see below. For $h\ge 5$ the following variant is
more flexible. It uses the hyperplanes of $S_{h+1}$ as in the proof of Lemma~\ref{lemma_moves_general}.(g).

\begin{Proposition}
  \label{prop_qh2_packing}
  Let $h\ge 3$, $1\le b\le q-1$, $\vartheta:=2q^2+(q^2-1)(q-b)$, $\beta:=(q-b)q^h+b$, and let $B^*$ be the dual space of $B=\F_q^{h+1}$. Suppose that there are $2q^3$ subspaces $W_1,\dots,W_{2q^3}$ of codimension $2$ of $B^*$,
  $\beta_1$ pairs $(\pi_s,P_s)$ of a $3$-dimensional subspace $\pi_s$ of $B^*$ and a point $P_s$ of $\pi_s$, and $\beta-\beta_1$ subspaces $L_t$ of dimension $2$ of $B^*$ such that
  $$
    \ell(x):=q\cdot\#\left\{j\,:\,x\le W_j\right\}+\#\left\{s\,:\,x\le\pi_s,\,x\neq P_s\right\}+q\cdot\#\left\{t\,:\,x\le L_t\right\}\le\vartheta
  $$
  for every point $x$ of $B^*$. Then the type of Proposition~\ref{thm_qh2_f2}, taken for this $h$, is $h$-partitionable for every $0\le a\le q-1$. Hence $2q^{h+3}+q^{h+2}+2q^3+aq+b$ is the length of a $q^{h+2}$-divisible faithful additive code over $\F_{q^h}$.
\end{Proposition}
The condition of Proposition~\ref{prop_qh2_packing} is a small packing problem in $\PG(h,q)$, which can only be solved for small $h$. Indeed, $\vartheta=q\left(2q+(q-b)q\right)-(q-b)<q\left(2q+(q-b)q\right)$, so every point lies in at most
$2q+(q-b)q-1$ of the subspaces $W_j$, while any $2q+(q-b)q$ subspaces of codimension $2$ of $B^*$ have a common point if $h+1>2\left(2q+(q-b)q\right)$.

For $(q,h,b)=(2,5,1)$ integer linear programming finds $16$ subspaces $W_j$, $18$ pairs $(\pi_s,P_s)$, and $15$ lines $L_t$
of $\PG(5,2)$ with $\ell(x)\le 11$ for all $63$ points $x$. This gives $657,659\in\cN_2^{128,5}$ with hyperplane weights $\{128^7,384^{24},512^{32},640^{4032}\}$ and $\{128^1,256^6,512^{24},640^{4064}\}$, respectively. For $(q,h,b)=(2,6,1)$ no configuration as in Proposition~\ref{prop_qh2_packing} is known, and
the dual approach of Proposition~\ref{thm_qh2_f2} cannot be applied for $(q,h)=(2,6)$ either. Hence $1297$ and $1299$ remain open.

\begin{Lemma}[E.g.~{\cite[Lemmas 4 and 6]{kurz2024additive}}]
  \label{lemma_dual_lines}
  Let $m\ge 3$, let $E=\F_q^m$ be equipped with a non-degenerate symmetric bilinear form, and let $\cK$ be a multiset of $2$-dimensional subspaces of $E$ with point multiset $\lambda=\sum_{L\in\cK}\chi_L$.
  Then every point $Q$ of $E$ is contained in exactly $\left(\lambda(Q^\perp)-\#\cK\right)/q$ elements of the multiset $\cK^\perp:=\left\{L^\perp\,:\,L\in\cK\right\}$ of $(m-2)$-dimensional subspaces of $E$.
\end{Lemma}
\begin{proof}
  We have $Q\le L^\perp$ if and only if $L\le Q^\perp$. Since $Q^\perp$ is a hyperplane of $E$, a line $L\le Q^\perp$ contributes $q+1$ to $\lambda(Q^\perp)$ and every other line of $\cK$ contributes exactly $1$.
\end{proof}

\begin{proof}[Proof of Proposition~\ref{thm_qh2_f2}]
  Let $S_1<\dots<S_{2h+2}=V$ be a chamber, $B:=S_{h+1}$, $E:=S_{h+2}$, $D:=S_{2h-1}$, and $F:=S_{2h+1}$. For $h=3$ we have $D=E$. Set $\beta:=(q-b)q^h+b$ and, for a statement $X$, let $\delta_X$ be $1$ if $X$ holds and $0$ otherwise.
  The type is $q^{2h+1}$-divisible by Lemma~\ref{lemma_type_divisible} with $r=2h+2$, since $q^3$ divides the coefficient $2q^3$ of $[2h-1]$. So the corresponding multiset of $(h-1)$-spaces is $q^{h+2}$-divisible by Lemma~\ref{lemma_exclusion}.
  Since $[2h-1]_q=[h]_q+q^h[h-1]_q$, the cardinality of the type exceeds the one of the type $\cT_2$ with $g=0$ of Lemma~\ref{lemma_cone_family} by $2q^3[2h-1]_q-2q^{h+3}[h-1]_q=2q^3[h]_q$, which gives the stated length. It remains to give the partition. Step~(1) reduces the task to $\cT'$, step~(2) to subspaces of $E$, and step~(3) dualizes in $E$, where the members become lines and Lemma~\ref{lemma_dual_lines} turns the prescribed point multiset into a condition on the point multiset of these lines.
  \begin{enumerate}
    \item[(1)] Take $b$ copies of a lifted MRD partition of the points of $V\backslash E$, $aq$ copies of a lifted MRD partition of the points of $F\backslash B$ (Lemma~\ref{lemma_moves}.(a) with $i=2h+2$ and $i=2h+1$), and $(q-1-a)(q-1)q$
               copies of each hyperplane of $B$. Since every point of $B$ lies in $[h]_q$ hyperplanes of $B$ and $(q-1)[h]_q=q^h-1$, the point multiset of these subspaces is $b\left(\chi_V-\chi_E\right)+aq\left(\chi_F-\chi_B\right)+(q-1-a)\left(q^{h+1}-q\right)\chi_B$.
               So, it remains to partition $\cT':=2q^3\chi_D+\beta\chi_E+(q-1)q\chi_B$.
    \item[(2)] Let $h=4$, so that $\dim(E)=6$ and $\dim(D)=7$. Let $T_0,\dots,T_{q^3}$ be a plane spread of $E$, i.e., $3$-dimensional subspaces that partition the points of $E$. Counting points shows that every hyperplane of $E$ contains
               exactly one of them, and we may assume $T_0\le B$. For $1\le j\le q^3$ Lemma~\ref{lemma_moves_general}.(e), applied to a chamber through $T_j<E<D$ with $(a,m,i)=(3,6,7)$, gives $q^3$ solids through $T_j$ with point multiset
               $\chi_D-\chi_E+q^3\chi_{T_j}$. Taking all of them twice gives $2q^3\left(\chi_D-\chi_E\right)+2q^3\left(\chi_E-\chi_{T_0}\right)$. It remains to partition $\cT'':=\beta\chi_E+(q-1)q\chi_B+2q^3\chi_{T_0}$ into
               $4$-dimensional subspaces of $E$. For $h=3$ let $\cT'':=\cT'=\left(2q^3+\beta\right)\chi_E+(q-1)q\chi_B$, which has to be partitioned into $3$-dimensional subspaces of $E=D$.
    \item[(3)] In both cases the members are the subspaces of codimension $2$ of $E$, and we apply Lemma~\ref{lemma_dual_lines} with $m=h+2$. Let $P:=B^\perp$, which is a point of $E$, and let $R$ be a $3$-dimensional subspace
               of $E$ containing $P$, where we choose $R:=T_0^\perp$ for $h=4$ (note that $T_0\le B$ implies $P\le T_0^\perp$). The $q+1$ lines of $R$ through $P$ have point multiset $\chi_R+q\chi_P$ and the $q^2$ other lines of $R$
               have point multiset $q\left(\chi_R-\chi_P\right)$. For a point $Q$ of $E$ we have $P\le Q^\perp$ if and only if $Q\le B$, and $R\le Q^\perp$ if and only if $Q\le R^\perp$.

               For $h=3$ let $\vartheta:=2q^2+(q^2-1)(q-b)$, $y:=q^2+(q-1)(q-b)$, and $z:=q+(q-1)(q-b)$. Let $\cK$ consist of $\vartheta$ copies of a partition of the points of $E\backslash R$ into $q^3$ lines
               (Lemma~\ref{lemma_moves}.(a) with $h=2$ and $i=5$), $y$ copies of each line of $R$ through $P$, and $z$ copies of each other line of $R$. Since $y+qz=\vartheta$ and $y-z=q^2-q$, the point multiset of $\cK$ is
               $\lambda=\vartheta\chi_E+\left(q^3-q^2\right)\chi_P$, and $\#\cK=\vartheta q^3+(q+1)y+q^2z$. By Lemma~\ref{lemma_dual_lines} a point $Q$ of $E$ is contained in exactly $\left(\vartheta[4]_q-\#\cK\right)/q+\left(q^2-q\right)\delta_{Q\le B}$
               elements of $\cK^\perp$. From $\vartheta[4]_q-\#\cK=\vartheta[3]_q-(q+1)y-q^2z=2q^4+q^2+(q-b)\left(q^4-q\right)$ we obtain $\left(\vartheta[4]_q-\#\cK\right)/q=2q^3+\beta$. Hence $\cK^\perp$ is a partition of $\cT''$.

               For $h=4$ let $\cK$ consist of $(q^2-1)(q-b)$ copies of a line spread of $E$, $q^2$ copies of each line of $R$ through $P$, and $q$ copies of each other line of $R$. Then
               $\lambda=(q^2-1)(q-b)\chi_E+2q^2\chi_R+\left(q^3-q^2\right)\chi_P$ and $\#\cK=(q^2-1)(q-b)\left(q^4+q^2+1\right)+2q^3+q^2$. Since
               $\lambda(Q^\perp)=(q^2-1)(q-b)[5]_q+2q^2\left([2]_q+q^2\delta_{Q\le T_0}\right)+\left(q^3-q^2\right)\delta_{Q\le B}$ and $[5]_q-\left(q^4+q^2+1\right)=q^3+q$, Lemma~\ref{lemma_dual_lines} shows that $Q$ is contained in exactly
               $(q-b)\left(q^4-1\right)+q+2q^3\delta_{Q\le T_0}+\left(q^2-q\right)\delta_{Q\le B}$ elements of $\cK^\perp$. As $(q-b)\left(q^4-1\right)+q=\beta$, also here $\cK^\perp$ is a partition of $\cT''$.
  \end{enumerate}
\end{proof}

For $h\ge 5$ and $b=q-1$ the analogous construction would use cones over $2q^3$ subspaces $T_1,\dots,T_{2q^3}$ of dimension $3$ of $E$ in step~(2). By Lemma~\ref{lemma_dual_lines} the dual lines would then have the point multiset $\left(3q^2-1\right)\chi_E+\left(q^3-q^2\right)\chi_{B^\perp}-q^2\sum_k\chi_{T_k^\perp}$. Hence the subspaces $T_k^\perp$ cover every point of $E$ other than $B^\perp$ at most twice and $B^\perp$ at most $q+1$ times. This is impossible for $h\ge 8$, since any three subspaces of codimension $3$ of $\F_q^{h+2}$ have a common point, which then has to be $B^\perp$ for all triples, while $B^\perp$ may be covered at most $q+1<2q^3$ times.
This is the reason for Proposition~\ref{prop_qh2_packing}.

\begin{proof}[Proof of Proposition~\ref{prop_qh2_packing}]
  We use the notation of the proof of Proposition~\ref{thm_qh2_f2} and its step~(1). It remains to partition $\cT'=2q^3\chi_D+\beta\chi_E+(q-1)q\chi_B$. Write $D=B\oplus U$ with $\dim(U)=h-2$ and choose $u_0\in U$ with $E=B+\langle u_0\rangle$.
  \begin{enumerate}
    \item[(2')] For each $j$ the subspace $W_j^\perp\le B$ has dimension $2$, and Lemma~\ref{lemma_moves_general}.(e) with $(a,m,i)=(2,h+1,2h-1)$ gives $q^{h-1}$ $(h-1)$-spaces through $W_j^\perp$ with point multiset
                $\chi_D-\chi_B+q^{h-1}\chi_{W_j^\perp}$.
    \item[(3')] For each $s$ take the $q^2$ subspaces $\mathcal{Q}\left(u_0,\pi_s,P_s\right)$ and for each $t$ take the $q^2$ subspaces $\cP\left(u_0,L_t^\perp\right)$, both from the proof of
                Lemma~\ref{lemma_moves_general}.(g). Each of these $\beta$ units covers every point of $E\backslash B$ exactly once and no point of $D\backslash E$.
    \item[(4')] Take each hyperplane $H$ of $B$ exactly $\vartheta-\ell\left(H^\perp\right)$ times.
  \end{enumerate}
  Every point of $D\backslash E$ is covered $2q^3$ times and every point of $E\backslash B$ is covered $2q^3+\beta$ times. Let $Q$ be a point of $B$. The hyperplanes of $B$ through $Q$ correspond to the points of the hyperplane $Q^\perp$ of $B^*$.
  So step~(4') covers $Q$ exactly $\vartheta[h]_q-\sum_{x\le Q^\perp}\ell(x)$ times. Since $\#\left\{x\le W_j\cap Q^\perp\right\}=[h-2]_q+q^{h-2}\delta_{Q\le W_j^\perp}$, step~(2') and the contribution of $W_j$ to $\ell$
  together contribute $-q[h-2]_q$ for each $j$. As in the proof of Lemma~\ref{lemma_moves_general}.(g), each unit of step~(3') together with its contribution to $\ell$ contributes $-q$. So, $Q$ is covered exactly
  $\vartheta[h]_q-2q^4[h-2]_q-q\beta$ times. Using $\vartheta[h]_q-(q+1)\beta=2q^2[h]_q-q^2-q$ and $2q^2[h]_q=2q^2+2q^3+2q^4[h-2]_q$ this number equals $2q^3+\beta+(q-1)q=\cT'(Q)$.
\end{proof}

The dual approach of Proposition~\ref{thm_qh2_f2} in $S_{h+2}$ cannot be applied for $(q,h)=(2,6)$ either. Indeed, $16$ subspaces of codimension $3$ of $\F_2^8$
pairwise meet in at least $3$ points. Hence $\sum_x\binom{c(x)}{2}\ge 360$ for the number $c(x)$ of these subspaces containing a point $x$, while $c(x)\le 2$ for $x\neq B^\perp$ and $c(B^\perp)\le 3$ only allow $\sum_x\binom{c(x)}{2}\le 257$.

\subsection{Constructions for \texorpdfstring{$\Delta=q^{h+3}$}{Delta=q\textasciicircum{}(h+3)}}\label{app_sec4_qh3}
\begin{proof}[Proof of Corollary~\ref{cor_qh3}]
  (a) Write $\rho=\sum_{i=0}^{j-1}c_iq^i$ and apply Lemma~\ref{lemma_cone_family} with $g=0$, $\sigma=\tau=\nu=0$, $\mu_k=0$ for $k\le j-2$, and
  $\mu_{j-1}:=\lfloor\rho/q\rfloor$ if $c_{j-1}\ge 1$ and $\mu_{j-1}:=0$ otherwise. In the first case $\mu_{j-1}<q^{j-1}\le c_{j-1}q^{j-1}$, $y_B=c_0$, and
  $y_{AB}=q^h\left(q^j-\rho\right)+\mu_{j-1}-[h]_qc_0\ge q^h-(q-1)[h]_q=1$. In the second case $j\ge 2$, $\rho<q^{j-1}$, $y_B=\rho$, and
  $y_{AB}=q^h\left(q^j-\rho\right)-[h]_q\rho>q^{j-1}\left((q-1)q^h-[h]_q\right)>0$.

  (b)  Let $z:=\min\{aq,(b-1)(q-1)\}$. We apply Lemma~\ref{lemma_cone_family} with $j=2$, $g=1$, $(c_0,c_1)=(b,a)$, $\mu_0=b$, $\mu_1=aq-z$, $\sigma=\nu=0$, and $\tau=1$, which is possible since $h\ge 3$.
  Then $\cT_2$ is the stated type, and the lemma gives the first length. Taking the resulting code $q$ times gives the second length. It remains to show $y_{AB},y_B\ge 0$. We have
  $y_B=(q-1)\left((q-a-b)q+q-b\right)-1+qz$ and, using $q^h-1=(q-1)[h]_q$ and $[h+1]_q=q[h]_q+1$, $y_{AB}=(b-1)q^{h+1}+q^2-qb-z[h+1]_q$. Since $(q-1)[h+1]_q=q^{h+1}-1$ and $z\le (b-1)(q-1)$, we have
  $z[h+1]_q\le (b-1)q^{h+1}-(b-1)$. Hence $y_{AB}\ge q^2-qb+b-1=(q-1)(q+1-b)\ge 0$. If $z=(b-1)(q-1)$, then $y_B=(q-1)\left(q^2-aq-b\right)-1\ge (q-1)(q-b)-1\ge 0$. Otherwise $z=aq$ and
  $y_B=\left(q^2-1\right)(q-b)+aq-1\ge 0$.

  (c) Let $j=3$ and $\rho=c+dq+eq^2$. Since $q^h-1=(q-1)[h]_q$, we have $y_{AB}=[h]_qZ+W$ with
  $Z:=(q-1)\left(q^3-\rho\right)-y_B-\tau-\nu$ and $W:=q^3-\rho+\mu_0+\mu_1+\mu_2+\nu-q\tau$. Note that $W\ge q^3-\rho+\mu_0-q\tau$, and that $y_{AB}\ge 0$ holds if $Z\ge 1$ and $W\ge -[h]_q$, or if $Z=0$ and $W\ge 0$.
  For $g=1$ let $\mu_0:=c$, $\tau:=c$, $\sigma:=0$, $T:=q^3-cq^2+eq+d-c$, and
  $$
    \mu_2:=\min\left\{eq^2,T\right\},\quad \mu_1:=\min\left\{dq,\left\lfloor\left(T-\mu_2\right)/q\right\rfloor\right\},\quad \nu:=T-\mu_2-q\mu_1.
  $$
  Since $T\ge q^3-q+1-(q-1)q^2>0$, these numbers are non-negative, and a direct computation gives $y_B=0$. So $Z=(q-1)\left(q^3-\rho\right)-c-\nu$ and $W\ge q^3-\rho-(q-1)c\ge 1-(q-1)^2>-[h]_q$.

  If $\nu\ge q$, then $\mu_2=eq^2$ and $\mu_1=dq$, and hence $\nu=T-q^2(d+e)$. In this case $(q-1-d)q^2+(q-c)q-\nu=c\left(q^2-q+1\right)+e\left(q^2-q\right)-d>0$ and
  $$
    Z=q^4-2q^3+c\left(q^2-q+1\right)+d(q-1)-eq(q-1)^2\ge q^4-2q^3-q(q-1)^3+q^2-q+1=q^3-2q^2+1\ge 1.
  $$

  If $\nu\le q-1$, then $\nu\le (q-c)q$ and $Z\ge (q-1)\left(q^3-\rho-1\right)-c$, which is at least $q-1$ if $\rho\le q^3-3$. If $\rho\ge q^3-2$, then $d=e=q-1$ and $c\ge q-2$. For $q\ge 4$ we then have
  $T\le 3q^2-q+1\le eq^2$. Hence $\nu=0$ and $Z=(q-1)\left(q^3-\rho\right)-c$, which equals $q$ for $c=q-2$ and $0$ for $c=q-1$. In the latter case $\mu_2=2q^2-q$ and $W=q^2+q$.
  The three remaining cases $(q,c,d,e)\in\{(2,1,1,1),(3,1,2,2),(3,2,2,2)\}$ give $(Z,W)=(0,5)$, $(2,21)$, and $(0,12)$.

  For $g=2$ the conditions of Lemma~\ref{lemma_cone_family} depend on $h$ only via $y_{AB}=[h]_qZ+W$. A computer search over the admissible $\mu_0,\mu_1,\mu_2,\sigma,\tau,\nu$ finds for each of the stated
  $(q,c,d,e)$ parameters with $y_B\ge 0$, $Z\ge 1$, and $[4]_qZ+W\ge 0$. Hence $y_{AB}\ge 0$ for all $h\ge 4$.
\end{proof}

\begin{proof}[Proof of Proposition~\ref{thm_qh3_q2}]
  The type is $2^{2h+2}$-divisible by Lemma~\ref{lemma_type_divisible} with $r=2h+3$. So the corresponding multiset of $(h-1)$-spaces is $2^{h+3}$-divisible by Lemma~\ref{lemma_exclusion}. With $[i]_2=2^i-1$ the number of points of the type is $5\cdot 2^{2h+3}-(39-2d-4e)2^h-(1+2d+4e)=\left(2^h-1\right)\left(5\cdot 2^{h+3}+1+2d+4e\right)$.

  Let $S_1<\dots<S_{2h+3}=V$ be a chamber, $A':=S_{h-2}$, $B:=S_{h+1}$, $E:=S_{h+2}$, $F:=S_{2h+1}$, and $F':=S_{2h+2}$. For $h+2\le i\le 2h+3$ Lemma~\ref{lemma_moves_general}.(e) with $(a,m)=(h-2,i-2)$ gives $2^{i-h}$
  $(h-1)$-spaces containing $A'$ with point multiset $\chi_{S_i}-\chi_{S_{i-2}}+2^{i-h}\chi_{A'}$. For $2\le k\le h+3$ let $\varepsilon_k\in\{0,1\}$ with $\varepsilon_k\equiv k\pmod 2$, and let $\Gamma_k$ be the union of these
  sets for all $h+2\le i\le h+k$ with $i\equiv h+k\pmod 2$, together with the subspace $S_h$ if $k$ is even. Telescoping shows that the point multiset of $\Gamma_k$ is $\chi_{S_{h+k}}-\varepsilon_k\chi_B+t_k\chi_{A'}$, where
  $t_k=\left(2^{k+2}-4-4\varepsilon_k\right)/3$. In particular, $t_2=4$, $t_3=8$, and $\varepsilon_3=1$.
  Let $\varepsilon:=\varepsilon_{h+1}$. Then $\varepsilon_{h+3}=\varepsilon$ and $\varepsilon_{h+2}=1-\varepsilon$. The numbers $\alpha$, $\beta$, $y_3$, and $y_4$ below depend on $(d,e)$ as follows.
  \begin{center}
    \begin{tabular}{c|cccc}
      $(d,e)$ & $\alpha$ & $\beta$ & $y_3$ & $y_4$\\
      \hline
      $(0,0)$ & $0$ & $0$ & $\left(2^{h+2}+1+\varepsilon\right)/3$ & $1$\\
      $(0,1)$ & $0$ & $2$ & $1+\varepsilon$ & $1$\\
      $(1,0)$ & $2$ & $0$ & $\left(2^{h+1}+5-\varepsilon\right)/3$ & $3$\\
      $(1,1)$ & $1$ & $3$ & $1+\varepsilon$ & $1$
    \end{tabular}
  \end{center}
  The numbers $y_3$ are positive integers, since $2^{h+1}\equiv 1+\varepsilon\pmod 3$. Take
  \begin{enumerate}
    \item[(1)] $\Gamma_{h+3}$, $\alpha$ copies of $\Gamma_{h+2}$, and $\beta$ copies of $\Gamma_{h+1}$,
    \item[(2)] $2d-\alpha$ copies of a lifted MRD partition of the points of $F'\backslash E$ and $4e-\beta$ copies of a lifted MRD partition of the points of $F\backslash B$ (Lemma~\ref{lemma_moves}.(a)),
    \item[(3)] $2^h$ copies of $\Gamma_3$ and $(1-d)2^{h+1}+2d-\alpha$ copies of $\Gamma_2$,
    \item[(4)] $y_3$ copies of each $(h-1)$-space $X$ with $A'<X<B$, which by Lemma~\ref{lemma_moves_general}.(d) with $(a,b)=(h-2,h+1)$ have the point multiset $3\chi_B+4\chi_{A'}$, and $y_4$ copies of each hyperplane of $B$.
  \end{enumerate}
  The coefficients of $\chi_V$, $\chi_{F'}$, $\chi_F$, $\chi_{S_{h+3}}$, and $\chi_E$ in the point multiset are $1$, $2d$, $4e$, $2^h$, and $(1-d)2^{h+1}$. The coefficient of $\chi_{A'}$ is
  $t_{h+3}+\alpha t_{h+2}+\beta t_{h+1}+2^{h+3}+4\left((1-d)2^{h+1}+2d-\alpha\right)+4y_3$ and the one of $\chi_B$ is
  $-\varepsilon-\alpha(1-\varepsilon)-\beta\varepsilon-(4e-\beta)-2^h+3y_3+\left(2^h-1\right)y_4$. Inserting $t_{h+3}=\left(2^{h+5}-4-4\varepsilon\right)/3$, $t_{h+2}=\left(2^{h+4}-8+4\varepsilon\right)/3$, and
  $t_{h+1}=\left(2^{h+3}-4-4\varepsilon\right)/3$ shows that these coefficients are $2^{h+5}$ and $(1-e)2^{h+2}$ in all four cases. For example, for $(d,e)=(0,0)$ they are $\left(2^{h+5}-4-4\varepsilon\right)/3+2^{h+4}+4y_3=2^{h+5}$ and
  $-\varepsilon-2^h+\left(2^{h+2}+1+\varepsilon\right)+2^h-1=2^{h+2}$.
\end{proof}

\section{Proofs for Section~\ref{sec_partitions}}\label{app_sec5}
The following lemma contains the first two standard equations for the hyperplane multiplicities and a second moment bound.
\begin{Lemma}
  Let $\cN$ assign an integer $\cN(P)$ to every point $P$ of $\PG(v-1,2)$, let $\#\cN:=\sum_P\cN(P)$, $\sigma:=\sum_P\cN(P)^2$, and $w(H):=\#\cN-\sum_{P\in H}\cN(P)$ for every
  hyperplane $H$. Then $\sum_H w(H)=2^{v-1}\#\cN$ and $\sum_H w(H)^2=2^{v-2}\left(\#\cN^2+\sigma\right)$. If $2^r$ divides $w(H)$ for all $H$, then
  $$
    2^{v-2}\left(\left(\#\cN-2^r(2i+1)\right)^2+\sigma-2^{2r}\right)\ge 2^{2r}i(i+1)\quad\text{for all }i\in\Z.
  $$
  In particular, if a $2^r$-divisible set of $n$ points contains $a$ pairwise disjoint subspaces of dimension $r$, then $\left(n-2a[r]_2-2^r(2i+1)\right)^2+n\ge 2^{2r}$
  for all $i\in\Z$, with equality only if $i\in\{-1,0\}$.
\end{Lemma}
\begin{proof}
  Since $w(H)=\sum_{P\notin H}\cN(P)$, a point lies outside $2^{v-1}$ hyperplanes, and two distinct points lie outside $2^{v-2}$ common hyperplanes, we obtain the two sums.
  If $2^r$ divides $w(H)$, then the inequality follows by summing $\left(w(H)-2^ri\right)\left(w(H)-2^r(i+1)\right)\ge 0$ over the $2^v-1$ hyperplanes and using
  $4i(i+1)=(2i+1)^2-1$.
  For the last part let $\cN:=\chi_\cS-2\sum_{j=1}^a\chi_{U_j}$ for the set $\cS$ and the subspaces $U_1,\dots,U_a$. Every hyperplane contains $[r]_2$ or $[r-1]_2$ points of
  $U_j$. Hence $2^r$ divides $w(H)$, and $\cN(P)\in\{0,\pm 1\}$ gives $\#\cN=n-2a[r]_2$ and $\sigma=n$.
\end{proof}

\subsection{Types of \texorpdfstring{$4$}{4}-divisible sets}\label{app_sec5_4div}
\begin{Lemma}
  \label{lemma_2_div_lines}
  If $\cS$ is a set of points with type $2^{m_2}$ in $\PG(v-1,2)$, then we have $m_2\ge 5$ and each of these cases occurs.
\end{Lemma}
\begin{proof}
  Lemma~\ref{lemma_4_div_gf2} implies $m_2\ge 5$ and Lemma~\ref{lemma_2_div_gf4} gives the corresponding constructions via field reduction.
\end{proof}
Due to Lemma~\ref{lemma_2_div_gf2} there are no sets of points of type $2^{m_2}1^1$ or of type $2^{m_2}1^2$ in $\PG(v-1,2)$. 
\begin{Lemma}
  \label{lemma_types_2_13t}
  If $\cS$ is a set of points with type $2^{m_2}1^{3t}$ in $\PG(v-1,2)$, then we have $m_2\ge \max\{5-t,0\}$ and each of these cases occurs.
\end{Lemma}
\begin{proof}
  Lemma~\ref{lemma_4_div_gf2} implies $m_2\ge 5-t$ and the corresponding constructions can be obtained from Lemma~\ref{lemma_2_div_lines} by replacing each of $t$ lines by its 
  three points.
\end{proof}
Each $2$-divisible set of three points over $\F_2$ is the characteristic function of a line, see, e.g., \cite{projective_divisible_binary_codes}.
\begin{Lemma}
  \label{lemma_2_div_lines_p4}
  If $\cS$ is a set of points with type $2^{m_2}1^4$ in $\PG(v-1,2)$, then we have $m_2=1$ or $m_2\ge 4$ and each of these cases occurs.
\end{Lemma}
\begin{proof}
  Lemma~\ref{lemma_4_div_gf2} implies $m_2=1$ or $m_2\ge 4$. For $m_2=1$ we can consider the characteristic function of a plane which we partition into a line and four points.
  For $m_2=4$ and $m_2=5$ examples in $\PG(4,2)$ and $\PG(6,2)$ are obtained from the classification of projective $4$-divisible binary codes \cite{projective_divisible_binary_codes} together with the partition integer linear program of Appendix~\ref{app_sec5_ilp}.\footnote{A geometric construction can be described as follows. Choose two affine solids $S_1\backslash E_1$ and $S_2\backslash E_2$ for solids $S_1,S_2$ and 
  planes $E_1,E_2$. Now choose $S_2$ such that it intersects $S_1$ in $E_1$ and $E_2$ such that it intersects $E_1$ in a line $L$. So, all points of $S_1\backslash L$ are contained in the support 
  of $\cS:=\chi_{S_1\backslash E_1}+\chi_{S_2\backslash E_2}$ and choosing a line spread of $S_1$ including line $L$ gives an example of type $2^4 1^4$ in $\PG(v-1,2)$ for $v\ge 5$.}
  For $m_2\ge 6$ we can start with an example consisting of $m_2-1$ lines, see Lemma~\ref{lemma_2_div_lines}, add a plane, and replace the plane 
  by a vector space partition of type $2^1 1^4$.
\end{proof}
Each $2$-divisible set of four points over $\F_2$ is the characteristic function of an affine plane, see, e.g., \cite{projective_divisible_binary_codes}.

\begin{Lemma}
  \label{lemma_types_2_143t}
  If $\cS$ is a set of points with type $2^{m_2}1^{4+3t}$ in $\PG(v-1,2)$, where $t\ge 1$, then we have $m_2\ge \max\{4-t,0\}$ or $\left(m_2,t\right)=(0,1)$ and each of these cases occurs.
\end{Lemma}
\begin{proof}
  If $t=1$, then $m_2=0$ can be realized by the characteristic function of a plane. Besides the case $\left(m_2,t\right)=(0,1)$, 
  Lemma~\ref{lemma_4_div_gf2} gives $m_2\ge 4-t$. The other necessary constructions can be obtained from Lemma~\ref{lemma_2_div_lines_p4} by replacing $t$ lines by the $3t$ points they contain.
\end{proof}

\begin{Lemma}
  \label{lemma_exclusion_1}
  There is no set of points in $\PG(v-1,2)$ with type $2^1 1^5$ or type $2^3 1^5$.
\end{Lemma}
\begin{proof}
  The projective $4$-divisible binary codes of the lengths $8$ and $14$ are classified, see \cite{bouyukliev2021computer,projective_divisible_binary_codes}. The partition integer linear program of Appendix~\ref{app_sec5_ilp} applied to each of them shows that the type does not occur.
\end{proof}
  
\begin{Lemma}[E.g.~{\cite[Lemma 12]{kurz2022vector}}]
  \label{lemma_exclusion_2}
  There does not exist a set of points in $\PG(v-1,2)$ of type $2^4 1^5$.
\end{Lemma}
\begin{proof}
  The projective $4$-divisible binary codes of length $17$ are classified, see \cite{bouyukliev2021computer,projective_divisible_binary_codes}. The partition integer linear program of Appendix~\ref{app_sec5_ilp} applied to each of them shows that the type does not occur.
\end{proof}

\begin{Lemma}
  \label{lemma_exclusion_3}
  There exists a set of points in $\PG(6,2)$ of type $2^5 1^5$.
\end{Lemma}  
\begin{proof}
  Such a set exists in $\PG(6,2)$ by the classification of projective $4$-divisible binary codes \cite{projective_divisible_binary_codes} together with the partition integer linear program of Appendix~\ref{app_sec5_ilp}.
\end{proof}  

\begin{Lemma}
  \label{lemma_types_2_15}
  If $\cS$ is a set of points with type $2^{m_2}1^5$ in $\PG(v-1,2)$, then we have $m_2\ge 5$ and each of these cases occurs.
\end{Lemma}
\begin{proof}
  The cases $m_2\in\{0,2\}$ are excluded by Lemma~\ref{lemma_4_div_gf2}. The cases $m_2\in\{1,3,4\}$ are excluded by Lemma~\ref{lemma_exclusion_1} and Lemma~\ref{lemma_exclusion_2}. For $m_2=5$ see Lemma~\ref{lemma_exclusion_3}. For $m_2\in\{6,\dots,10\}$ there exist examples in $\PG(4,2)$ and $\PG(5,2)$ by the classification of projective $4$-divisible binary codes \cite{projective_divisible_binary_codes} together with the partition integer linear program of Appendix~\ref{app_sec5_ilp}.
  For $m_2>10$ we can argue recursively, using 
  a line spread of $\PG(3,2)$ consisting of five lines.
\end{proof}

\begin{Lemma}
  \label{lemma_types_2_18}
  If $\cS$ is a set of points with type $2^{m_2}1^8$ in $\PG(v-1,2)$, then we have $m_2=0$ or $m_2\ge 2$ and each of these cases occurs.
\end{Lemma}
\begin{proof}
  An affine solid is an example for $m_2=0$ (indeed it is the unique example). Since there is no $4$-divisible set of cardinality $11$, see Lemma~\ref{lemma_4_div_gf2}, $m_2=1$ is impossible. 
  For $m_2=2$ we can consider the union of two disjoint planes and replace both planes by a vector space partition of type $2^1 1^4$.
  Examples of types $2^3 1^8$ in $\PG(5,2)$ and $2^4 1^8$ in $\PG(6,2)$ are obtained from the classification of projective $4$-divisible binary codes \cite{projective_divisible_binary_codes} together with the partition integer linear program of Appendix~\ref{app_sec5_ilp}.
  For $m_2\ge 5$ take an affine solid and a set of type $2^{m_2}$ (Lemma~\ref{lemma_2_div_lines}) in subspaces that intersect trivially.
\end{proof}

\begin{Lemma}
  \label{lemma_types_2_111}
  If $\cS$ is a set of points with type $2^{m_2}1^{11}$ in $\PG(v-1,2)$, then we have $m_2\ge 1$ and each of these cases occurs. 
  If $\cS$ is a set of points with type $2^{m_2}1^{11+3t}$ in $\PG(v-1,2)$, where $t\ge 1$, then we have $m_2\ge 0$ and each of these cases occurs.
\end{Lemma}
\begin{proof}
  Consider a set of points with type $2^{m_2'} 1^8$ in $\PG(v-1,2)$ with the same number of points and replace a suitable number of lines by their contained three points.
\end{proof}

\begin{Proposition}
  \label{prop_types_2_1_large}
  If $m_1\ge 14$, then there exists a set of points with type $2^{m_2}1^{m_1}$ in $\PG(v-1,2)$ for all non-negative integers $m_2$ and a suitably large dimension $v$.  
\end{Proposition}
\begin{proof}
  Let $m_1\ge 14$. If $m_1=3t$, then $t\ge 5$ and Lemma~\ref{lemma_types_2_13t} gives all $m_2\ge\max\{5-t,0\}=0$. If $m_1=4+3t$, then $t\ge 4$ and Lemma~\ref{lemma_types_2_143t} gives all $m_2\ge\max\{4-t,0\}=0$.
  If $m_1=11+3t$, then $t\ge 1$ and Lemma~\ref{lemma_types_2_111} gives all $m_2\ge 0$. These three cases cover all $m_1\ge 14$.
\end{proof}

\subsection{Types of \texorpdfstring{$8$}{8}-divisible sets with \texorpdfstring{$n\le 51$}{n<= 51}}\label{app_sec5_8div}

\begin{Lemma}[E.g.~{\cite[Lemma 11.2]{kurz2021divisible}}, {\cite[Proposition 23]{kurz2022vector}}]
  \label{lemma_n15}
  Let $\cS$ be an $8$-divisible set of points in $\PG(v-1,2)$ with cardinality $15$. Then $\cS$ can be of type $3^1 1^8$ or $2^5$ but not of type $3^2 1^1$, $3^1 2^1 1^5$, or $3^1 2^2 1^2$.
\end{Lemma}
\begin{proof}
  The projective $8$-divisible binary codes of length $15$ are classified, see \cite{betsumiya2012triply,projective_divisible_binary_codes}. The partition integer linear program of Appendix~\ref{app_sec5_ilp} applied to each of them yields exactly the stated types.
\end{proof} 
We can also replace each $(i-1)$-space by a vector space 
partition of $\PG(i-1,q)$ (working for general fields $\F_q$). We call these types derived types and usually do not mention them. In our case we have types $2^1 1^{12}$ and $2^i 1^{3(5-i)}$ 
for $0\le i\le 4$.   

\begin{Lemma}
  \label{lemma_n16}
  Let $\cS$ be an $8$-divisible set of points in $\PG(v-1,2)$ with cardinality $16$. Then $\cS$ can be of type $1^{16}$ but not of type $2^1 1^{13}$.
\end{Lemma}
\begin{proof}
  The projective $8$-divisible binary codes of length $16$ are classified, see \cite{betsumiya2012triply,projective_divisible_binary_codes}. The partition integer linear program of Appendix~\ref{app_sec5_ilp} applied to each of them yields exactly the stated types.
\end{proof}
Consequently, types from which $2^1 1^{13}$ would be derived, such as $3^1 1^9$ or $2^2 1^{10}$, cannot be realized either.

\begin{Lemma}[E.g.~{\cite[Lemma 11.6]{kurz2021divisible}}]
  \label{lemma_n30}
  Let $\cS$ be an $8$-divisible set of points in $\PG(v-1,2)$ with cardinality $30$. Then $\cS$ can be of type $3^2 1^{16}$, $3^1 2^5 1^8$, or $2^{10}$ but not of type $3^3 1^9$, $3^2 2^1 1^{13}$, 
  or $3^1 2^6 1^5$.
\end{Lemma}
\begin{proof}
  The projective $8$-divisible binary codes of length $30$ are classified, see \cite{betsumiya2012triply,projective_divisible_binary_codes}. The partition integer linear program of Appendix~\ref{app_sec5_ilp} applied to each of them yields exactly the stated types.
\end{proof}
By the classification of the projective $8$-divisible binary codes of length $30$, see \cite{betsumiya2012triply,projective_divisible_binary_codes}, each $8$-divisible set of points in $\PG(v-1,2)$ with cardinality $30$ is the sum of the characteristic functions of two disjoint solids.

\begin{Lemma}
  \label{lemma_n31}
  Let $\cS$ be an $8$-divisible set of points in $\PG(v-1,2)$ with cardinality $31$. Then $\cS$ can be of type $3^1 2^{8}$ or $2^9 1^4$ but not of type $3^2 1^{17}$ or $2^{10} 1^1$.
\end{Lemma}
\begin{proof}
  The projective $8$-divisible binary codes of length $31$ are classified, see \cite{betsumiya2012triply,projective_divisible_binary_codes}. The partition integer linear program of Appendix~\ref{app_sec5_ilp} applied to each of them yields exactly the stated types.
\end{proof}
  
\begin{Lemma}
  \label{lemma_n32}
  Let $\cS$ be an $8$-divisible set of points in $\PG(v-1,2)$ with cardinality $32$. Then $\cS$ can be of type $2^8 1^8$ but not of type $3^1 1^{25}$ or $2^9 1^5$.
\end{Lemma}  
\begin{proof}
  The projective $8$-divisible binary codes of length $32$ are classified, see \cite{betsumiya2012triply,projective_divisible_binary_codes}. The partition integer linear program of Appendix~\ref{app_sec5_ilp} applied to each of them yields exactly the stated types.
\end{proof}  

\begin{Lemma}[E.g.~the proof of {\cite[Lemma 15]{kurz2022vector}}]
  \label{lemma_n45}
  Let $\cS$ be an $8$-divisible set of points in $\PG(v-1,2)$ with cardinality $45$. Then $\cS$ can be of type $3^3 2^8 $, $3^2 2^9 1^{4}$, $3^1 2^{10} 1^8$ or $2^{15}$ but not of type 
  $3^4 1^{17}$, $3^2 2^{10} 1^{1}$, or $3^1 2^{11} 1^5$.
\end{Lemma}
\begin{proof}
  The projective $8$-divisible binary codes of length $45$ are classified, see \cite{betsumiya2012triply,projective_divisible_binary_codes}. The partition integer linear program of Appendix~\ref{app_sec5_ilp} applied to each of them yields exactly the stated types.
\end{proof}

\begin{Lemma}
  \label{lemma_n46}
  Let $\cS$ be an $8$-divisible set of points in $\PG(v-1,2)$ with cardinality $46$. Then $\cS$ can be of type $3^2 2^8 1^{8}$, $3^1 2^{13} $, or $3^0 2^{14} 1^4$ but not of type 
  $3^3 1^{25}$, $3^2 2^9 1^{5}$, or $3^0 2^{15} 1^1$.
\end{Lemma}  
\begin{proof}
  The projective $8$-divisible binary codes of length $46$ are classified, see \cite{betsumiya2012triply,projective_divisible_binary_codes}. The partition integer linear program of Appendix~\ref{app_sec5_ilp} applied to each of them yields exactly the stated types.
\end{proof}

\begin{Lemma}
  \label{lemma_n47}
  Let $\cS$ be an $8$-divisible set of points in $\PG(v-1,2)$ with cardinality $47$. Then $\cS$ can be of type $3^2 2^9 1^6$, $3^1 2^{12} 1^4$, or $3^0 2^{14} 1^5$  but not of type 
  $3^3 1^{26}$, $3^2 2^{10} 1^{3}$, $3^1 2^{13} 1^1$, or $3^0 2^{15} 1^2$.
\end{Lemma}  
\begin{proof}
  The projective $8$-divisible binary codes of length $47$ are classified, see \cite{betsumiya2012triply,projective_divisible_binary_codes}. The partition integer linear program of Appendix~\ref{app_sec5_ilp} applied to each of them yields exactly the stated types.
\end{proof}  

\begin{Lemma}
  \label{lemma_n48}
  Let $\cS$ be an $8$-divisible set of points in $\PG(v-1,2)$ with cardinality $48$. Then $\cS$ can be of type $3^1 2^{11} 1^8$ or $3^0 2^{16}$  but not of type 
  $3^2 2^{0} 1^{34}$ or $3^1 2^{12} 1^5$.
\end{Lemma}  
\begin{proof}
  The projective $8$-divisible binary codes of length $48$ are classified, see \cite{betsumiya2012triply,projective_divisible_binary_codes}. The partition integer linear program of Appendix~\ref{app_sec5_ilp} applied to each of them yields exactly the stated types.
\end{proof}

\begin{Lemma}
  \label{lemma_n49}
  Let $\cS$ be an $8$-divisible set of points in $\PG(v-1,2)$ with cardinality $49$. Then $\cS$ can be of type $3^0 2^{14} 1^7$ but not of type 
  $3^1 1^{42}$ or $3^0 2^{15} 1^4$.
\end{Lemma}  
\begin{proof}
  A generator matrix of an $8$-divisible set of points that can be partitioned into type $3^0 2^{14} 1^7$ is given in Appendix~\ref{app_examples}. For the non-existence results all projective $8$-divisible $[49,k]_2$ codes are checked, see Appendix~\ref{app_sec5_ilp}.
\end{proof}  
  
\begin{Lemma}
  \label{lemma_n50}
  Let $\cS$ be an $8$-divisible set of points in $\PG(v-1,2)$ with cardinality $50$. Then $\cS$ can be of type $3^0 2^{10} 1^{20}$ but not of type 
  $3^1 1^{43}$ or $3^0 2^{11} 1^{17}$.
\end{Lemma}  
\begin{proof}
  The projective $8$-divisible binary codes of length $50$ are classified, see \cite{bouyukliev2021computer,kurz2021divisible}. The partition integer linear program of Appendix~\ref{app_sec5_ilp} applied to each of them yields exactly the stated types.
\end{proof}  
  
\begin{Lemma}[E.g.~{\cite[Lemma 24]{honold2018partial}}]
  \label{lemma_n51}
  Let $\cS$ be an $8$-divisible set of points in $\PG(v-1,2)$ with cardinality $51$. Then $\cS$ can be of type $3^0 2^{17}$ but not of type 
  $3^1 1^{44}$.
\end{Lemma}
\begin{proof}
  The projective $8$-divisible binary codes of length $51$ are classified, see \cite{bouyukliev2021computer,kurz2021divisible}. The partition integer linear program of Appendix~\ref{app_sec5_ilp} applied to each of them yields exactly the stated types.
\end{proof}  

\section{Computational details}\label{app_comp}
This appendix collects explicit examples, tables, and computer results that support statements of the main text.

\subsection{Classification data for Section~\ref{sec_isomorphism_types}}\label{app_sec3}
The $3$-divisible $[n,k]_3$ codes with $n\le 13$ have been enumerated with \texttt{LinCode}, see Table~\ref{table_3_divisible}.

For $q=3$ and $n=9$ the generator matrix $\left(\begin{smallmatrix}111110100\\011221010\\122220001
\end{smallmatrix}\right)$ gives the multiset $2\left(\chi_L-\chi_P\right)+\left(\chi_{L'}-\chi_P\right)$ for two lines $L\neq L'$ through a point $P$, which is an example of the former type with a doubled line. For $n=10$ there exist three non-canonical examples: the projective $[10,4,\{6,9\}]_3$ two-weight code, which geometrically corresponds to an ovoid, and two 
non-projective codes given by generator matrices $\left(\begin{smallmatrix}1111100100\\0022211010\\0211201001\end{smallmatrix}\right)$ and $\left(\begin{smallmatrix}1111101000\\0112210100\\
1222200010\\0221110001\end{smallmatrix}\right)$. For $n=11$ another projective two-weight code appears: the $[11,5,\{6,9\}]_3$ from family RT6 \cite{calderbank1986geometry}, which is the dual of the 
$[11,6,5]_3$ (ternary) Golay code. Besides this code there are four further non-canonical examples for $n=11$, all of them non-projective. Hence there are five non-canonical examples for $n=11$ in total.

\begin{table}[htp]
  \begin{center}
    \begin{tabular}{rrrrrrrrrrrr}
      \hline 
      k/n & 3 & 4 & 5 & 6 & 7 & 8 & 9 & 10 & 11 & 12 & 13 \\ 
      \hline
      1 & 1 &   & & 1 & & & 1 & & & 1 & \\
      2 &   & 1 & & 1 & 1 & 1 & 2 & 2 & 1 &  4 &  3 \\
      3 &   &   & &   & 1 & 1 & 3 & 4 & 5 & 10 & 15 \\
      4 &   &   & &   &   & 1 & 1 & 3 & 5 & 16 & 26 \\ 
      5 &   &   & &   &   &   &   &   & 3 &  7 & 17 \\ 
      6 &   &   & &   &   &   &   &   &   &  3 &  4 \\ 
      \hline
      $\sum$ & 1 & 1 & 0 & 2 & 2 & 3 & 7 & 9 & 14 & 41 & 65 \\ 
      \hline
    \end{tabular}
    \caption{Number of $3$-divisible $[n,k]_3$ codes for $n\le 13$.} 
    \label{table_3_divisible} 
  \end{center}
\end{table} 

\begin{table}[htp]
  \begin{center}
    \scriptsize\setlength{\tabcolsep}{3pt}
    \begin{tabular}{rrrrrrrrrrrrrrrrrrrrrrrr}
      \hline 
      k/n & 9 & 12 & 13 & 18 & 21 & 22 & 24 & 25 & 26 & 27 & 30 & 31 & 33 & 34 & 35 & 36 & 37 & 38 &  39 &  40 &  42 & 43 & 44 \\  
      \hline
      1 &   1 &    &    &  1 &    &    &    &    &    &  1 &    &    &    &    &    &  1 &    &    &     &     &     &    &    \\
      2 &     &  1 &    &  1 &  1 &    &  1 &    &    &  2 &  2 &    &  1 &    &    &  4 &    &    &   3 &     &   2 &    &    \\
      3 &     &    &  1 &    &  1 &  1 &  1 &  1 &  1 &  3 &  4 &  2 &  5 &  2 &  1 & 10 &  2 &  1 &  15 &   6 &  17 &  6 &  2 \\
      4 &     &    &    &    &    &  1 &  1 &  1 &  1 &  3 &  3 &  3 &  5 &  5 &  4 & 22 &  7 &  6 &  34 &  25 &  52 & 32 & 14 \\
      5 &     &    &    &    &    &    &    &  1 &  1 &  1 &    &  1 &  3 &  4 &  4 & 13 & 10 & 12 &  41 &  40 &  44 & 40 & 22 \\
      6 &     &    &    &    &    &    &    &    &  1 &    &    &    &    &  1 &  3 &  4 &  3 & 10 &  23 &  30 &  15 & 16 & 17 \\ 
      7 &     &    &    &    &    &    &    &    &    &    &    &    &    &    &  1 &    &  1 &  3 &   8 &  10 &     &  3 &  6 \\ 
      8 &     &    &    &    &    &    &    &    &    &    &    &    &    &    &    &    &    &  1 &   2 &   1 &     &    &  1 \\
      9 &     &    &    &    &    &    &    &    &    &    &    &    &    &    &    &    &    &    &   1 &     &     &    &    \\
      \hline
      $\sum$ & 1 & 1 & 1 & 2 &  2 &  2 &  3 &  3 &  4 & 10 &  9 &  6 & 14 & 12 & 13 & 54 & 23 & 33 & 127 & 112 & 130 & 97 & 62 \\ 
      \hline
    \end{tabular}
    \caption{Number of $9$-divisible $[n,k]_3$ codes for $n\le 44$.} 
    \label{table_9_divisible} 
  \end{center}
\end{table} 

Among the $127$ codes with $n=39$ in Table~\ref{table_9_divisible}, exactly $46$ are canonical, i.e., three planes or a $3$-fold line plus three $9$-fold points.

\subsection{Known sets \texorpdfstring{$\cN_q^{\Delta,h}$}{N\textunderscore{}q\textasciicircum{}(Delta,h)} and supporting computations}\label{app_N_sets}
For $1\le l\le 12$ the sets of Theorem~\ref{thm_N_explicit} agree with the minimal generators of $\{n\,:\,3n\in \mathcal{L}_2^{l+1}\}$ computed directly from the generators of $\mathcal{L}_2^{l+1}$.
The non-existence of a $3$-divisible multiset of $11$ lines in $\PG(v-1,2)$, see Theorem~\ref{thm_N_t1}, follows from a computer search.
For the line spreads of Proposition~\ref{prop_line_spreads_PG52} and $q=2$ the implication (ii)$\Rightarrow$(i) can be checked by computer.
Up to the action of $\operatorname{GL}(6,2)$ we may fix two lines $L_1$, $L_2$ of $\cS$. A computer search shows that exactly $64$ normal line spreads contain $L_1$ and $L_2$,
and these coincide with the $64$ Desarguesian line spreads containing $L_1$ and $L_2$.

The lengths $73$ for $(q,h,\Delta)=(4,3,24)$ and $142$ for $(4,3,28)$ are realized by multisets of lines in $\PG(4,4)$ found as follows. Using Lemma~\ref{lemma_dual_cover} we searched for multisets of $n$ lines in $\PG(4,4)$ such that every point is covered a number of times congruent to $n$ modulo $\Delta$ and that are invariant under the cyclic group generated by a random matrix of order $3$ resp.\ $5$. For $73$ (group order $3$) the weights are $48^1\,72^{340}$, and for $142$ (group order $5$, one point on $30$ lines and all others on two lines) the weights are $112^1\,140^{340}$. So $73$ is a length for every divisor $\Delta$ of $24$, and $142$ for every divisor $\Delta$ of $28$.

\begin{Remark}
  \label{rem_apery_nonminimal}
  For $j=q$ and $h=j+2$ the set $\mathcal{G}_q^{h,j}$ is not minimal, since its element $(q^2+1)q^{2q+2}+1+q^{q+2}$ with $(i,u,f)=(0,1,q)$ equals $q^{2h}+q^h+1+q^{2q+2}$, e.g., $337=273+64$ for $(q,h,j)=(2,4,2)$, although Equation~(\ref{eq_apery_qhj}) still holds by Theorem~\ref{thm_apery_qhj}. A computation for $(q,h,j)\in\{(2,4,2),(3,5,3)\}$ shows that this element ($337$ resp.\ $65854$) is the only element of $\mathcal{G}_q^{h,j}$ that is not a minimal generator.
  For $j\ge q+1$ Equation~(\ref{eq_apery_qhj}) fails in general, e.g., for $16$ residues if $(q,j)=(2,3)$ and $5\le h\le 15$, and the proof of Theorem~\ref{thm_apery_qhj} breaks down, see Appendix~\ref{app_sec4_apery}. Lemma~\ref{lemma_cone_family} with $j=3$ does not apply to the parameters excluded in Corollary~\ref{cor_qh3}.(c), for any $h\ge 4$, since $[h]_qZ+W$ decreases in $h$ for $Z<0$.
\end{Remark}

\begin{Remark}
  \label{rem_N_qh3}
  For $\Delta=q^{h+j}$ with $j\in\{2,3\}$, $h\ge 4$ for $j=2$, and $h\ge 5$ for $j=3$, the minimal generators of the semigroup of necessary conditions are the elements of $\mathcal{G}_q^{h,j}$ of Corollary~\ref{cor_gens_qhj}, which was checked by computations for $q\in\{3,4,5\}$, $j=3$ (for $q=3$, $h=5$ except $65854$, see Remark~\ref{rem_apery_nonminimal}). For $q=2$, $j=3$, and $6\le h\le 8$ there are $19$ generators, namely $2^{h+3}$, $2^{h+3}+8$, $2^{2h}+2^h+1$, $3\cdot 2^{h+3}+4$, $3\cdot 2^{h+3}+20$, $5\cdot 2^{h+3}+2+4e+16g$ with $e\in\{0,1\}$, $0\le g\le 2$, and $5\cdot 2^{h+3}+1+2d+4e+32g$ with $d,e,g\in\{0,1\}$.
  The generators of layers $i\ge 1$ are lengths, and in layer $i=0$ so are those with $f\le 1$ for $j=2$ (Corollary~\ref{cor_qh3}.(a),(b)) and those covered by Corollary~\ref{cor_qh3}.(a),(c) for $j=3$, $q\ge 3$. For $j=2$ the generators with $f=2$ are lengths for $h=4$ (Proposition~\ref{thm_qh2_f2}) and for $(q,h)=(2,5)$ (Proposition~\ref{prop_qh2_packing}), and open otherwise. For $q=2$, $j=3$, the lengths $5\cdot 2^{h+3}+1+2d+4e$, $5\cdot 2^{h+3}+18$, and $5\cdot 2^{h+3}+22$ are realized (Proposition~\ref{thm_qh3_q2}, Corollary~\ref{cor_qh3}.(b)), and $5\cdot 2^{h+3}+34$, $5\cdot 2^{h+3}+38$, $5\cdot 2^{h+3}+33+2d+4e$ are open. For $j=3$, $q\ge 3$, the remaining generators are open.
\end{Remark}

The method of Subsection~\ref{sub_qh23} gives
\begin{equation}\label{eq_N_243}
\begin{aligned}
\cN_3^{243,3}&=\{243,252,757,975,978,1002,1003,1005,1006,\\
&\phantom{=\{}1702,1703,1705,1706,1708,1709,1730,1733,1736\},\\
\cN_3^{243,2}&=\{243,270,273,820,981,982,990,991,992,1000,1001,\\
&\phantom{=\{}1071,1072,1073,1081,1082,1091,1704,1707,1713,1714,\\
&\phantom{=\{}1716,1717,1722,1723,1724,1725,1726,1727\},
\end{aligned}
\end{equation}
where $820$ is realized by the field reduction of all points of $\PG(3,9)$.

\begin{Remark}
  \label{rem_N_qh2_small}
  For $(q,h)=(4,3)$ and $\Delta=q^{h+2}=1024$ the computation of the minimal generators of the semigroup of necessary conditions gives the $36$ numbers
  $q^{h+2}$, $q^{h+2}+q^2$, $q^{2h}+q^h+1$, $(q+1)q^{h+2}+cq+\eta q^3$, $(q+1)q^{h+2}+q^3+cq+1$ for $1\le c\le 3$ and $\eta\in\{0,1\}$, and $(2q+1)q^{h+2}+fq^3+aq+b$ with $0\le a\le 3$ and
  $(f,b)\in\{(0,1),(0,2),(0,3),(1,2),(1,3),(2,3)\}$. All of them are lengths, so they form $\cN_4^{1024,3}$.
\end{Remark}

{\small
\begin{longtable}{rrr>{\raggedright\arraybackslash}p{6.0cm}>{\raggedright\arraybackslash}p{5.2cm}}
  \caption{Known sets $\cN_q^{\Delta,h}$. The marks in the status column are explained below the table.}
  \label{table_N_sets}\\
    \hline
    $q$ & $h$ & $\Delta$ & $\cN_q^{\Delta,h}$ & status\\
    \hline
  \endfirsthead
    \hline
    $q$ & $h$ & $\Delta$ & $\cN_q^{\Delta,h}$ & status\\
    \hline
  \endhead
    \hline
  \endfoot
    $q$ & $h$ & $\Delta\mid q^h$ & $\{\Delta,q^h+1\}$ & Theorem~\ref{thm_delta_divides_qh}\\
    2 & 2 & $2^l$, $l\ge 1$ & see Theorem~\ref{thm_N_explicit} & Theorem~\ref{thm_N_explicit}\\
    2 & 2 & 8 & $\{8,10,21,25\}$ & Theorem~\ref{thm_N_explicit}\\
    2 & 2 & 16 & $\{16,20,21,50\}$ & Theorem~\ref{thm_N_explicit}\\
    2 & 3 & 16 & $\{16,18,49,53,73\}$ & Theorem~\ref{thm_N_qh1}\\
    2 & 3 & 32 & $\{32,36,73,98,106,107,$ $161,163\}$ & Proposition~\ref{prop_N_qh2}\\
    2 & 3 & 64 & $\{64,72,73,196,212,214,$ $321,322,323,325,326,327\}$ & computer\\
    2 & 4 & 64 & $\{64,68,194,202,273,321,323,$ $329,331,339\}$ & computer, $329,331$ by Corollary~\ref{cor_qh3}.(b), $339$ by Proposition~\ref{thm_qh2_f2}\\
    2 & 5 & 128 & $\{128,132,386,394,641,643,$ $649,651,657,659,1057\}$ & computer, $649,651$ by Corollary~\ref{cor_qh3}.(b), $657,659$ by Proposition~\ref{prop_qh2_packing}, completeness by Corollary~\ref{cor_gens_qhj}\\
    2 & 4 & 32 & $\{32,34,97,101,273\}$ & Theorem~\ref{thm_N_qh1}\\
    2 & 5 & 64 & $\{64,66,193,197,1057\}$ & Theorem~\ref{thm_N_qh1}\\
    3 & 2 & 27 & $\{27,30,91,109,110,119\}$ & Theorem~\ref{thm_N_qh1}\\
    3 & 2 & 81 & $\{81,90,91,327,330,357,$ $568,569,571,572,574,575\}$ & Proposition~\ref{prop_N_qh2}\\
    3 & 2 & 243 & see Equation~(\ref{eq_N_243}) & computer\\
    3 & 3 & 243 & see Equation~(\ref{eq_N_243}) & computer\\
    3 & 3 & 81 & $\{81,84,325,326,334,335,757\}$ & Theorem~\ref{thm_N_qh1}\\
    3 & 4 & 243 & $\{243,246,973,974,982,983,6643\}$ & Theorem~\ref{thm_N_qh1}\\
    3 & 4 & 729 & $\{729,738,2919,2922,2946,2949,$ $5104,5105,5107,5108,5110,5111,$ $5131,5132,5134,5135,5137,5138,$ $5158,5159,5161,5162,5164,5165,$ $6643\}$ & computer, $5131,\dots,5138$ by Corollary~\ref{cor_qh3}.(b), $5158,\dots,5165$ by Proposition~\ref{thm_qh2_f2}, completeness by Corollary~\ref{cor_gens_qhj}\\
    4 & 3 & 1024 & see Remark~\ref{rem_N_qh2_small} & computer, $9359$ by Proposition~\ref{thm_qh2_f2}\\
    4 & 2 & 32 & $\{32,34,161,169,273\}$ & Subsection~\ref{sub_qh23}\\
    4 & 2 & 64 & $\{64,68,273,321,322,323,338,339\}$ & Theorem~\ref{thm_N_qh1}\\
    5 & 2 & 125 & $\{125,130,651,751,752,$ $753,754,777,778,779\}$ & Theorem~\ref{thm_N_qh1}\\
    $q$ & 2 & $q^3$ & $\cG_q^2\backslash\{q^4+q^3+q^2+1\}$ & Theorem~\ref{thm_N_qh1}\\
    $q$ & $\ge 3$ & $q^{h+1}$ & $\cG_q^h$ & Theorem~\ref{thm_N_qh1}\\
    $q$ & 2 & $q^4$ & see Proposition~\ref{prop_N_qh2}.(a) & Proposition~\ref{prop_N_qh2} for $q\in\{3,4,5,7\}$\\
    2 & 2 & 6 & $\{6,7,15\}$ & Proposition~\ref{prop_N_coprime}\\
    2 & 2 & 12 & $\{12,14,15,31,37\}$ & computer\\
    2 & 2 & 24 & $\{24,28,30,31,63,74,75,77\}$ & computer\\
    2 & 2 & 48 & $\{48,56,60,62,63,127,148,$ $150,151,154,155,157\}$ & computer\\
    3 & 2 & 2 & $\{2,13\}$ & Proposition~\ref{prop_N_coprime}\\
    3 & 2 & 6 & $\{6,13,20\}$ & Proposition~\ref{prop_N_coprime}\\
    3 & 2 & 18 & $\{18,20,39,73,121\}$ & computer\\
    2 & 3 & 7 & $\{7,15,31\}$ & Proposition~\ref{prop_N_coprime_dual}\\
    2 & 3 & 14 & $\{14,15,31,63\}$ & Proposition~\ref{prop_N_coprime_dual}\\
    2 & 3 & 28 & $\{28,30,31,63,85,127\}$ & computer\\
    2 & 4 & 3 & $\{3,25,38\}$ & Propositions~\ref{prop_N_coprime_dual} and~\ref{prop_25_38}\\
    2 & 4 & 6 & $\{6,25,38,51\}$ & Propositions~\ref{prop_N_coprime_dual} and~\ref{prop_25_38}\\
    2 & 4 & 12 & $\{12,25,38,51\}$ & Propositions~\ref{prop_N_coprime_dual} and~\ref{prop_25_38}\\
    4 & 3 & 3 & $\{3,73\}$ & computer\\
    4 & 3 & 7 & $\{7,85,142,341\}$ & computer\\
    $q$ & $h$ & $p^a[h]_q$ & see Theorem~\ref{thm_N_t1}, $p^a\mid q$ & Theorem~\ref{thm_N_t1}\\
    2 & 4 & 5 & $\{5,21\}$ & Proposition~\ref{prop_N_coprime_dual}\\
    3 & 2 & 4 & $\{4,13\}$ & Proposition~\ref{prop_N_coprime_dual}\\
    3 & 2 & 12 & $\{12,13,40\}$ & Proposition~\ref{prop_N_coprime_dual}\\
    3 & 2 & 36 & $\{36,39,40,121,145,146\}$ & computer\\
    4 & 2 & 5 & $\{5,21\}$ & Proposition~\ref{prop_N_coprime_dual}\\
    4 & 2 & 10 & $\{10,21,85\}$ & Proposition~\ref{prop_N_coprime_dual}\\
    4 & 2 & 20 & $\{20,21,85\}$ & Proposition~\ref{prop_N_coprime_dual}\\
    4 & 2 & 40 & $\{40,42,85,201,341\}$ & computer\\
    5 & 2 & 2 & $\{2,31\}$ & Proposition~\ref{prop_N_coprime_dual}\\
    5 & 2 & 3 & $\{3,31\}$ & Proposition~\ref{prop_N_coprime_dual}\\
    5 & 2 & 6 & $\{6,31\}$ & Proposition~\ref{prop_N_coprime_dual}\\
\end{longtable}

}

\smallskip\noindent
\begin{Proposition}
  \label{prop_N_qh2}
  \begin{enumerate}
    \item[(a)] For $q\in\{3,4,5,7\}$ we have
               \begin{eqnarray*}
                 \cN_q^{q^4,2}&=&\left\{q^4,q^4+q^2,q^4+q^2+1\right\}\cup\left\{q^5+q^4+cq\,:\,1\le c\le q-1\right\}\\
                 & &\cup\left\{q^5+q^4+q^3+cq\,:\,2\le c\le q-1\right\}\\
                 & &\cup\left\{2q^5+q^4+aq+b\,:\,0\le a\le q-1,\,1\le b\le q-1\right\}.
               \end{eqnarray*}
    \item[(b)] $\cN_2^{32,3}=\{32,36,73,98,106,107,161,163\}$.
  \end{enumerate}
\end{Proposition}
\begin{proof}
  (a) The lengths are realized by a $q^4$-fold line, $q^2$ copies of a line spread of $\PG(3,q)$, the field reduction of all points of $\PG(2,q^2)$, $q$ copies of the codes of Theorem~\ref{thm_C} and
  Proposition~\ref{prop_ring_family}, and Corollary~\ref{cor_qh3}.(a).

  (b) The lengths are realized by a $32$-fold plane, four copies of a plane spread of $\PG(5,2)$, the field reduction of $\PG(2,8)$, two copies of
  the codes of length $49$ and $53$ from Theorems~\ref{thm_C} and~\ref{thm_plus_q2}, and Corollary~\ref{cor_qh3}.(a). For $107$ let $S_4<S_5<S_7<S_8=V=\F_2^8$ be subspaces with $\dim(S_i)=i$, $K<S_4$ with $\dim(K)=2$,
  $W$ a complement of $K$ in $S_4$, and $u\in S_5\backslash S_4$. Take
  \begin{itemize}
    \item each $3$-dimensional subspace of $S_4$ containing $K$ once and each other one three times,
    \item the four subspaces $K+\langle u+w\rangle$, $w\in W$,
    \item two lifted MRD partitions of $S_7\backslash S_4$, and
    \item one lifted MRD partition of $S_8\backslash S_5$.
  \end{itemize}
  This gives $3+36+4+32+32=107$ subspaces with point multiset of type $[8]+2[7]+16[4]$
  (Lemma~\ref{lemma_moves_general}.(g) with $(k,u,m_0)=(1,1,3)$), which is $128$-divisible by Lemma~\ref{lemma_type_divisible} with $r=8$. The multiset of subspaces is $32$-divisible by Lemma~\ref{lemma_exclusion}.

  Completeness follows as in the proof of Theorem~\ref{thm_N_qh1}.(b), using the lengths of $q^{2h+1}$-divisible multisets of points, by a computer calculation of the minimal generators of
  the semigroup given by the necessary conditions.
\end{proof}

\begin{table}[htp]
  \begin{center}
  \small
  \begin{tabular}{llllr}
    \hline
    construction & $(q,h)$ & $[n,r/h,d]_q^h$ & weight enumerator $\prod w^{A_w}$ & $d_G(n)$\\
    \hline
    Theorem~\ref{thm_C}, $c=1$ & $(2,2)$ & $[25,5/2,8]$ & $8^{3}\,16^{12}\,24^{16}$ & $19$\\
    Proposition~\ref{prop_ring_family}, $c=1$ & $(2,2)$ & $[29,5/2,8]$ & $8^{3}\,24^{28}$ & $22$\\
    Proposition~\ref{prop_ring_family}, $c=2$ & $(2,2)$ & $[30,5/2,16]$ & $16^{3}\,24^{28}$ & $23$\\
    Corollary~\ref{cor_qh3}.(a), $j=2$, $(a,b)=(0,1)$ & $(2,2)$ & $[81,3,16]$ & $16^{3}\,32^{4}\,48^{24}\,80^{32}$ & $61$\\
    Corollary~\ref{cor_qh3}.(a), $j=2$, $(a,b)=(1,1)$ & $(2,2)$ & $[83,3,16]$ & $16^{1}\,32^{2}\,48^{28}\,80^{32}$ & $63$\\
    Theorem~\ref{thm_C}, $c=1$ & $(2,3)$ & $[49,7/3,16]$ & $16^{7}\,32^{24}\,48^{96}$ & $43$\\
    Theorem~\ref{thm_plus_q2}, $c=1$ & $(2,3)$ & $[53,7/3,16]$ & $16^{3}\,32^{4}\,48^{120}$ & $46$\\
    Proposition~\ref{prop_ring_family}, $1\le c\le 6$ & $(2,3)$ & $[120+c,7/3,16c]$ & $(16c)^{7}\,112^{120}$ & $105+c$\\
    Proposition~\ref{prop_N_qh2}.(b), $n=107$ & $(2,3)$ & $[107,8/3,32]$ & $32^{1}\,64^{14}\,96^{240}$ & $93$\\
    Corollary~\ref{cor_qh3}.(a), $j=2$, $(a,b)=(0,1)$ & $(2,3)$ & $[161,8/3,32]$ & $32^{7}\,64^{8}\,96^{48}\,160^{192}$ & $141$\\
    Corollary~\ref{cor_qh3}.(a), $j=2$, $(a,b)=(1,1)$ & $(2,3)$ & $[163,8/3,32]$ & $32^{1}\,64^{6}\,96^{56}\,160^{192}$ & $143$\\
    Theorem~\ref{thm_C}, $1\le c\le 2$ & $(3,2)$ & $[108+c,5/2,27c]$ & $(27c)^{8}\,81^{72}\,108^{162}$ & $96+c$\\
    Proposition~\ref{prop_ring_family}, $1\le c\le 3$ & $(3,2)$ & $[117+c,5/2,27c]$ & $(27c)^{8}\,108^{234}$ & $104+c$\\
    Theorem~\ref{thm_C}, $1\le c\le 3$ & $(4,2)$ & $[320+c,5/2,64c]$ & $(64c)^{15}\,256^{240}\,320^{768}$ & $300+c$\\
    Proposition~\ref{prop_ring_family}, $1\le c\le 4$ & $(4,2)$ & $[336+c,5/2,64c]$ & $(64c)^{15}\,320^{1008}$ & $315+c$\\
    $\Delta=32$, $n=161$, Subsection~\ref{sub_qh23} & $(4,2)$ & $[161,5/2,64]$ & $64^{15}\,128^{240}\,160^{768}$ & $150$\\
    $\Delta=32$, $n=169$, Subsection~\ref{sub_qh23} & $(4,2)$ & $[169,5/2,64]$ & $64^{15}\,160^{1008}$ & $158$\\
    Theorem~\ref{thm_C}, $c=1$ & $(2,4)$ & $[97,9/4,32]$ & $32^{15}\,64^{48}\,96^{448}$ & $91$\\
    Theorem~\ref{thm_plus_q2}, $c=1$ & $(2,4)$ & $[101,9/4,32]$ & $32^{3}\,64^{12}\,96^{496}$ & $94$\\
    \hline
  \end{tabular}
  \caption{Parameters and weight enumerators of some of the additive codes of Subsections~\ref{sub_constructions} and~\ref{sub_qh23}. $d_G(n)$ is the largest minimum distance compatible with the Griesmer type bound
           for additive $[n,r/h]_q^h$ codes \cite{ball2025griesmer,kurz2024additive}. The weights of all hyperplanes were computed for every entry. The exponents in the weight enumerators are numbers of codewords.}
  \label{table_few_weight}
  \end{center}
\end{table}

\subsection{Explicit partitions for Section~\ref{sec_partitions}}\label{app_examples}
For $n\ge 46$ we give a generator matrix of the $8$-divisible set of points, whose columns are ordered such that a partition into planes, lines, and points is visible: the first $7a$ columns form $a$ planes with $7$ consecutive columns each, the next $3b$ columns form $b$ lines with $3$ consecutive columns each, and the remaining columns are points.
\begin{itemize}
  \item $n=46$ (Lemma~\ref{lemma_n46}), type $3^22^81^8$:
  $$
    \left(\begin{smallmatrix}
1110010111100000011011000011011011011011110000\\
0111100001111010100001111011001101111000011000\\
0111100010101100000001111011000010101101111110\\
1010101011010111000000000000011010101100000000\\
0000000001111011001111011011011001111011110000\\
0000000101001100011000001100001100010100111100\\
0000000001111000011000000001100011010111110000\\
0000000111100000011001100000001111000001010101
\end{smallmatrix}\right)
  $$
  type $3^12^{13}$:
  $$
    \left(\begin{smallmatrix}
1110001000110000110000000000110110110110000000\\
0111100110110000011110000000000110000011000000\\
0000000000000110000000000110000000000000110110\\
1010110000011000101110000000000011011110000000\\
0000000000000000000000110011000000000000011110\\
1101010000000000011011000000011011110011000000\\
0000000000000000000000101101000000000000110011\\
1010110101000000011000000000101110101101000000\\
0000000000000101000000000101000000000000011101
\end{smallmatrix}\right)
  $$
  \item $n=47$ (Lemma~\ref{lemma_n47}), types $3^22^91^6$, $3^12^{12}1^4$, and $3^02^{14}1^5$:
  $$
  \left(\begin{smallmatrix}
11110000000000110110110110110110000110000111111\\
00111101111000000011000000110110110011110001111\\
01011010011110000011011110000000011110011010111\\
11110000101011000110101011000011110101011001001\\
00000001100110101011110101110011110110011100010\\
11110001111000000011000101011110011011101100100
\end{smallmatrix}\right)
  $$
  The columns of this matrix are ordered according to the first type. The other two types are obtained from the following partitions, given by column numbers. Type $3^{1}2^{12}1^{4}$: planes (2,15,28,39,40,41,47), lines (16,23,25) (6,9,26) (7,29,35) (3,11,43) (5,13,33) (14,22,36) (8,18,31) (27,32,44) (4,17,19) (12,30,37) (10,42,46) (21,34,45). Type $3^{0}2^{14}1^{5}$: lines (12,16,22) (7,33,40) (5,26,34) (18,19,20) (3,29,44) (10,23,38) (17,25,36) (9,24,47) (8,15,30) (4,27,41) (28,35,46) (1,32,42) (11,21,37) (13,31,45).
  \item $n=48$ (Lemma~\ref{lemma_n48}), type $3^02^{16}$:
  $$
  \left(\begin{smallmatrix}
101110110110110110110110110110110110110110110110\\
110011110000011110110011011011000000011011000110\\
110000000011101011110000011011011110000110110011\\
000110000110000101011011101011000011110011110011\\
110110011000000101110000000101101101101101110110\\
000110000000011000110101110101011011011110011101
\end{smallmatrix}\right)
  $$
  type $3^12^{11}1^8$:
  $$
  \left(\begin{smallmatrix}
000000011011000011011011011011011011011011110000\\
111001000000011000001101111001111001111000111100\\
011110000001101100000010101110101111011011010100\\
101010100000010111011010100000010111001101111110\\
011110010111010100010111000011011000010100011110\\
011110000011000001100010101101101110101110111000\\
000000000000000000000000011011011001111001111011
\end{smallmatrix}\right)
  $$
  \item $n=49$ (Lemma~\ref{lemma_n49}), type $3^02^{14}1^7$:
  $$
  \left(\begin{smallmatrix}
1100001101101101100001101101100001101100001111000\\
0001011100001100111100110000111101100111100011000\\
0001100111100001100000110110111100110000111101100\\
1011100000001100110110001100110110001101011101100\\
0000001100110000110001011010000110000110000100100\\
0000001100000111101011101011100000111101011011100\\
0001100001100000110000111100111011011101011001110\\
0001101100001100111100110001101101100111101000001
\end{smallmatrix}\right)
  $$
  \item $n=50$ (Lemma~\ref{lemma_n50}), type $3^02^{10}1^{20}$:
  $$
  \left(\begin{smallmatrix}
11011011011011011011011011011011111111111100000000\\
00001111011011001100001100001100000011111111110000\\
00000001101100001111011011010100001100111100111100\\
10100011001111010100010101111000110100001101010100\\
00011011000001111000010110110101011001010110011100\\
00011000010111000011010100000000111110011111101110\\
00011011000000000000000001110111011111101111011101\\
00011000010100010101101110111010100101110101011000
\end{smallmatrix}\right)
  $$
\end{itemize}

\subsection{Enumeration and integer linear programs for Section~\ref{sec_partitions}}\label{app_sec5_ilp}
\begin{table}
  \begin{center}
  \begin{tabular}{rrrrrrrrrrrrrrr}
  \hline
  $n$ & $\sum$ & $k=$ & 4 & 5 & 6 & 7 & 8 & 9 & 10 & 11 & 12 & 13 & 14 & 15 \\ 
  \hline
  15 &    1 & & 1 &   &   &    &     &     &      &     &     &    &   &   \\
  16 &    1 & &   & 1 &   &    &     &     &      &     &     &    &   &   \\
  30 &    1 & &   &   &   &    &   1 &     &      &     &     &    &   &   \\
  31 &    6 & &   & 1 & 1 &  1 &   2 &   1 &      &     &     &    &   &   \\
  32 &   11 & &   &   & 2 &  2 &   3 &   3 &    1 &     &     &    &   &   \\
  45 &    6 & &   &   &   &    &   2 &   1 &    1 &   1 &   1 &    &   &   \\
  46 &   51 & &   &   &   &    &   6 &  18 &   14 &   8 &   4 &  1 &   &   \\
  47 &  856 & &   &   & 1 & 11 & 100 & 299 &  274 & 122 &  40 &  8 & 1 &   \\
  48 & 2973 & &   &   & 1 & 15 & 211 & 921 & 1071 & 529 & 173 & 44 & 7 & 1 \\
  49 &  120 & &   &   &   &    &   9 &  38 &   44 &  21 &   7 &  1 &   &  \\
  50 &    1 & &   &   &   &    &   1 &     &      &     &     &    &   &   \\
  51 &    1 & &   &   &   &    &   1 &     &      &     &     &    &   &   \\
  \hline
  \end{tabular}
  \caption{Enumeration of projective $8$-divisible $[n,k]_2$ codes for $n<60$. The rows $n\le 48$ are taken from \cite{projective_divisible_binary_codes}, where
  they are derived from the classification of the triply even codes of length $48$ in \cite{betsumiya2012triply}. The rows $n=49,50,51$ are taken from \cite[Table 5.1]{kurz2021divisible}, where they were enumerated with the software package \texttt{LinCode}.}
  \label{table_enumeration_8_div}
  \end{center}
\end{table}  

\noindent\emph{Partition ILP.} Let $\cS$ be the point set of a projective $8$-divisible code from Table~\ref{table_enumeration_8_div}, and let $\mathcal{E}$ and $\mathcal{K}$ be the sets of planes and lines of $\PG(v-1,2)$ that are contained in $\cS$.
Then $\cS$ has a partition of type $3^a2^b1^c$ with $7a+3b+c=\#\cS$ if and only if there are binary variables $x_E$ for $E\in\mathcal{E}$ and $y_L$ for $L\in\mathcal{K}$ with $\sum_Ex_E=a$, $\sum_Ly_L=b$, such that every point of $\cS$ lies in at most one selected plane or line.
The remaining $c$ points of $\cS$ are the points of the partition. The lemmas for $n\in\{45,\dots,50\}$ use this integer linear program for all codes of the respective row of the table.

\begin{Remark}\label{rem_hyperplane_system}
  For a set $\cS$ of points in $\PG(v-1,2)$ of type $3^a2^b1^c$ and a hyperplane $H$ let $p_H$, $l_H$, $r_H$ denote the numbers of planes, lines, and points of the
  partition that are contained in $H$. Then $\cS(H)=7p_H+3(a-p_H)+3l_H+(b-l_H)+r_H$, $c-r_H$ is even, $2(b-l_H)+(c-r_H)\equiv 0\pmod 4$ if $\cS$ is $8$-divisible, and
  $\cS|_H$ minus the $p_H$ planes is a $4$-divisible set of type $2^{a-p_H+l_H}1^{b-l_H+r_H}$, which has to be realizable according to Corollary~\ref{cor_types_2_1}. Counting the incidences of hyperplanes with single parts and with
  pairs of parts gives ten linear equations for the numbers of hyperplanes of each kind $(p_H,l_H,r_H)$. If for no $v$ with $\log_2(n+1)\le v\le 3a+2b+c$ this integer linear system has a non-negative solution,
  then type $3^a2^b1^c$ does not exist. This excludes the $4$-divisible types $2^11^5$, $2^31^5$, $2^41^5$ and all non-existence statements of Lemmas~\ref{lemma_exclusion_1}, \ref{lemma_exclusion_2}, and \ref{lemma_n15}--\ref{lemma_n51}
  except for the types $3^12^{11}1^5$ ($n=45$), $3^22^91^5$ ($n=46$), $3^12^{12}1^5$ ($n=48$), $3^11^{42}$, $2^{15}1^4$ ($n=49$), $3^11^{43}$, and $2^{11}1^{17}$ ($n=50$).
  Adding the incidence counts with triples and with $4$-subsets of $\cS$, with the number of lines and of $4$-subsets with sum zero in $\cS$ as further unknowns, excludes $3^11^{43}$ for all $v$. For the
  six other types the extended system has no integer solution for $v\ge 13$ ($3^12^{11}1^5$, $3^22^91^5$) or $v\ge 17$ (the other four), so that these statements still rest on the enumeration of the codes in Table~\ref{table_enumeration_8_div}.
\end{Remark}

\subsection{Types for \texorpdfstring{$60\le n\le 70$}{60<= n<= 70}}\label{app_types_60_70}
Further types were found by a computer search. Starting from an $8$-divisible set $\cS$ in $\PG(v-1,2)$, $v\in\{7,8\}$, we repeatedly replace $\cS$ by
$\cS\triangle\cB$, where $\cB$ is a solid, a $4$-space, or an affine $4$-space such that $\cS\cap\cB$ is $4$-divisible, which implies that $\cS\triangle\cB$ is
$8$-divisible. For the sets with $64\le\#\cS\le 69$ we maximize the number of lines for a given number of planes in a packing of $\cS$ by an integer linear program.
Twelve of the resulting sets suffice to generate all new types, e.g., types $3^42^{12}$ ($n=64$, $v=7$), $3^62^8$ and $3^32^{15}$ ($n=66$, $v=8$), and $3^12^{20}$ ($n=67$, $v=8$).
The exclusions in the table include all admissible types for which the linear system of Remark~\ref{rem_hyperplane_system} has no solution. For $n=68$ these are the six types $3^72^b1^{19-3b}$ with $0\le b\le 5$.
For $n=64$ the type $3^62^51^7$ is excluded by Lemma~\ref{lemma_n64}, and hence also $3^62^61^4$, whose derived type it is.

\begin{Lemma}
  \label{lemma_n64}
  Let $k\in\N_{>0}$. The point set $\cS\subseteq\PG(k-1,2)$ of a projective $8$-divisible $[64,k]_2$ code does not admit a partition of type $3^62^51^7$.
\end{Lemma}
\begin{proof}
  The proof is by enumeration of the six planes up to symmetry and ILP computations afterwards.
\end{proof}

\begin{table}[htp]
  \begin{center}
  \begin{tabular}{rrrrrrr}
  \hline
  $n$ & admissible & realized & excluded & open & max.\ $m_3$ realized & max.\ $m_3$ not excluded\\
  \hline
  60 &  92 &  75 &  7 &  10 &  4 &  5 \\
  61 &  93 &  62 & 14 &  17 &  3 &  4 \\
  62 &  95 &  65 & 16 &  14 &  3 &  4 \\
  63 & 105 & 105 &  0 &   0 &  9 &  9 \\
  64 & 107 &  98 &  6 &   3 &  8 &  8 \\
  65 & 107 &  83 &  4 &  20 &  7 &  7 \\
  66 & 110 & 105 &  4 &   1 &  6 &  6 \\
  67 & 114 &  89 &  7 &  18 &  5 &  6 \\
  68 & 115 &  68 &  6 &  41 &  4 &  6 \\
  69 & 119 &  54 &  6 &  59 &  3 &  6 \\
  70 & 128 &  66 &  1 &  61 & 10 & 10 \\
  \hline
  \end{tabular}
  \caption{Partition types $3^{m_3}2^{m_2}1^{m_1}$ of $8$-divisible sets of points in $\PG(v-1,2)$ for $60\le n\le 70$. The admissible types are those with $7m_3+3m_2+m_1=n$ for which the $4$-divisible part $2^{m_2}1^{m_1}$ is realizable by Corollary~\ref{cor_types_2_1}.}
  \label{table_types_60_70}
  \end{center}
\end{table}

\section*{Declaration of AI use}
The author used Claude Sonnet 5.5 (Anthropic) to improve the language and wording of the manuscript and to check proofs for errors, inconsistencies, and notation clashes. All mathematical statements and proofs were obtained and verified by the author, who takes full responsibility for the content of this paper.


\begin{thebibliography}{10}

\bibitem{assmus1992designs}
E.~F. Assmus, Jr. and J.~D. Key.
\newblock {\em Designs and their Codes}, volume 103 of {\em Cambridge Tracts in
  Mathematics}.
\newblock Cambridge University Press, 1992.

\bibitem{bader2011desarguesian}
L.~Bader and G.~Lunardon.
\newblock Desarguesian spreads.
\newblock {\em Ricerche di Matematica}, 60(1):15--37, 2011.

\bibitem{ball2025griesmer}
S.~Ball, M.~Lavrauw, and T.~Popatia.
\newblock Griesmer type bounds for additive codes over finite fields, integral
  and fractional {MDS} codes.
\newblock {\em Designs, Codes and Cryptography}, 93(1):175--196, 2025.

\bibitem{barlotti1974finite}
A.~Barlotti and J.~Cofman.
\newblock Finite {S}perner spaces constructed from projective and affine
  spaces.
\newblock {\em Abhandlungen aus dem Mathematischen Seminar der Universit{\"a}t
  Hamburg}, 40:231--241, 1974.

\bibitem{betsumiya2012triply}
K.~Betsumiya and A.~Munemasa.
\newblock On triply even binary codes.
\newblock {\em Journal of the London Mathematical Society (2)}, 86(1):1--16,
  2012.

\bibitem{blokhuis1994size}
A.~Blokhuis.
\newblock On the size of a blocking set in {PG}(2,{$p$}).
\newblock {\em Combinatorica}, 14(1):111--114, 1994.

\bibitem{blokhuis2004small}
A.~Blokhuis and A.~E. Brouwer.
\newblock Small additive quaternary codes.
\newblock {\em European Journal of Combinatorics}, 25(2):161--167, 2004.

\bibitem{bouyukliev2021computer}
I.~Bouyukliev, S.~Bouyuklieva, and S.~Kurz.
\newblock Computer classification of linear codes.
\newblock {\em IEEE Transactions on Information Theory}, 67(12):7807--7814,
  2021.

\bibitem{calderbank1998quantum}
A.~R. Calderbank, E.~M. Rains, P.~W. Shor, and N.~J.~A. Sloane.
\newblock Quantum error correction via codes over {GF}(4).
\newblock {\em IEEE Transactions on Information Theory}, 44(4):1369--1387,
  1998.

\bibitem{calderbank1986geometry}
R.~Calderbank and W.~M. Kantor.
\newblock The geometry of two-weight codes.
\newblock {\em Bulletin of the London Mathematical Society}, 18(2):97--122,
  1986.

\bibitem{delsarte1970generalized}
P.~Delsarte, J.-M. Goethals, and F.~J. MacWilliams.
\newblock On generalized {R}eed--{M}uller codes and their relatives.
\newblock {\em Information and Control}, 16:403--442, 1970.

\bibitem{dhaeseleer2025generalized}
J.~D'haeseleer and S.~Kurz.
\newblock Generalized {H}amming weights of additive codes and geometric
  counterparts.
\newblock {\em Designs, Codes and Cryptography}, 94(7):142, 2026.

\bibitem{heden2012survey}
O.~Heden.
\newblock A survey of the different types of vector space partitions.
\newblock {\em Discrete Mathematics, Algorithms and Applications},
  4(1):1250001, 2012.

\bibitem{projective_divisible_binary_codes}
D.~Heinlein, T.~Honold, M.~Kiermaier, S.~Kurz, and A.~Wassermann.
\newblock Projective divisible binary codes.
\newblock In {\em The Tenth International Workshop on Coding and Cryptography},
  pages 1--10, 2017.
\newblock arXiv:1703.08291.

\bibitem{honold2018partial}
T.~Honold, M.~Kiermaier, and S.~Kurz.
\newblock Partial spreads and vector space partitions.
\newblock In M.~Greferath, M.~O. Pav{\v{c}}evi{\'c}, N.~Silberstein, and
  M.~{\'A}. V{\'a}zquez-Castro, editors, {\em Network Coding and Subspace
  Designs}, Signals and Communication Technology, pages 131--170. Springer,
  Cham, 2018.

\bibitem{honold2019lengths}
T.~Honold, M.~Kiermaier, S.~Kurz, and A.~Wassermann.
\newblock The lengths of projective triply-even binary codes.
\newblock {\em IEEE Transactions on Information Theory}, 66(5):2713--2716,
  2020.

\bibitem{kasami1970weight}
T.~Kasami and N.~Tokura.
\newblock On the weight structure of {R}eed--{M}uller codes.
\newblock {\em IEEE Transactions on Information Theory}, 16(6):752--759, 1970.

\bibitem{kiermaier2020lengths}
M.~Kiermaier and S.~Kurz.
\newblock On the lengths of divisible codes.
\newblock {\em IEEE Transactions on Information Theory}, 66(7):4051--4060,
  2020.

\bibitem{kim2009generalized}
J.-L. Kim and X.~Liu.
\newblock A generalized {G}leason--{P}ierce--{W}ard theorem.
\newblock {\em Designs, Codes and Cryptography}, 52(3):363--380, 2009.

\bibitem{korner2024lengths}
T.~K{\"o}rner and S.~Kurz.
\newblock Lengths of divisible codes with restricted column multiplicities.
\newblock {\em Advances in Mathematics of Communications}, 18(2):505--534,
  2024.

\bibitem{kurz2021divisible}
S.~Kurz.
\newblock Divisible codes, 2021.
\newblock arXiv:2112.11763.

\bibitem{kurz2022vector}
S.~Kurz.
\newblock Vector space partitions of {$\mathrm{GF}(2)^8$}.
\newblock {\em Serdica Journal of Computing}, 16(2):71--100, 2022.

\bibitem{kurz2024additive}
S.~Kurz.
\newblock Additive codes attaining the {Griesmer} bound, 2024.
\newblock arXiv:2412.14615.

\bibitem{kurz2024lengths}
S.~Kurz.
\newblock Lengths of divisible codes: the missing cases.
\newblock {\em Designs, Codes and Cryptography}, 92(8):2367--2378, 2024.

\bibitem{kurz2021generalization}
S.~Kurz and S.~Mattheus.
\newblock A generalization of the cylinder conjecture for divisible codes.
\newblock {\em IEEE Transactions on Information Theory}, 68(4):2281--2289,
  2022.

\bibitem{macwilliams1977theory}
F.~J. MacWilliams and N.~J.~A. Sloane.
\newblock {\em The Theory of Error-Correcting Codes}.
\newblock North-Holland, Amsterdam, 1977.

\bibitem{mateva2009line}
Z.~T. Mateva and S.~T. Topalova.
\newblock Line spreads of {$\mathrm{PG}(5,2)$}.
\newblock {\em Journal of Combinatorial Designs}, 17(1):90--102, 2009.

\bibitem{rolland2010second}
R.~Rolland.
\newblock The second weight of generalized {R}eed--{M}uller codes in most
  cases.
\newblock {\em Cryptography and Communications}, 2:19--40, 2010.
\newblock arXiv:0902.0058.

\bibitem{ward1981divisible}
H.~N. Ward.
\newblock Divisible codes.
\newblock {\em Archiv der Mathematik}, 36(1):485--494, 1981.

\end{thebibliography}
\end{document}